\documentclass[a4paper,reqno]{amsart}

\usepackage{amsmath,amssymb,amsthm,mathtools,mathrsfs,bm,empheq}
\usepackage{enumitem}
\usepackage{array}
\usepackage{xcolor}
\usepackage{tikz}
\usetikzlibrary{arrows.meta}
\usepackage{aliascnt}
\usepackage{etoolbox}
\makeatletter
\renewenvironment{proof}[1][\proofname]{\par
  \pushQED{\qed}
  \normalfont\topsep6\p@\@plus6\p@\relax
  \trivlist
  \item[\hskip\labelsep\bfseries\itshape #1\@addpunct{.}]\ignorespaces
}{\popQED\endtrivlist\@endpefalse}
\makeatother
\usepackage[mathlines]{lineno}


\usepackage[
colorlinks,
citecolor=magenta,
urlcolor=blue,
linkcolor=blue,
allcolors=black,
backref=page
]{hyperref}
\numberwithin{equation}{section}
\makeatletter
\patchcmd{\@settitle}{\uppercasenonmath\@title}{}{}{}
\patchcmd{\section}{\scshape}{}{}{}
\makeatother

\allowdisplaybreaks
\newcommand{\iu}{\mathrm{i}}
\newcommand{\N}{\mathbb N}
\newcommand{\Z}{\mathbb Z}
\newcommand{\R}{\mathbb R}
\newcommand{\C}{\mathbb C}
\newcommand{\T}{\mathbb T}
\newcommand{\Pp}{\mathbb P}
\newcommand{\Ee}{\mathbb E}

\newcommand{\Ree}{\operatorname{Re}}

\newcommand{\norm}[1]{\left\|#1\right\|}
\newcommand{\set}[1]{\left\{#1\right\}}
\newcommand{\starpow}[2]{\left\{#1\right\}^{\ast^{[#2]}}}

\theoremstyle{definition}
\newtheorem{definition}{Definition}[section]
\newaliascnt{example}{definition}
\newtheorem{example}[example]{Example}
\aliascntresetthe{example}
\newaliascnt{problem}{definition}

\aliascntresetthe{problem}
\theoremstyle{plain}
\newaliascnt{lemma}{definition}
\newtheorem{lemma}[lemma]{Lemma}
\aliascntresetthe{lemma}
\newaliascnt{proposition}{definition}
\newtheorem{proposition}[proposition]{Proposition}
\aliascntresetthe{proposition}
\newaliascnt{theorem}{definition}
\newtheorem{theorem}[theorem]{Theorem}
\aliascntresetthe{theorem}
\newtheorem{introtheorem}{Theorem}
\renewcommand{\theintrotheorem}{\ifcase\value{introtheorem}\or P\or N\else\Alph{introtheorem}\fi}
\newaliascnt{corollary}{definition}
\newtheorem{corollary}[corollary]{Corollary}
\aliascntresetthe{corollary}
\theoremstyle{definition}
\newaliascnt{remark}{definition}
\newtheorem{remark}[remark]{Remark}
\aliascntresetthe{remark}

\title[]
{Persistence and Breakdown of the Initial Gaussian Rate Function for 1D Quasi-Periodic Cubic NLS}

\author{Fei Xu}
\address{\scriptsize (F. Xu)~Institute of Mathematics, Jilin University, Changchun 130012, P.R. China.}
\email{\color{magenta}stuxuf@outlook.com}
\thanks{The first author (F.~Xu) was supported in part by 
the National Natural Science Foundation of China (Grant No.~12501235), 
the China Postdoctoral Science Foundation (Grant Nos.~BX20240138 and 2025M783104) 
and the Jilin Provincial Postdoctoral Science Foundation (2025). 
}

\author{Yong Li}
\address{\scriptsize  (Y. Li)~Institute of Mathematics, Jilin University, Changchun 130012, P.R. China. 
School of Mathematics and Statistics, Center for Mathematics and Interdisciplinary Sciences, Northeast Normal University, Changchun, Jilin 130024, P.R. China.}
\email{\color{magenta}liyong@jlu.edu.cn}
\thanks{The second author (Y.~Li) is the corresponding author. He was supported in part by the National Natural Science Foundation of China (Grant Nos. 12531009 and 12471183).}

\date{\today}



\begin{document}

\begin{abstract}
We study rare events for the spatial supremum of the solution amplitude
in the 1D cubic nonlinear Schr\"odinger equation with weak nonlinearity
and Gaussian quasi-periodic initial data.
Under suitable polynomial decay assumptions in Fourier space and a Diophantine condition controlling small divisors arising from spatial quasi-periodicity, we prove that the initial Gaussian rate function persists uniformly on the \emph{critical time scale} with a \emph{small enough} logarithmic correction prefactor.
Using the \emph{Akhmediev breather}, we prove that the same Gaussian rate function \emph{fails} to hold uniformly on this time scale with a \emph{large} prefactor.
This also implies that the initial Gaussian rate function cannot, in general, persist uniformly on arbitrarily long polynomial time scale beyond the critical case.
\end{abstract}

\nolinenumbers
\maketitle
\begingroup
\makeatletter
\settowidth{\@tempdima}{\bfseries1.\quad}
\edef\tocsectionnumberwidth{\the\@tempdima}
\edef\tocleveltwoindent{\the\dimexpr1pc+\@tempdima\relax}
\renewcommand{\l@section}{\@tocline{1}{6pt}{1pc}{\tocsectionnumberwidth}{\bfseries}}
\renewcommand{\l@subsection}{\@tocline{2}{0pt}{\tocleveltwoindent}{}{}}
\renewcommand{\l@subsubsection}{\@tocline{3}{0pt}{\tocleveltwoindent}{}{}}
\patchcmd{\@tocline}{\hfil}
  {\nobreak\leaders\hbox{\normalfont\kern.25em.\kern.25em}\hfill}{}{}
\makeatother
\tableofcontents
\endgroup

\nolinenumbers
\section{Introduction}
\subsection{The Problem and Main Results}
\label{subsec:intro-main-results}
\emph{Does the initial Gaussian rate function persist beyond the critical
time scale?}
We address this question for the spatial supremum of the solution amplitude
in the 1D cubic nonlinear Schr\"odinger equation with Gaussian
quasi-periodic initial data.

We consider the 1D cubic NLS
\begin{equation}\label{eq:qp-cubic-NLS}
 \iu\partial_tu+\partial_{xx}u+\varepsilon^\alpha|u|^2u=0,
 \qquad (t,x)\in\R\times\R,
\end{equation}
which describes the evolution of a complex-valued field
$u(t,\cdot):\R_x\to\C$ on its interval of existence.
We work in a \emph{weakly nonlinear setting}, with $0<\varepsilon\ll1$
and fixed $\alpha>1$. The prefactor $\varepsilon^\alpha$ measures the
nonlinear strength.

The initial data are Gaussian and quasi-periodic with
$\nu\ge2$ basic frequencies and admit the Fourier expansion
\begin{equation}\label{eq:qp-randomized-initial-data}
 u(0,x)=\sum_{n\in\Z^\nu}c(n)g_ne^{\iu\langle n,\omega\rangle x}.
\end{equation}
Here $\omega\in\R^\nu$ is a fixed nonresonant frequency vector;
see \autoref{subsec:quasi-periodic-functions}.
The deterministic Fourier amplitudes $c(n)$ satisfy product polynomial
decay; see \autoref{subsubsec:deterministic-fourier-amplitudes}.
The $g_n$ are independent and identically distributed (i.i.d.) standard
complex Gaussian random variables on
$(\Omega,\mathcal F,\Pp)$;
see \autoref{subsubsec:complex-gaussian-exceptional-set}.
The total variance is
$S:=\sum_{n\in\Z^\nu}|c(n)|^2\in(0,\infty)$.

Motivated by the large deviation approach to the study of \emph{rogue waves}
in \cite{DGV18,GGKS23}, we fix the observation scale at
$\varepsilon^{-1/2}$. The quadratic exponent in the initial Gaussian tails
gives large deviations with speed
$\varepsilon^{-1}$; see \eqref{eq:nf-gaussian-tail-lower-bound}.
At this amplitude, the effective nonlinear strength is
$\mu_\varepsilon:=\varepsilon^{\alpha-1}$.
Following \cite[Remark~1.7]{GGKS23}, \emph{criticality} refers to \emph{the onset of
the first nonlinear effects comparable in size to the observation amplitude}.
The cubic scaling gives the accumulated nonlinear strength
$t\varepsilon^\alpha(\varepsilon^{-1/2})^2=t\mu_\varepsilon$.
Thus the \emph{subcritical}, \emph{critical} and \emph{supercritical} time regimes correspond
to $t\ll\mu_\varepsilon^{-1}$, $t\asymp\mu_\varepsilon^{-1}$ and
$t\gg\mu_\varepsilon^{-1}$, respectively.

To describe events of large amplitude at this observation scale, fix a
threshold $z>0$. For a real random variable $f$, define the
exceedance event
\begin{equation}\label{eq:nonlinear-spatial-exceedance-event}
 \mathscr A^f_{\varepsilon,z}:=\{f>z\varepsilon^{-1/2}\}.
\end{equation}
Its complement is denoted by
\begin{equation}\label{eq:spatial-exceedance-complement}
 \mathscr B^f_{\varepsilon,z}:=(\mathscr A^f_{\varepsilon,z})^c
 =\{f\le z\varepsilon^{-1/2}\}.
\end{equation}
In what follows, the observable is the spatial supremum of the solution
amplitude for each flow under consideration.

We state our main positive and negative results informally below,
with references to their precise formulations.

\begin{introtheorem}[Persistence of the Initial Gaussian Rate Function]
\label{thm:intro-persistence}
Let $u$ solve \eqref{eq:qp-cubic-NLS} with initial data
\eqref{eq:qp-randomized-initial-data}. Assume that the frequency vector is Diophantine and that the deterministic
Fourier amplitudes have sufficiently fast polynomial decay.
For each fixed $z>0$, there is a small enough constant $c>0$ such that
\begin{equation}\label{eq:intro-gaussian-rate}
 \lim_{\varepsilon\to0^+}
 \varepsilon\log\Pp\!\left(
 \mathscr A_{\varepsilon,z}^{\|u(t)\|_\infty}\right)
 =-\frac{z^2}{S},
\end{equation}
uniformly for $0\le t\le c\mu_\varepsilon^{-1}|\log\mu_\varepsilon|$; see \autoref{thm:log-critical-ldp}.
In particular, this includes the critical time scale
$t\asymp\mu_\varepsilon^{-1}$; see \autoref{thm:first-critical-ldp-local}.
\end{introtheorem}

\begin{introtheorem}[Failure at the Critical Time Scale with a Large Logarithmic Correction Prefactor]
\label{thm:intro-failure}
For every fixed Diophantine frequency vector, one can choose
product polynomial Fourier amplitudes and a fixed $\alpha>1$
sufficiently close to $1$ such that the following holds.
Let $u$ be the corresponding solution of
\eqref{eq:qp-cubic-NLS}--\eqref{eq:qp-randomized-initial-data}.
There exist a fixed prefactor $\mathsf C_{\mathrm{fail}}>0$ and a sequence
$\varepsilon_j\to0^+$ such that,
at the times
$t_j=\mathsf C_{\mathrm{fail}}
\mu_{\varepsilon_j}^{-1}|\log\mu_{\varepsilon_j}|$,
\begin{equation}\label{eq:intro-gaussian-rate-failure}
 \liminf_{j\to\infty}\varepsilon_j\log\Pp\!\left(
 \mathscr A_{\varepsilon_j,2}^{\|u(t_j)\|_\infty}\right)
 \ge-1>-\frac4S.
\end{equation}
The prefactor $\mathsf C_{\mathrm{fail}}$ exceeds the small prefactor in
\autoref{thm:intro-persistence}; see \autoref{rem:counterexample-time-range}.
The bound above rules out uniform persistence of the initial Gaussian
rate function over this interval; see \autoref{thm:general-failure}.
Thus uniform persistence also fails on any polynomial time scale beyond the
critical scale; see
\autoref{cor:counterexample-polynomial-obstruction}.
\end{introtheorem}

\subsection{Background and Related Work}
\label{subsec:intro-background}

\noindent\emph{Rogue Waves and Large Deviations.}
In oceanography, a rogue wave is commonly defined as a \emph{surface wave}
whose height exceeds twice the \emph{significant wave height},
a measure of the surrounding sea state; see \cite{ORBMA13}.
This raises a probabilistic question: \emph{how likely is such an extreme
event under nonlinear wave evolution?}
Dematteis, Grafke and Vanden-Eijnden~\cite{DGV18} studied this question
for a modified NLS model with Gaussian initial data, using
\emph{large deviation optimization} (minimizing the \emph{Gaussian cost} of
initial data that produce a target large wave).
Their approach links the probability of a large wave to the Gaussian
cost of the corresponding initial configurations; see also \cite{DGOV19,DGV19}.
Garrido, Grande, Kurianski and Staffilani~\cite{GGKS23} subsequently
gave a rigorous justification for the weakly nonlinear
1D periodic cubic NLS.

\par\medskip\noindent\emph{Large Deviations for Periodic NLS and Related Models.}
The result of \cite{GGKS23} establishes the Gaussian rate function
through the critical time scale for exponentially decaying
initial Fourier amplitudes.
Subsequent work weakens this decay assumption: Liang and
Wang~\cite{LW25} treated absolutely summable amplitudes, and Fan and
Ye~\cite{FY25} treated square-summable polynomial amplitudes with decay
exponents arbitrarily close to $1/2$ from above. 
An extension with a logarithmic correction and a small enough prefactor is discussed in
\cite[Remark~1.2]{FY25}.
In contrast to these persistence results, Grande~\cite{Gra25} showed that \emph{resonant energy exchange} in a spatially modulated \emph{beating} NLS with two-mode random initial data can \emph{lower} the initial Gaussian rate function.

\par\medskip\noindent\emph{Large Deviations for Other Models.}
Similar questions about time scales and rate functions have been
studied for other dispersive equations.
Berti, Grande, Maspero and Staffilani~\cite{BGMS25} established
a large deviations principle for 2D pure gravity deep water waves
with a periodic horizontal coordinate, combining \emph{normal form theory}
with probabilistic estimates to reach \emph{long time scales}.
Building on their approach, Berforini D'Aquino and Grande~\cite{BDG26} used KdV integrability
to prove persistence of the initial Gaussian rate function
over arbitrarily long polynomial time scales.

\par\medskip\noindent\emph{From Periodic to Quasi-Periodic Spatial Settings.}
Spatially quasi-periodic patterns arise through modulational
instabilities of periodic wavetrains, a mechanism believed to
contribute to rogue wave formation; see \cite{WZ21}.
The evolution of almost periodic data already presents difficult
deterministic problems; see \cite[Problem~1]{Dei07}.
For randomized initial data, we ask whether the initial Gaussian
rate function persists over longer times.
In \cite{XL26}, the present authors established this persistence in the subcritical
regime for higher-dimensional NLS under product polynomial Fourier
decay, without imposing a Diophantine condition.
Here we study \emph{persistence over long times} for the 1D cubic NLS.

\subsection{Main Contributions}
\label{subsec:intro-contributions}

\begin{enumerate}[label=\textup{(\roman*)},leftmargin=*,itemsep=\medskipamount,parsep=\smallskipamount]
\item \emph{Failure over the Critical Time Scale with a Large Logarithmic Correction Prefactor.}

We prove that the initial Gaussian rate function fails to persist uniformly on the critical time scale with a logarithmic correction and a suitable large prefactor; see \autoref{thm:general-failure}.
  
The key is a strict \emph{cost gap}: the breather data have Gaussian cost tending to one and produce normalized amplitudes above two, whereas the initial Gaussian rate function predicts a higher cost \(4/S>3\); see \eqref{eq:general-initial-cost} and \eqref{eq:general-reference-bounds}.

This also rules out uniform persistence on every polynomial time scale beyond the critical scale; see \autoref{cor:counterexample-polynomial-obstruction}.

\item
\emph{Small Frequency Embedding of the Akhmediev Breather and Tube Stability.}

We embed the frequencies of periodic Akhmediev breathers into a fixed
quasi-periodic module using arbitrarily small nonzero physical
frequencies, constructed via the \emph{Pell equation} in the special case
and \emph{Dirichlet-type approximation} in the general case.
This matches the breather's dispersive scale to the
vanishing nonlinear strength; see \autoref{lem:general-small-frequencies}.
Such a construction is unavailable in a fixed frequency module of rank one,
where nonzero frequencies are bounded away from zero, so their squares
cannot match a nonlinear coefficient tending to zero.

We also construct initial neighborhoods for which \emph{tube stability} keeps all corresponding
solutions, including normalized random solutions, in a deterministic
breather tube up to focusing; see \eqref{eq:general-tube-stability}.
The neighborhood probability estimate then converts the cost gap in (i)
into a probability lower bound incompatible with the initial Gaussian rate
function; see \eqref{eq:general-failure-lower-bound}.

\begin{remark}
\label{rem:deterministic-references-random-dynamics}
The second contribution reflects a common philosophy: \emph{controlling
random dynamics through deterministic reference solutions}.
For \emph{typical behavior} and \emph{fluctuations},
Girotti, Grava, McLaughlin and Najnudel~\cite{GGMN26} used a deterministic \emph{soliton
gas} generated by averaged scattering data to establish a Law of
Large Numbers and a Central Limit Theorem. For \emph{rare events}, we use the deterministic
\emph{Akhmediev breather} as a reference solution to establish
\emph{tube stability} and prove failure of the initial Gaussian
rate function at the critical time scale with a large logarithmic correction prefactor. 
\end{remark}

\item \emph{Persistence over the Critical Time Scale with a Small Enough Logarithmic Correction Prefactor.}

For product polynomial Fourier data and Diophantine frequency vectors, we prove persistence of the initial
Gaussian rate function through the critical time scale and its logarithmic
correction with a small enough prefactor; see \autoref{thm:log-critical-ldp}.
This extends the subcritical result of
\cite{XL26} to longer times for 1D quasi-periodic cubic NLS.
\end{enumerate}

\par\addvspace{\medskipamount}

\subsection{Main Difficulties and Proof Strategy}
\label{subsec:intro-difficulties}

\begin{enumerate}[label=\textup{(\roman*)},leftmargin=*,itemsep=\medskipamount,parsep=\smallskipamount]
\item \emph{From Deterministic Breather Focusing to a Probability Lower Bound.}

We seek an initial neighborhood that preserves breather focusing without
introducing an additional leading exponential cost. This is essential
to retain the cost gap underlying the failure result.

On the deterministic side, tube stability keeps solutions close to the
breather until focusing, so their normalized amplitudes still exceed two.
It suffices to impose
$r_j e^{150\widetilde\mu_j\widetilde t_j}\le1/8$; see
\eqref{eq:general-tube-stability}. 

On the probabilistic side, the neighborhood must be large enough to
avoid an additional leading exponential cost. Requiring
$\sqrt{\widetilde\varepsilon_j}/r_j\to0$ ensures that the centered
neighborhood has probability tending to one; see
\eqref{eq:general-centered-probability}.

To meet both requirements, we choose a radius satisfying
\[
 \sqrt{\widetilde\varepsilon_j}\ll r_j
 \le\tfrac18e^{-150\widetilde\mu_j\widetilde t_j}.
\]
Under \eqref{eq:general-parameter-compatibility}, we therefore choose
$r_j=\tfrac18e^{-150\widetilde\mu_j\widetilde t_j}$, as in
\eqref{eq:general-tube-radius}.
The translation estimate \eqref{eq:general-translated-probability} then gives
a probability lower bound with the Gaussian cost of the center,
preserving the cost gap.

\item \emph{Logarithmic Extension: Small Divisors and Weighted Estimates.}

Spatial quasi-periodicity introduces a \emph{small divisor problem} absent in the
periodic setting when the spatial period is fixed. Integration by parts in
the nonresonant cubic term produces $1/\Sigma_3$, where $\Sigma_3$
is the resonance function of \autoref{def:resonance-data}.
For a fixed spatial period $L$, nonzero phases satisfy
$|\Sigma_3|\gtrsim L^{-2}$, so no small divisor problem arises.
In the quasi-periodic setting, however, nonzero phases can be arbitrarily
small, and their reciprocals are therefore not uniformly bounded.

To overcome this small divisor problem, we impose the \emph{Diophantine
condition} \eqref{eq:HDC}, which yields
\[
 \frac{1}{|\Sigma_3|}\le C_{\gamma,\tau}
 \max_{1\le j\le3}\langle n_j\rangle^{2\tau},
\]
where $\tau$ is the Diophantine exponent and
$\langle n\rangle=1+|n|_1$. The \emph{tame estimates} of
\autoref{lem:log-pointwise-tame} absorb this $2\tau$-power loss into
one higher weight norm by choosing $p=q+2\tau$ with $q>\nu$, while
the other two factors require only the lower weight norm.
Assuming a bound on the low weight norm, we obtain the high weight
bound \eqref{eq:log-high-envelope} and the comparison estimate
$E_{\log}(t)\le C\mu_\varepsilon e^{\Lambda\mu_\varepsilon t}$
from \eqref{eq:log-error-growth}. Hence, for
$0\le t\le c\mu_\varepsilon^{-1}|\log\mu_\varepsilon|$, we have
$E_{\log}(t)=\mathcal{O}(\mu_\varepsilon^{1-\Lambda c})$, which tends
to zero uniformly on this interval when $\Lambda c<1$.
This improves the assumed low weight bound and closes the bootstrap argument.
\autoref{lem:log-continuation} then permits extension.

\item \emph{Relative Probability Errors at the Large Deviation Scale.}

To transfer the resonant LDP to the nonlinear solution, we need the
exceptional probability to be negligible relative to the resonant tail.
Since this tail is exponentially small, convergence to zero alone does
not suffice. For each fixed $z>0$, choose
$\delta_\varepsilon=\varepsilon/(1+z^2/S)$, and let $\Omega_\delta$
be the Gaussian exceptional set of \autoref{def:exceptional}.
Let $\mathscr T_\varepsilon=c\mu_\varepsilon^{-1}|\log\mu_\varepsilon|$
and let $\vartheta\in(0,1)$ be the exponent in the comparison estimate.
For the resonant reference field $u_{\mathrm{res}}$ defined in
\eqref{eq:nf-resonant-reconstruction}, we obtain, for both choices of
sign and all sufficiently small $\varepsilon$,
\[
 \sup_{0\le t\le\mathscr T_\varepsilon}
 \frac{\Pp(\Omega_{\delta_\varepsilon})}
 {\Pp\!\left(
   \mathscr A_{\varepsilon,z\pm C\mu_\varepsilon^\vartheta}
   ^{\|u_{\mathrm{res}}(t)\|_\infty}\right)}
 \le e^{-1/(2\varepsilon)}\longrightarrow0;
\]
see \autoref{prop:nf-resonant-law} and
\eqref{eq:log-relative-exception-small}.
The uniform comparison between the two flows shifts the threshold by
a quantity tending to zero. Together with this relative probability
estimate, it transfers the initial Gaussian rate function to the
nonlinear solution uniformly in time.
\end{enumerate}

\par\addvspace{\medskipamount}

\subsection{Organization of the Paper and Proof Routes}
\label{subsec:intro-organization}
\autoref{sec:structures} introduces the quasi-periodic Fourier framework,
random initial data and exceptional set estimates.
\autoref{sec:normal-forms} constructs the resonant flow, establishes its
LDP and transfers it to the nonlinear solution at critical time.
\autoref{sec:log-critical} develops weighted estimates and a bootstrap
argument to obtain logarithmic extension and persistence of the initial
Gaussian rate function. For failure at a larger prefactor,
\autoref{sec:longtime-obstruction} constructs a breather tube in the
special case using the Pell equation. Finally,
\autoref{sec:general-failure} extends the construction to general
Diophantine frequencies through Dirichlet-type approximation.

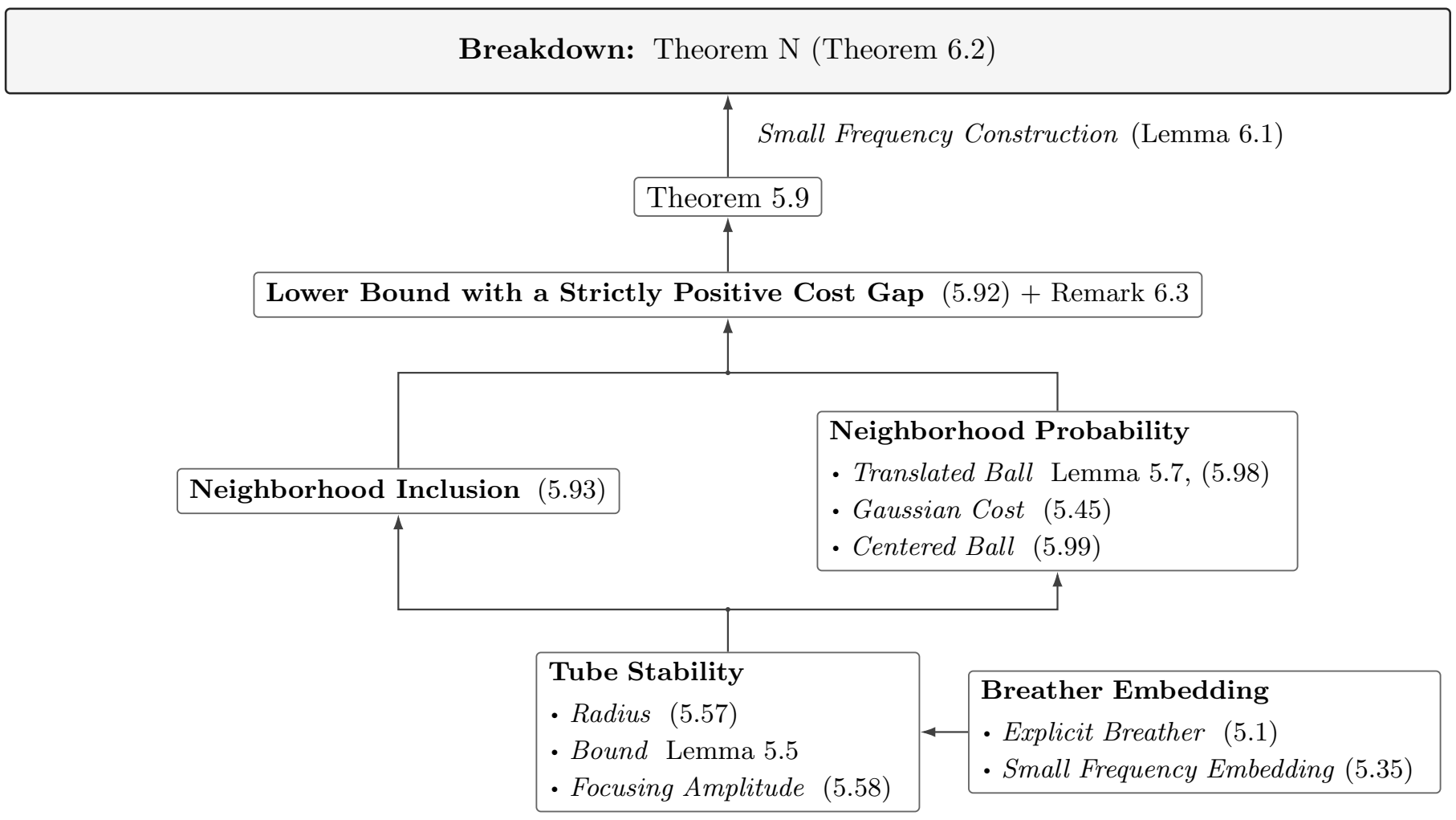
\begin{figure}[!htbp]
\centering
\setlength{\abovecaptionskip}{10pt}
\resizebox{1.01\linewidth}{!}{
\begingroup
\let\breakdownitemlabel\labelitemi
\renewcommand{\labelitemi}{\makebox[.5em][c]{\raisebox{.25ex}{\scalebox{.5}{\breakdownitemlabel}}}}
\begin{tikzpicture}[
 x=1cm,y=.92cm,font=\fontsize{9.5}{11}\selectfont,
 route box/.style={draw=black!60,fill=white,line width=.5pt,
   rounded corners=1.5pt,align=center,text width=4.6cm,
   minimum height=0pt,inner xsep=4pt,inner ysep=3pt,
   execute at begin node={\hyphenpenalty=10000\exhyphenpenalty=10000}},
 argument box/.style={route box,align=left,text width=4cm,
   font=\fontsize{9.5}{11}\selectfont,minimum height=0pt},
 input box/.style={route box,align=left,text width=3.2cm,
   font=\fontsize{9.5}{11}\selectfont,minimum height=0pt},
 route arrow/.style={-{Latex[length=1.8mm,width=1.2mm]},
   draw=black!75,line width=.6pt}
]
 \node[input box,text width=5.3cm,inner xsep=4pt,anchor=east] (breather-embedding) at (13.18,0)
   {\textbf{Breather Embedding}\par
    {\fontsize{9.5}{11}\selectfont
     \begin{itemize}[leftmargin=7pt,labelsep=2.5pt,topsep=3pt,
       itemsep=1.5pt,parsep=0pt,partopsep=0pt]
      \item \mbox{\emph{Explicit Breather}\hspace{.6em}\eqref{eq:explicit-breather}}
      \item \mbox{\emph{Small Frequency Embedding}\hspace{.3em}\eqref{eq:embedded-breather}}
     \end{itemize}}};
 \node[argument box,text width=4.25cm] (tube-stability) at (4.74,0)
   {\textbf{Tube Stability}\par
    {\fontsize{9.5}{11}\selectfont
     \begin{itemize}[leftmargin=7pt,labelsep=2.5pt,topsep=3pt,
       itemsep=1.5pt,parsep=0pt,partopsep=0pt]
      \item \emph{Radius}\hspace{.6em}\eqref{eq:tube-radius}
      \item \emph{Bound}\hspace{.6em}\autoref{lem:breather-tube}
      \item \emph{Focusing Amplitude}\hspace{.6em}\eqref{eq:breather-tube-stability}
     \end{itemize}}};
 \coordinate (neighborhood-fork) at ([yshift=.5cm]tube-stability.north);
 \node[route box,text width={},minimum height=0pt,font=\fontsize{9.5}{11}\selectfont]
   (event-inclusion) at ([xshift=-3.9cm,yshift=1.4cm]neighborhood-fork)
   {\mbox{\textbf{Neighborhood Inclusion}\hspace{.6em}{\fontsize{9.5}{11}\selectfont\eqref{eq:counterexample-tube-inclusion}}}};
 \node[argument box,text width=5.4cm]
   (neighborhood-probability) at ([xshift=3.9cm,yshift=1.4cm]neighborhood-fork)
   {\textbf{\mbox{Neighborhood Probability}}\par
    {\fontsize{9.5}{11}\selectfont
     \begin{itemize}[leftmargin=7pt,labelsep=2.5pt,topsep=3pt,
       itemsep=1.5pt,parsep=0pt,partopsep=0pt]
      \item \emph{Translated Ball}\hspace{.6em}\autoref{lem:gaussian-tube},
        \eqref{eq:counterexample-tube-probability}
      \item \emph{Gaussian Cost}\hspace{.6em}\eqref{eq:precursor-gaussian-cost}
      \item \emph{Centered Ball}\hspace{.6em}\eqref{eq:counterexample-centered-tube}
     \end{itemize}}};
 \coordinate (lower-bound-join) at ([xshift=-3.9cm,yshift=.45cm]neighborhood-probability.north);
 \node[route box,text width={},minimum height=0pt,font=\fontsize{9.5}{11}\selectfont,anchor=south]
   (tail-bound) at ([yshift=.65cm]lower-bound-join)
   {\mbox{\textbf{Lower Bound with a Strictly Positive Cost Gap}\hspace{.6em}{\fontsize{9.5}{11}\selectfont\eqref{eq:counterexample-LDP-bound} +
     \autoref{rem:counterexample-time-range}}}};
 \node[route box,text width={},minimum height=0pt,font=\normalsize,anchor=south]
   (special-failure) at ([yshift=.65cm]tail-bound.north)
   {\autoref{thm:longtime-rate-obstruction}};
 \coordinate (breakdown-level) at ([yshift=.975cm]special-failure.north);
 \node[route box,text width=16.6cm,minimum height=1cm,inner xsep=7pt,inner ysep=5pt,font=\normalsize,
   draw=black!85,line width=.7pt,fill=black!4,anchor=south]
   (breakdown) at (4.74,0 |- breakdown-level)
   {\textbf{Breakdown:}\hspace{.6em}\autoref{thm:intro-failure}
    (\autoref{thm:general-failure})};
 \draw[route arrow] (breather-embedding.west) -- (tube-stability.east);
 \draw[draw=black!75,line width=.6pt]
   (tube-stability.north) -- (neighborhood-fork);
 \fill[black!75] (neighborhood-fork) circle (.8pt);
 \draw[route arrow] (neighborhood-fork) -| (event-inclusion.south);
 \draw[route arrow] (neighborhood-fork) -| (neighborhood-probability.south);
 \draw[draw=black!75,line width=.6pt]
   (event-inclusion.north) |- (lower-bound-join);
 \draw[draw=black!75,line width=.6pt]
   (neighborhood-probability.north) |- (lower-bound-join);
 \fill[black!75] (lower-bound-join) circle (.8pt);
 \draw[route arrow] (lower-bound-join) -- (tail-bound.south);
 \draw[route arrow] (tail-bound.north) -- (special-failure.south);
 \draw[route arrow] (special-failure.north) --
   node[midway,right=6pt,align=left,font=\fontsize{9.5}{11}\selectfont]
     {\emph{Small Frequency Construction}\hspace{.4em}(\autoref{lem:general-small-frequencies})}
   (breakdown.south);
\path[use as bounding box] (-4,0) -- (14,0);
\end{tikzpicture}
\endgroup
}
\par\nointerlineskip
\begingroup
\setlength{\captionindent}{0pt}
\makeatletter\def\@captionfont{\normalfont\footnotesize}\makeatother
\caption[Proof Route for Breakdown]{The proof route for breakdown begins with an explicit Akhmediev
breather and its small frequency embedding into a fixed quasi-periodic
module. Tube stability controls solutions starting near the embedded
breather data up to focusing. The argument then separates into two
branches. Neighborhood inclusion connects these initial data to the target
exceedance event. Neighborhood probability combines a Gaussian translation
estimate, the cost of the neighborhood center and a centered ball estimate.
The branches merge to give a probability lower bound with a strictly
positive cost gap. The special case yields
\autoref{thm:longtime-rate-obstruction}, and the general 
construction extends the conclusion to \autoref{thm:intro-failure}.}\label{fig:proof-routes}
\endgroup
\end{figure}

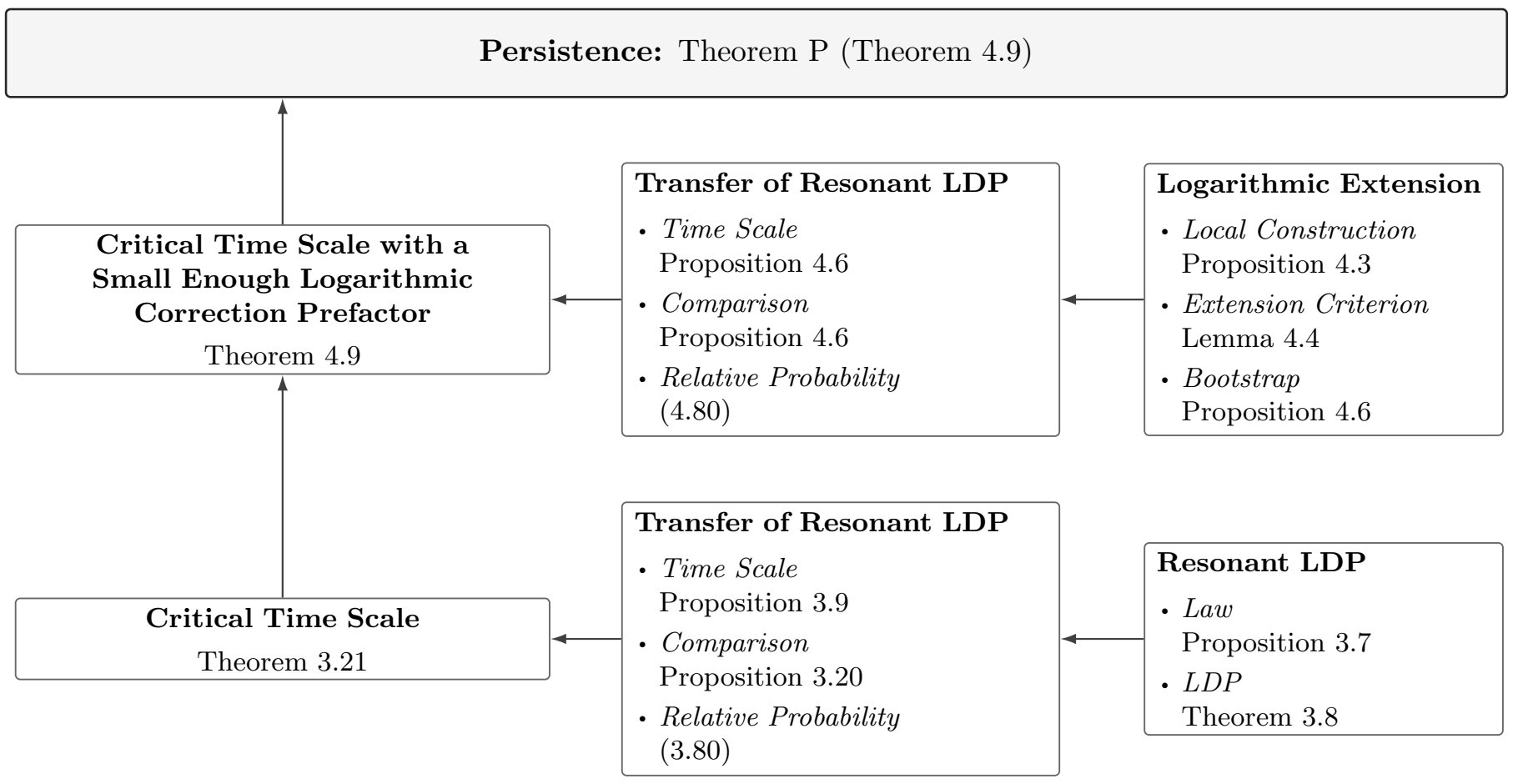
\begin{figure}[!htbp]
\centering
\vspace*{6pt}
\setlength{\abovecaptionskip}{10pt}
\resizebox{1.01\linewidth}{!}{
\begingroup
\let\persistenceitemlabel\labelitemi
\renewcommand{\labelitemi}{\makebox[.5em][c]{\raisebox{.25ex}{\scalebox{.5}{\persistenceitemlabel}}}}
\begin{tikzpicture}[
 x=1cm,y=.92cm,font=\fontsize{9.5}{11}\selectfont,
 route box/.style={draw=black!60,fill=white,line width=.5pt,
   rounded corners=1.5pt,align=center,text width=4.6cm,
   minimum height=0pt,inner xsep=4pt,inner ysep=3pt,
   execute at begin node={\hyphenpenalty=10000\exhyphenpenalty=10000}},
 transfer box/.style={route box,align=left,text width=4.7cm,
   font=\fontsize{9.5}{11}\selectfont,minimum height=0pt,inner xsep=4pt,inner ysep=3pt},
 input box/.style={route box,align=left,text width=3.8cm,
   font=\fontsize{9.5}{11}\selectfont,minimum height=0pt,inner xsep=4pt,inner ysep=3pt},
 route arrow/.style={-{Latex[length=1.8mm,width=1.2mm]},
   draw=black!75,line width=.6pt}
]
 \node[route box,text width=16.6cm,minimum height=1cm,inner xsep=7pt,inner ysep=5pt,font=\normalsize,
   draw=black!85,line width=.7pt,fill=black!4]
   (persistence) at (4.74,3.05)
   {\textbf{Persistence:} \autoref{thm:intro-persistence}
    (\autoref{thm:log-critical-ldp})};
 \pgfmathsetlengthmacro{\leftcolumnx}{-.55cm-3pt}
 \node[route box,text width=5.8cm,font=\fontsize{9.5}{11}\selectfont] (log-extension) at (\leftcolumnx,0)
   {\textbf{Critical Time Scale with a Small Enough Logarithmic Correction Prefactor}\\[3pt]
    {\fontsize{9.5}{11}\selectfont\autoref{thm:log-critical-ldp}}};
 \node[route box,text width=5.8cm,font=\fontsize{9.5}{11}\selectfont] (critical-case) at (\leftcolumnx,-4.2)
   {\textbf{Critical Time Scale}\\[3pt]
    {\fontsize{9.5}{11}\selectfont\autoref{thm:first-critical-ldp-local}}};
 \node[transfer box] (log-transfer) at (5.7,0)
   {\textbf{\mbox{Transfer of Resonant LDP}}\par
    {\fontsize{9.5}{11}\selectfont
     \begin{itemize}[leftmargin=8pt,labelsep=3pt,topsep=4pt,
       itemsep=2pt,parsep=0pt,partopsep=0pt,listparindent=0pt]
      \item \emph{Time Scale}\par\noindent
        \autoref{prop:log-pointwise-comparison}
      \item \emph{Comparison}\par\noindent
        \autoref{prop:log-pointwise-comparison}
      \item \emph{Relative Probability}\par\noindent
        \eqref{eq:log-relative-exception-small}
     \end{itemize}}};
 \node[input box] (bootstrap) at (11.2,0)
   {\textbf{\mbox{Logarithmic Extension}}\par
    {\fontsize{9.5}{11}\selectfont
     \begin{itemize}[leftmargin=8pt,labelsep=3pt,topsep=4pt,
       itemsep=2pt,parsep=0pt,partopsep=0pt,listparindent=0pt]
      \item \emph{Local Construction}\par\noindent
        \autoref{prop:log-local-wellposedness}
      \item \emph{\mbox{Extension Criterion}}\par\noindent
        \autoref{lem:log-continuation}
      \item \emph{Bootstrap}\par\noindent
        \autoref{prop:log-pointwise-comparison}
     \end{itemize}}};
 \node[transfer box] (critical-transfer) at (5.7,-4.2)
   {\textbf{\mbox{Transfer of Resonant LDP}}\par
    {\fontsize{9.5}{11}\selectfont
     \begin{itemize}[leftmargin=8pt,labelsep=3pt,topsep=4pt,
       itemsep=2pt,parsep=0pt,partopsep=0pt,listparindent=0pt]
      \item \emph{Time Scale}\par\noindent
        \autoref{prop:critical-pointwise-lifespan}
      \item \emph{Comparison}\par\noindent
        \autoref{prop:resonant-comparison}
      \item \emph{Relative Probability}\par\noindent
        \eqref{eq:critical-relative-error-small}
     \end{itemize}}};
 \node[input box] (resonant-reference) at (11.2,-4.2)
   {\textbf{\mbox{Resonant LDP}}\par
    {\fontsize{9.5}{11}\selectfont
     \begin{itemize}[leftmargin=8pt,labelsep=3pt,topsep=4pt,
       itemsep=2pt,parsep=0pt,partopsep=0pt,listparindent=0pt]
      \item \emph{Law}\par\noindent
        \autoref{prop:nf-resonant-law}
      \item \emph{LDP}\par\noindent
        \autoref{thm:nf-resonant-LDP}
     \end{itemize}}};
 \draw[route arrow] (bootstrap.west) -- (log-transfer.east);
 \draw[route arrow] (log-transfer.west) -- (log-extension.east);
 \draw[route arrow] (resonant-reference.west) -- (critical-transfer.east);
 \draw[route arrow] (critical-transfer.west) -- (critical-case.east);
 \draw[route arrow] (critical-case.north) -- (log-extension.south);
 \draw[route arrow] (log-extension.north)
   -- (log-extension.north |- persistence.south);
\path[use as bounding box] (-4,0) -- (14,0);
\end{tikzpicture}
\endgroup
}
\par\nointerlineskip
\begingroup
\setlength{\captionindent}{0pt}
\makeatletter\def\@captionfont{\normalfont\footnotesize}\makeatother
\caption[Proof Route for Persistence]{The proof route for persistence begins
with the resonant flow, whose law agrees with that of the initial data
and yields a reference LDP. At critical time, existence and comparison
estimates control the nonlinear solution, while relative probability
bounds make the exceptional set negligible at the large deviation scale.
Together, these estimates transfer the resonant LDP to the nonlinear flow.
For logarithmic extension, local construction and an extension criterion
provide the starting point for a bootstrap argument. The bootstrap controls
the solution and its difference from the resonant flow on the longer
interval. The same probability transfer then gives persistence with a
small enough logarithmic prefactor, yielding \autoref{thm:intro-persistence}.}\label{fig:proof-routes-persistence}
\endgroup
\end{figure}

\autoref{fig:proof-routes} and \autoref{fig:proof-routes-persistence}
outline the proof routes for breakdown and persistence, respectively;
see their captions for explanations.

\section{Functional Setting}\label{sec:structures}

Throughout the paper, fix an integer $\nu\ge2$, which denotes the
number of basic frequencies and the dimension of the Fourier lattice
$\Z^\nu$. For any integer $d\ge1$ and vector
$v=(v_1,\ldots,v_d)\in\R^d$, we use the $\ell^1$ norm
\begin{align}
    |v|_1:=\sum_{j=1}^d|v_j|
\end{align}
and the associated polynomial weight
\begin{align}
    \langle v\rangle:=1+|v|_1.
\end{align}
For $d=1$, these reduce to $|v|_1=|v|$ and
$\langle v\rangle=1+|v|$.

\subsection{Label for Complex Conjugation}\label{subsec:complex-conjugation}

For any \(z\in\C\), we use \(z^{\ast^{[0]}}\) and \(z^{\ast^{[1]}}\) to denote \(z\) and its complex conjugate \(\bar z\), respectively; that is,
\begin{align*}
z\triangleq z^{\ast^{[0]}}
\quad\text{and}\quad
\bar z\triangleq z^{\ast^{[1]}}.
\end{align*}
In general, applying complex conjugation \(m\) times to \(z\) yields \(z^{\ast^{[m]}}\), where $[\cdot]=[\hat{\cdot}]$ and $\hat{\cdot}\in\{0,1\}$ is determined by the congruence $\hat{m}\equiv m~(\text{mod}~2)$.
This notation satisfies 
\begin{itemize}
    \item \((z_1+z_2)^{\ast^{[m]}}=z_1^{\ast^{[m]}}+z_2^{\ast^{[m]}}\);
    \item \((z_1z_2)^{\ast^{[m]}}=z_1^{\ast^{[m]}}z_2^{\ast^{[m]}}\);
    \item \(\overline{z^{\ast^{[m]}}}=z^{\ast^{[m+1]}}\);
    \item \(\bigl(z^{\ast^{[m]}}\bigr)^{\ast^{[m']}}=z^{\ast^{[m+m']}}\),
\end{itemize}
for any \(z_1,z_2\in\C\) and \(m,m'\in\N\); see \cite{DLX24I}.

\subsection{Quasi-Periodic Functions}
\label{subsec:quasi-periodic-functions}

We work with continuous quasi-periodic functions obtained by restricting
hull functions on a finite-dimensional torus to a linear orbit. These functions are uniformly almost periodic, and their Fourier
frequencies lie in a finitely generated subgroup of $\R$; see
\cite{Bes55,Tam35}. We first specify the frequency vector and then
develop the Fourier representation and mean
used below.

\begin{definition}[Rational Independence and Diophantine Condition]
\label{def:dio}
Let $\omega=(\omega_1,\ldots,\omega_\nu)\in\R^\nu$ be a fixed \emph{frequency vector}.
\begin{enumerate}[label=\textup{(\roman*)},leftmargin=*,
    topsep=3pt,itemsep=4pt,parsep=0pt]
\item The vector $\omega$ is \emph{rationally independent}
(or \emph{nonresonant}) if
\[
\langle n,\omega\rangle\ne0, 
\qquad\text{for every }n\in\Z^\nu\setminus\{0\}.
\]
\item The vector $\omega$ is said to satisfy the \emph{Diophantine condition}
if there exist constants $\gamma>0$ and $\tau>\nu-1$ such that
\begin{equation}\label{eq:HDC}
|\langle n,\omega\rangle|\ge\gamma|n|_1^{-\tau}, 
\qquad\text{for every }n\in\Z^\nu\setminus\{0\}.
\end{equation}
For fixed $\gamma$ and $\tau$, we denote by $\mathrm{DC}(\gamma,\tau)$
the set of all vectors satisfying \eqref{eq:HDC}.
\end{enumerate}
\end{definition}
\begin{remark}\label{rem:standing-frequency}
Throughout the paper, the basic frequency vector $\omega$ is assumed
to be rationally independent. The Diophantine condition is a
quantitative strengthening of this assumption and will be imposed
explicitly when needed.
\end{remark}

\begin{definition}[Quasi-Periodic Function]\label{def:qp}

Let $\T^\nu:=(\R/2\pi\Z)^\nu$ be the $\nu$-dimensional torus.
For $y\in\R^\nu$, write $\theta=[y]\in\T^\nu$ for its \emph{equivalence
class} modulo $2\pi\Z^\nu$, and call $\theta$ the \emph{angle variable}.
A function $u:\R\to\C$ is called \emph{quasi-periodic with frequency
vector $\omega$} if there exists a continuous function $U:\T^\nu\to\C; \theta\mapsto U(\theta)$ such that
\begin{equation}\label{eq:qp-representation}
u(x)=U([\omega x]),\qquad x\in\R.
\end{equation}
We call \eqref{eq:qp-representation} a \emph{torus representation} of $u$
and $U$ its \emph{hull function}.
\end{definition}

\begin{remark}
(a)~For convenience, we write $U(\omega x)$ for $U([\omega x])$. 

(b)~\phantomsection\label{sym:def:qp-space}Let $\mathrm{QP}_{\omega}(\R)$ denote the space of continuous
quasi-periodic functions with nonresonant frequency $\omega$. 
\end{remark}

\begin{definition}[Fourier Modes and Frequency Module]
\label{def:frequency-module}
Let $n\in\Z^\nu$ be a (dual) \emph{lattice vector} of $\T^\nu$. The restriction of the \emph{torus character}, or
\emph{Fourier mode}, $e^{\iu\langle n,\theta\rangle}$ to the \emph{linear orbit}
$x\mapsto[\omega x]$ is called the \emph{spatial Fourier mode} $
e^{\iu\langle n,\omega\rangle x}$,
with \emph{spatial frequency} $\langle n,\omega\rangle$.

The set of all
such spatial frequencies is called the \emph{frequency module} generated by
the frequency vector $\omega$ and is denoted by
\begin{align*}
\mathcal M_\omega
:=\Z\omega_1+\cdots+\Z\omega_\nu=\set{\langle n,\omega\rangle:n\in\Z^\nu}.
\end{align*}

\end{definition}

\begin{lemma}[Uniqueness of Spatial Frequency]\label{lem:frequency-label}
The map
\[
\iota_\omega:\Z^\nu\longrightarrow\mathcal M_\omega,
\qquad n\longmapsto\langle n,\omega\rangle
\]
is a bijection. In particular, each spatial frequency in
$\mathcal M_\omega$ has a unique label in $\Z^\nu$; equivalently,
for $n,m\in\Z^\nu$,
\[
\langle n,\omega\rangle=\langle m,\omega\rangle
\quad\Longleftrightarrow\quad
n=m.
\]
\end{lemma}

\begin{proof}[Proof of \autoref{lem:frequency-label}]
Surjectivity follows from the definition of $\mathcal M_\omega$.
If $\langle n,\omega\rangle=\langle m,\omega\rangle$, then
$\langle n-m,\omega\rangle=0$. Since $n-m\in\Z^\nu$ and
$\omega$ is rationally independent, we have $n-m=0$.
Thus the map is also injective.
\end{proof}

\begin{remark}
In what follows, the term \emph{mode} refers either to a lattice vector $n\in\Z^\nu$ or to the corresponding spatial Fourier mode $e^{\iu\langle n,\omega\rangle x}$.
\end{remark}

\begin{proposition}[Fourier Series of the Hull Function]
\label{prop:torus-fourier-convergence}

Let $U\in C(\T^\nu;\C)$ and define its Fourier coefficients by
\begin{equation*}
\widehat U(n)
:=
\int_{\T^\nu} U(\theta)e^{-\iu\langle n,\theta\rangle}\,dm_\nu(\theta),
\qquad n\in\Z^\nu,
\end{equation*}
where $m_\nu$ is the normalized Haar measure on $\T^\nu$, that is,  $dm_\nu(\theta)=(2\pi)^{-\nu}d\theta$. 
If $\widehat U\in\ell^1(\Z^\nu)$, then
\[
U(\theta)
=
\sum_{n\in\Z^\nu}\widehat U(n)e^{\iu\langle n,\theta\rangle},
\qquad\theta\in\T^\nu,
\]
with absolute and uniform convergence.
\end{proposition}

\begin{proof}[Proof of \autoref{prop:torus-fourier-convergence}]
By the \emph{Weierstrass test}, the series converges absolutely and uniformly
to some $V\in C(\T^\nu;\C)$. Termwise integration yields
$\widehat V(n)=\widehat U(n)$ for all $n\in\Z^\nu$.
Uniqueness of Fourier coefficients gives $V=U$ almost everywhere,
and continuity implies equality everywhere.
\end{proof}

\phantomsection\label{sym:def:fourier-coefficient}
For $u\in\mathrm{QP}_{\omega}(\R)$ with hull function $U$, set
$\widehat u(n):=\widehat U(n)$. If $\widehat u\in\ell^1(\Z^\nu)$,
the corresponding Fourier
expansion is
\begin{equation}\label{eq:qp-fourier}
u(x)
=
\sum_{n\in\Z^\nu}
\widehat u(n)e^{\iu\langle n,\omega\rangle x},\qquad x\in\R.
\end{equation}
Proposition~\ref{prop:torus-fourier-convergence}
shows that \eqref{eq:qp-fourier} converges absolutely and uniformly to
$u$. 

We now introduce a spatial average that gives an intrinsic description of these coefficients
and transfers the \emph{torus orthogonality} and \emph{Parseval identity} to $\R$.

\begin{definition}[Spatial Mean]\label{def:spatial-mean}

Let $f\in L^1_{\mathrm{loc}}(\R)$. If the limit exists, its
\emph{spatial mean} is defined by
\begin{equation*}
\operatorname{Mean}(f)
:=
\lim_{L\to\infty}
\frac{1}{2L}\int_{-L}^{L}f(x)\,dx.
\end{equation*}
\end{definition}

\begin{lemma}[Spatial Mean of a Fourier Mode]\label{lem:mean-mode}
For every $\lambda\in\R$,
\[
\operatorname{Mean}(e^{\iu\lambda x})
=
\begin{cases}
1, & \lambda=0,\\
0, & \lambda\neq0.
\end{cases}
\]
\end{lemma}

\begin{proof}[Proof of \autoref{lem:mean-mode}]
The assertion is immediate for $\lambda=0$. If $\lambda\neq0$, then
\[
\frac{1}{2L}\int_{-L}^{L}e^{\iu\lambda x}\,dx
=
\frac{\sin(\lambda L)}{\lambda L},
\]
which tends to $0$ as $L\to\infty$.
\end{proof}

\begin{corollary}[Orthogonality of Modes with Respect to the Spatial Mean]
\label{cor:qp-orthogonality} For every
$n,m\in\Z^\nu$,
\begin{equation*}
\operatorname{Mean}\left(
e^{\iu\langle n,\omega\rangle x}
e^{-\iu\langle m,\omega\rangle x}
\right)
=
\delta_{n,m},
\end{equation*}
where $\delta_{n,m}$ is the Kronecker delta.
\end{corollary}

\begin{proof}[Proof of \autoref{cor:qp-orthogonality}]
Apply \autoref{lem:mean-mode} with
$\lambda=\langle n-m,\omega\rangle$ and use
\autoref{lem:frequency-label}.
\end{proof}

\begin{proposition}[Spatial and Torus Means]\label{prop:spatial-torus-mean}
For every $F\in C(\T^\nu;\C)$, the spatial mean of $F(\omega\,\cdot)$ exists and
\begin{equation*}
\operatorname{Mean}\bigl(F(\omega\,\cdot)\bigr)
=\int_{\T^\nu}F(\theta)\,dm_\nu(\theta),
\end{equation*}
where the right-hand side is the torus mean of $F$.
\end{proposition}

\begin{proof}[Proof of \autoref{prop:spatial-torus-mean}]
For a trigonometric polynomial $P$, the identity follows from
\autoref{cor:qp-orthogonality}.
Given $\varepsilon>0$, choose such a $P$ with
$\norm{F-P}_{L^\infty(\T^\nu)}<\varepsilon$. For every $L>0$,
\begin{align*}
&\left|\frac{1}{2L}\int_{-L}^{L}F(\omega x)\,dx
-\int_{\T^\nu}F(\theta)\,dm_\nu(\theta)\right|\\
\le~& \left|\frac{1}{2L}\int_{-L}^{L}P(\omega x)\,dx
-\int_{\T^\nu}P(\theta)\,dm_\nu(\theta)\right|+2\varepsilon.
\end{align*}
Taking the limit superior as $L\to\infty$ and then letting
$\varepsilon\to0^+$ prove the claim.
\end{proof}

\begin{proposition}[Recovery of Fourier Coefficients]
\label{prop:fourier-recovery}
Let $u\in\mathrm{QP}_{\omega}(\R)$. For every $n\in\Z^\nu$,
\begin{equation}\label{eq:qp-fourier-coefficient}
\widehat u(n)
=\operatorname{Mean}\left(u(x)e^{-\iu\langle n,\omega\rangle x}\right).
\end{equation}
More generally, for $\lambda\in\R\setminus\mathcal M_\omega$,
\[
\operatorname{Mean}\left(u(x)e^{-\iu\lambda x}\right)=0.
\]
Thus the Fourier--Bohr spectrum of $u$ is
$\{\langle n,\omega\rangle:\widehat u(n)\ne0\}\subseteq\mathcal M_\omega$.
\end{proposition}

\begin{proof}[Proof of \autoref{prop:fourier-recovery}]
If $U$ is the hull function of $u$, apply
\autoref{prop:spatial-torus-mean} to
$F(\theta)=U(\theta)e^{-\iu\langle n,\theta\rangle}$ to obtain
\eqref{eq:qp-fourier-coefficient}.
For $\lambda\notin\mathcal M_\omega$, the asserted mean vanishes for
every torus trigonometric polynomial by Lemma~\ref{lem:mean-mode}.
Uniform approximation of $U$ by such polynomials gives the result.
\end{proof}

\begin{proposition}[Parseval Identity]
\label{prop:qp-parseval}
Let $u,v\in\mathrm{QP}_{\omega}(\R)$ have hull functions $U,V$.
Then
\begin{align*}
\operatorname{Mean}(u\overline v)
&=\int_{\T^\nu}U\overline {V}\,dm_\nu\\
&=\sum_{n\in\Z^\nu}\widehat u(n)\overline{\widehat v(n)}.
\end{align*}
In particular,
\begin{equation*}
\operatorname{Mean}(|u|^2)=\sum_{n\in\Z^\nu}|\widehat u(n)|^2.
\end{equation*}
\end{proposition}

\begin{proof}[Proof of \autoref{prop:qp-parseval}]
Apply Proposition~\ref{prop:spatial-torus-mean} to $U\overline V$
and use \emph{Parseval's identity} on $(\T^\nu,m_\nu)$. 
\end{proof}

\begin{proposition}[Fourier Calculus]\label{prop:qp-fourier-calculus}
Let $u,v\in\mathrm{QP}_{\omega}(\R)$ with
$\widehat u,\widehat v\in\ell^1(\Z^\nu)$.

\begin{enumerate}[label=\textup{(\roman*)}]
\item Complex conjugation acts by
\begin{equation*}
\widehat{\overline{u}}(n)
=
\overline{\widehat u(-n)}.
\end{equation*}

\item Pointwise multiplication corresponds to discrete convolution, that is, 
\begin{equation*}
\widehat{uv}(n)
=
\sum_{\substack{n_1,n_2\in\Z^\nu\\n_1+n_2=n}}
\widehat u(n_1)\widehat v(n_2).
\end{equation*}

\item In particular,
\begin{equation*}
\widehat{|u|^2u}(n)
=
\sum_{\substack{n_1,n_2,n_3\in\Z^\nu\\n_1-n_2+n_3=n}}
\prod_{j=1}^3\starpow{\widehat u(n_j)}{j-1}.
\end{equation*}

\item Let $j\ge1$ be an integer. If
\[
\sum_{n\in\Z^\nu}
|\langle n,\omega\rangle|^j|\widehat u(n)|<\infty,
\]
then $u\in C^j(\R)$, $\partial_x^j u\in\mathrm{QP}_{\omega}(\R)$ and
\begin{equation}\label{eq:qp-differentiation}
\widehat{\partial_x^j u}(n)
=
\bigl(\iu\langle n,\omega\rangle\bigr)^j
\widehat u(n).
\end{equation}
\end{enumerate}
\end{proposition}

\begin{proof}[Proof of \autoref{prop:qp-fourier-calculus}]
Absolute convergence permits conjugation, multiplication and
regrouping of the Fourier series, giving (i)--(iii).
For (iv), the bound $|\lambda|^k\le1+|\lambda|^j$ for $0\le k\le j$
ensures absolute and uniform convergence of every derivative series
up to order $j$. Termwise differentiation gives
\[
\partial_x^j u(x)
=
\sum_{n\in\Z^\nu}
\bigl(\iu\langle n,\omega\rangle\bigr)^j
\widehat u(n)e^{\iu\langle n,\omega\rangle x},
\]
and the same absolutely summable coefficients define a continuous hull
function for $\partial_x^j u$, proving \eqref{eq:qp-differentiation}.
\end{proof}

\subsection{Fourier Coefficients and Randomization}\label{subsec:fourier-randomization}
The representation \eqref{eq:qp-representation} identifies a spatially
quasi-periodic function on $\R$ with the restriction of a function on
$\T^\nu$ to the linear orbit $x\longmapsto[\omega x]$.
By the rational independence of $\omega$, this orbit is dense in
$\T^\nu$, and the map
$\Z^\nu\ni n\longmapsto\langle n,\omega\rangle\in\mathcal M_\omega$
is bijective. Hence each quasi-periodic Fourier mode is uniquely
indexed by $n\in\Z^\nu$. Accordingly, a function
$u(x)=\sum_{n\in\Z^\nu}\widehat u(n)e^{\iu\langle n,\omega\rangle x}$
is naturally encoded by its Fourier coefficient sequence
$\widehat u=(\widehat u(n))_{n\in\Z^\nu}$.

We randomize the initial Fourier coefficients by setting
\[
 \widehat u(n)=c(n)g_n,\qquad n\in\Z^\nu,
\]
where $c=(c(n))_{n\in\Z^\nu}$ is the deterministic amplitude sequence
and $(g_n)_{n\in\Z^\nu}$ is the complex Gaussian family. Their assumptions
are stated separately below.

\subsubsection{Deterministic Amplitudes and Weight}
\label{subsubsec:deterministic-fourier-amplitudes}
Fix an exponent vector
$\rho=(\rho_1,\ldots,\rho_\nu)\in\R^\nu$. 
For every mode $n=(n_1,\ldots,n_\nu)\in\Z^\nu$, define the
product polynomial weight
\begin{align}\label{eq:product-polynomial-weight}
 w_\rho(n):=\prod_{j=1}^{\nu}\langle n_j\rangle^{\rho_j}.
\end{align}
Assume that $c\not\equiv0$
and, for some $A>0$ and $\rho_j>1$ for every $j$,
\begin{equation}\label{eq:deterministic-decay}
 |c(n)|\le A w_{-\rho}(n),
 \qquad n\in\Z^\nu.
\end{equation}\label{eq:nf-polynomial-decay}

\subsubsection{Complex Gaussian and Exceptional Set}
\label{subsubsec:complex-gaussian-exceptional-set}

The following definitions specify the complex Gaussian and the
associated exceptional set.
For further background on complex Gaussian random variables, we refer the reader to Appendix A of \cite{XL26}.

\begin{definition}[Standard Complex Gaussian Family]\label{def:gaussian}

On a probability space $(\Omega,\mathcal F,\Pp)$, let $(g_n)_{n\in\Z^\nu}$ be independent and identically distributed (i.i.d.) standard complex Gaussian random variables. Equivalently, $g_n=X_n+\iu Y_n$, where $X_n,Y_n$ are i.i.d. real Gaussian random variables with law $\mathcal N_{\R}(0,1/2)$. Then
\[
 \Ee [g_n]=0,\qquad \Ee[g_ng_m]=0,\qquad
 \Ee\bigl[g_n\overline{g_m}\bigr]=\delta_{n,m}.
\]
In polar coordinates $g_n=R_ne^{\iu\Theta_n}$, the squared radii satisfy $R_n^2\sim\mathscr E(1)$\footnote{Since $X_n,Y_n\sim\mathcal N_{\R}(0,1/2)$, we have
$2(X_n^2+Y_n^2)\sim\chi^2_2=\mathscr E(1/2)=\operatorname{Gamma}(1,2)$. Hence
$R_n^2=X_n^2+Y_n^2\sim\mathscr E(1)=\operatorname{Gamma}(1,1)$.}, and the angles satisfy $\Theta_n\sim\mathscr U([0,2\pi))$. The radii and angles are mutually independent. 
\end{definition}

\begin{definition}[Exceptional Set]\label{def:exceptional}

Let $\delta>0$ and $\kappa=(\kappa_1,\cdots,\kappa_\nu)\in(0,\infty)^\nu$. Define
\begin{equation}\label{eq:Omega-delta}
 \Omega_\delta
 :=\bigcup_{n\in\Z^\nu}
 \set{|g_n|>\delta^{-1/2}w_\kappa(n)}.
\end{equation}
\end{definition}

\begin{proposition}[Probability of the Exceptional Set]\label{prop:exception-prob}
There exists a constant
$C_{\kappa,\nu}>0$, depending only on $\kappa$ and $\nu$, such that
\begin{equation}\label{eq:exception-prob}
 \Pp(\Omega_\delta)\le C_{\kappa,\nu} e^{-1/\delta},
 \qquad 0<\delta\le1.
\end{equation}
\end{proposition}
\begin{proof}[Proof of \autoref{prop:exception-prob}]
Set $w_n:=w_\kappa(n)\ge1$. Since $|g_n|^2\sim\mathscr{E}(1)$,
\begin{align*}
 \Pp\bigl(|g_n|>\delta^{-1/2}w_n\bigr)
 &=\Pp\bigl(|g_n|^2>\delta^{-1}w_n^2\bigr)\\
 &=e^{-w_n^2/\delta}.
\end{align*}
The union bound gives
\begin{align*}
 \Pp(\Omega_\delta)
 &\le\sum_{n\in\Z^\nu}e^{-w_n^2/\delta}\\
 &=e^{-1/\delta}\sum_{n\in\Z^\nu}e^{-(w_n^2-1)/\delta}\\
 &\le e^{-1/\delta}\sum_{n\in\Z^\nu}e^{-(w_n^2-1)},
\end{align*}
where the last inequality uses $0<\delta\le1$ and $w_n^2-1\ge0$.
Thus it suffices to show that
\[
 C_{\kappa,\nu}:=\sum_{n\in\Z^\nu}e^{-(w_n^2-1)}<\infty.
\]
Set
$\kappa_{\min}:=\min_{1\le j\le\nu}\kappa_j>0$ and
$|n|_\infty:=\max_{1\le j\le\nu}|n_j|$. Then
\begin{align*}
 w_n^2
 &=\prod_{j=1}^\nu\langle n_j\rangle^{2\kappa_j}\\
 &\ge\prod_{j=1}^\nu\langle n_j\rangle^{2\kappa_{\min}}\\
 &\ge\langle n\rangle^{2\kappa_{\min}}.
\end{align*}
Grouping lattice points according to $m=|n|_\infty$ and using
$\#\{n\in\Z^\nu:|n|_\infty=m\}\le(2m+1)^\nu$, we obtain
\begin{align*}
 C_{\kappa,\nu}&\le e\sum_{m=0}^\infty
 (2m+1)^\nu e^{-(1+m)^{2\kappa_{\min}}}\\
 &<\infty,
\end{align*}
since $\kappa_{\min}>0$ and the exponential decay dominates the
polynomial growth. This proves
\eqref{eq:exception-prob}.
\end{proof}

\begin{corollary}[Almost Sure Decay of the Randomized Fourier Coefficients]
\label{cor:randomized-decay}
Assume that $0<\kappa_j<\rho_j$ for every $j=1,\ldots,\nu$, and set
$\beta:=\rho-\kappa\in(0,\infty)^\nu$. Then, on $G_\delta:=\Omega\setminus\Omega_\delta$, the randomized Fourier coefficients satisfy
the uniform bound
\begin{equation}\label{eq:randomized-decay}
 |c(n)g_n|\le A\delta^{-1/2}w_{-\beta}(n),
 \qquad n\in\Z^\nu.
\end{equation}
\end{corollary}
\begin{proof}[Proof of \autoref{cor:randomized-decay}]
On $G_\delta$, definition \eqref{eq:Omega-delta} gives
$|g_n|\le\delta^{-1/2}w_\kappa(n)$ for all $n\in\Z^\nu$.
Combining this with \eqref{eq:deterministic-decay} and
$w_{-\rho}(n)w_\kappa(n)=w_{-\beta}(n)$ yields
\eqref{eq:randomized-decay}.
Since \eqref{eq:exception-prob} gives $\Pp(\Omega_\delta)\to0$
as $\delta\downarrow0$, the bound holds with a finite random constant
almost surely.
\end{proof}

\section{LDP at the Critical Time and Cubic Resonance}\label{sec:normal-forms}

In this section, we first construct the cubic resonant approximation,
following \emph{the resonant decomposition technique} of
\cite{JM21}; see \autoref{subsec:cubic-resonance} through
\autoref{sec:nf-resonant-gauge}.
We then prove that the law of $u_{\mathrm{res}}(t,\cdot)$ is independent
of $t$ and coincides with that of the initial Gaussian field; see
\autoref{prop:nf-resonant-law}.
Using this equality in law, we combine a one-point lower bound with
\emph{the Fourier truncation argument} of \cite{XL26} to prove the
resonant LDP; see \autoref{thm:nf-resonant-LDP}.
Next, we use \emph{the tree expansion method} of \cite{DG16,DLX24I} to establish existence, uniqueness and uniform pointwise Fourier decay on the critical time scale;
see \autoref{prop:critical-pointwise-lifespan}.
Combining these bounds with \emph{integration by parts} in time and the
Diophantine condition, we obtain the uniform resonant comparison; see
\autoref{prop:resonant-comparison}.
Finally, we control the exceptional probability relative to the resonant
tail and transfer the LDP to the nonlinear solution, uniformly on the
critical interval; see \eqref{eq:critical-relative-error-small} and
\autoref{thm:first-critical-ldp-local}.

\subsection{Dispersion, Momentum and Cubic Resonance}
\label{subsec:cubic-resonance}

\begin{definition}[Dispersion and Cubic Resonance]\label{def:resonance-data}

For $n,n_1,n_2,n_3\in\Z^\nu$, define the \emph{dispersion relation} $Q$, the \emph{momentum function} $\mu_3$ and the \emph{resonance function} $\Sigma_3$ by
\begin{align*}
 Q(n)&:=\langle n,\omega\rangle^2,\\
 \mu_3(n_1,n_2,n_3;n)&:=n_1-n_2+n_3-n,\\
 \Sigma_3(n_1,n_2,n_3;n)&:=Q(n_1)-Q(n_2)+Q(n_3)-Q(n).
\end{align*}

The triple $(n_1,n_2,n_3)$ is \emph{resonant relative to $n$} if it satisfies both \emph{momentum conservation}, $\mu_3(n_1,n_2,n_3;n)=0$, and the \emph{resonance condition}, $\Sigma_3(n_1,n_2,n_3;n)=0$.

\end{definition}

\begin{remark}
When the input and output indices are clear from the context, we
abbreviate the momentum and resonance functions as $\mu_3$ and $\Sigma_3$.
\end{remark}

\begin{lemma}[Factorization of the Resonance Function]\label{lem:factorization}
If $\mu_3(n_1,n_2,n_3;n)=0$, then
\begin{equation}\label{eq:Sigma-factorization}
 \Sigma_3(n_1,n_2,n_3;n)
 =2\langle n_1-n_2,\omega\rangle
   \langle n_2-n_3,\omega\rangle.
\end{equation}
\end{lemma}

\begin{proof}[Proof of \autoref{lem:factorization}]
Set $a=\langle n_1,\omega\rangle$, $b=\langle n_2,\omega\rangle$ and $c=\langle n_3,\omega\rangle$.  Momentum conservation gives $\langle n,\omega\rangle=a-b+c$.  Therefore
\begin{align*}
 \Sigma_3&=a^2-b^2+c^2-(a-b+c)^2\\
 &=2(a-b)(b-c),
\end{align*}
which is \eqref{eq:Sigma-factorization}.
\end{proof}

\begin{corollary}[Classification of Cubic Resonance]\label{cor:cubic-resonances}
For $n,n_1,n_2,n_3\in\Z^\nu$, if
$\mu_3(n_1,n_2,n_3;n)=\Sigma_3(n_1,n_2,n_3;n)=0$, then
\[
 n_1=n_2\quad\text{or}\quad n_2=n_3.
\]
Equivalently,
\[
 \{n_1,n_3\}=\{n_2,n\}.
\]
\end{corollary}

\begin{proof}[Proof of \autoref{cor:cubic-resonances}]
By Lemma \ref{lem:factorization}, $\Sigma_3=0$ implies
\[
 \langle n_1-n_2,\omega\rangle
 \langle n_2-n_3,\omega\rangle=0.
\]
Rational independence gives $n_1=n_2$ or $n_2=n_3$.  Momentum conservation then yields respectively $n_3=n$ or $n_1=n$.
\end{proof}
\subsection{Removing the Linear Flow}\label{sec:normal-form-integral}

Let $u$ solve \eqref{eq:qp-cubic-NLS} with initial data
\eqref{eq:qp-randomized-initial-data}. Expand $u$ in a Fourier series as
\begin{equation}\label{eq:nf-fourier-expansion}
 u(t,x)=\sum_{n\in\Z^\nu}c(t,n)e^{\iu \langle n,\omega\rangle x},
\end{equation}
  with randomized initial data $c(0,n)=c(n)g_n$. The Fourier
coefficients in \eqref{eq:nf-fourier-expansion} satisfy
\begin{equation}\label{eq:nf-fourier-ode}
 \partial_t c(t,n)=-\iu Q(n)c(t,n)
 +\iu\varepsilon^\alpha\sum_{\substack{n_1,n_2,n_3\in\Z^\nu\\n_1-n_2+n_3=n}}\prod_{j=1}^3\starpow{c(t,n_j)}{j-1}.
\end{equation}
Introduce the interaction profile
\begin{equation}\label{eq:nf-interaction-transform}
 v_n(t)=e^{\iu tQ(n)}c(t,n),
 \qquad c(t,n)=e^{-\iu tQ(n)}v_n(t).
\end{equation}
Then $v$ has initial data $v_n(0)=c(n)g_n$ and satisfies
\begin{align}
 \dot v_n(t)&=\iu\varepsilon^\alpha\mathcal N^3[v](t,n),
\label{eq:nf-profile-ode}
\end{align}
where
\begin{align}
 \mathcal N^3[v](t,n)&:=
 \sum_{\substack{n_1,n_2,n_3\in\Z^\nu\\n_1-n_2+n_3=n}}
 e^{-\iu t\Sigma_3(n_1,n_2,n_3;n)}
 \prod_{j=1}^3\starpow{v_{n_j}(t)}{j-1}.
\label{eq:nf-cubic-interaction}
\end{align}
Equation \eqref{eq:nf-profile-ode} removes the linear flow while preserving
its effect through the cubic oscillatory phase. It is the starting
point for the resonant decomposition below.

\subsection{Resonant Decomposition and Mass Conservation}\label{sec:nf-resonant-gauge}
\subsubsection{Splitting the Cubic Interaction}
Fix an output mode $n\in\Z^\nu$. The nonlinear term in
\eqref{eq:nf-profile-ode} sums over all input triples
$(n_1,n_2,n_3)\in(\Z^\nu)^3$ satisfying the momentum conservation
condition $n_1-n_2+n_3=n$. We partition these triples according to
whether $\Sigma_3(n_1,n_2,n_3;n)$ vanishes. The triples with
$\Sigma_3=0$ form the \emph{resonant contribution}, for which
$e^{-\iu s\Sigma_3}=1$ at every time $s$; those with
$\Sigma_3\ne0$ form the \emph{nonresonant contribution}, which retains
the time-dependent oscillatory factor. Thus, at each time $s$,
\begin{align}
 \mathcal N^3[v](s,n)
 =\mathcal{N}_0^3[v](s,n)+\mathcal{N}_1^3[v](s,n),
\label{eq:nf-cubic-splitting}
\end{align}
where $\Sigma_3$ abbreviates $\Sigma_3(n_1,n_2,n_3;n)$, and 
\begin{align}
\mathcal{N}_0^3[v](t,n)&=\sum_{\substack{n_1,n_2,n_3\in\Z^\nu\\n_1-n_2+n_3=n\\\Sigma_3=0}}
 \prod_{j=1}^3\starpow{v_{n_j}(t)}{j-1},\nonumber\\
\mathcal{N}_1^3[v](t,n)&=\sum_{\substack{n_1,n_2,n_3\in\Z^\nu\\n_1-n_2+n_3=n\\\Sigma_3\ne0}}
 e^{-\iu t\Sigma_3}\prod_{j=1}^3\starpow{v_{n_j}(t)}{j-1}.\label{eq:nf-nonresonant-interaction}
\end{align}

\subsubsection{Counting the Resonant Modes}
To compute $\mathcal{N}_0^3[v](t,n)$, apply \autoref{cor:cubic-resonances}. Rational independence of
$\omega$ implies that, under momentum conservation, $\Sigma_3=0$
precisely when $n_1=n_2$ or $n_2=n_3$. In the first case,
momentum conservation forces $n_3=n$, so the triples are $(m,m,n)$;
in the second case it forces $n_1=n$, so the triples are $(n,m,m)$,
with $m\in\Z^\nu$ in both cases. Their respective contributions are
\begin{align*}
 \sum_{m\in\Z^\nu}v_m(t)\overline{v_m(t)}v_n(t)
 &=\left(\sum_{m\in\Z^\nu}|v_m(t)|^2\right)v_n(t),\\
 \sum_{m\in\Z^\nu}v_n(t)\overline{v_m(t)}v_m(t)
 &=\left(\sum_{m\in\Z^\nu}|v_m(t)|^2\right)v_n(t).
\end{align*}

\begin{proposition}[Mass Conservation]\label{prop:nf-mass-conservation}
Let $v=(v_n)_n$ solve \eqref{eq:nf-profile-ode} with initial data
$v_n(0)=c(n)g_n$, and assume that its coefficients are absolutely
summable locally uniformly in time.
If $\beta_j>1$ for every $j=1,\ldots,\nu$, then, on $G_\delta$,
throughout its interval of existence,
\begin{equation}\label{eq:nf-mass-conservation}
 M:=\sum_{m\in\Z^\nu}|v_m(t)|^2
 =\sum_{m\in\Z^\nu}|c(m)|^2|g_m|^2.
\end{equation}
\end{proposition}
\begin{proof}[Proof of \autoref{prop:nf-mass-conservation}]
On $G_\delta$, \autoref{cor:randomized-decay} and the assumption
$\beta_j>1$ for every $j=1,\ldots,\nu$ imply that the initial Fourier
coefficients are absolutely summable, that is, 
\begin{align*}
 \sum_{n\in\Z^\nu}|c(n)g_n|
 &\le A\delta^{-1/2}\sum_{n\in\Z^\nu}w_{-\beta}(n)\\
 &\le A\delta^{-1/2}\prod_{j=1}^\nu
\sum_{n_j\in\Z}\langle n_j\rangle^{-\beta_j}\\
 &<\infty.
\end{align*}

Define the quartic sum
\begin{align*}
 \mathscr S(t)
 &:=\sum_{n\in\Z^\nu}\overline{v_n(t)}\mathcal N^3[v](t,n)\\
 &=\sum_{\substack{n,n_1,n_2,n_3\in\Z^\nu\\
                   n_1-n_2+n_3=n}}
 e^{-\iu t\Sigma_3}\overline{v_n(t)}v_{n_1}(t)
 \overline{v_{n_2}(t)}v_{n_3}(t). 
\end{align*}
The involution
$(n,n_1,n_2,n_3)\mapsto(n_1,n,n_3,n_2)$ preserves momentum
conservation, reverses the sign of $\Sigma_3$ and conjugates the
product of Fourier coefficients. Hence
$\mathscr S(t)=\overline{\mathscr S(t)}$.

Differentiating the mass termwise and using \eqref{eq:nf-profile-ode}, we obtain
\begin{align*}
 \frac{d}{dt}\sum_n|v_n(t)|^2
 &=2\operatorname{Re}\sum_n\overline{v_n(t)}\dot v_n(t)\\
 &=2\operatorname{Re}\sum_n\overline{v_n(t)}(\iu\varepsilon^\alpha\mathcal N^3[v](t,n))\\
 &=2\operatorname{Re}\bigl(\iu\varepsilon^\alpha\mathscr S(t)\bigr)\\
 &=0.
\end{align*}
Evaluating the constant mass at $t=0$ proves
\eqref{eq:nf-mass-conservation}.
\end{proof}

By 
Proposition~\ref{prop:nf-mass-conservation}, each of the two sums displayed above
contributes $Mv_n(t)$.
Their intersection consists only of $(n,n,n)$, whose contribution
is $|v_n(t)|^2v_n(t)$. Adding these two sums counts this term
twice, so it must be subtracted once. Consequently,
\begin{equation}\label{eq:nf-cubic-resonant-sum}
\mathcal{N}_0^3[v](t,n)
 =2Mv_n(t)-|v_n(t)|^2v_n(t).
\end{equation}
Substituting \eqref{eq:nf-cubic-resonant-sum} into \eqref{eq:nf-cubic-splitting},
and then substituting \eqref{eq:nf-cubic-splitting} into
\eqref{eq:nf-profile-ode}, we obtain
\begin{align}
 \dot v_n(t)=2\iu\varepsilon^\alpha Mv_n(t)
 -\iu\varepsilon^\alpha|v_n(t)|^2v_n(t)+\iu\varepsilon^\alpha\mathcal{N}_1^3[v](t,n).
\label{eq:nf-profile-resonant-ode}
\end{align}
The nonlinear contribution is now separated into the
\emph{mass-dependent term}, the \emph{diagonal cubic term} and the \emph{nonresonant
oscillatory term}. We next remove the mass-dependent term
by a (random) gauge transform, while retaining the other two contributions.

\subsubsection{Removing the Mass-Dependent Term}
The conserved mass $M$ generates a common phase rotation in
\eqref{eq:nf-profile-resonant-ode}. To remove it, define the (random)
gauge transform
\begin{equation}\label{eq:nf-mass-gauge}
 a_n(t)=e^{-2\iu\varepsilon^\alpha Mt}v_n(t),
 \qquad v_n(t)=e^{2\iu\varepsilon^\alpha Mt}a_n(t).
\end{equation}
The gauge factor has modulus one, so $|a_n(t)|=|v_n(t)|$ and
$\sum_{n\in\Z^\nu}|a_n(t)|^2=M$. Differentiating
\eqref{eq:nf-mass-gauge} gives
\[
 \dot a_n(t)=e^{-\iu\psi(t)}
 \bigl(\dot v_n(t)-2\iu\varepsilon^\alpha Mv_n(t)\bigr),
\]
where $\psi(t)=2\varepsilon^\alpha Mt$. Substitution of \eqref{eq:nf-profile-resonant-ode} cancels the
mass-dependent term. With our conjugation convention, the remaining cubic terms take the following form in $a_n$:
\begin{align*}
 \prod_{j=1}^3\starpow{v_{n_j}(t)}{j-1}
 &=\prod_{j=1}^3\starpow{e^{\iu\psi(t)}a_{n_j}(t)}{j-1}\\
 &=e^{\iu\psi(t)\sum_{j=1}^3(-1)^{j-1}}
   \prod_{j=1}^3\starpow{a_{n_j}(t)}{j-1}\\
 &=e^{\iu\psi(t)}\prod_{j=1}^3\starpow{a_{n_j}(t)}{j-1}.
\end{align*}
In particular, when $n_1=n_2=n_3=n$, the diagonal cubic term satisfies
$|v_n|^2v_n=e^{\iu\psi(t)}|a_n|^2a_n$, so multiplication by
$e^{-\iu\psi(t)}$ gives $|a_n|^2a_n$. Since the gauge
leaves the linear oscillatory factor $e^{-\iu t\Sigma_3}$ unchanged,
we obtain
\begin{equation}\label{eq:nf-gauged-ode}
 \dot a_n(t)=-\iu\varepsilon^\alpha|a_n(t)|^2a_n(t)
 +\iu\varepsilon^\alpha\mathcal{N}_1^3[a](t,n).
\end{equation}
The initial data are $a_n(0)=c(n)g_n$.
\begin{remark}
The original solution $u$ is recovered from
the gauged interaction coefficients $a_n$ by
\begin{equation}\label{eq:nf-solution-reconstruction}
 u(t,x)=e^{2\iu\varepsilon^\alpha Mt}
 \sum_{n\in\Z^\nu}a_n(t)e^{\iu\langle n,\omega\rangle x-\iu tQ(n)}.
\end{equation}
\end{remark}
\subsection{LDP for the Cubic Resonant Flow}
\label{sec:cubic-resonant-ldp}

Dropping the
nonresonant term from \eqref{eq:nf-gauged-ode} gives the decoupled system
\begin{equation}\label{eq:nf-resonant-ode}
 \dot b_n(t)=-\iu\varepsilon^\alpha|b_n(t)|^2b_n(t),
 \qquad b_n(0)=c(n)g_n.
\end{equation}
For each $n$, the \emph{product rule} and \eqref{eq:nf-resonant-ode} give
\begin{align*}
 \frac{d}{dt}|b_n(t)|^2
 &=2\Ree\bigl(\dot b_n(t)\overline{b_n(t)}\bigr)\\
  &=2\Ree\bigl(-\iu\varepsilon^\alpha|b_n(t)|^4\bigr)\\
  &=0.
\end{align*}
Thus $b_n$ has constant modulus, and
\begin{align*}
 |b_n(t)|^2&=|b_n(0)|^2\\
 &=|c(n)|^2|g_n|^2.
\end{align*}
Substituting this identity into the equation for $b_n$ reduces it
to a linear initial value problem with constant coefficients
\[
\begin{cases}
 \dot b_n(t)=-\iu\varepsilon^\alpha|c(n)|^2|g_n|^2b_n(t),\\
 b_n(0)=c(n)g_n.
\end{cases}
\]
Multiplying by the \emph{integrating factor} gives
\[
 \frac{d}{dt}\left(
 e^{\iu\varepsilon^\alpha t|c(n)|^2|g_n|^2}b_n(t)\right)=0.
\]
Using $b_n(0)=c(n)g_n$, we obtain
\begin{equation}\label{eq:nf-resonant-profile}
 b_n(t)=c(n)g_n
 \exp\bigl(-\iu\varepsilon^\alpha t|c(n)|^2|g_n|^2\bigr).
\end{equation}
Substituting $b_n$ for $a_n$ in \eqref{eq:nf-solution-reconstruction}
and using \eqref{eq:nf-resonant-profile}, we obtain the cubic resonant flow
\begin{align}
 u_{\mathrm{res}}(t,x)
 &=e^{2\iu\varepsilon^\alpha Mt}
   \sum_{n\in\Z^\nu}b_n(t)
   e^{\iu\langle n,\omega\rangle x-\iu tQ(n)}
   \label{eq:nf-resonant-reconstruction}\\
 &=e^{2\iu\varepsilon^\alpha Mt}
   \sum_{n\in\Z^\nu}c(n)g_n
   e^{\iu\langle n,\omega\rangle x-\iu tQ(n)-\iu\varepsilon^\alpha t|c(n)|^2|g_n|^2}.
   \label{eq:nf-resonant-reference}
\end{align}

\begin{proposition}[Law of the Cubic Resonant Flow]
\label{prop:nf-resonant-law}
Assume $c\in\ell^1(\Z^\nu)$ and set
\[
 Z(x)\triangleq u(0,x)=\sum_{n\in\Z^\nu}c(n)g_ne^{\iu\langle n,\omega\rangle x}.
\]
For every $t\ge0$,
$u_{\mathrm{res}}(t,\cdot)$ and $Z$ have the same law as bounded
continuous random functions. In particular, for every $r>0$,
\[
 \Pp\bigl(\|u_{\mathrm{res}}(t)\|_\infty>r\bigr)
 =\Pp\bigl(\|Z\|_\infty>r\bigr).
\]
\end{proposition}

\begin{proof}[Proof of \autoref{prop:nf-resonant-law}]
Write $g_n=r_ne^{\iu\theta_n}$, where $r_n^2\sim\mathscr{E}(1)$ and
$\theta_n\sim\mathscr U([0,2\pi))$, with all radii and angles
 mutually independent. In polar coordinates, the two fields are
  \begin{align*}
   Z(x)&=\sum_{n\in\Z^\nu}c(n)r_n
    e^{\iu(\theta_n+\langle n,\omega\rangle x)},\\
   u_{\mathrm{res}}(t,x)&=\sum_{n\in\Z^\nu}c(n)r_n
    e^{\iu(\theta_n+\langle n,\omega\rangle x +\phi_n)},
  \end{align*}
  where the additional phase in $u_{\mathrm{res}}$ is
\[
 \phi_n=2\varepsilon^\alpha Mt-tQ(n)
  -\varepsilon^\alpha t|c(n)|^2r_n^2.
\]
Here $M=\sum_n|c(n)|^2r_n^2$, so every $\phi_n$ depends only on
the radii. Conditioning on $(r_n)_n$
leaves the angles independent and uniformly distributed on
$[0,2\pi)$, with joint distribution
\[
 dm(\theta)=\bigotimes_{n\in\Z^\nu}\frac{d\theta_n}{2\pi}.
\] 

For any Borel set $A\subset C_b(\R)$, let $\mathbf 1_A$ denote its
indicator function. Independence and translation invariance of the product measure give,
almost surely,
\begin{align*}
 \Ee\!\left[\mathbf 1_A(u_{\mathrm{res}}(t,\cdot))\mid(r_n)_n\right]
 &=\int\mathbf 1_A\!\left(
   \sum_n c(n)r_ne^{\iu(\theta_n+\phi_n)}
                e^{\iu\langle n,\omega\rangle x}\right)dm(\theta)\\
 &=\int\mathbf 1_A\!\left(
   \sum_n c(n)r_ne^{\iu\theta_n}
                e^{\iu\langle n,\omega\rangle x}\right)dm(\theta)\\
 &=\Ee\!\left[\mathbf 1_A(Z)\mid(r_n)_n\right].
\end{align*}
Taking expectations over the radii and using the \emph{tower property},
we obtain
\begin{align*}
 \Ee\!\left[\mathbf 1_A(u_{\mathrm{res}}(t,\cdot))\right]
 &=\Ee\!\left[\Ee\!\left[
   \mathbf 1_A(u_{\mathrm{res}}(t,\cdot))\mid(r_n)_n
   \right]\right]\\
 &=\Ee\!\left[\Ee\!\left[\mathbf 1_A(Z)\mid(r_n)_n\right]\right]\\
 &=\Ee\!\left[\mathbf 1_A(Z)\right].
\end{align*}
Thus $\Pp(u_{\mathrm{res}}(t,\cdot)\in A)=\Pp(Z\in A)$ for every
Borel set $A$, proving equality in law. Taking
$A=\{f\in C_b(\R):\|f\|_\infty>r\}$ gives the tail identity.
\end{proof}

\begin{theorem}[Large Deviations for the Cubic Resonant Flow]
\label{thm:nf-resonant-LDP}
Under the assumptions of \autoref{cor:randomized-decay}, suppose that
$\beta_j=\rho_j-\kappa_j>1$ for every $j=1,\ldots,\nu$. 
Then, for every $z>0$ and $t\ge0$,
\begin{equation}\label{eq:nf-resonant-LDP}
 \lim_{\varepsilon\to0^+}
 \varepsilon\log\Pp\left(
 \|u_{\mathrm{res}}(t)\|_\infty
 >z\varepsilon^{-1/2}
 \right)=-\frac{z^2}{S}.
\end{equation}
\end{theorem}

\begin{proof}[Proof of \autoref{thm:nf-resonant-LDP}]
The decay condition \eqref{eq:deterministic-decay} implies
$c\in\ell^1(\Z^\nu)$. By \autoref{prop:nf-resonant-law},
$u_{\mathrm{res}}(t,\cdot)$ has the same law as the Gaussian field
$Z$ defined there. It therefore suffices to estimate
\begin{equation}\label{eq:gaussian-tail-probability}
 p_\varepsilon(z):=\Pp\bigl(\|Z\|_\infty>z\varepsilon^{-1/2}\bigr).
\end{equation}
Following
\cite{XL26}, we obtain the lower bound
\eqref{eq:nf-gaussian-tail-lower-bound} at one spatial point and
the upper bound by Fourier truncation.

\medskip
\noindent\textbf{Step 1. Lower Bound.}
Since $\|Z\|_\infty\ge |Z(0)|$ and $|Z(0)|^2$ is \emph{exponentially
distributed} with rate $1/S$, we have
\begin{align}
 p_\varepsilon(z)&=\Pp\bigl(\|Z\|_\infty>z\varepsilon^{-1/2}\bigr)\nonumber\\
 &\ge\Pp\left(|Z(0)|>z\varepsilon^{-1/2}\right)\nonumber\\
 &=\Pp\left(|Z(0)|^2>\frac{z^2}{\varepsilon}\right)\nonumber\\
 &=\exp\left(-\frac{z^2}{S\varepsilon}\right).\label{eq:nf-gaussian-tail-lower-bound}
\end{align}
Consequently,
\[
 \liminf_{\varepsilon\to0^+}\varepsilon\log p_\varepsilon(z)
 \ge-\frac{z^2}{S}.
\]

\medskip
\noindent\textbf{Step 2. Upper Bound for the Truncated Field.}
For $N\in\N$, use the coordinate norm $|n|_\infty$ defined in
\autoref{prop:exception-prob} and set
\begin{equation}\label{eq:nf-truncation-set}
 \Lambda_N=\{n\in\Z^\nu:|n|_\infty\le N\},\qquad
 Z_N(x)=\sum_{n\in\Lambda_N}c(n)g_ne^{\iu\langle n,\omega\rangle x}.
\end{equation}
Set $S_N:=\sum_{n\in\Lambda_N}|c(n)|^2\le S$ and
$d_N:=|\Lambda_N|=(2N+1)^\nu$. By the \emph{Cauchy--Schwarz inequality},
\begin{align*}
 \|Z_N\|_\infty^2
 &\le S_N\sum_{n\in\Lambda_N}|g_n|^2\\
 &\le S\sum_{n\in\Lambda_N}|g_n|^2.
\end{align*}
By \autoref{def:gaussian}, the variables $(|g_n|^2)_{n\in\Lambda_N}$
are independent, with
$|g_n|^2\sim\mathscr E(1)=\operatorname{Gamma}(1,1)$. Hence
\[
Y_N:=\sum_{n\in\Lambda_N}|g_n|^2\sim\operatorname{Gamma}(d_N,1),
\]
with density
$\mathbf 1_{(0,\infty)}(x)x^{d_N-1}e^{-x}/(d_N-1)!$. 
So for $q\ge0$,
\begin{align*}
 \Pp(Y_N>q)
 & =\frac{1}{(d_N-1)!}\int_q^\infty r^{d_N-1}e^{-r}\,dr\\
 & =e^{-q}\sum_{k=0}^{d_N-1}\frac{q^k}{k!}.
\end{align*}
For a threshold $a>0$, set
$q_\varepsilon(a):=a^2/(S\varepsilon)$.
For $\varepsilon$ sufficiently small that $q_\varepsilon(a)\ge1$,
we have
\begin{align}
 \Pp\bigl(\|Z_N\|_\infty>a\varepsilon^{-1/2}\bigr)
 &=\Pp\bigl(\|Z_N\|_\infty^2>a^2\varepsilon^{-1}\bigr)\nonumber\\
 &\le\Pp\bigl(SY_N>a^2\varepsilon^{-1}\bigr)\nonumber\\
 &=\Pp\bigl(Y_N>q_\varepsilon(a)\bigr)\nonumber\\
 &=e^{-q_\varepsilon(a)}
   \sum_{k=0}^{d_N-1}\frac{q_\varepsilon(a)^k}{k!}.\label{eq:nf-truncated-gamma-tail}
\end{align}
Since $q_\varepsilon(a)\ge1$, the finite sum satisfies
\begin{align}
 \sum_{k=0}^{d_N-1}\frac{q_\varepsilon(a)^k}{k!}
 &\le q_\varepsilon(a)^{d_N-1}
   \sum_{k=0}^{d_N-1}\frac{1}{k!}\nonumber\\
 &\le e\,q_\varepsilon(a)^{d_N-1}.\label{eq:nf-gamma-polynomial-bound}
\end{align}
Combining \eqref{eq:nf-truncated-gamma-tail} and
\eqref{eq:nf-gamma-polynomial-bound} gives
\begin{equation}\label{eq:nf-truncated-tail}
 \Pp\bigl(\|Z_N\|_\infty>a\varepsilon^{-1/2}\bigr)
 \le e^{-q_\varepsilon(a)+1}q_\varepsilon(a)^{d_N-1}.
\end{equation}
Taking logarithms in \eqref{eq:nf-truncated-tail} and multiplying
by $\varepsilon>0$ give, with $d_N=(2N+1)^\nu$,
\begin{align}
 \varepsilon\log\Pp\bigl(\|Z_N\|_\infty
          >a\varepsilon^{-1/2}\bigr)
 \le{}&d_N\varepsilon|\log\varepsilon|
 +(2\log a-\log S)d_N\varepsilon\nonumber\\
 &-\varepsilon|\log\varepsilon|
 -(2\log a-\log S-1)\varepsilon-a^2/S.\label{eq:nf-truncation-polynomial}
\end{align}
Fix $\mu\in(0,1/\nu)$ and choose the integer-valued truncation
\begin{align}
  N(\varepsilon)
  &=\bigl\lfloor\varepsilon^{-\frac{1}{\nu}+\mu}\bigr\rfloor\nonumber\\
  &\asymp\varepsilon^{-\frac{1}{\nu}+\mu}\longrightarrow\infty.\label{eq:nf-truncation-scale}
\end{align}
For this truncation, $d_{N(\varepsilon)}=(2N(\varepsilon)+1)^\nu
=\mathcal{O}(\varepsilon^{-1+\nu\mu})$. Since $\nu\mu>0$, it follows that
\begin{align*}
  \varepsilon d_{N(\varepsilon)}
  &=\mathcal{O}(\varepsilon^{\nu\mu})=o(1),\\
  \varepsilon d_{N(\varepsilon)}|\log\varepsilon|
  &=\mathcal{O}\bigl(\varepsilon^{\nu\mu}|\log\varepsilon|\bigr)=o(1).
\end{align*}
Hence, for fixed $a>0$, all terms other than $-a^2/S$ in
\eqref{eq:nf-truncation-polynomial} vanish as $\varepsilon\to0^+$.
Taking the upper limit gives
\begin{align*}
  \limsup_{\varepsilon\to0^+}\varepsilon\log
 \Pp\bigl(\|Z_{N(\varepsilon)}\|_\infty>a\varepsilon^{-1/2}\bigr)
 \le-\frac{a^2}{S}.
\end{align*}
\medskip
\noindent\textbf{Step 3. Bound for the Remainder.}
With $Z_N$ and $\Lambda_N$ as in \eqref{eq:nf-truncation-set},
write the truncation remainder $Z-Z_N$ as
\begin{equation}\label{eq:nf-high-mode-definition}
 \tilde{Z}_N(x)=\sum_{n\notin\Lambda_N}
 c(n)g_ne^{\iu\langle n,\omega\rangle x},\qquad x\in\R.
\end{equation}
For a fixed threshold $z>0$, choose the cutoff from
\autoref{cor:randomized-decay} by setting
\begin{equation}\label{eq:nf-high-mode-delta}
 \delta_\varepsilon:=\frac{\varepsilon}{1+z^2/S}.
\end{equation}
We explain the choice in \eqref{eq:nf-high-mode-delta}.
The target logarithmic probability has rate $-z^2/S$, whereas
\eqref{eq:exception-prob} bounds the probability of the exceptional set by $C_{\kappa,\nu}e^{-1/\delta}$. Thus we choose
$\delta_\varepsilon$ proportional to $\varepsilon$ so that
$\varepsilon/\delta_\varepsilon>z^2/S$, with a fixed positive gap.
The choice \eqref{eq:nf-high-mode-delta} makes this gap equal to one.
Indeed, for $0<\varepsilon\le1$, we have
$0<\delta_\varepsilon\le\varepsilon\le1$, so
\eqref{eq:exception-prob} applies and gives
\begin{align*}
 \Pp(G_{\delta_\varepsilon}^c)
 &=\Pp(\Omega_{\delta_\varepsilon})\\
 &\le C_{\kappa,\nu}
       \exp\left(-\frac{1+z^2/S}{\varepsilon}\right).
\end{align*}
Thus
\begin{align*}
 \limsup_{\varepsilon\to0^+}
 \varepsilon\log\Pp(G_{\delta_\varepsilon}^c)
 &\le-\left(1+\frac{z^2}{S}\right)
 <-\frac{z^2}{S}.
\end{align*}
At the same time,
\[
 \delta_\varepsilon^{-1/2}
 =(1+z^2/S)^{1/2}\varepsilon^{-1/2}.
\]
Since $z$ and $S$ are fixed, the cutoff therefore changes the
coefficient bound only by a constant factor at the
$\varepsilon^{-1/2}$ scale. The decay of the Fourier tail will
provide the additional small factor needed for the remainder.
By \autoref{cor:randomized-decay}, on $G_{\delta_\varepsilon}$ we have
\begin{equation}\label{eq:nf-high-mode-majorant}
 |c(n)g_n|\le A\delta_\varepsilon^{-1/2}w_{-\beta}(n)
 =A\delta_\varepsilon^{-1/2}
 \prod_{j=1}^\nu(1+|n_j|)^{-\beta_j},
 \qquad n\in\Z^\nu.
\end{equation}
We next sum the majorant in \eqref{eq:nf-high-mode-majorant} over
$\Lambda_N^c$ to bound the truncation remainder. For each
coordinate $\ell\in\{1,\ldots,\nu\}$, define
\[
 E_{\ell,N}:=\{n\in\Z^\nu:|n_\ell|>N\}.
\]
By \eqref{eq:nf-truncation-set},
$\Lambda_N^c=\bigcup_{\ell=1}^\nu E_{\ell,N}$. For each
$E_{\ell,N}$, the coordinates $k\ne\ell$ are unrestricted and
their 1D sums are
\[
 B_k:=\sum_{m\in\Z}(1+|m|)^{-\beta_k}<\infty.
\]
The $\ell$th coordinate contributes a tail sum, for which \emph{integral
comparison} gives
\[
 \sum_{|m|>N}(1+|m|)^{-\beta_\ell}
 \le \frac{2}{\beta_\ell-1}(1+N)^{1-\beta_\ell}.
\]
Set
\begin{equation}\label{eq:nf-high-mode-beta}
 \beta_*:=\min_{1\le j\le\nu}(\beta_j-1)>0.
\end{equation}
Since the summands are nonnegative, the sum over
$\Lambda_N^c=\bigcup_\ell E_{\ell,N}$ is bounded by the sum over
the individual sets, even if they overlap. Factoring the sums on
each $E_{\ell,N}$ and applying the tail estimate give
\begin{align}
 \sum_{n\notin\Lambda_N}
 \prod_{j=1}^\nu(1+|n_j|)^{-\beta_j}
 &\le \sum_{\ell=1}^\nu\sum_{n\in E_{\ell,N}}
 \prod_{j=1}^\nu(1+|n_j|)^{-\beta_j}\nonumber\\
 &= \sum_{\ell=1}^\nu
 \left(\sum_{|m|>N}(1+|m|)^{-\beta_\ell}\right)
 \prod_{k\ne\ell}B_k\nonumber\\
 &\le C_{\beta,\nu}(1+N)^{-\beta_*},\label{eq:nf-polynomial-tail}
\end{align}
where $C_{\beta,\nu}:=\sum_{\ell=1}^\nu
\frac{2}{\beta_\ell-1}\prod_{k\ne\ell}B_k<\infty$ depends only on
$\beta$ and $\nu$. The last inequality uses
$\beta_\ell-1\ge\beta_*$ from \eqref{eq:nf-high-mode-beta}.
Now take $N=N(\varepsilon)$ from \eqref{eq:nf-truncation-scale}
and set
\begin{equation}\label{eq:nf-high-mode-eta}
 \eta:=\beta_*\left(\frac1\nu-\mu\right)>0.
\end{equation}
Apply the triangle inequality to \eqref{eq:nf-high-mode-definition},
then use \eqref{eq:nf-high-mode-majorant} and
\eqref{eq:nf-polynomial-tail} for the first inequality below.
For the second, substitute $\delta_\varepsilon$ from
\eqref{eq:nf-high-mode-delta} and use
$1+N(\varepsilon)>\varepsilon^{-1/\nu+\mu}$ from
\eqref{eq:nf-truncation-scale}. Thus, on $G_{\delta_\varepsilon}$,
\begin{align}
 \|\tilde Z_{N(\varepsilon)}\|_\infty
 &\le A\delta_\varepsilon^{-1/2}C_{\beta,\nu}
       (1+N(\varepsilon))^{-\beta_*}\nonumber\\
 &\le A(1+z^2/S)^{1/2}C_{\beta,\nu}
       \varepsilon^{-1/2+\eta}.\label{eq:nf-high-mode-remainder}
\end{align}

\medskip
\noindent\textbf{Step 4. Relative Probability and Conclusion.}
Write the remainder bound \eqref{eq:nf-high-mode-remainder} as
$\|\tilde Z_{N(\varepsilon)}\|_\infty
\le h_\varepsilon\varepsilon^{-1/2}$ on $G_{\delta_\varepsilon}$, where
\[
 h_\varepsilon:=A(1+z^2/S)^{1/2}C_{\beta,\nu}\varepsilon^\eta.
\]
By \eqref{eq:nf-high-mode-eta}, $h_\varepsilon\to0$; hence
$0<h_\varepsilon<z/2$ for all sufficiently small $\varepsilon$.
On $G_{\delta_\varepsilon}$, the decomposition
$Z=Z_{N(\varepsilon)}+\tilde Z_{N(\varepsilon)}$ and the reverse
triangle inequality give
\[
 \|Z_{N(\varepsilon)}\|_\infty
 \ge \|Z\|_\infty-\|\tilde Z_{N(\varepsilon)}\|_\infty
 \ge \|Z\|_\infty-h_\varepsilon\varepsilon^{-1/2}.
\]
Consequently, if $\|Z\|_\infty>z\varepsilon^{-1/2}$ holds on
$G_{\delta_\varepsilon}$, then
$\|Z_{N(\varepsilon)}\|_\infty
>(z-h_\varepsilon)\varepsilon^{-1/2}$. Equivalently,
\[
 \{\|Z\|_\infty>z\varepsilon^{-1/2}\}\cap G_{\delta_\varepsilon}
 \subseteq
 \{\|Z_{N(\varepsilon)}\|_\infty
     >(z-h_\varepsilon)\varepsilon^{-1/2}\}.
\]

For the threshold $a_\varepsilon:=z-h_\varepsilon>0$, define
\[
  p_\varepsilon^{\mathrm{tr}}:=\Pp\bigl(\|Z_{N(\varepsilon)}\|_\infty
        >a_\varepsilon\varepsilon^{-1/2}\bigr),
 \qquad p_\varepsilon^{\mathrm{exc}}:=\Pp(G_{\delta_\varepsilon}^c).
\]
Splitting $p_\varepsilon(z)$ according to the
partition $G_{\delta_\varepsilon}\cup G_{\delta_\varepsilon}^c$ and applying the event
inclusion on $G_{\delta_\varepsilon}$ yield
\begin{align}
 p_\varepsilon(z)
 &=\Pp\bigl(\{\|Z\|_\infty>z\varepsilon^{-1/2}\}
                  \cap G_{\delta_\varepsilon}\bigr)
   +\Pp\bigl(\{\|Z\|_\infty>z\varepsilon^{-1/2}\}
                  \cap G_{\delta_\varepsilon}^c\bigr)\nonumber\\
 &\le\Pp\bigl(\{\|Z_{N(\varepsilon)}\|_\infty
                  >a_\varepsilon\varepsilon^{-1/2}\}\bigr)
       +\Pp(G_{\delta_\varepsilon}^c)\nonumber\\
 &\le p_\varepsilon^{\mathrm{tr}}
       +p_\varepsilon^{\mathrm{exc}}.\label{eq:nf-final-union}
\end{align}
For sufficiently small $\varepsilon$, the cutoff $N(\varepsilon)$ retains a mode
with $c(n)\ne0$, so $p_\varepsilon^{\mathrm{tr}}>0$.
Taking logarithms in \eqref{eq:nf-final-union}, multiplying by
$\varepsilon>0$, and factoring out $p_\varepsilon^{\mathrm{tr}}$
give
\begin{equation}\label{eq:nf-final-log-bound}
 \varepsilon\log p_\varepsilon(z)
 \le\varepsilon\log p_\varepsilon^{\mathrm{tr}}
 +\varepsilon\log\left(1+
 \frac{p_\varepsilon^{\mathrm{exc}}}{p_\varepsilon^{\mathrm{tr}}}\right).
\end{equation}
To estimate $p_\varepsilon^{\mathrm{tr}}$, apply
\eqref{eq:nf-truncation-polynomial} with $a=a_\varepsilon$ and
$q_\varepsilon(a_\varepsilon)=a_\varepsilon^2/(S\varepsilon)$. Since
$z/2<a_\varepsilon<z$ for small $\varepsilon$,
$\log q_\varepsilon(a_\varepsilon)
=\log(1/\varepsilon)+\mathcal{O}(1)$ and
$q_\varepsilon(a_\varepsilon)\ge1$ eventually.
The scale \eqref{eq:nf-truncation-scale} gives
$\varepsilon d_{N(\varepsilon)}\log(1/\varepsilon)\to0$.
Thus all terms in the substituted estimate except its main term
$-a_\varepsilon^2/S$ tend to zero. As
$a_\varepsilon\to z$, we obtain
\[
 \limsup_{\varepsilon\to0^+}\varepsilon\log
 p_\varepsilon^{\mathrm{tr}}
 \le-\frac{z^2}{S}.
\]
For small $\varepsilon$, apply \eqref{eq:exception-prob} with
$\delta=\delta_\varepsilon$ from \eqref{eq:nf-high-mode-delta}
to get
\[
 p_\varepsilon^{\mathrm{exc}}
 \le C_{\kappa,\nu}e^{-(1+z^2/S)/\varepsilon}.
\]
Combining the lower bound
$p_\varepsilon(z)\ge e^{-z^2/(S\varepsilon)}$ from Step 1
with \eqref{eq:nf-final-union} and the preceding bound for the exceptional probability gives
\begin{align*}
 p_\varepsilon^{\mathrm{tr}}
 &\ge p_\varepsilon(z)-p_\varepsilon^{\mathrm{exc}}\\
 &\ge e^{-z^2/(S\varepsilon)}
       \bigl(1-C_{\kappa,\nu}e^{-1/\varepsilon}\bigr).
\end{align*}
Dividing the upper bound for
$p_\varepsilon^{\mathrm{exc}}$ by this lower bound for
$p_\varepsilon^{\mathrm{tr}}$ yields the relative probability
needed in \eqref{eq:nf-final-log-bound}:
\begin{equation*}
 0\le\frac{p_\varepsilon^{\mathrm{exc}}}
             {p_\varepsilon^{\mathrm{tr}}}
 \le\frac{C_{\kappa,\nu}e^{-1/\varepsilon}}
          {1-C_{\kappa,\nu}e^{-1/\varepsilon}}
 \le2C_{\kappa,\nu}e^{-1/\varepsilon}\longrightarrow0.
\end{equation*}
Using $0\le\log(1+r)\le r$ for $r\ge0$, the correction term in
\eqref{eq:nf-final-log-bound} lies between zero and
$2\varepsilon C_{\kappa,\nu}e^{-1/\varepsilon}$ and therefore
tends to zero. Combining this with the truncated-field limit gives
\[
 \limsup_{\varepsilon\to0^+}\varepsilon\log p_\varepsilon(z)
 \le-\frac{z^2}{S}.
\]
The lower bound in Step 1 also gives
$\liminf_{\varepsilon\to0^+}\varepsilon\log p_\varepsilon(z)
\ge-z^2/S$. The upper and lower limits agree, proving
\eqref{eq:nf-resonant-LDP}.

\end{proof}

\subsection{Time Scale: Existence, Uniqueness and Uniform Decay}
\label{sec:critical-existence}
\begin{proposition}[Time Scale: Existence, Uniqueness and Uniform Decay]
\label{prop:critical-pointwise-lifespan}
Fix $\delta>0$ and set $d_n:=c(n)g_n$.
On $G_\delta$, the coefficients $d_n$ satisfy the decay estimate
\eqref{eq:randomized-decay}, where
$\beta_j=\rho_j-\kappa_j>2$ for every $j=1,\ldots,\nu$; see \autoref{cor:randomized-decay}.
Define
\begin{equation}\label{eq:critical-pointwise-lifespan}
 T_\delta:=\frac{\delta}{32\varepsilon^\alpha A^2
 \|w_{-\beta/2}\|_{\ell^1(\Z^\nu)}^2}.
\end{equation}
Then \eqref{eq:nf-gauged-ode} with initial data $a_n(0)=d_n$ admits
a $C^1$ solution on $[0,T_\delta]$ satisfying, for
every $n\in\Z^\nu$,
\begin{equation}\label{eq:critical-pointwise-decay}
 \sup_{0\le t\le T_\delta}|a_n(t)|
 \le2A\delta^{-1/2}\|w_{-\beta/2}\|_{\ell^1(\Z^\nu)}
 w_{-\beta/2}(n).
\end{equation}
This solution is unique among $C^1$ solutions for which
$|a_n(t)|\le Cw_{-\beta/2}(n)$ for all $n\in\Z^\nu$ and
$t\in[0,T_\delta]$, with some finite constant $C$.

\end{proposition}

\begin{proof}[Proof of \autoref{prop:critical-pointwise-lifespan}]
\medskip
\noindent\textbf{Step 1. The Cubic Coefficient Equation.}
For each output index $n$, define the discrete cubic interaction
kernel on triples satisfying $n_1-n_2+n_3=n$ by
\[
 \mathcal K_{t,n}(n_1,n_2,n_3)
 :=-\mathbf 1_{\{n_1=n_2=n_3=n\}}
   +e^{-\iu t\Sigma_3}\mathbf 1_{\{\Sigma_3\ne0\}},
\]
where $\Sigma_3=\Sigma_3(n_1,n_2,n_3;n)$. The first indicator gives
$-|a_n|^2a_n$; the second retains the nonresonant interactions in
\eqref{eq:nf-gauged-ode}. Their supports are disjoint because
$\Sigma_3(n,n,n;n)=0$, so $|\mathcal K_{t,n}|\le1$.
Define the cubic convolution map on sequences $a=(a_n)_{n\in\Z^\nu}$ by
\begin{equation}\label{eq:critical-cubic-map}
 H_{t,n}[a]:=
 \sum_{n_1-n_2+n_3=n}\mathcal K_{t,n}(n_1,n_2,n_3)
 \prod_{j=1}^3\starpow{a_{n_j}}{j-1}.
\end{equation}
We write $H_t[a]:=(H_{t,n}[a])_{n\in\Z^\nu}$ for the corresponding
sequence-valued map.
Thus
equation \eqref{eq:nf-gauged-ode} becomes
\begin{equation}\label{eq:critical-gauged-cubic}
 \dot a_n(t)=\iu\varepsilon^\alpha H_{t,n}[a(t)],
 \qquad a_n(0)=d_n.
\end{equation}

\medskip
\noindent\textbf{Step 2. Picard Iteration and Tree Expansion.}
We solve \eqref{eq:critical-gauged-cubic} by \emph{Picard successive
approximation}, using the \emph{tree expansion method} developed in
\cite{DLX24I}, following \cite{C07,DG16}. Starting from the initial data $d$, define
\begin{equation}\label{eq:critical-pointwise-picard}
 a_n^{(0)}(t)=d_n,\qquad
 a_n^{(k)}(t)=d_n+\iu\varepsilon^\alpha
        \int_0^t H_{r,n}[a^{(k-1)}(r)]\,dr,
 \quad k\ge1.
\end{equation}
\par\noindent
\begin{definition}[Branch Sets]
The branch sets $\Gamma^{(k)}$ are defined recursively by
\[
 \Gamma^{(1)}=\{0,1\}
 =\left\{
  \scalebox{0.75}{\begin{tikzpicture}[baseline=(current bounding box.center),x=.55cm,y=.55cm,
   line width=.55pt,line cap=round]
   \draw (0,0)--(0,1);
   \fill (0,0) circle (1.9pt);
   \fill (0,1) circle (1.6pt);
   \node[below=1pt] at (0,0) {$n$};
   \node[above=1pt] at (0,1) {$n$};
   \end{tikzpicture}},\quad
  \scalebox{0.75}{\begin{tikzpicture}[baseline=(current bounding box.center),x=.55cm,y=.55cm,
   line width=.55pt,line cap=round]
   \coordinate (root) at (0,0);
   \draw (root)--(-.8,1);
   \draw[red!75!black] (root)--(0,1);
   \draw (root)--(.8,1);
   \fill (root) circle (1.9pt);
   \fill (-.8,1) circle (1.6pt);
   \fill[red!75!black] (0,1) circle (1.6pt);
   \fill (.8,1) circle (1.6pt);
   \node[below=1pt] at (root) {$n$};
   \node[above=1pt] at (-.8,1) {$n_1$};
   \node[above=1pt] at (0,1) {$n_2$};
   \node[above=1pt] at (.8,1) {$n_3$};
   \end{tikzpicture}}
 \right\},
 \qquad 
 \Gamma^{(k)}=\{0\}\cup(\Gamma^{(k-1)})^3,\qquad k\ge2.
\]
The diagrams represent the labels $0$ and $1$. In $0$, the edge
carries $n$ unchanged to its terminal leaf. In a diagram at level
$k$, a $0$ branch continues by identity transmission to the terminal
level; this extension is fixed by $k$ and introduces no further
cubic vertex or branch choice. In $1$,
the root has label $n$ and its three
ordered leaves have labels $n_1,n_2,n_3$, where
$n=n_1-n_2+n_3$. The red middle edge marks the minus sign and the
conjugated factor. At every level $k\ge2$, the separate label $0$
represents the initial data.
\end{definition}
\begin{definition}[Dual Tree]
The dual of a colored ordered ternary tree preserves its root,
ordered branches, leaves and frequency labels, while exchanging red
and black on every vertex and edge. Taking the dual twice recovers
the original tree. For example, the dual of $1$ is
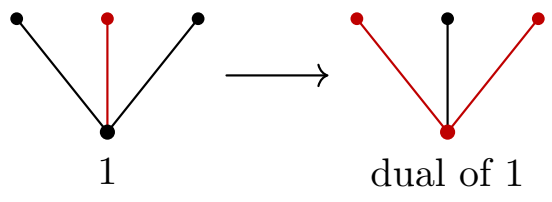
\begin{center}
\scalebox{0.75}{\begin{tikzpicture}[x=1cm,y=1cm,line width=.55pt,line cap=round]
  \coordinate (leftroot) at (0,0);
  \draw (leftroot)--(-.8,1);
  \draw[red!75!black] (leftroot)--(0,1);
  \draw (leftroot)--(.8,1);
  \fill (leftroot) circle (1.9pt);
  \fill (-.8,1) circle (1.6pt);
  \fill[red!75!black] (0,1) circle (1.6pt);
  \fill (.8,1) circle (1.6pt);
  \node[below=3pt] at (leftroot) {$1$};
  \draw[->] (1.05,.5)--(1.95,.5);
  \coordinate (rightroot) at (3,0);
  \draw[red!75!black] (rightroot)--(2.2,1);
  \draw (rightroot)--(3,1);
  \draw[red!75!black] (rightroot)--(3.8,1);
  \fill[red!75!black] (rightroot) circle (1.9pt);
  \fill[red!75!black] (2.2,1) circle (1.6pt);
  \fill (3,1) circle (1.6pt);
  \fill[red!75!black] (3.8,1) circle (1.6pt);
  \node[below=3pt] at (rightroot) {dual of $1$};
\end{tikzpicture}}
\end{center}
\end{definition}
\begin{definition}[Cubic Tree Product Rule]
\label{def:critical-cubic-tree-product}
For $k\ge2$, graft the ordered triple
$(\gamma_1,\gamma_2,\gamma_3)\in(\Gamma^{(k-1)})^3$ onto a new cubic
root, from left to right. The middle red edge replaces $\gamma_2$ by
its dual; the black outer edges preserve $\gamma_1$ and $\gamma_3$:
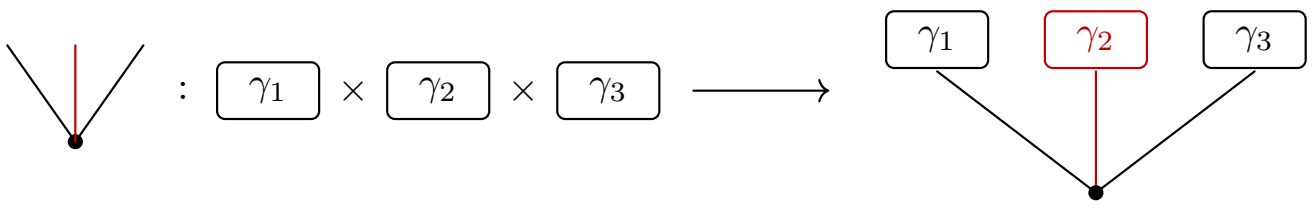
\begin{center}
\scalebox{0.75}{\begin{tikzpicture}[x=1cm,y=1cm,line width=.55pt,line cap=round,
  every node/.style={inner sep=2pt}]
  \fill (-6.3,-.1) circle (1.9pt);
  \foreach \x in {-6.9,-5.7} {
    \draw (-6.3,-.1)--(\x,.75);
  }
  \draw[red!75!black] (-6.3,-.1)--(-6.3,.75);
  \node at (-5.35,.35) {$:$};
  \foreach \x/\j in {-4.6/1,-3.1/2,-1.6/3}
    \node[draw,rounded corners=2pt,minimum width=.9cm,minimum height=.5cm]
      at (\x,.35) {$\gamma_{\j}$};
  \node at (-3.85,.35) {$\times$};
  \node at (-2.35,.35) {$\times$};
  \draw[->] (-.85,.35)--(.35,.35);
  \foreach \x/\j in {1.3/1,4.1/3} {
    \node[draw,rounded corners=2pt,minimum width=.9cm,minimum height=.5cm]
      at (\x,.8) {$\gamma_{\j}$};
    \draw (\x,.52)--(2.7,-.55);
  }
  \node[draw=red!75!black,text=red!75!black,rounded corners=2pt,
    minimum width=.9cm,minimum height=.5cm] at (2.7,.8) {$\gamma_2$};
  \draw[red!75!black] (2.7,.52)--(2.7,-.55);
  \fill (2.7,-.55) circle (1.9pt);
\end{tikzpicture}}
\end{center}
Correspondingly, the dual cubic map exchanges red and black in the
whole diagram: its outer edges dualize $\gamma_1$ and $\gamma_3$,
while its middle edge preserves $\gamma_2$:
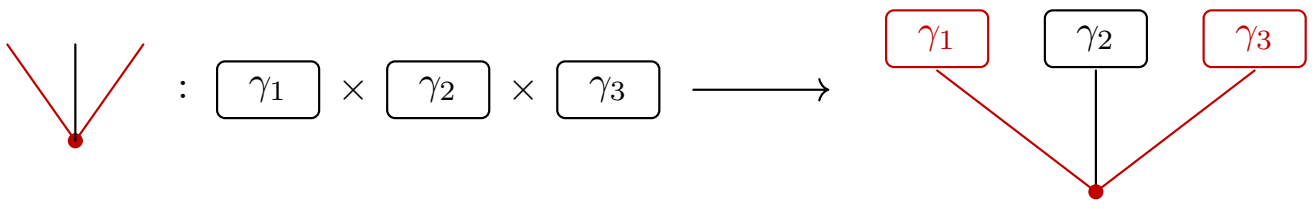
\begin{center}
\scalebox{0.75}{\begin{tikzpicture}[x=1cm,y=1cm,line width=.55pt,line cap=round,
  every node/.style={inner sep=2pt}]
  \fill[red!75!black] (-6.3,-.1) circle (1.9pt);
  \foreach \x in {-6.9,-5.7}
    \draw[red!75!black] (-6.3,-.1)--(\x,.75);
  \draw (-6.3,-.1)--(-6.3,.75);
  \node at (-5.35,.35) {$:$};
  \foreach \x/\j in {-4.6/1,-3.1/2,-1.6/3}
    \node[draw,rounded corners=2pt,minimum width=.9cm,minimum height=.5cm]
      at (\x,.35) {$\gamma_{\j}$};
  \node at (-3.85,.35) {$\times$};
  \node at (-2.35,.35) {$\times$};
  \draw[->] (-.85,.35)--(.35,.35);
  \foreach \x/\j in {1.3/1,4.1/3} {
    \node[draw=red!75!black,text=red!75!black,rounded corners=2pt,
      minimum width=.9cm,minimum height=.5cm]
      at (\x,.8) {$\gamma_{\j}$};
    \draw[red!75!black] (\x,.52)--(2.7,-.55);
  }
  \node[draw,rounded corners=2pt,minimum width=.9cm,minimum height=.5cm]
    at (2.7,.8) {$\gamma_2$};
  \draw (2.7,.52)--(2.7,-.55);
  \fill[red!75!black] (2.7,-.55) circle (1.9pt);
\end{tikzpicture}}
\end{center}
\end{definition}\par\noindent

Set $\operatorname{Term}_0(t,n)=d_n$. The tree $1$ has three leaves,
so its contribution is
\[
 \operatorname{Term}_1(t,n)
 =\iu\varepsilon^\alpha\int_0^t
 \sum_{n_1-n_2+n_3=n}\mathcal K_{r,n}(n_1,n_2,n_3)
 \prod_{j=1}^3\starpow{d_{n_j}}{j-1}\,dr.
\]
For $k\ge2$ and
$\gamma=(\gamma_1,\gamma_2,\gamma_3)\in(\Gamma^{(k-1)})^3$,
define recursively
\[
 \begin{aligned}
 \operatorname{Term}_\gamma(t,n)
 =\iu\varepsilon^\alpha\int_0^t
 \sum_{n_1-n_2+n_3=n}\mathcal K_{r,n}(n_1,n_2,n_3)\prod_{j=1}^3
 \starpow{\operatorname{Term}_{\gamma_j}(r,n_j)}{j-1}\,dr.
 \end{aligned}
\]
\begin{example}[$\Gamma^{(2)}\setminus\{0\}$]
\label{ex:critical-second-level-trees}
Since $\Gamma^{(1)}=\{0,1\}$, the eight ordered triples in
$\Gamma^{(2)}\setminus\{0\}=\{0,1\}^3$ contribute the following
cubic products at time $r$:
\[
\begin{array}{c@{\quad\longmapsto\quad}l}
(0,0,0)&d_{n_1}\overline{d_{n_2}}d_{n_3}\\
(0,0,1)&d_{n_1}\overline{d_{n_2}}\operatorname{Term}_1(r,n_3)\\
(0,1,0)&d_{n_1}\overline{\operatorname{Term}_1(r,n_2)}d_{n_3}\\
(0,1,1)&d_{n_1}\overline{\operatorname{Term}_1(r,n_2)}\operatorname{Term}_1(r,n_3)\\
(1,0,0)&\operatorname{Term}_1(r,n_1)\overline{d_{n_2}}d_{n_3}\\
(1,0,1)&\operatorname{Term}_1(r,n_1)\overline{d_{n_2}}\operatorname{Term}_1(r,n_3)\\
(1,1,0)&\operatorname{Term}_1(r,n_1)\overline{\operatorname{Term}_1(r,n_2)}d_{n_3}\\
(1,1,1)&\operatorname{Term}_1(r,n_1)\overline{\operatorname{Term}_1(r,n_2)}\operatorname{Term}_1(r,n_3)
\end{array}
\]
The colored diagrams in \autoref{fig:critical-second-level-trees}
represent these eight products in the same order; each red branch
records the conjugation in the corresponding factor.
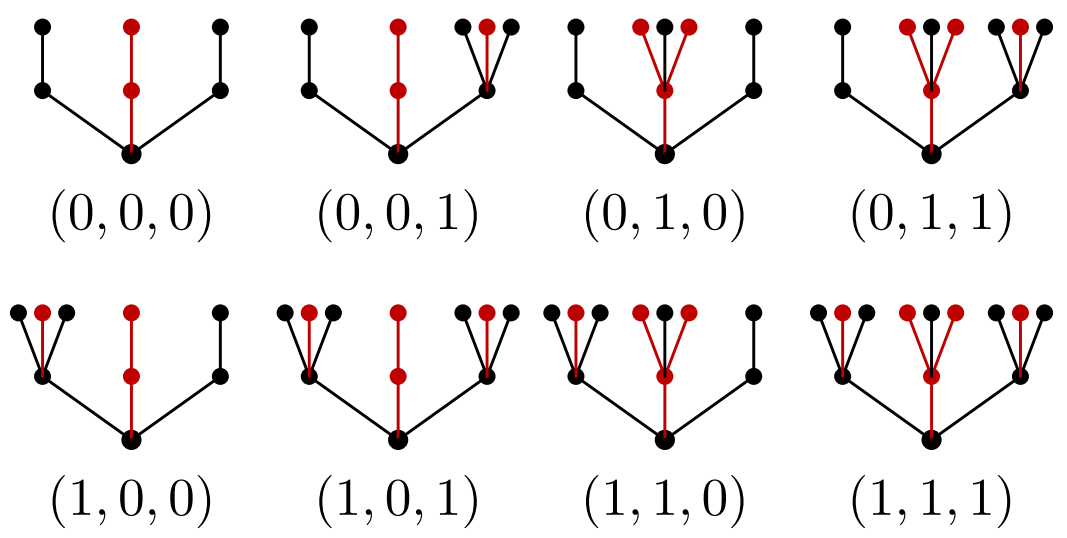
\begin{figure}[htbp]
\centering
\begin{tikzpicture}[x=.42cm,y=.42cm,line width=.55pt,line cap=round]
\foreach \a/\b/\c/\col/\row in {
  0/0/0/0/0,0/0/1/1/0,0/1/0/2/0,0/1/1/3/0,
  1/0/0/0/1,1/0/1/1/1,1/1/0/2/1,1/1/1/3/1} {
  \begin{scope}[shift={({\col*4.2},{-\row*4.5})}]
    \fill (0,0) circle (1.9pt);
    \foreach \x/\v in {-1.4/\a,0/\b,1.4/\c} {
      \ifdim\x pt=0pt
        \draw[red!75!black] (0,0)--(\x,1);
        \fill[red!75!black] (\x,1) circle (1.6pt);
        \ifnum\v=0
          \draw[red!75!black] (\x,1)--(\x,2);
          \fill[red!75!black] (\x,2) circle (1.6pt);
        \else
          \foreach \s in {-.38,.38} {
            \draw[red!75!black] (\x,1)--({\x+\s},2);
            \fill[red!75!black] ({\x+\s},2) circle (1.6pt);
          }
          \draw (\x,1)--(\x,2);
          \fill (\x,2) circle (1.6pt);
        \fi
      \else
        \draw (0,0)--(\x,1);
        \fill (\x,1) circle (1.6pt);
        \ifnum\v=0
          \draw (\x,1)--(\x,2);
          \fill (\x,2) circle (1.6pt);
        \else
          \foreach \s in {-.38,.38} {
            \draw (\x,1)--({\x+\s},2);
            \fill ({\x+\s},2) circle (1.6pt);
          }
          \draw[red!75!black] (\x,1)--(\x,2);
          \fill[red!75!black] (\x,2) circle (1.6pt);
        \fi
      \fi
    }
    \node[below=3pt] at (0,0) {$(\a,\b,\c)$};
  \end{scope}
}
\end{tikzpicture}
\caption{The eight colored level-$2$ tree diagrams, in the order of
the cubic products above. The red middle edge dualizes its child.}
\label{fig:critical-second-level-trees}
\end{figure}
\end{example}

\begin{lemma}[Tree Expansion of Picard Iterates]
\label{lem:critical-picard-tree-expansion}
Every Picard iterate has the expansion
\begin{equation}\label{eq:critical-picard-tree-expansion}
 a_n^{(k)}(t)=\sum_{\gamma\in\Gamma^{(k)}}
 \operatorname{Term}_\gamma(t,n),\qquad k\ge1.
\end{equation}
\end{lemma}
\begin{proof}[Proof of \autoref{lem:critical-picard-tree-expansion}]
We prove this by induction on $k$. For $k=1$, the three children of
$1$ are leaves, so the definition of $\operatorname{Term}_1$ gives
\begin{align*}
 a_n^{(1)}(t)
 &=d_n+\iu\varepsilon^\alpha\int_0^t H_{r,n}[d]\,dr\\
 &=\operatorname{Term}_0(t,n)+\operatorname{Term}_1(t,n)
 =\sum_{\gamma\in\Gamma^{(1)}}\operatorname{Term}_\gamma(t,n).
 \end{align*}
Suppose \eqref{eq:critical-picard-tree-expansion} holds at level
$k-1\ge1$ for every time and index. Substitute the three finite
\emph{tree expansions} into the cubic term of
\eqref{eq:critical-pointwise-picard} to get
\begin{align*}
 a_n^{(k)}(t)
 &=d_n+\iu\varepsilon^\alpha\int_0^t
   \sum_{n_1-n_2+n_3=n}\mathcal K_{r,n}(n_1,n_2,n_3)\cdot\prod_{j=1}^3
   \starpow{a_{n_j}^{(k-1)}(r)}{j-1}\,dr\\
 &=\operatorname{Term}_0(t,n)
   +\sum_{(\gamma_1,\gamma_2,\gamma_3)
             \in(\Gamma^{(k-1)})^3}
    \operatorname{Term}_{(\gamma_1,\gamma_2,\gamma_3)}(t,n)\\
 &=\sum_{\gamma\in\Gamma^{(k)}}
    \operatorname{Term}_\gamma(t,n).
 \end{align*}
\end{proof}
\medskip
\noindent\textbf{Step 3. Uniform Bound and the Time Scale.}
Define
\begin{equation}\label{eq:critical-tree-vertex-count}
 \ell(\gamma)=
 \begin{cases}
  0, & \gamma=0,\\
  1, & \gamma=1,\\
  1+\displaystyle\sum_{j=1}^3\ell(\gamma_j),
     & \gamma=(\gamma_1,\gamma_2,\gamma_3).
 \end{cases}
\end{equation}
Geometrically, $\ell(\gamma)$ is the number of cubic vertices of
$\gamma$. Analytically, it is the total number of time integration
operators in the recursive expression for
$\operatorname{Term}_\gamma$. These counts agree because $0$ has
neither a cubic vertex nor a time integral, $1$ has one of each,
and grafting three child trees adds one cubic vertex and one time
integral at the new root.

\begin{lemma}[Number of Leaves]\label{lem:critical-tree-leaf-count}
Every tree $\gamma$ has exactly $2\ell(\gamma)+1$ leaves.
\end{lemma}
\begin{proof}[Proof of \autoref{lem:critical-tree-leaf-count}]
The tree $0$ has one leaf and no cubic vertex, while $1$ has three
leaves and one cubic vertex. Suppose
$\gamma=(\gamma_1,\gamma_2,\gamma_3)$ and the formula holds for its
three child trees. Grafting them introduces one cubic vertex but no
new leaf. Hence their leaf counts add to
\[
 \sum_{j=1}^3\bigl(2\ell(\gamma_j)+1\bigr)
 =2\Bigl(1+\sum_{j=1}^3\ell(\gamma_j)\Bigr)+1
 =2\ell(\gamma)+1.
\]
Induction on the construction of $\gamma$ completes the proof.
\end{proof}

\begin{lemma}[Momentum from Leaf Indices]\label{lem:critical-tree-momentum}
Let $\gamma$ be a tree with root $n$ and ordered leaf
indices $(m_j)_{1\le j\le 2\ell(\gamma)+1}$. Then
\[
 n=\sum_{j=1}^{2\ell(\gamma)+1}(-1)^{j-1}m_j.
\]
\end{lemma}
\begin{proof}[Proof of \autoref{lem:critical-tree-momentum}]
For $\gamma=0$, the sole leaf carries the root index, so the formula
holds. For $\gamma=1$, it is the cubic momentum relation
$n=m_1-m_2+m_3$.

Suppose the formula holds for three child trees and write
$\gamma_{123}:=(\gamma_1,\gamma_2,\gamma_3)$
for their grafting, with root index
$n_{123}$. For $i=1,2,3$, let $n_i$ be
the root index of $\gamma_i$ and
$(m_{i,j})_{1\le j\le
2\ell(\gamma_i)+1}$ its ordered leaf indices. The induction
hypothesis gives
\[
 n_i=\sum_{j=1}^{2\ell(\gamma_i)+1}(-1)^{j-1}m_{i,j},
 \qquad i=1,2,3.
\]
Concatenate the three leaf lists, in order, and relabel them by
\begin{align}
 m_{1,j}&=m_j,
 &&1\le j\le 2\ell(\gamma_1)+1,\nonumber\\
 m_{2,j}&=m_{j+2\ell(\gamma_1)+1},
 &&1\le j\le 2\ell(\gamma_2)+1,\nonumber\\
 m_{3,j}&=m_{j+2\ell(\gamma_1)+2\ell(\gamma_2)+2},
 &&1\le j\le 2\ell(\gamma_3)+1.\label{eq:critical-tree-leaf-relabel}
\end{align}
Each child list has odd
length by \autoref{lem:critical-tree-leaf-count}, so the alternating
sign changes at each boundary. The induction hypothesis and momentum
conservation at the new root give
\begin{align*}
 n_{123}
 &=n_1-n_2+n_3\\
 &=\sum_{i=1}^3\sum_{j=1}^{2\ell(\gamma_i)+1}
   (-1)^{i+j-2}m_{i,j}\\
 &=\sum_{j=1}^{2\ell(\gamma_{123})+1}(-1)^{j-1}m_j.
 \end{align*}
This completes the induction.
\end{proof}

\begin{lemma}[Pointwise Bound for One Tree]
\label{lem:critical-single-tree-bound}
Suppose $|d_q|\le A\delta^{-1/2}w_{-\beta}(q)$ for every
$q\in\Z^\nu$, with $\beta_j>2$ for each coordinate. Let $\gamma$ be
a tree. Then, for every fixed $n\in\Z^\nu$ and $t\ge0$,
\begin{align}
 &|\operatorname{Term}_\gamma(t,n)|\nonumber\\
  &\qquad\le A\delta^{-1/2}\|w_{-\beta/2}\|_{\ell^1(\Z^\nu)}
  w_{-\beta/2}(n)\bigl(\varepsilon^\alpha tA^2\delta^{-1}
  \|w_{-\beta/2}\|_{\ell^1(\Z^\nu)}^2\bigr)^{\ell(\gamma)}
  .\label{eq:critical-pointwise-tree-bound}
\end{align}
\end{lemma}
\begin{proof}[Proof of \autoref{lem:critical-single-tree-bound}]
Fix a tree $\gamma$ with the root index $n$. We prove by induction
that
\begin{equation}\label{eq:critical-tree-leaf-sum}
 |\operatorname{Term}_\gamma(t,n)|
 \le(\varepsilon^\alpha t)^{\ell(\gamma)}
 \sum_{\substack{m_1,\ldots,m_{2\ell(\gamma)+1}\in\Z^\nu\\
                   \sum_{j=1}^{2\ell(\gamma)+1}(-1)^{j-1}m_j=n}}
 \prod_{j=1}^{2\ell(\gamma)+1}|d_{m_j}|.
\end{equation}
For the two base cases, $\ell(0)=0$ and $\ell(1)=1$.
The only leaf of $0$ has index $n$, and
$\operatorname{Term}_0(t,n)=d_n$. For $1$, its defining integral
and $|\mathcal K_{r,n}|\le1$ yield
\[
 |\operatorname{Term}_1(t,n)|
 \le\varepsilon^\alpha t
 \sum_{\substack{m_1,m_2,m_3\in\Z^\nu\\m_1-m_2+m_3=n}}
 \prod_{j=1}^3|d_{m_j}|.
\]
Let
$\gamma=(\gamma_1,\gamma_2,\gamma_3)$ and assume the estimate \eqref{eq:critical-tree-leaf-sum} holds for
each child. Write its ordered leaf indices as
$(m_{i,j})_{1\le j\le 2\ell(\gamma_i)+1}$, so the induction
hypothesis applies with
$n_i=\sum_{j=1}^{2\ell(\gamma_i)+1}(-1)^{j-1}m_{i,j}$.
Substituting the three bounds into the defining cubic integral and
using $|\mathcal K_{r,n}|\le1$ give
\begin{align*}
 |\operatorname{Term}_\gamma(t,n)|
 &\le\varepsilon^\alpha\int_0^t
  \sum_{n_1-n_2+n_3=n}
  \prod_{i=1}^3|\operatorname{Term}_{\gamma_i}(r,n_i)|\,dr\\
&\le\varepsilon^\alpha\int_0^t
  (\varepsilon^\alpha r)^{\sum_{i=1}^3\ell(\gamma_i)}
  \sum_{n_1-n_2+n_3=n}\\
&\qquad\times
  \prod_{i=1}^3\left(
    \sum_{\substack{m_{i,1},\ldots,m_{i,2\ell(\gamma_i)+1}\in\Z^\nu\\
      \sum_{j=1}^{2\ell(\gamma_i)+1}(-1)^{j-1}m_{i,j}=n_i}}
    \prod_{j=1}^{2\ell(\gamma_i)+1}|d_{m_{i,j}}|
  \right)dr.
\end{align*}
Applying the relabeling \eqref{eq:critical-tree-leaf-relabel} gives
\[
 n_1-n_2+n_3
 =\sum_{i=1}^3\sum_{j=1}^{2\ell(\gamma_i)+1}
   (-1)^{i+j-2}m_{i,j}
 =\sum_{j=1}^{2\ell(\gamma)+1}(-1)^{j-1}m_j
\]
and then
\begin{align*}
 |\operatorname{Term}_\gamma(t,n)|
 &\le\varepsilon^\alpha\int_0^t
   (\varepsilon^\alpha r)^{\ell(\gamma)-1}\,dr
   \sum_{\substack{m_1,\ldots,m_{2\ell(\gamma)+1}\in\Z^\nu\\
                   \sum_{j=1}^{2\ell(\gamma)+1}(-1)^{j-1}m_j=n}}
   \prod_{j=1}^{2\ell(\gamma)+1}|d_{m_j}|\\
 &\le(\varepsilon^\alpha t)^{\ell(\gamma)}
   \sum_{\substack{m_1,\ldots,m_{2\ell(\gamma)+1}\in\Z^\nu\\
                   \sum_{j=1}^{2\ell(\gamma)+1}(-1)^{j-1}m_j=n}}
   \prod_{j=1}^{2\ell(\gamma)+1}|d_{m_j}|.
 \end{align*}
This completes the proof of \eqref{eq:critical-tree-leaf-sum}.

To retain decay at the root, split the initial weight as
$w_{-\beta}=w_{-\beta/2}\,w_{-\beta/2}$.
Use one half of the exponent $\beta$ to transfer the product of
leaf weights to the root through momentum conservation; the other
half remains on the leaves and is summed using
$\beta_h/2>1$ in every coordinate. To justify the transfer to the
root explicitly, \autoref{lem:critical-tree-momentum} gives
\[
 n_h=\sum_{j=1}^{2\ell(\gamma)+1}(-1)^{j-1}(m_j)_h,
 \qquad h=1,\ldots,\nu.
\]
The triangle inequality therefore implies
\[
 1+|n_h|
 \le 1+\sum_{j=1}^{2\ell(\gamma)+1}|(m_j)_h|
 \le\prod_{j=1}^{2\ell(\gamma)+1}(1+|(m_j)_h|).
\]
For the second inequality, expand the product: its constant term
is one, its terms of degree one sum to
$\sum_j|(m_j)_h|$, and every higher-degree term is nonnegative.
Since $\beta_h/2>0$, the function
$s\mapsto s^{-\beta_h/2}$ is decreasing on $(0,\infty)$.
Taking this negative power reverses the inequality and yields
\[
 \prod_{j=1}^{2\ell(\gamma)+1}
 (1+|(m_j)_h|)^{-\beta_h/2}
 \le(1+|n_h|)^{-\beta_h/2}.
\]
Multiply over all coordinates $h$ and regroup the finite products.
By the definition of the product weight, this gives
\begin{align*}
 \prod_{j=1}^{2\ell(\gamma)+1}w_{-\beta/2}(m_j)
 &=\prod_{h=1}^{\nu}\prod_{j=1}^{2\ell(\gamma)+1}
   (1+|(m_j)_h|)^{-\beta_h/2}\\
 &\le\prod_{h=1}^{\nu}(1+|n_h|)^{-\beta_h/2}
 =w_{-\beta/2}(n).
 \end{align*}
Thus, for every leaf configuration satisfying the root momentum
constraint,
\begin{equation}\label{eq:critical-product-decay}
 \prod_{j=1}^{2\ell(\gamma)+1}w_{-\beta/2}(m_j)
 \le w_{-\beta/2}(n).
\end{equation}
Since \autoref{lem:critical-tree-leaf-count} gives
$2\ell(\gamma)+1$ leaves,
\eqref{eq:critical-product-decay} yields the pointwise product bound
\[
 \prod_{j=1}^{2\ell(\gamma)+1}|d_{m_j}|
 \le (A\delta^{-1/2})^{2\ell(\gamma)+1}w_{-\beta/2}(n)
       \prod_{j=1}^{2\ell(\gamma)+1}w_{-\beta/2}(m_j).
\]
For this fixed root $n$, drop the root constraint only after
extracting $w_{-\beta/2}(n)$. The sum of the remaining leaf weights is
$\|w_{-\beta/2}\|_{\ell^1(\Z^\nu)}^{2\ell(\gamma)+1}$,
which is finite because $\beta_j/2>1$ for every $j$.
Combining this with the time integral bound proves
\eqref{eq:critical-pointwise-tree-bound}.
\end{proof}
We now sum the bound over the complete diagrams in $\Gamma^{(k)}$.
The recursion for $\Gamma^{(k)}$ gives
$\ell(\gamma)\le(3^k-1)/2$. For each fixed $n$, the
factor outside the parentheses in
\eqref{eq:critical-pointwise-tree-bound} is independent of $\gamma$;
only the exponent $\ell(\gamma)$ depends on the tree. Grouping the
trees by their number of cubic vertices gives
\begin{align*}
 \sum_{\gamma\in\Gamma^{(k)}}|\operatorname{Term}_\gamma(t,n)|
 &\le A\delta^{-1/2}\|w_{-\beta/2}\|_{\ell^1(\Z^\nu)}
 w_{-\beta/2}(n)\\
&\quad\times\sum_{\gamma\in\Gamma^{(k)}}
 \bigl(\varepsilon^\alpha tA^2\delta^{-1}
 \|w_{-\beta/2}\|_{\ell^1(\Z^\nu)}^2\bigr)^{\ell(\gamma)}.
\end{align*}
Grouping this sum by vertex count gives
\begin{align*}
 &\sum_{\gamma\in\Gamma^{(k)}}
 \bigl(\varepsilon^\alpha tA^2\delta^{-1}
 \|w_{-\beta/2}\|_{\ell^1(\Z^\nu)}^2\bigr)^{\ell(\gamma)}\\
 =~& \sum_{m=0}^{(3^k-1)/2}
 \#\{\gamma\in\Gamma^{(k)}:\ell(\gamma)=m\}
 \bigl(\varepsilon^\alpha tA^2\delta^{-1}
 \|w_{-\beta/2}\|_{\ell^1(\Z^\nu)}^2\bigr)^m.
\end{align*}
The remaining task is to bound $\#\{\gamma\in\Gamma^{(k)}:\ell(\gamma)=m\}$.
\begin{lemma}
\label{lem:critical-ternary-tree-count}
At each fixed iteration level $k\ge1$, for any $0\le m\le \frac{3^k-1}{2}$, we have
\begin{equation}\label{eq:critical-fixed-level-tree-count}
 \#\{\gamma\in\Gamma^{(k)}:\ell(\gamma)=m\}
 \le\frac{1}{2m+1}\binom{3m}{m}\le8^m.
\end{equation}
\end{lemma}
\begin{proof}[Proof of \autoref{lem:critical-ternary-tree-count}]
\noindent\emph{Compression of Identity Paths.}
Fix $k$ and contract every nonbranching identity path represented by
$0$ to a single leaf. Retain the cubic vertices, the order of their
three children and the edge colors. The result is an ordered full
ternary tree with exactly $m$ cubic vertices. Conversely, the depth
of each leaf in this reduced tree determines how many identity links
must be restored to reach level $k$. Restoring those links with their
inherited colors recovers the complete diagram uniquely.

\begin{figure}[htbp]
\centering
\begin{tikzpicture}[x=.42cm,y=.42cm,line width=.55pt,line cap=round]
\foreach \a/\b/\c/\col/\row in {
  0/0/0/0/0,0/0/1/1/0,0/1/0/2/0,0/1/1/3/0,
  1/0/0/0/1,1/0/1/1/1,1/1/0/2/1,1/1/1/3/1} {
  \begin{scope}[shift={({\col*4.2},{-\row*4.5})}]
    \fill (0,0) circle (1.9pt);
    \foreach \x/\v in {-1.4/\a,0/\b,1.4/\c} {
      \ifdim\x pt=0pt
        \draw[red!75!black] (0,0)--(\x,1);
        \fill[red!75!black] (\x,1) circle (1.6pt);
        \ifnum\v=1
          \foreach \s in {-.38,.38} {
            \draw[red!75!black] (\x,1)--({\x+\s},2);
            \fill[red!75!black] ({\x+\s},2) circle (1.6pt);
          }
          \draw (\x,1)--(\x,2);
          \fill (\x,2) circle (1.6pt);
        \fi
      \else
        \draw (0,0)--(\x,1);
        \fill (\x,1) circle (1.6pt);
        \ifnum\v=1
          \foreach \s in {-.38,.38} {
            \draw (\x,1)--({\x+\s},2);
            \fill ({\x+\s},2) circle (1.6pt);
          }
          \draw[red!75!black] (\x,1)--(\x,2);
          \fill[red!75!black] (\x,2) circle (1.6pt);
        \fi
      \fi
    }
    \node[below=3pt] at (0,0) {$(\a,\b,\c)$};
  \end{scope}
}
\end{tikzpicture}
\caption{The compressed versions of the eight diagrams in
\autoref{ex:critical-second-level-trees}, in the same order.}
\label{fig:critical-compressed-second-level-trees}
\end{figure}
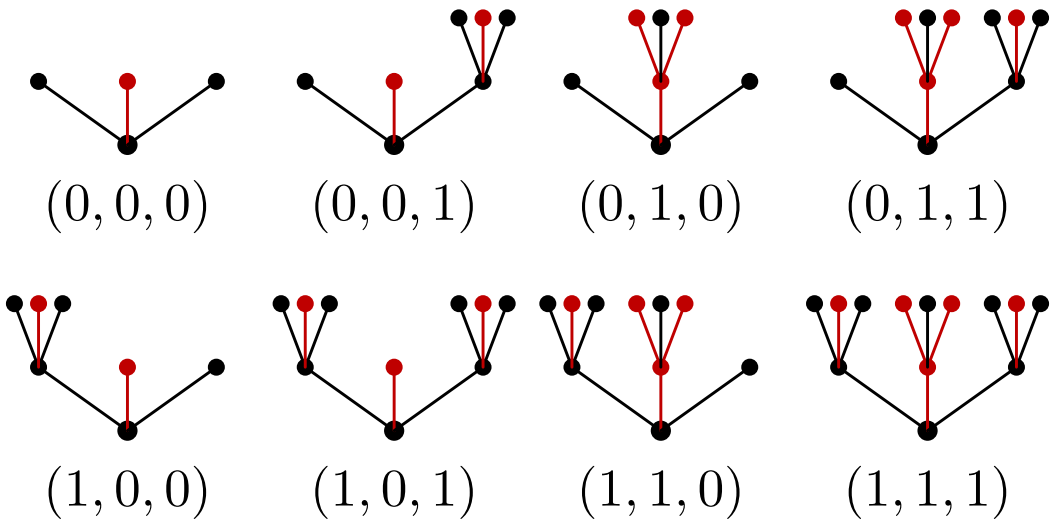

\medskip
\noindent\emph{Counting the Compressed Trees.}
 \emph{The Fuss--Catalan formula} gives
$\frac{1}{2m+1}\binom{3m}{m}$ unrestricted ordered full ternary
trees with $m$ cubic vertices; see  \cite{S99}. The compression injects the complete diagrams
in $\Gamma^{(k)}$ with $\ell(\gamma)=m$ into this family, proving the
first inequality of \eqref{eq:critical-fixed-level-tree-count}.
Finally,
$\binom{3m}{m}\le\sum_{j=0}^{3m}\binom{3m}{j}=2^{3m}$ and
$(2m+1)^{-1}\le1$, so the count is at most $8^m$.
\end{proof}
At each fixed level $k$, \eqref{eq:critical-fixed-level-tree-count}
bounds the number of diagrams with $m$ cubic vertices by $8^m$.
For $0\le t\le T_\delta$, the definition of $T_\delta$ gives
\[
 8\varepsilon^\alpha tA^2\delta^{-1}
 \|w_{-\beta/2}\|_{\ell^1(\Z^\nu)}^2\le\frac14.
\]
Grouping \eqref{eq:critical-pointwise-tree-bound} by $m$ therefore
gives, for every fixed $n$ and every $k$,
\begin{align*}
 \sum_{\gamma\in\Gamma^{(k)}}|\operatorname{Term}_\gamma(t,n)|
 &\le A\delta^{-1/2}\|w_{-\beta/2}\|_{\ell^1(\Z^\nu)}
 w_{-\beta/2}(n)\sum_{m=0}^{(3^k-1)/2}4^{-m}\\
 &\le\frac43A\delta^{-1/2}\|w_{-\beta/2}\|_{\ell^1(\Z^\nu)}
 w_{-\beta/2}(n).
 \end{align*}
By \autoref{lem:critical-picard-tree-expansion}, the \emph{Picard
iterates} satisfy, for each fixed $n$,
\begin{equation}\label{eq:critical-picard-uniform-decay}
 \sup_{0\le t\le T_\delta}|a_n^{(k)}(t)|
 \le\frac43A\delta^{-1/2}\|w_{-\beta/2}\|_{\ell^1(\Z^\nu)}
 w_{-\beta/2}(n),\qquad k\ge1.
\end{equation}

\begin{remark}[Comparison with Earlier Tree Expansions]
\label{rem:critical-tree-comparison}
The \emph{tree expansion} follows \cite{DG16,DLX24I}.
For the pointwise bounds needed here, the simpler estimate
\eqref{eq:critical-tree-leaf-sum} and the \emph{Fuss--Catalan count}
\eqref{eq:critical-fixed-level-tree-count} reduce the tree sum to a
\emph{geometric series}, simplifying the combinatorial analysis without using
finer time integral estimates.
\end{remark}

\medskip
\noindent\textbf{Step 4. Construction of the Solution.}
Set
$
 C=\frac43A\delta^{-1/2}\|w_{-\beta/2}\|_{\ell^1(\Z^\nu)}$. 
By \eqref{eq:critical-picard-uniform-decay},
$|a_n^{(k)}(t)|\le Cw_{-\beta/2}(n)$ for every $k\ge0$.
For $k\ge1$, subtract \eqref{eq:critical-pointwise-picard} at
levels $k+1$ and $k$ to obtain the increment identity
\begin{equation}\label{eq:critical-picard-increment-identity}
 a_n^{(k+1)}(t)-a_n^{(k)}(t)
 =\iu\varepsilon^\alpha\int_0^t
 \bigl(H_{r,n}[a^{(k)}(r)]-H_{r,n}[a^{(k-1)}(r)]\bigr)\,dr.
\end{equation}
By \eqref{eq:critical-cubic-map}, the cubic difference expands as
\[
\begin{aligned}
&H_{r,n}[a^{(k)}(r)]-H_{r,n}[a^{(k-1)}(r)]\\
&\quad= \sum_{n_1-n_2+n_3=n}\mathcal K_{r,n}(n_1,n_2,n_3)
 \left[\begin{aligned}
 &
 (a_{n_1}^{(k)}-a_{n_1}^{(k-1)})
   \overline{a_{n_2}^{(k)}}a_{n_3}^{(k)}+{}\\
 &a_{n_1}^{(k-1)}
   \overline{(a_{n_2}^{(k)}-a_{n_2}^{(k-1)})}a_{n_3}^{(k)}+{}\\
 &a_{n_1}^{(k-1)}\overline{a_{n_2}^{(k-1)}}
   (a_{n_3}^{(k)}-a_{n_3}^{(k-1)})
 \end{aligned}\right].
\end{aligned}
\]
Define 
$F_k(t):=\sup_{n\in\Z^\nu}|a_n^{(k)}(t)-a_n^{(k-1)}(t)|$.
Bound the difference factor in
each of the three summands by $F_k(r)$ and the other two factors by $Cw_{-\beta/2}$ at their respective indices. Under
$n_1-n_2+n_3=n$, those two indices range freely, so their weights
sum to $\|w_{-\beta/2}\|_{\ell^1(\Z^\nu)}^2$. Hence
\begin{equation}\label{eq:critical-cubic-increment-bound}
 \bigl|H_{r,n}[a^{(k)}(r)]-H_{r,n}[a^{(k-1)}(r)]\bigr|
 \le3C^2\|w_{-\beta/2}\|_{\ell^1(\Z^\nu)}^2F_k(r).
\end{equation}
Inserting this estimate into the time integral and taking the supremum over $n$  give
\[
 F_{k+1}(t)
 \le3\varepsilon^\alpha C^2
 \|w_{-\beta/2}\|_{\ell^1(\Z^\nu)}^2
 \int_0^t F_k(r)\,dr.
\]
For the first increment, use $w_{-\beta/2}(n)\le1$ and the
uniform bound for both $a^{(1)}$ and $a^{(0)}$ to obtain
\[
 F_1(t)\le\sup_n\bigl(|a_n^{(1)}(t)|+|a_n^{(0)}(t)|\bigr)
 \le2C,\qquad 0\le t\le T_\delta.
\]
Substituting this bound into the integral inequality with $k=1$,
and then substituting the resulting bound with $k=2$, give
\begin{align*}
 F_2(t)
 &\le3\varepsilon^\alpha C^2
 \|w_{-\beta/2}\|_{\ell^1(\Z^\nu)}^2\int_0^t2C\,dr
 =2C\bigl(3\varepsilon^\alpha C^2
 \|w_{-\beta/2}\|_{\ell^1(\Z^\nu)}^2\bigr)t,\\
 F_3(t)
 &\le3\varepsilon^\alpha C^2
 \|w_{-\beta/2}\|_{\ell^1(\Z^\nu)}^2
 \int_0^t2C\bigl(3\varepsilon^\alpha C^2
 \|w_{-\beta/2}\|_{\ell^1(\Z^\nu)}^2\bigr)r\,dr\\
 &=2C\bigl(3\varepsilon^\alpha C^2
 \|w_{-\beta/2}\|_{\ell^1(\Z^\nu)}^2\bigr)^2\frac{t^2}{2!}.
\end{align*}
In general, suppose that
\[
 F_k(r)\le2C\bigl(3\varepsilon^\alpha C^2
 \|w_{-\beta/2}\|_{\ell^1(\Z^\nu)}^2\bigr)^{k-1}
 \frac{r^{k-1}}{(k-1)!}.
\]
The same integral inequality then yields
\begin{align*}
 F_{k+1}(t)
 &\le2C\bigl(3\varepsilon^\alpha C^2
 \|w_{-\beta/2}\|_{\ell^1(\Z^\nu)}^2\bigr)^k
 \int_0^t\frac{r^{k-1}}{(k-1)!}\,dr\\
 &=2C\bigl(3\varepsilon^\alpha C^2
 \|w_{-\beta/2}\|_{\ell^1(\Z^\nu)}^2\bigr)^k
 \frac{t^k}{k!}.
\end{align*}
Thus each additional iteration integrates one more time, producing
the factor $1/k!$. Induction starting from $F_1\le2C$ proves
\begin{equation}\label{eq:critical-picard-factorial-bound}
 F_{k+1}(t)
 \le 2C\frac{\bigl(3\varepsilon^\alpha C^2
 \|w_{-\beta/2}\|_{\ell^1(\Z^\nu)}^2t\bigr)^k}{k!},
 \qquad 0\le t\le T_\delta.
\end{equation}
We now substitute \eqref{eq:critical-picard-factorial-bound}, with
$k$ replaced by $k-1$ and $t$ by $r$, into
\eqref{eq:critical-cubic-increment-bound}. For $k\ge1$, this gives
\begin{align*}
 &\bigl|H_{r,n}[a^{(k)}(r)]-H_{r,n}[a^{(k-1)}(r)]\bigr|\\
 &\qquad\le 
 6C^3\|w_{-\beta/2}\|_{\ell^1(\Z^\nu)}^2
 \bigl(3\varepsilon^\alpha C^2
 \|w_{-\beta/2}\|_{\ell^1(\Z^\nu)}^2\bigr)^{k-1}
 \frac{r^{k-1}}{(k-1)!}.
\end{align*}
Inserting this nonlinear difference bound into
\eqref{eq:critical-picard-increment-identity} and taking absolute
values yield the desired estimate for each Fourier coefficient:
\begin{align*}
 |a_n^{(k+1)}(t)-a_n^{(k)}(t)|
 &\le\varepsilon^\alpha\int_0^t
 \bigl|H_{r,n}[a^{(k)}(r)]-H_{r,n}[a^{(k-1)}(r)]\bigr|\,dr\\
 &\le2C\bigl(3\varepsilon^\alpha C^2
 \|w_{-\beta/2}\|_{\ell^1(\Z^\nu)}^2\bigr)^k
 \int_0^t\frac{r^{k-1}}{(k-1)!}\,dr\\
 &=2C\frac{\bigl(3\varepsilon^\alpha C^2
 \|w_{-\beta/2}\|_{\ell^1(\Z^\nu)}^2t\bigr)^k}{k!}.
\end{align*}
To compare any two iterates, let $k\ge1$ and $m\ge1$ be integers.
Telescoping their difference and applying the preceding estimate give
\begin{align*}
 |a_n^{(k+m)}(t)-a_n^{(k)}(t)|
 =~&\left|\sum_{j=k}^{k+m-1}
     \bigl(a_n^{(j+1)}(t)-a_n^{(j)}(t)\bigr)\right|\\
 \le~&\sum_{j=k}^{k+m-1}
     |a_n^{(j+1)}(t)-a_n^{(j)}(t)|\\
 \le~&2C\sum_{j=k}^{k+m-1}
 \frac{\bigl(3\varepsilon^\alpha C^2
 \|w_{-\beta/2}\|_{\ell^1(\Z^\nu)}^2t\bigr)^j}{j!}\\
 \le~&2C\sum_{j=k}^{\infty}
 \frac{\bigl(3\varepsilon^\alpha C^2
 \|w_{-\beta/2}\|_{\ell^1(\Z^\nu)}^2t\bigr)^j}{j!}.
\end{align*}
For fixed $\varepsilon$ and $\delta$, the full series on the last
line is an exponential series with a finite argument. Its tail
therefore tends to zero as $k\to\infty$. This bound is
independent of $m$, $n$ and $t\in[0,T_\delta]$. Consequently,
\[
 \lim_{k\to\infty}\sup_{m\ge1}\sup_{n\in\Z^\nu}
 \sup_{0\le t\le T_\delta}
 |a_n^{(k+m)}(t)-a_n^{(k)}(t)|=0.
\]
Hence $a_n^{(k)}$ converges uniformly to a continuous
function $a_n$ with $a_n(0)=d_n$ and the bound
\eqref{eq:critical-pointwise-decay}.
Expanding the difference of the two cubic products as above and
using $|\mathcal K_{r,n}|\le1$ give the quantitative estimate
\begin{align*}
&\sup_{0\le r\le T_\delta}
 |H_{r,n}[a^{(k)}(r)]-H_{r,n}[a(r)]|\\
&\qquad\le 3C^2\|w_{-\beta/2}\|_{\ell^1(\Z^\nu)}^2
 \sup_{n\in\Z^\nu}\sup_{0\le r\le T_\delta}
 |a_n^{(k)}(r)-a_n(r)|\longrightarrow0.
\end{align*}
The displayed estimate shows that $H_{r,n}[a^{(k)}(r)]$ converges
uniformly to $H_{r,n}[a(r)]$. 
Passing to the limit in \eqref{eq:critical-pointwise-picard} yields
\[
 a_n(t)=d_n+\iu\varepsilon^\alpha
 \int_0^t H_{r,n}[a(r)]\,dr.
\]
Uniform convergence of $H_{r,n}[a^{(k)}(r)]$ to $H_{r,n}[a(r)]$ and continuity of $H_{r,n}[a^{(k)}(r)]$ imply that $H_{r,n}[a(r)]$ is continuous. The \emph{fundamental theorem of calculus} yields
$a_n\in C^1([0,T_\delta])$, and $a_n$ solves
\eqref{eq:nf-gauged-ode} at each $n$. 

\medskip
\noindent\textbf{Step 5. Uniqueness.}
Let $a$ and $\widetilde a$ be two coefficientwise $C^1$ solutions
with the same initial data and pointwise bounds
$|a_n(t)|\le Cw_{-\beta/2}(n)$ and
$|\widetilde a_n(t)|\le\widetilde Cw_{-\beta/2}(n)$.
The initial data agree, so subtraction of the two integral equations
gives, for each $n$,
\begin{equation}\label{eq:critical-uniqueness-integral-difference}
 a_n(t)-\widetilde a_n(t)
 =\iu\varepsilon^\alpha\int_0^t
 \bigl(H_{r,n}[a(r)]-H_{r,n}[\widetilde a(r)]\bigr)\,dr.
\end{equation}
Expanding the difference of the cubic products as in Step 4 gives
\begin{equation}\label{eq:critical-uniqueness-cubic-difference-bound}
 \begin{aligned}[t]
 |H_{t,n}[a(t)]-H_{t,n}[\widetilde a(t)]|
 \le F(t)(C^2+C\widetilde C+\widetilde C^2)
 \|w_{-\beta/2}\|_{\ell^1(\Z^\nu)}^2,
 \end{aligned}
\end{equation}
where $F(t):=\sup_n|a_n(t)-\widetilde a_n(t)|$.
Substituting \eqref{eq:critical-uniqueness-cubic-difference-bound}
into \eqref{eq:critical-uniqueness-integral-difference}, and then taking the supremum over output indices
yield
\[
 F(t)\le\varepsilon^\alpha(C^2+C\widetilde C+\widetilde C^2)
 \|w_{-\beta/2}\|_{\ell^1(\Z^\nu)}^2
 \int_0^t F(r)\,dr.
\]
Gronwall's inequality gives $F(t)=0$ throughout
$[0,T_\delta]$, proving uniqueness.
\end{proof}

\subsection{Comparison with the Cubic Resonant Flow}
\label{sec:resonant-comparison}
\begin{proposition}[Uniform Resonant Approximation on the Critical Interval]
\label{prop:resonant-comparison}
Assume that $c$ satisfies \eqref{eq:deterministic-decay}, with
$0<\kappa_j<\rho_j$ and $\beta_j:=\rho_j-\kappa_j>4\tau+2$
for every $j=1,\ldots,\nu$, and that
$\omega\in\mathrm{DC}(\gamma,\tau)$.
On $G_\delta$, let $u$ and $u_{\mathrm{res}}$ be defined by
\eqref{eq:nf-solution-reconstruction} and
\eqref{eq:nf-resonant-reconstruction}, respectively, with the same
initial coefficients $d_n=c(n)g_n$, where the nonlinear solution
is supplied by \autoref{prop:critical-pointwise-lifespan}. Set
\[
 K_\delta:=2A\delta^{-1/2}
 \|w_{-\beta/2}\|_{\ell^1(\Z^\nu)}
 \|\langle\cdot\rangle^{2\tau}w_{-\beta/2}\|_{\ell^1(\Z^\nu)}.
\]
Then
\begin{align}
 \sup_{0\le t\le T_\delta}
 \|u(t)-u_{\mathrm{res}}(t)\|_\infty
 \le C_\gamma\varepsilon^\alpha K_\delta^3
       (1+\varepsilon^\alpha T_\delta K_\delta^2)
       e^{3\varepsilon^\alpha K_\delta^2T_\delta}.
 \label{eq:resonant-comparison-general}
\end{align}
In particular, for $\delta=\delta_\varepsilon$ chosen in
\eqref{eq:nf-high-mode-delta}, the error is negligible relative to
the observation scale $\varepsilon^{-1/2}$, in the sense that
\begin{equation}\label{eq:resonant-comparison-critical}
 \sup_{0\le t\le T_{\delta_\varepsilon}}
 \|u(t)-u_{\mathrm{res}}(t)\|_\infty
 \le r_\varepsilon\varepsilon^{-1/2},
\end{equation}
where $r_\varepsilon:=C_{\gamma,A,\beta,\tau,z,S}
 \varepsilon^{\alpha-1}$.
\end{proposition}
\begin{proof}[Proof of \autoref{prop:resonant-comparison}]
Recalling \eqref{eq:nf-solution-reconstruction} and
\eqref{eq:nf-resonant-reconstruction} for $u$ and $u_{\mathrm{res}}$,
respectively, we have
\begin{equation}\label{eq:resonant-comparison-reconstruction-bound}
 \|u(t)-u_{\mathrm{res}}(t)\|_\infty
 \le\sum_n|a_n(t)-b_n(t)|\triangleq E(t).
\end{equation}
Since $a_n(0)=b_n(0)=d_n$, integrating
\eqref{eq:nf-gauged-ode} and \eqref{eq:nf-resonant-ode} and subtracting give
\begin{align}
 a_n(t)-b_n(t)
 &=-\iu\varepsilon^\alpha\int_0^t
   \bigl(|a_n(r)|^2a_n(r)-|b_n(r)|^2b_n(r)\bigr)\,dr\nonumber\\
   &\quad+\iu\varepsilon^\alpha\int_0^t
   \mathcal N_1^3[a](r,n)\,dr.\label{eq:resonant-comparison-coefficient-difference}
\end{align}
Substituting \eqref{eq:resonant-comparison-coefficient-difference} into
\eqref{eq:resonant-comparison-reconstruction-bound} and applying the
triangle inequality yield the two contributions to be estimated:
\begin{align}
E(t)
 &\le\varepsilon^\alpha\int_0^t\sum_n
 \bigl||a_n(r)|^2a_n(r)-|b_n(r)|^2b_n(r)\bigr|\,dr
 \label{eq:resonant-comparison-two-contributions}\\
 &\qquad+\varepsilon^\alpha\sum_n\left|\int_0^t
 \mathcal N_1^3[a](r,n)\,dr\right|.
 \label{eq:resonant-comparison-nonresonant-contribution}
\end{align}

\medskip
\noindent\textbf{Step 1. Preliminary Bounds.}
By \eqref{eq:critical-pointwise-decay} and
\autoref{cor:randomized-decay}, respectively,
\begin{align}
 |a_n(r)|&\le 2A\delta^{-1/2}
 \|w_{-\beta/2}\|_{\ell^1(\Z^\nu)}w_{-\beta/2}(n),
 \label{eq:resonant-comparison-pointwise-a}\\
 |b_n(r)|&=|d_n|\le A\delta^{-1/2}w_{-\beta}(n).
 \label{eq:resonant-comparison-pointwise-profiles}
\end{align}
Using $\langle n\rangle\le\prod_{j=1}^\nu\langle n_j\rangle$
and the weight definition \eqref{eq:product-polynomial-weight}, we obtain
\begin{align}
 \sum_{n\in\Z^\nu}\langle n\rangle^{2\tau}w_{-\beta/2}(n)
 &\le\sum_{n\in\Z^\nu}
   \prod_{j=1}^\nu\langle n_j\rangle^{-(\beta_j/2-2\tau)}\nonumber\\
 &=\prod_{j=1}^\nu
   \sum_{m\in\Z}\langle m\rangle^{-(\beta_j/2-2\tau)}<\infty.\label{eq:resonant-comparison-weight-summability}
\end{align}
The final bound holds since $\beta_j/2-2\tau>1$ for every $j$.
Multiplying \eqref{eq:resonant-comparison-pointwise-a} and
\eqref{eq:resonant-comparison-pointwise-profiles} by
$\langle n\rangle^{2\tau}$ and summing, with
\eqref{eq:resonant-comparison-weight-summability}, give for
$0\le r\le T_\delta$,
\begin{align}
 \sum_n\langle n\rangle^{2\tau}|a_n(r)|
 &\le 2A\delta^{-1/2}\|w_{-\beta/2}\|_{\ell^1(\Z^\nu)}
       \|\langle\cdot\rangle^{2\tau}w_{-\beta/2}\|_{\ell^1(\Z^\nu)}
  =K_\delta,
  \label{eq:resonant-weighted-a}\\
 \sum_n\langle n\rangle^{2\tau}|b_n(r)|
 &\le A\delta^{-1/2}
       \|\langle\cdot\rangle^{2\tau}w_{-\beta}\|_{\ell^1(\Z^\nu)}
  \le K_\delta.
  \label{eq:resonant-weighted-b}
\end{align}
For the derivative estimate, we use the evolution equation
\eqref{eq:critical-gauged-cubic} with the cubic map
\eqref{eq:critical-cubic-map}, $|\mathcal K_{r,n}|\le1$,
the momentum conservation $n=n_1-n_2+n_3$ and the inequality
$\langle n_1-n_2+n_3\rangle\le\prod_{j=1}^3\langle n_j\rangle$.
Applying \eqref{eq:resonant-weighted-a} at the last step, we have
\begin{align}
 \sum_n\langle n\rangle^{2\tau}|\dot a_n(r)|
 =~& \varepsilon^\alpha\sum_n\langle n\rangle^{2\tau}
   |H_{r,n}[a(r)]|\nonumber\\
 =~& \varepsilon^\alpha\sum_n\langle n\rangle^{2\tau}
   \left|\sum_{n_1-n_2+n_3=n}\mathcal K_{r,n}(n_1,n_2,n_3)
   \prod_{j=1}^3\starpow{a_{n_j}(r)}{j-1}\right|\nonumber\\
 \le~& \varepsilon^\alpha\sum_n\sum_{n_1-n_2+n_3=n}
   \langle n\rangle^{2\tau}\prod_{j=1}^3|a_{n_j}(r)|\nonumber\\
 =~& \varepsilon^\alpha
 \sum_{n_1,n_2,n_3}
 \langle n_1-n_2+n_3\rangle^{2\tau}
  \prod_{j=1}^3|a_{n_j}(r)|\nonumber\\
 \le~& \varepsilon^\alpha
 \sum_{n_1,n_2,n_3}
  \prod_{j=1}^3\langle n_j\rangle^{2\tau}|a_{n_j}(r)|\nonumber\\
 =~& \varepsilon^\alpha\prod_{j=1}^3
  \left(\sum_{n_j\in\Z^\nu}
   \langle n_j\rangle^{2\tau}|a_{n_j}(r)|\right)\nonumber\\
 \le~& \varepsilon^\alpha K_\delta^3.\label{eq:resonant-derivative-bound}
\end{align}

\medskip
\noindent\textbf{Step 2. Extracting an Oscillatory Gain from the Nonresonant Term \eqref{eq:resonant-comparison-nonresonant-contribution}.}

\smallskip
\noindent\emph{Integration by Parts.}
For a triple in \eqref{eq:resonant-comparison-nonresonant-contribution}
with $n_1-n_2+n_3=n$ and $\Sigma_3\ne0$, the factor
$e^{-\iu r\Sigma_3}$ is a nonconstant time oscillation. Since
\[
 e^{-\iu r\Sigma_3}
 =\frac{1}{-\iu\Sigma_3}\frac{d}{dr}e^{-\iu r\Sigma_3},
\]
we \emph{integrate by parts}, moving the time derivative from the oscillatory
factor to the cubic product. In
\eqref{eq:resonant-comparison-nonresonant-contribution}, this gives
\begin{align*}
&\varepsilon^\alpha\sum_n \left|\int_0^t\mathcal N_1^3[a](r,n)\,dr\right|\\
&\quad=\varepsilon^\alpha\sum_n\left|
 \sum_{\substack{n_1-n_2+n_3=n\\\Sigma_3\ne0}}
 \int_0^t e^{-\iu r\Sigma_3}
 a_{n_1}(r)\overline{a_{n_2}(r)}a_{n_3}(r)\,dr\right|\\
&\quad\le\varepsilon^\alpha\sum_n
 \sum_{\substack{n_1-n_2+n_3=n\\\Sigma_3\ne0}}
 \frac1{|\Sigma_3|}
 \left|\Bigl[e^{-\iu r\Sigma_3}
 a_{n_1}(r)\overline{a_{n_2}(r)}a_{n_3}(r)\Bigr]_{r=0}^{r=t}\right|\\
&\qquad+\varepsilon^\alpha\sum_n
 \sum_{\substack{n_1-n_2+n_3=n\\\Sigma_3\ne0}}
 \frac1{|\Sigma_3|}\int_0^t
 \left|\frac{d}{dr}\bigl(a_{n_1}\overline{a_{n_2}}a_{n_3}\bigr)(r)\right|\,dr.
\end{align*}

\smallskip
\noindent\emph{Small Divisors and the Diophantine Condition.}
Since $\Sigma_3\ne0$, the factorization
\eqref{eq:Sigma-factorization} implies that both
$\langle n_1-n_2,\omega\rangle$ and
$\langle n_2-n_3,\omega\rangle$ are nonzero. In the
quasi-periodic setting, these scalar products need not be uniformly
separated from zero as the integer differences vary. This is the
\emph{small divisor problem} absent from the 1D periodic setting
of \cite{GGKS23}. Using the Diophantine condition \eqref{eq:HDC}, $|n_i-n_j|_1\le C\langle n_i\rangle
\langle n_j\rangle$ and $\langle n_j\rangle\ge1$, we obtain
\begin{align}
 \frac1{|\Sigma_3|}
 &=\frac1{2|\langle n_1-n_2,\omega\rangle|
             |\langle n_2-n_3,\omega\rangle|}\nonumber\\
 &\le\frac{|n_1-n_2|_1^\tau|n_2-n_3|_1^\tau}{2\gamma^2}\nonumber\\
 &\le C_\gamma\langle n_1\rangle^\tau
   \langle n_2\rangle^{2\tau}\langle n_3\rangle^\tau\nonumber\\
 &\le C_\gamma\prod_{j=1}^3\langle n_j\rangle^{2\tau}.\label{eq:resonant-reciprocal-phase}
\end{align}

\smallskip
\noindent\emph{Boundary Terms.}
At either endpoint $r\in\{0,t\}$, \eqref{eq:resonant-reciprocal-phase}
and \eqref{eq:resonant-weighted-a} give
\begin{align*}
 &\sum_n\sum_{\substack{n_1-n_2+n_3=n\\\Sigma_3\ne0}}
 \frac{|a_{n_1}(r)a_{n_2}(r)a_{n_3}(r)|}{|\Sigma_3|}\\
 &\quad\le  C_\gamma
 \sum_n\sum_{\substack{n_1-n_2+n_3=n\\\Sigma_3\ne0}}
 \prod_{j=1}^3\langle n_j\rangle^{2\tau}|a_{n_j}(r)|\\
 &\quad\le  C_\gamma\sum_{n_1,n_2,n_3}
 \prod_{j=1}^3\langle n_j\rangle^{2\tau}|a_{n_j}(r)|\\
 &\quad= C_\gamma
 \left(\sum_m\langle m\rangle^{2\tau}|a_m(r)|\right)^3
 \le C_\gamma K_\delta^3.
\end{align*}
The first inequality uses \eqref{eq:resonant-reciprocal-phase}.

\smallskip
\noindent\emph{Derivative Terms.}
Differentiating the product gives three terms, one for each
differentiated factor. For example, the term with $\dot a_{n_1}$
is bounded by
\begin{align*}
 &\sum_n\sum_{\substack{n_1-n_2+n_3=n\\\Sigma_3\ne0}}
 \frac{|\dot a_{n_1}(r)|\,|a_{n_2}(r)|\,|a_{n_3}(r)|}
 {|\Sigma_3|}\\
&\quad\le C_\gamma\sum_n
 \sum_{\substack{n_1-n_2+n_3=n\\\Sigma_3\ne0}}
 \langle n_1\rangle^{2\tau}|\dot a_{n_1}(r)|
 \prod_{j=2}^3\langle n_j\rangle^{2\tau}|a_{n_j}(r)|\\
&\quad\le  C_\gamma\sum_{n_1,n_2,n_3}
 \langle n_1\rangle^{2\tau}|\dot a_{n_1}(r)|
 \prod_{j=2}^3\langle n_j\rangle^{2\tau}|a_{n_j}(r)|\\
&\quad= C_\gamma
 \left(\sum_m\langle m\rangle^{2\tau}|\dot a_m(r)|\right)
 \left(\sum_m\langle m\rangle^{2\tau}|a_m(r)|\right)^2
 \\ 
&\quad\le  C_\gamma(\varepsilon^\alpha K_\delta^3)K_\delta^2
 =C_\gamma\varepsilon^\alpha K_\delta^5.
\end{align*}
The last inequality uses \eqref{eq:resonant-weighted-a} and
\eqref{eq:resonant-derivative-bound}.
The other two terms obey the same bound. 

\smallskip
\noindent\emph{Nonresonant Contribution.}
Hence the two boundary
terms contribute at most $C_\gamma\varepsilon^\alpha K_\delta^3$,
and the derivative terms contribute at most
$C_\gamma\varepsilon^{2\alpha}tK_\delta^5\le C_\gamma\varepsilon^{2\alpha}T_\delta K_\delta^5$.
Thus the
nonresonant contribution in
\eqref{eq:resonant-comparison-nonresonant-contribution} satisfies, for any $0\le t\le T_\delta$, 
\begin{align}
 \varepsilon^\alpha\sum_n\left|\int_0^t
 \mathcal N_1^3[a](r,n)\,dr\right|
\le C_\gamma\varepsilon^\alpha K_\delta^3
 (1+\varepsilon^\alpha T_\delta K_\delta^2).
 \label{eq:resonant-integral-bound}
\end{align}

\medskip
\noindent\textbf{Step 3. The Diagonal Cubic Term \eqref{eq:resonant-comparison-two-contributions}.}
Since $\langle n\rangle^{2\tau}\ge1$, the bounds
\eqref{eq:resonant-weighted-a} and \eqref{eq:resonant-weighted-b}
give $|a_n(r)|,|b_n(r)|\le K_\delta$ for every $n$ and
$0\le r\le T_\delta$. Suppressing the common time argument, we write
\[
 |a_n|^2a_n-|b_n|^2b_n
 =|a_n|^2(a_n-b_n)+(|a_n|^2-|b_n|^2)b_n.
\]
The reverse triangle inequality then gives
\begin{align*}
 \bigl||a_n|^2-|b_n|^2\bigr|
 &=(|a_n|+|b_n|)\bigl||a_n|-|b_n|\bigr|\\
 &\le(|a_n|+|b_n|)|a_n-b_n|.
\end{align*}
Thus, pointwise in $n$ and for $0\le r\le T_\delta$, we have
\begin{align}
 \bigl||a_n(r)|^2a_n(r)-|b_n(r)|^2b_n(r)\bigr|
 &\le 3\max\{|a_n(r)|,|b_n(r)|\}^2
 |a_n(r)-b_n(r)|\nonumber\\
 &\le 3K_\delta^2|a_n(r)-b_n(r)|.\label{eq:resonant-diagonal-pointwise-difference}
\end{align}
Substituting \eqref{eq:resonant-diagonal-pointwise-difference} into
the diagonal term in \eqref{eq:resonant-comparison-two-contributions} gives
\begin{equation}\label{eq:resonant-diagonal-integral-bound}
\begin{aligned}[t]
 \varepsilon^\alpha\int_0^t\sum_n
 \bigl||a_n(r)|^2a_n(r)-|b_n(r)|^2b_n(r)\bigr|\,dr
 \le3\varepsilon^\alpha K_\delta^2
 \int_0^tE(r)\,dr.
\end{aligned}
\end{equation}

\medskip
\noindent\textbf{Step 4. Gronwall Argument.}
Substituting \eqref{eq:resonant-diagonal-integral-bound} and
\eqref{eq:resonant-integral-bound} into 
\eqref{eq:resonant-comparison-two-contributions} and
\eqref{eq:resonant-comparison-nonresonant-contribution}, respectively,
gives the integral inequality for $E$:
\begin{equation}\label{eq:resonant-error-growth}
\begin{aligned}[t]
 E(t)
 \le C_\gamma\varepsilon^\alpha K_\delta^3
 (1+\varepsilon^\alpha T_\delta K_\delta^2)+3\varepsilon^\alpha K_\delta^2
 \int_0^t E(r)\,dr, 
 \qquad 0\le t\le T_\delta.
\end{aligned}
\end{equation}
The first term is constant in $t$, so Gronwall's inequality yields
\begin{align}
 E(t)
 &\le C_\gamma\varepsilon^\alpha K_\delta^3
 (1+\varepsilon^\alpha T_\delta K_\delta^2)
 e^{3\varepsilon^\alpha K_\delta^2t}\\
 &\le C_\gamma\varepsilon^\alpha K_\delta^3
 (1+\varepsilon^\alpha T_\delta K_\delta^2)
 e^{3\varepsilon^\alpha T_\delta K_\delta^2}.
\nonumber\end{align}
Together with \eqref{eq:resonant-comparison-reconstruction-bound},
this proves \eqref{eq:resonant-comparison-general}.
Finally, $\delta_\varepsilon=\varepsilon/(1+z^2/S)$ gives
$K_{\delta_\varepsilon}=\mathcal{O}(\varepsilon^{-1/2})$, while
\eqref{eq:critical-pointwise-lifespan} gives
$
 \varepsilon^\alpha T_{\delta_\varepsilon}
 K_{\delta_\varepsilon}^2
 =1/8
 \|\langle\cdot\rangle^{2\tau}w_{-\beta/2}\|_{\ell^1(\Z^\nu)}^2$.
This is independent of $\varepsilon$, so
\eqref{eq:resonant-comparison-general} yields
\eqref{eq:resonant-comparison-critical}.
\end{proof}

\subsection{LDP at the Critical Time}
\label{sec:first-critical-ldp-local}
We now transfer the resonant large deviations principle of
Theorem~\ref{thm:nf-resonant-LDP} to the nonlinear solution using
Proposition~\ref{prop:resonant-comparison}. 

\begin{theorem}[LDP for the Nonlinear Solution at the Critical Time]
\label{thm:first-critical-ldp-local}
Consider the weakly nonlinear Schr\"odinger equation \eqref{eq:qp-cubic-NLS} with random initial data
\eqref{eq:qp-randomized-initial-data}, under the hypotheses of
\autoref{prop:resonant-comparison}. 
Choose $\delta=\delta_\varepsilon$ as in
\eqref{eq:nf-high-mode-delta} for
\eqref{eq:critical-pointwise-lifespan} and set
\[
 T_\varepsilon:=T_{\delta_\varepsilon}
 =\frac{\varepsilon^{-(\alpha-1)}}
 {32A^2(1+z^2/S)\|w_{-\beta/2}\|_{\ell^1(\Z^\nu)}^2}.
\]
Then
\begin{equation}\label{eq:nonlinear-first-critical-ldp}
\begin{aligned}[t]
 \lim_{\varepsilon\to0^+}
 \sup_{0\le t\le T_\varepsilon}
 \varepsilon\log\Pp\!\left(
 \mathscr A_{\varepsilon,z}^{\|u(t)\|_\infty}\right)
 =\lim_{\varepsilon\to0^+}
 \inf_{0\le t\le T_\varepsilon}
 \varepsilon\log\Pp\!\left(
 \mathscr A_{\varepsilon,z}^{\|u(t)\|_\infty}\right)
 =-\frac{z^2}{S}.
\end{aligned}
\end{equation}
In particular, for every deterministic choice
$t_\varepsilon\in[0,T_\varepsilon]$,
\[
 \lim_{\varepsilon\to0^+}\varepsilon\log\Pp\!\left(
 \mathscr A_{\varepsilon,z}^{\|u(t_\varepsilon)\|_\infty}\right)
 =-\frac{z^2}{S}.
\]
\end{theorem}

\begin{proof}[Proof of \autoref{thm:first-critical-ldp-local}]
We use the event notation in
\eqref{eq:nonlinear-spatial-exceedance-event} and
\eqref{eq:spatial-exceedance-complement} throughout.
By \autoref{prop:resonant-comparison}, we have
\begin{equation}\label{eq:critical-comparison-event-inclusion}
\begin{aligned}[t]
 G_{\delta_\varepsilon}
 \subseteq
 \mathscr B_{\varepsilon,r_\varepsilon}^{\|u(t)-u_{\mathrm{res}}(t)\|_\infty},\qquad
 \mathscr A_{\varepsilon,r_\varepsilon}^{\|u(t)-u_{\mathrm{res}}(t)\|_\infty}
 \subseteq\Omega_{\delta_\varepsilon}.
\end{aligned}
\end{equation}
\medskip
\noindent\textbf{Step 1. Upper Bound.}
The triangle inequality gives the event inclusion
\begin{align}
 \mathscr A_{\varepsilon,z}^{\|u(t)\|_\infty}
 &\subseteq
 \mathscr A_{\varepsilon,z}^{\|u_{\mathrm{res}}(t)\|_\infty+\|u(t)-u_{\mathrm{res}}(t)\|_\infty},
 \label{eq:critical-triangle-event-inclusion}
\end{align}
where
\begin{align}
 \mathscr A_{\varepsilon,z}^{\|u_{\mathrm{res}}(t)\|_\infty+\|u(t)-u_{\mathrm{res}}(t)\|_\infty}
 \cap\mathscr B_{\varepsilon,r_\varepsilon}^{\|u(t)-u_{\mathrm{res}}(t)\|_\infty}
 &\subseteq
 \mathscr A_{\varepsilon,z-r_\varepsilon}^{\|u_{\mathrm{res}}(t)\|_\infty},
 \label{eq:critical-small-error-event-inclusion}\\[4pt]
 \mathscr A_{\varepsilon,z}^{\|u_{\mathrm{res}}(t)\|_\infty+\|u(t)-u_{\mathrm{res}}(t)\|_\infty}
 \cap\mathscr A_{\varepsilon,r_\varepsilon}^{\|u(t)-u_{\mathrm{res}}(t)\|_\infty}
 &\subseteq
 \mathscr A_{\varepsilon,r_\varepsilon}^{\|u(t)-u_{\mathrm{res}}(t)\|_\infty}.
 \label{eq:critical-large-error-event-inclusion}
\end{align}
For $0\le t\le T_\varepsilon$, we first apply
\eqref{eq:critical-triangle-event-inclusion} and then split the resulting
event into its intersections with $\mathscr B_{\varepsilon,r_\varepsilon}^{\|u(t)-u_{\mathrm{res}}(t)\|_\infty}$ and its complement $\mathscr A_{\varepsilon,r_\varepsilon}^{\|u(t)-u_{\mathrm{res}}(t)\|_\infty}$. By the \emph{total probability
formula}, followed by \eqref{eq:critical-small-error-event-inclusion},
\eqref{eq:critical-large-error-event-inclusion} and
\eqref{eq:critical-comparison-event-inclusion}, we obtain
\begin{align*}
 \Pp\!\left(\mathscr A_{\varepsilon,z}^{\|u(t)\|_\infty}\right)
 \le~&\Pp\!\left(
 \mathscr A_{\varepsilon,z}^{\|u_{\mathrm{res}}(t)\|_\infty+\|u(t)-u_{\mathrm{res}}(t)\|_\infty}
 \right)\\
 =~&\Pp\!\left(
 \mathscr A_{\varepsilon,z}^{\|u_{\mathrm{res}}(t)\|_\infty+\|u(t)-u_{\mathrm{res}}(t)\|_\infty}
 \cap\mathscr B_{\varepsilon,r_\varepsilon}^{\|u(t)-u_{\mathrm{res}}(t)\|_\infty}
 \right)\\
 &\quad+\Pp\!\left(
 \mathscr A_{\varepsilon,z}^{\|u_{\mathrm{res}}(t)\|_\infty+\|u(t)-u_{\mathrm{res}}(t)\|_\infty}
 \cap\mathscr A_{\varepsilon,r_\varepsilon}^{\|u(t)-u_{\mathrm{res}}(t)\|_\infty}
 \right)\\
 \le~&\Pp\!\left(
 \mathscr A_{\varepsilon,z-r_\varepsilon}^{\|u_{\mathrm{res}}(t)\|_\infty}
 \right)
 +\Pp\!\left(
 \mathscr A_{\varepsilon,r_\varepsilon}^{\|u(t)-u_{\mathrm{res}}(t)\|_\infty}
 \right)\\
 \le~&\Pp\!\left(
 \mathscr A_{\varepsilon,z-r_\varepsilon}^{\|u_{\mathrm{res}}(t)\|_\infty}
 \right)+\Pp(\Omega_{\delta_\varepsilon}).
\end{align*}
Taking logarithms of the preceding bound and multiplying by
$\varepsilon$, we factor out the resonant tail probability inside
the logarithm and use $\log(ab)=\log a+\log b$:
\begin{align*}
 \varepsilon\log\Pp\!\left(
 \mathscr A_{\varepsilon,z}^{\|u(t)\|_\infty}\right)
 \le~&\varepsilon\log\left[
 \Pp\!\left(\mathscr A_{\varepsilon,z-r_\varepsilon}^{\|u_{\mathrm{res}}(t)\|_\infty}\right)
 +\Pp(\Omega_{\delta_\varepsilon})\right]\\
 =~&\varepsilon\log\left\{
 \Pp\!\left(\mathscr A_{\varepsilon,z-r_\varepsilon}^{\|u_{\mathrm{res}}(t)\|_\infty}\right)
 \left[1+\frac{\Pp(\Omega_{\delta_\varepsilon})}
 {\Pp\!\left(\mathscr A_{\varepsilon,z-r_\varepsilon}^{\|u_{\mathrm{res}}(t)\|_\infty}\right)}\right]
 \right\}\\
 =~&\varepsilon\log
 \Pp\!\left(\mathscr A_{\varepsilon,z-r_\varepsilon}^{\|u_{\mathrm{res}}(t)\|_\infty}\right)+\varepsilon\log\left[1+
 \frac{\Pp(\Omega_{\delta_\varepsilon})}
 {\Pp\!\left(\mathscr A_{\varepsilon,z-r_\varepsilon}^{\|u_{\mathrm{res}}(t)\|_\infty}\right)}\right].
\end{align*}

\medskip
\noindent\textbf{Step 2. Lower Bound.}
The reverse triangle inequality gives the event inclusions
\begin{equation}\label{eq:critical-good-tail-inclusions}
 \mathscr A_{\varepsilon,z+r_\varepsilon}^{\|u_{\mathrm{res}}(t)\|_\infty}
 \cap\mathscr B_{\varepsilon,r_\varepsilon}^{\|u(t)-u_{\mathrm{res}}(t)\|_\infty}
 \subseteq\mathscr A_{\varepsilon,z}^{\|u_{\mathrm{res}}(t)\|_\infty-\|u(t)-u_{\mathrm{res}}(t)\|_\infty}
 \subseteq\mathscr A_{\varepsilon,z}^{\|u(t)\|_\infty}.
\end{equation}
By \eqref{eq:critical-good-tail-inclusions}, the \emph{total probability
formula} and
\eqref{eq:critical-comparison-event-inclusion}, we obtain
\begin{align*}
 \Pp\!\left(\mathscr A_{\varepsilon,z}^{\|u(t)\|_\infty}\right)
 &\ge\Pp\!\left(
 \mathscr A_{\varepsilon,z}^{\|u_{\mathrm{res}}(t)\|_\infty-\|u(t)-u_{\mathrm{res}}(t)\|_\infty}\right)\\
 &\ge\Pp\!\left(
 \mathscr A_{\varepsilon,z+r_\varepsilon}^{\|u_{\mathrm{res}}(t)\|_\infty}
 \cap\mathscr B_{\varepsilon,r_\varepsilon}^{\|u(t)-u_{\mathrm{res}}(t)\|_\infty}\right)\\
 &=\Pp\!\left(
 \mathscr A_{\varepsilon,z+r_\varepsilon}^{\|u_{\mathrm{res}}(t)\|_\infty}\right)
 -\Pp\!\left(
 \mathscr A_{\varepsilon,z+r_\varepsilon}^{\|u_{\mathrm{res}}(t)\|_\infty}
 \cap\mathscr A_{\varepsilon,r_\varepsilon}^{\|u(t)-u_{\mathrm{res}}(t)\|_\infty}\right)\\
 &\ge\Pp\!\left(
 \mathscr A_{\varepsilon,z+r_\varepsilon}^{\|u_{\mathrm{res}}(t)\|_\infty}\right)
 -\Pp\!\left(
 \mathscr A_{\varepsilon,r_\varepsilon}^{\|u(t)-u_{\mathrm{res}}(t)\|_\infty}\right)\\
 &\ge\Pp\!\left(
 \mathscr A_{\varepsilon,z+r_\varepsilon}^{\|u_{\mathrm{res}}(t)\|_\infty}\right)
 -\Pp(\Omega_{\delta_\varepsilon}).
\end{align*}
Taking logarithms, multiplying by $\varepsilon$ and factoring out
the resonant probability give
\begin{align*}
 \varepsilon\log\Pp\!\left(
 \mathscr A_{\varepsilon,z}^{\|u(t)\|_\infty}\right)
 &\ge\varepsilon\log\left[
 \Pp\!\left(\mathscr A_{\varepsilon,z+r_\varepsilon}^{\|u_{\mathrm{res}}(t)\|_\infty}\right)
 -\Pp(\Omega_{\delta_\varepsilon})\right]\\
 &=\varepsilon\log\left\{
 \Pp\!\left(\mathscr A_{\varepsilon,z+r_\varepsilon}^{\|u_{\mathrm{res}}(t)\|_\infty}\right)
 \left[1-\frac{\Pp(\Omega_{\delta_\varepsilon})}
 {\Pp\!\left(\mathscr A_{\varepsilon,z+r_\varepsilon}^{\|u_{\mathrm{res}}(t)\|_\infty}\right)}\right]
 \right\}\\
 &=\varepsilon\log
 \Pp\!\left(\mathscr A_{\varepsilon,z+r_\varepsilon}^{\|u_{\mathrm{res}}(t)\|_\infty}\right)+\varepsilon\log\left[1-
 \frac{\Pp(\Omega_{\delta_\varepsilon})}
 {\Pp\!\left(\mathscr A_{\varepsilon,z+r_\varepsilon}^{\|u_{\mathrm{res}}(t)\|_\infty}\right)}\right].
\end{align*}

\medskip
\noindent\textbf{Step 3. Resonant Main Terms.}
Set
\[
 p_\varepsilon(q):=
 \Pp\!\left(\mathscr A_{\varepsilon,q}^{\|u_{\mathrm{res}}(t)\|_\infty}\right),
 \qquad q>0.
\]
Fix $0<h<z$. Since $r_\varepsilon\to0$,
for sufficiently small $\varepsilon$ we have
\[
 p_\varepsilon(z+h)\le p_\varepsilon(z\pm r_\varepsilon)
 \le p_\varepsilon(z-h).
\]
Taking logarithms, multiplying by $\varepsilon$ and applying
the cubic resonant LDP, \autoref{thm:nf-resonant-LDP}, give
\begin{align*}
 -\frac{(z+h)^2}{S}
 &\le\liminf_{\varepsilon\to0^+}
 \varepsilon\log p_\varepsilon(z\pm r_\varepsilon)\\
 &\le\limsup_{\varepsilon\to0^+}
 \varepsilon\log p_\varepsilon(z\pm r_\varepsilon)
 \le-\frac{(z-h)^2}{S}.
\end{align*}
Letting $h\downarrow0$, we obtain the common limit of the main terms:
\begin{equation}\label{eq:critical-moving-threshold-ldp}
 \varepsilon\log p_\varepsilon(z\pm r_\varepsilon)
 \longrightarrow-\frac{z^2}{S}.
\end{equation}

\medskip
\noindent\textbf{Step 4. Relative Probabilities.}
For these ratios, \autoref{prop:exception-prob}
and \eqref{eq:nf-high-mode-delta} imply
\[
 \Pp(\Omega_{\delta_\varepsilon})
 \le C_{\kappa,\nu}\exp\!\left(-\frac{1+z^2/S}{\varepsilon}\right).
\]
Hence, for either choice of sign,
\[
\begin{aligned}
 \varepsilon\log
 \frac{\Pp(\Omega_{\delta_\varepsilon})}
 {p_\varepsilon(z\pm r_\varepsilon)}
 \le\varepsilon\log C_{\kappa,\nu}-1-\frac{z^2}{S}
 -\varepsilon\log p_\varepsilon(z\pm r_\varepsilon).
\end{aligned}
\]
By \eqref{eq:critical-moving-threshold-ldp}, the right-hand side
converges to $-1$. It is therefore at most $-1/2$ for sufficiently
small $\varepsilon$, so that
\[
 \varepsilon\log
 \frac{\Pp(\Omega_{\delta_\varepsilon})}
 {p_\varepsilon(z\pm r_\varepsilon)}\le-\frac12.
\]
Dividing by $\varepsilon$ and exponentiating give
\begin{equation}\label{eq:critical-relative-error-small}
 0\le\frac{\Pp(\Omega_{\delta_\varepsilon})}
 {p_\varepsilon(z\pm r_\varepsilon)}\le e^{-1/(2\varepsilon)}.
\end{equation}
Using $0\le\log(1+x)\le x$ and
$-2x\le\log(1-x)\le0$ for $0\le x\le1/2$, we find
\begin{align}
 0&\le\varepsilon\log\left[1+
 \frac{\Pp(\Omega_{\delta_\varepsilon})}
 {p_\varepsilon(z-r_\varepsilon)}\right]
 \le\varepsilon e^{-1/(2\varepsilon)},\nonumber\\
 -2\varepsilon e^{-1/(2\varepsilon)}
 &\le\varepsilon\log\left[1-
 \frac{\Pp(\Omega_{\delta_\varepsilon})}
 {p_\varepsilon(z+r_\varepsilon)}\right]\le0.\label{eq:critical-logarithmic-corrections}
\end{align}
Both corrections therefore tend to zero. 

\medskip
\noindent\textbf{Step 5. Matching Upper and Lower Limits.}
All the bounds in Steps 3 and 4 are independent of $t$.
Combining the upper logarithmic bound of Step 1, obtained from
\eqref{eq:critical-triangle-event-inclusion}--\eqref{eq:critical-large-error-event-inclusion},
with \eqref{eq:critical-moving-threshold-ldp} and
\eqref{eq:critical-logarithmic-corrections}, we obtain
\begin{equation}\label{eq:critical-uniform-upper-limit}
 \limsup_{\varepsilon\to0^+}\sup_{0\le t\le T_\varepsilon}
 \varepsilon\log\Pp\!\left(
 \mathscr A_{\varepsilon,z}^{\|u(t)\|_\infty}\right)
 \le-\frac{z^2}{S}.
\end{equation}
Similarly, the lower logarithmic bound of Step 2, based on
\eqref{eq:critical-good-tail-inclusions}, together with
\eqref{eq:critical-moving-threshold-ldp} and
\eqref{eq:critical-logarithmic-corrections}, gives
\begin{equation}\label{eq:critical-uniform-lower-limit}
 \liminf_{\varepsilon\to0^+}\inf_{0\le t\le T_\varepsilon}
 \varepsilon\log\Pp\!\left(
 \mathscr A_{\varepsilon,z}^{\|u(t)\|_\infty}\right)
 \ge-\frac{z^2}{S}.
\end{equation}
Since the infimum is at most the supremum for each $\varepsilon$,
\eqref{eq:critical-uniform-upper-limit} and
\eqref{eq:critical-uniform-lower-limit} show that both converge to
$-z^2/S$. Equivalently,
\[
 \lim_{\varepsilon\to0^+}\sup_{0\le t\le T_\varepsilon}
 \varepsilon\log\Pp\!\left(
 \mathscr A_{\varepsilon,z}^{\|u(t)\|_\infty}\right)
 =\lim_{\varepsilon\to0^+}\inf_{0\le t\le T_\varepsilon}
 \varepsilon\log\Pp\!\left(
 \mathscr A_{\varepsilon,z}^{\|u(t)\|_\infty}\right)
 =-\frac{z^2}{S},
\]
which is \eqref{eq:nonlinear-first-critical-ldp}. 

\end{proof}

\section{LDP at the Critical Time Scale with a Small Enough Logarithmic Correction Prefactor}
\label{sec:log-critical}

Motivated by the periodic 1D results \cite{GGKS23,LW25,FY25}, we extend
the critical time LDP of \autoref{thm:first-critical-ldp-local} by a
logarithmic factor with a small enough prefactor.
To this end, we first construct a local solution by
\emph{Picard iteration in the high weight space} and establish
an \emph{extension criterion}; see
\autoref{prop:log-local-wellposedness} and \autoref{lem:log-continuation}.
Starting from this local solution, we combine \emph{the tame estimates}
with \emph{integration by parts in time} to compare the nonlinear and
resonant flows under an assumed low weight bound;
see \autoref{lem:log-bootstrap-comparison}.
These estimates close \emph{the low weight bootstrap} while controlling
high weight growth, yielding the logarithmic extension and a vanishing
normalized error; see \autoref{prop:log-pointwise-comparison}.
Finally, we use this vanishing error together with
\emph{a relative probability estimate} to transfer the resonant LDP to
the nonlinear solution on the logarithmic interval; see
\eqref{eq:log-relative-exception-small} and \autoref{thm:log-critical-ldp}.
\subsection{Motivation for the Logarithmic Extension}
\label{subsec:logarithmic-time-mechanism}
\label{app:logarithmic-time-mechanism}
\subsubsection{Observation}
The critical time comparison \eqref{eq:resonant-comparison-critical} gives
\[
 \sup_{0\le t\le T_\varepsilon}
 \sqrt\varepsilon\,\|u(t)-u_{\mathrm{res}}(t)\|_\infty
 \le r_\varepsilon=\mathcal{O}(\varepsilon^{\alpha-1})\longrightarrow0,
 \qquad \alpha>1.
\]
Transferring the LDP from the resonant approximation to the nonlinear
solution requires only vanishing normalized error. This leaves room
to extend the observation interval with a weaker bound that still
tends to zero.

\subsubsection{Model Analysis}

Let $0<\mu\ll1$. In the slow time $t':=\mu t$, the critical scale is $t'\asymp1$;
we write $t$ for $t'$ below. Consider a nonnegative, absolutely continuous normalized error $e(t)$ satisfying
\begin{equation}\label{eq:appendix-model-error}
 e'(t)\le Ce(t)+\mu,\qquad e(0)=0,
\end{equation}
where $C>0$ is independent of $\mu$. We call $Ce(t)$ the
\emph{stability term} and $\mu$ the \emph{error source}.
Without the source, Gronwall's inequality gives $e\equiv0$.
Multiplying \eqref{eq:appendix-model-error} by $\exp(-Ct)$ and
integrating give
\begin{equation}\label{eq:appendix-model-solution}
 e(t)\le\mu\int_0^t\exp(C(t-r))\,dr
 =\frac{\mu}{C}(\exp(Ct)-1).
\end{equation}
The error need only vanish, rather than remain $\mathcal{O}(\mu)$.
For $0<c<1/C$, this allows
\begin{equation}\label{eq:appendix-logarithmic-gain}
 \sup_{0\le t\le c\log(1/\mu)}e(t)
 \le C^{-1}\mu^{1-Cc}\longrightarrow0.
\end{equation}
The logarithmic extension follows from balancing the small source
$\mu$ against the exponential growth $\exp(Ct)$ in the stability estimate.

\subsection{Normalization}
\label{subsec:log-normalization}
Put
$\mu_\varepsilon:=\varepsilon^{\alpha-1}$ and set
\begin{equation}\label{eq:log-normalization}
 \widetilde a_n(t):=\sqrt\varepsilon\,a_n(t),\qquad
 \widetilde b_n(t):=\sqrt\varepsilon\,b_n(t).
\end{equation}
Substituting \eqref{eq:log-normalization} into the equations
\eqref{eq:critical-gauged-cubic} for $a_n$ and
\eqref{eq:nf-resonant-ode} for $b_n$ and using cubic homogeneity, we obtain
\begin{equation}\label{eq:log-normalized-equations}
 \dot{\widetilde a}_n=\iu\mu_\varepsilon H_{t,n}[\widetilde a],\qquad
 \dot{\widetilde b}_n=-\iu\mu_\varepsilon|\widetilde b_n|^2\widetilde b_n.
\end{equation}

Define the normalized initial sequence by
\begin{equation}\label{eq:log-normalized-initial-data}
 \widetilde d_n:=\sqrt\varepsilon\,c(n)g_n,
 \qquad \widetilde a_n(0)=\widetilde b_n(0)=\widetilde d_n.
\end{equation}
By \eqref{eq:randomized-decay} and \eqref{eq:nf-high-mode-delta},
\begin{equation}\label{eq:log-normalized-initial-decay}
 |\widetilde d_n|\le A(1+z^2/S)^{1/2}w_{-\beta}(n)\quad\text{on}~G_{\delta_\varepsilon}.
\end{equation}
The explicit solution \eqref{eq:nf-resonant-profile} gives
\[
 \widetilde b_n(t)=\widetilde d_n
 \exp\bigl(-\iu\mu_\varepsilon t|\widetilde d_n|^2\bigr),
 \qquad |\widetilde b_n(t)|=|\widetilde d_n|.
\]

\subsection{Local Theory and Tame Estimates}
\label{subsec:log-local-theory}

For a coefficient sequence $f=(f_n)$ and $s\ge1$, define its
\emph{weighted supremum norm} by
\[
 \|f\|_{\ell^\infty_s(\Z^\nu)}:=\sup_{n\in\Z^\nu}\langle n\rangle^s|f_n|.
\]
We use two weighted supremum norms: the lower weight makes the
cubic sums summable, while the higher weight compensates for the
small divisors. 

\subsubsection{Tame Estimates}
The following estimates control the high weight norm of the time
derivative of $\widetilde a$ in terms of the high and low weight
norms of $\widetilde a$, and bound the boundary terms produced by
\emph{integration by parts}.
\begin{lemma}[Tame Estimates]\label{lem:log-pointwise-tame}
Let $\omega\in\mathrm{DC}(\gamma,\tau)$, $q>\nu$ and
$p=q+2\tau$. For $f\in\ell^\infty_p(\Z^\nu)$, define the nonresonant boundary operator by
\begin{equation}\label{eq:log-primitive}
 \mathcal B_{t,n}[f]
 :=\sum_{\substack{n_1-n_2+n_3=n\\\Sigma_3\ne0}}
 \frac{e^{-\iu t\Sigma_3}}{\Sigma_3}
 f_{n_1}\overline{f_{n_2}}f_{n_3}.
\end{equation}
Set $\mathcal B_t[f]:=(\mathcal B_{t,n}[f])_{n\in\Z^\nu}$.
Under these hypotheses, $H_t$ maps $\ell^\infty_p$ to $\ell^\infty_p$,
and $\mathcal B_t$ maps $\ell^\infty_p$ to $\ell^\infty_q$.
Uniformly in $t$,
\begin{align}
 \|H_t[f]\|_{\ell^\infty_p(\Z^\nu)}&\le C_{p,q}\|f\|_{\ell^\infty_p(\Z^\nu)}\|f\|_{\ell^\infty_q(\Z^\nu)}^2,
 \label{eq:log-high-cubic-tame}\\
 \|H_t[f]\|_{\ell^\infty_q(\Z^\nu)}&\le C_q\|f\|_{\ell^\infty_q(\Z^\nu)}^3,
 \label{eq:log-cubic-tame}\\
 \|\mathcal B_t[f]\|_{\ell^\infty_q(\Z^\nu)}
 &\le C_{p,q,\gamma,\tau}\|f\|_{\ell^\infty_p(\Z^\nu)}\|f\|_{\ell^\infty_q(\Z^\nu)}^2.
 \label{eq:log-primitive-tame}
\end{align}
\end{lemma}

\begin{proof}[Proof of \autoref{lem:log-pointwise-tame}]
Since $q>\nu$, 
\begin{align}
 \|\langle\cdot\rangle^{-q}\|_{\ell^1(\Z^\nu)}
 =\sum_{n\in\Z^\nu}\langle n\rangle^{-q}<\infty.
\end{align}
\par\medskip\noindent\textbf{Step 1. Weighted Convolution Bound.}
Momentum conservation gives
$\langle n\rangle\le\sum_{j=1}^3\langle n_j\rangle$.
Convexity therefore yields, for $s\ge1$,
\begin{equation}\label{eq:log-momentum-weight}
 \langle n\rangle^s\le C_s\sum_{j=1}^3\langle n_j\rangle^s,
 \quad C_s=3^{s-1}.
\end{equation}
For
$s=p$ or $s=q$, we have
\begin{align}
 \sum_{n_1-n_2+n_3=n}\langle n_1\rangle^s
       |f_{n_1}f_{n_2}f_{n_3}|
       &= \sum_{n_2,n_3}\langle n+n_2-n_3\rangle^s
             |f_{n+n_2-n_3}|\,|f_{n_2}|\,|f_{n_3}|\nonumber\\
 &\le\|f\|_{\ell^\infty_s(\Z^\nu)}
       \sum_{n_2,n_3}|f_{n_2}|\,|f_{n_3}|\nonumber\\
 &\le\|f\|_{\ell^\infty_s(\Z^\nu)}
       \|f\|_{\ell^\infty_q(\Z^\nu)}^2
       \sum_{n_2,n_3}\langle n_2\rangle^{-q}\langle n_3\rangle^{-q}\nonumber\\
 &=\|\langle\cdot\rangle^{-q}\|_{\ell^1(\Z^\nu)}^2
       \|f\|_{\ell^\infty_s(\Z^\nu)}
       \|f\|_{\ell^\infty_q(\Z^\nu)}^2.\label{eq:log-weighted-convolution-bound}
\end{align}
The other two weighted sums with $\langle n_2\rangle^s$ or $\langle n_3\rangle^s$ have the same bound by solving the momentum
relation for $n_2$ or $n_3$, respectively. 

\par\medskip\noindent\textbf{Step 2. Cubic Nonlinear Bounds.}
It follows from the definition \eqref{eq:critical-cubic-map} of $H_{t,n}$, $|\mathcal K_{t,n}|\le1$ and
\eqref{eq:log-momentum-weight}--\eqref{eq:log-weighted-convolution-bound} that
\begin{align*}
 \langle n\rangle^s|H_{t,n}[f]|
 &\le\sum_{n_1-n_2+n_3=n}\langle n\rangle^s
             \prod_{j=1}^3|f_{n_j}|\\
 &\le C_s\sum_{j=1}^3\sum_{n_1-n_2+n_3=n}
             \langle n_j\rangle^s\prod_{i=1}^3|f_{n_i}|\\
 &\le C_{s,q}\|f\|_{\ell^\infty_s(\Z^\nu)}
             \|f\|_{\ell^\infty_q(\Z^\nu)}^2.
\end{align*}
Taking the supremum over $n$, first with $s=p$ and then with $s=q$,
proves \eqref{eq:log-high-cubic-tame} and \eqref{eq:log-cubic-tame}.

\par\medskip\noindent\textbf{Step 3. Boundary Operator Bound.}
For $\Sigma_3\ne0$, both $n_1-n_2$ and $n_2-n_3$ are nonzero.
The factorization \eqref{eq:Sigma-factorization} and the Diophantine
condition \eqref{eq:HDC} give
\begin{align}
 \frac1{|\Sigma_3|}
 &=\frac1{2|\langle n_1-n_2,\omega\rangle|
              |\langle n_2-n_3,\omega\rangle|}\nonumber\\
 &\le\frac{|n_1-n_2|_1^\tau|n_2-n_3|_1^\tau}{2\gamma^2}\nonumber\\
 &\le\frac{(|n_1|_1+|n_2|_1)^\tau(|n_2|_1+|n_3|_1)^\tau}{2\gamma^2}\nonumber\\
 &\le C_{\gamma,\tau}\max_{1\le j\le3}
                                  \langle n_j\rangle^{2\tau}.\label{eq:log-small-divisor-loss}
\end{align}
Since $p=q+2\tau$, multiplying the preceding bound by
$\langle n\rangle^q$ yields
\begin{equation}\label{eq:log-pointwise-divisor}
 \frac{\langle n\rangle^q}{|\Sigma_3|}
 \le C_{p,q,\gamma,\tau}
\max_j\langle n_j\rangle^p
 \le C_{p,q,\gamma,\tau}\sum_{j=1}^3\langle n_j\rangle^p.
\end{equation}
Substituting \eqref{eq:log-pointwise-divisor} into
\eqref{eq:log-primitive}, and then applying the weighted convolution estimates \eqref{eq:log-weighted-convolution-bound}, give
\begin{align*}
 \langle n\rangle^q|\mathcal B_{t,n}[f]|
 &\le\sum_{\substack{n_1-n_2+n_3=n\\\Sigma_3\ne0}}
       \frac{\langle n\rangle^q}{|\Sigma_3|}
       \prod_{j=1}^3|f_{n_j}|\\
 &\le C_{p,q,\gamma,\tau}\sum_{j=1}^3
       \sum_{n_1-n_2+n_3=n}
       \langle n_j\rangle^p\prod_{i=1}^3|f_{n_i}|\\
 &\le C_{p,q,\gamma,\tau}
       \|f\|_{\ell^\infty_p(\Z^\nu)}
       \|f\|_{\ell^\infty_q(\Z^\nu)}^2.
\end{align*}
Taking the supremum over $n$ shows \eqref{eq:log-primitive-tame}. 
\end{proof}

\subsubsection{Local Construction}
\label{subsubsec:log-local-construction}

As in \eqref{eq:critical-pointwise-picard}, set
$\widetilde a_n^{(0)}(t)=\widetilde d_n$ and define
\begin{equation}\label{eq:log-pointwise-picard}
 \widetilde a_n^{(k+1)}(t)=\widetilde d_n+
 \iu\mu_\varepsilon\int_0^t H_{r,n}[\widetilde a^{(k)}(r)]\,dr,
 \qquad k\ge0.
\end{equation}
We construct $\widetilde a$ in the weighted space
\begin{equation}\label{eq:log-weighted-tail-space}
 \ell^\infty_{p,0}(\Z^\nu)
 :=\left\{f\in\ell^\infty_p(\Z^\nu):
       \lim_{|n|\to\infty}\langle n\rangle^p|f_n|=0\right\}.
\end{equation}
This is a closed Banach subspace of $\ell^\infty_p$, in which finite
truncations converge in norm. Together with the cubic estimates below,
this approximation gives joint continuity of $(t,f)\mapsto H_t[f]$.
Closedness ensures that the Picard limit retains the vanishing weighted tail.
We then identify $\widetilde a$ with $\sqrt\varepsilon\,a$ on their
common interval, where both decay estimates hold.

\begin{lemma}[Cubic and Lipschitz Estimates]\label{lem:log-local-cubic}
For the iterates \eqref{eq:log-pointwise-picard}, the following bounds
hold uniformly in $k$ and $t$, with $k\ge0$ for the cubic bound
and $k\ge1$ for the Lipschitz bound:
\begin{align}
 \text{Cubic Bound:}&\quad
 \|H_t[\widetilde a^{(k)}]\|_{\ell^\infty_p}
 \le C\|\widetilde a^{(k)}\|_{\ell^\infty_p}^3.
 \label{eq:log-local-cubic-bound}\\[5pt]
 \text{Lipschitz Bound:}&\quad
 \|H_t[\widetilde a^{(k)}]-H_t[\widetilde a^{(k-1)}]\|_{\ell^\infty_p}
 \nonumber\\
 &\quad\le C\bigl(\|\widetilde a^{(k)}\|_{\ell^\infty_p}
             +\|\widetilde a^{(k-1)}\|_{\ell^\infty_p}\bigr)^2
 \|\widetilde a^{(k)}-\widetilde a^{(k-1)}\|_{\ell^\infty_p}.
 \label{eq:log-local-lipschitz}
\end{align}
The same estimates hold for arbitrary inputs in $\ell^\infty_p$.
\end{lemma}
\begin{proof}[Proof of \autoref{lem:log-local-cubic}]
For \eqref{eq:log-local-cubic-bound}, the kernel bound
$|\mathcal K_{t,n}|\le1$, \eqref{eq:log-momentum-weight} and
\eqref{eq:log-weighted-convolution-bound} give
\begin{align*}
 \langle n\rangle^p|H_{t,n}[\widetilde a^{(k)}]|
 &\le C_p\sum_{j=1}^3\sum_{n_1-n_2+n_3=n}
 \langle n_j\rangle^p\prod_{i=1}^3|\widetilde a_{n_i}^{(k)}|\\
 &\le C_{p,q}\|\widetilde a^{(k)}\|_{\ell^\infty_p}
             \|\widetilde a^{(k)}\|_{\ell^\infty_q}^2\\
 &\le C_{p,q}\|\widetilde a^{(k)}\|_{\ell^\infty_p}^3.
\end{align*}
The last step uses $p\ge q$. Taking the supremum over $n$ proves
the cubic bound.

For the difference estimate, expand the cubic map at a fixed time, 
\begin{equation}\label{eq:log-cubic-increment-expansion}
\begin{aligned}[b]
&H_{t,n}[\widetilde a^{(k)}]-H_{t,n}[\widetilde a^{(k-1)}]\\
&\quad=\sum_{n_1-n_2+n_3=n}\mathcal K_{t,n}(n_1,n_2,n_3)
\left[\begin{aligned}
 &(\widetilde a_{n_1}^{(k)}-\widetilde a_{n_1}^{(k-1)})
 \overline{\widetilde a_{n_2}^{(k)}}\widetilde a_{n_3}^{(k)}+{}\\
 &\widetilde a_{n_1}^{(k-1)}
 \overline{(\widetilde a_{n_2}^{(k)}-\widetilde a_{n_2}^{(k-1)})}
 \widetilde a_{n_3}^{(k)}+{}\\
 &\widetilde a_{n_1}^{(k-1)}\overline{\widetilde a_{n_2}^{(k-1)}}
 (\widetilde a_{n_3}^{(k)}-\widetilde a_{n_3}^{(k-1)})
\end{aligned}\right].
\end{aligned}
\end{equation}
By \eqref{eq:log-momentum-weight}, distribute the output weight
$\langle n\rangle^p$ among the three inputs. We first estimate
the first summand of \eqref{eq:log-cubic-increment-expansion};
all coefficients below are evaluated at the same time.

\par\medskip\textit{Weight on the Difference Factor.}
Solve $n_1=n+n_2-n_3$. Extract the weighted supremum of the
difference, then sum the fixed weights of the other two factors:
\begin{align}
 &\sum_{n_1-n_2+n_3=n}\langle n_1\rangle^p
 |\widetilde a_{n_1}^{(k)}-\widetilde a_{n_1}^{(k-1)}|
 |\widetilde a_{n_2}^{(k)}|\,|\widetilde a_{n_3}^{(k)}|\nonumber\\
 &\qquad=\sum_{n_2,n_3}\langle n+n_2-n_3\rangle^p
 |\widetilde a_{n+n_2-n_3}^{(k)}-\widetilde a_{n+n_2-n_3}^{(k-1)}|
 |\widetilde a_{n_2}^{(k)}|\,|\widetilde a_{n_3}^{(k)}|\nonumber\\
 &\qquad\le\|\widetilde a^{(k)}-\widetilde a^{(k-1)}\|_{\ell^\infty_p}
 \|\widetilde a^{(k)}\|_{\ell^\infty_q}^2
 \sum_{n_2,n_3}\langle n_2\rangle^{-q}\langle n_3\rangle^{-q}\nonumber\\
 &\qquad=\|\langle\cdot\rangle^{-q}\|_{\ell^1}^2
 \|\widetilde a^{(k)}-\widetilde a^{(k-1)}\|_{\ell^\infty_p}
 \|\widetilde a^{(k)}\|_{\ell^\infty_q}^2.
 \label{eq:log-increment-weight-on-difference}
\end{align}

\par\medskip\textit{Weight on the Second or Third Factor.}
Solve $n_2=n_1+n_3-n$ or $n_3=n-n_1+n_2$, respectively.
The weighted factor supplies the high weight norm, while the other
two supply the summable low weights. We detail the second-factor
case; the third-factor case has the same bound.
\begin{align}
 &\sum_{n_1-n_2+n_3=n}\langle n_2\rangle^p
 |\widetilde a_{n_1}^{(k)}-\widetilde a_{n_1}^{(k-1)}|
 |\widetilde a_{n_2}^{(k)}|\,|\widetilde a_{n_3}^{(k)}|\nonumber\\
 &\qquad=\sum_{n_1,n_3}
 |\widetilde a_{n_1}^{(k)}-\widetilde a_{n_1}^{(k-1)}|
 \langle n_1+n_3-n\rangle^p|\widetilde a_{n_1+n_3-n}^{(k)}|
 |\widetilde a_{n_3}^{(k)}|\nonumber\\
 &\qquad\le\|\widetilde a^{(k)}\|_{\ell^\infty_p}
 \|\widetilde a^{(k)}-\widetilde a^{(k-1)}\|_{\ell^\infty_q}
 \|\widetilde a^{(k)}\|_{\ell^\infty_q}
 \sum_{n_1,n_3}\langle n_1\rangle^{-q}\langle n_3\rangle^{-q}\nonumber\\
 &\qquad=\|\langle\cdot\rangle^{-q}\|_{\ell^1}^2
 \|\widetilde a^{(k)}\|_{\ell^\infty_p}
 \|\widetilde a^{(k)}-\widetilde a^{(k-1)}\|_{\ell^\infty_q}
 \|\widetilde a^{(k)}\|_{\ell^\infty_q}.
 \label{eq:log-increment-weight-on-amplitude}
\end{align}

Since $p\ge q$, all low weight norms are bounded by their
high weight counterparts. The second and third summands of
\eqref{eq:log-cubic-increment-expansion} are treated identically;
the last factor becomes, respectively,
$\|\widetilde a^{(k-1)}\|_{\ell^\infty_p}
 \allowbreak\|\widetilde a^{(k)}\|_{\ell^\infty_p}$ and
$\|\widetilde a^{(k-1)}\|_{\ell^\infty_p}^2$.
Adding the three bounds proves \eqref{eq:log-local-lipschitz},
since the sum of these three products is bounded by
$(\|\widetilde a^{(k)}\|_{\ell^\infty_p}
 +\|\widetilde a^{(k-1)}\|_{\ell^\infty_p})^2$.
The same argument applies to arbitrary inputs in $\ell^\infty_p$.
\end{proof}

\begin{proposition}[Local Existence, Uniqueness and Identification]
\label{prop:log-local-wellposedness}
Let $\omega\in\mathrm{DC}(\gamma,\tau)$, $q>\nu$ and
$p=q+2\tau<\min_j\beta_j$.
Choose a deterministic $R\ge1$ such that
$\|\widetilde d\|_{\ell^\infty_p(\Z^\nu)}\le R$ on
$G_{\delta_\varepsilon}$. There is a constant $C>0$, independent of
$\varepsilon$ and the sample point, such that, with
\begin{equation}\label{eq:log-local-time}
 \mathcal T_0:=\frac{1}{32C\mu_\varepsilon R^2},
\end{equation}
the nonlinear initial-value problem
\eqref{eq:log-normalized-equations}--\eqref{eq:log-normalized-initial-data}
has a unique solution
$\widetilde a\in C([0,\mathcal T_0];\ell^\infty_{p,0}(\Z^\nu))$
with each component continuously differentiable in time and
$\sup_{t\le\mathcal T_0}\|\widetilde a(t)\|_{\ell^\infty_p}\le2R$.
This solution agrees with $\sqrt\varepsilon\,a$, the normalized critical
solution of \autoref{prop:critical-pointwise-lifespan}, on their common interval.
\end{proposition}
\begin{proof}[Proof of \autoref{prop:log-local-wellposedness}]
By \eqref{eq:log-normalized-initial-decay},
\begin{equation}\label{eq:log-initial-pointwise}
 |\widetilde d_n|\le Cw_{-\beta}(n)
 \le C\langle n\rangle^{-\min_j\beta_j},\quad p<\min_j\beta_j.
\end{equation}

\par\medskip\noindent\textbf{Step 1. Convergence of the Iteration.}
Choose $C$ in \eqref{eq:log-local-time} large enough for the estimates in
\autoref{lem:log-local-cubic}. As in the proof of
\autoref{prop:critical-pointwise-lifespan}, the Picard map preserves
the ball of radius $2R$ in
$C([0,\mathcal T_0];\ell^\infty_{p,0})$,
with $\ell^\infty_{p,0}$ defined in \eqref{eq:log-weighted-tail-space}, since
\begin{equation}\label{eq:log-local-iterate-bound}
 R+8C\mu_\varepsilon\mathcal T_0R^3=\tfrac54R<2R.
\end{equation}
By \eqref{eq:log-local-lipschitz}, successive differences satisfy
\begin{equation}\label{eq:log-local-increment-bound}
 \sup_{t\le\mathcal T_0}
 \|\widetilde a^{(k+1)}(t)-\widetilde a^{(k)}(t)\|_{\ell^\infty_p}
 \le\tfrac12\sup_{t\le\mathcal T_0}
 \|\widetilde a^{(k)}(t)-\widetilde a^{(k-1)}(t)\|_{\ell^\infty_p}.
\end{equation}
The first increment is bounded by
$C\mu_\varepsilon \mathcal T_0R^3\le R$. Iterating
\eqref{eq:log-local-increment-bound} gives
\begin{equation}\label{eq:log-local-geometric-bound}
 \sup_{0\le t\le \mathcal T_0}
 \|\widetilde a^{(k+1)}(t)-\widetilde a^{(k)}(t)\|_{\ell^\infty_p(\Z^\nu)}
 \le R\,2^{-k},\qquad k\ge0.
\end{equation}
For $k\ge0$ and $m\ge1$, telescoping and
\eqref{eq:log-local-geometric-bound} give
\begin{align}
 \sup_{0\le t\le \mathcal T_0}
 \|\widetilde a^{(k+m)}(t)-\widetilde a^{(k)}(t)\|_{\ell^\infty_p(\Z^\nu)}\le R\sum_{j=k}^{k+m-1}2^{-j}
 \le2R\,2^{-k}\longrightarrow0.\label{eq:log-local-cauchy}
\end{align}
Thus the iterates converge in $C([0,\mathcal T_0];\ell^\infty_{p,0})$
to $\widetilde a$, with
$\widetilde a(0)=\widetilde d$ and
$\sup_{0\le t\le \mathcal T_0}\|\widetilde a(t)\|_{\ell^\infty_p(\Z^\nu)}\le2R$.

\par\medskip\noindent\textbf{Step 2. Passage to the Integral Equation.}
To verify that the limit solves the equation, we pass the
nonlinearity through the limit. By \eqref{eq:log-local-lipschitz},
\begin{align}
 &\sup_{0\le t\le \mathcal T_0}
 \|H_t[\widetilde a^{(k)}(t)]-H_t[\widetilde a(t)]\|_{\ell^\infty_p(\Z^\nu)}\nonumber\\
 &\qquad\le16CR^2\sup_{0\le t\le \mathcal T_0}
 \|\widetilde a^{(k)}(t)-\widetilde a(t)\|_{\ell^\infty_p(\Z^\nu)}
 \longrightarrow0.\label{eq:log-local-nonlinearity-limit}
\end{align}
In particular,
\begin{align}
 &\sup_{0\le t\le\mathcal T_0}
 \left\|\int_0^t\bigl(H_r[\widetilde a^{(k)}(r)]
       -H_r[\widetilde a(r)]\bigr)\,dr\right\|_{\ell^\infty_p}
 \nonumber\\
 &\qquad\le\mathcal T_0\sup_{0\le r\le\mathcal T_0}
 \|H_r[\widetilde a^{(k)}(r)]-H_r[\widetilde a(r)]\|_{\ell^\infty_p}
 \longrightarrow0.
 \label{eq:log-local-integral-limit}
\end{align}
Thus both sides of \eqref{eq:log-pointwise-picard} converge
uniformly in the high weight norm, giving
\begin{equation}\label{eq:log-local-integral}
 \widetilde a_n(t)=\widetilde d_n+
 \iu\mu_\varepsilon\int_0^t H_{r,n}[\widetilde a(r)]\,dr.
\end{equation}
Joint continuity of the cubic map and the \emph{fundamental theorem of calculus} give
$\widetilde a_n\in C^1([0,\mathcal T_0])$ for every $n$, with the equation
and initial data in \eqref{eq:log-normalized-equations} and
\eqref{eq:log-normalized-initial-data}.

\par\medskip\noindent\textbf{Step 3. Uniqueness.}
The local Lipschitz estimate \eqref{eq:log-local-lipschitz} applies
to any two inputs with bounded high weight norms. Subtracting
their integral equations \eqref{eq:log-local-integral} and applying
\emph{Gronwall} show that two solutions with the same initial data
coincide on every compact common interval.

\par\medskip\noindent\textbf{Step 4. Identification on the Common Interval.}
Both $\widetilde a$ and $\sqrt\varepsilon\,a$ satisfy
\eqref{eq:log-local-integral} with the same initial data.
These solutions satisfy different decay bounds. On any admissible compact common interval $[0,T]$,
\eqref{eq:critical-pointwise-decay} and the high weight bound give
\[
 |\widetilde a_n(r)|\le C\langle n\rangle^{-p},\qquad
 |\sqrt\varepsilon\,a_n(r)|\le Cw_{-\beta/2}(n),
 \qquad 0\le r\le T.
\]
By \eqref{eq:log-cubic-increment-expansion} and
$|\mathcal K_{r,n}|\le1$, with all coefficients evaluated at $r$,
\begin{equation}\label{eq:log-identification-cubic-bound}
\begin{aligned}[b]
 &|H_{r,n}[\widetilde a]-H_{r,n}[\sqrt\varepsilon\,a]|\\
 &\qquad\le C^2\sup_m|\widetilde a_m-\sqrt\varepsilon\,a_m|
 \left[\begin{aligned}
 &\sum_{n_1-n_2+n_3=n}
 \langle n_2\rangle^{-p}\langle n_3\rangle^{-p}+{}\\
 &\sum_{n_1-n_2+n_3=n}
 w_{-\beta/2}(n_1)\langle n_3\rangle^{-p}+{}\\
 &\sum_{n_1-n_2+n_3=n}
 w_{-\beta/2}(n_1)w_{-\beta/2}(n_2)
 \end{aligned}\right].
\end{aligned}
\end{equation}
Solving the momentum relation $n_1-n_2+n_3=n$ for $n_1$, $n_2$ and $n_3$,
respectively, gives
\begin{align*}
 \sum_{n_1-n_2+n_3=n}
 \langle n_2\rangle^{-p}\langle n_3\rangle^{-p}
 &=\sum_{n_2,n_3}\langle n_2\rangle^{-p}\langle n_3\rangle^{-p}
 =\|\langle\cdot\rangle^{-p}\|_{\ell^1}^2,\\
 \sum_{n_1-n_2+n_3=n}
 w_{-\beta/2}(n_1)\langle n_3\rangle^{-p}
 &=\sum_{n_1,n_3}w_{-\beta/2}(n_1)\langle n_3\rangle^{-p}
 =\|w_{-\beta/2}\|_{\ell^1}\|\langle\cdot\rangle^{-p}\|_{\ell^1},\\
 \sum_{n_1-n_2+n_3=n}
 w_{-\beta/2}(n_1)w_{-\beta/2}(n_2)
 &=\sum_{n_1,n_2}w_{-\beta/2}(n_1)w_{-\beta/2}(n_2)
 =\|w_{-\beta/2}\|_{\ell^1}^2.
\end{align*}
These quantities are finite because $p>\nu$ and $\beta_j/2>1$.
Thus the bracket in \eqref{eq:log-identification-cubic-bound} is
independent of $n$ and bounded by
\[
 \left(\|\langle\cdot\rangle^{-p}\|_{\ell^1}
       +\|w_{-\beta/2}\|_{\ell^1}\right)^2.
\]
Subtracting \eqref{eq:log-local-integral} for the two solutions
therefore gives
\begin{align}
 \sup_n|\widetilde a_n(t)-\sqrt\varepsilon\,a_n(t)|\le C\mu_\varepsilon\int_0^t
        \sup_n|\widetilde a_n(r)-\sqrt\varepsilon\,a_n(r)|\,dr,\label{eq:log-critical-uniqueness}
\end{align}
Gronwall's inequality applied to \eqref{eq:log-critical-uniqueness} gives, for every $n\in\Z^\nu$ and throughout the common interval,
\begin{equation}\label{eq:log-critical-identification}
 \widetilde a_n(t)=\sqrt\varepsilon\,a_n(t),\qquad n\in\Z^\nu.
\end{equation}
\end{proof}

\subsection{Extension Criterion}
To extend the local solution to the logarithmic time scale, we need
the following extension criterion: \emph{a finite high weight bound
prevents the solution from terminating at a finite time}.
\begin{lemma}[Extension Criterion]\label{lem:log-continuation}
Let $\widetilde a$ be the solution of
\autoref{prop:log-local-wellposedness}, extended to $[0,T)$
with $0<T<\infty$. If
\begin{equation}\label{eq:log-continuation-assumption}
 \sup_{0\le t<T}\|\widetilde a(t)\|_{\ell^\infty_p(\Z^\nu)}<\infty,
\end{equation}
then $\widetilde a$ extends uniquely beyond $T$. The extension
is continuous in time with values in $\ell^\infty_p(\Z^\nu)$,
and each component is continuously differentiable in time.
\end{lemma}
\begin{proof}[Proof of \autoref{lem:log-continuation}]
\par\medskip\noindent\textbf{Step 1. Construction of the Endpoint Value.} For $0\le s<t<T$,
\eqref{eq:log-local-integral}, \eqref{eq:log-high-cubic-tame}
and $p\ge q$ give
\begin{align}
 \|\widetilde a(t)-\widetilde a(s)\|_{\ell^\infty_p}
 &\le\mu_\varepsilon\int_s^t
       \|H_r[\widetilde a(r)]\|_{\ell^\infty_p}\,dr\nonumber\\
 &\le C\mu_\varepsilon\int_s^t
       \|\widetilde a(r)\|_{\ell^\infty_p}
       \|\widetilde a(r)\|_{\ell^\infty_q}^{2}\,dr\nonumber\\
 &\le C\mu_\varepsilon\int_s^t
       \|\widetilde a(r)\|_{\ell^\infty_p}^{3}\,dr\nonumber\\
 &\le C\mu_\varepsilon(t-s)
       \left(\sup_{0\le r<T}
       \|\widetilde a(r)\|_{\ell^\infty_p}\right)^3.
 \label{eq:log-continuation-increments}
\end{align}
The supremum is finite by \eqref{eq:log-continuation-assumption}.
Hence $\widetilde a(t)$ is Cauchy as $t\uparrow T$ in the Banach
space $\ell^\infty_p$, and there is an endpoint value
$\widetilde a(T)$ satisfying
\[
 \lim_{t\uparrow T}
 \|\widetilde a(t)-\widetilde a(T)\|_{\ell^\infty_p}=0.
\]

\par\medskip\noindent\textbf{Step 2. Vanishing Weighted Tail.}
For every $t<T$, the constructed solution has vanishing weighted
tail. Therefore
\begin{align}
 \limsup_{|n|\to\infty}\langle n\rangle^p|\widetilde a_n(T)|
 &\le\|\widetilde a(T)-\widetilde a(t)\|_{\ell^\infty_p}
   +\limsup_{|n|\to\infty}\langle n\rangle^p|\widetilde a_n(t)|
 \nonumber\\
 &=\|\widetilde a(T)-\widetilde a(t)\|_{\ell^\infty_p}
 \longrightarrow0\qquad(t\uparrow T).
 \label{eq:log-continuation-tail}
\end{align}
Thus $\widetilde a(T)$ belongs to the working space
\eqref{eq:log-weighted-tail-space} of the Picard construction.

\par\medskip\noindent\textbf{Step 3. Extension.}
Restarting at $T$ with initial value
$\widetilde a(T)$ gives a solution on a positive interval beyond
$T$, satisfying
\begin{equation}\label{eq:log-restarted-integral}
 \widetilde a_n(T+s)=\widetilde a_n(T)+
 \iu\mu_\varepsilon\int_0^s
 H_{T+r,n}[\widetilde a(T+r)]\,dr.
\end{equation}
Letting $t\uparrow T$ in \eqref{eq:log-local-integral} gives
\begin{equation}\label{eq:log-endpoint-integral}
 \widetilde a_n(T)=\widetilde d_n+
 \iu\mu_\varepsilon\int_0^T H_{r,n}[\widetilde a(r)]\,dr.
\end{equation}
Substituting \eqref{eq:log-endpoint-integral} into
\eqref{eq:log-restarted-integral} and \emph{changing variables} yield
\begin{align}
 \widetilde a_n(T+s)
 &=\widetilde d_n+
 \iu\mu_\varepsilon\int_0^T H_{r,n}[\widetilde a(r)]\,dr
 +\iu\mu_\varepsilon\int_0^s
 H_{T+r,n}[\widetilde a(T+r)]\,dr\nonumber\\
 &=\widetilde d_n+
 \iu\mu_\varepsilon\int_0^T H_{r,n}[\widetilde a(r)]\,dr
 +\iu\mu_\varepsilon\int_T^{T+s}
 H_{r,n}[\widetilde a(r)]\,dr\nonumber\\
 &=\widetilde d_n+
 \iu\mu_\varepsilon\int_0^{T+s}
 H_{r,n}[\widetilde a(r)]\,dr.
 \label{eq:log-joined-integral}
\end{align}
Thus the extension satisfies the original integral equation
across $T$. Continuity of the cubic map yields componentwise
$C^1$ regularity there, and the local Lipschitz estimate
\eqref{eq:log-local-lipschitz} gives uniqueness of the extension.
\end{proof}
\begin{remark}[Blow-Up Alternative]\label{rem:log-blowup-alternative}
Let
$T_{\max}$ denote the maximal lifespan
of the high weight solution. Then
\begin{equation}\label{eq:log-blowup-alternative}
 T_{\max}<\infty\quad\Longrightarrow\quad
 \limsup_{t\uparrow T_{\max}}
 \|\widetilde a(t)\|_{\ell^\infty_p(\Z^\nu)}=\infty.
\end{equation}
Indeed, boundedness of the high weight norm up to a finite
$T_{\max}$ would allow extension by \autoref{lem:log-continuation},
contradicting maximality. 
\end{remark}

\subsection{Logarithmic Extension and Resonant Comparison}
\label{subsec:log-continuation}
We now prove the longer time result. Comparison with the resonant
flow keeps the low weight norm bounded; the tame estimate then
controls high weight growth and allows the solution to be extended.
The final estimate below is the input for the LDP.

\begin{proposition}[Logarithmic Extension and Resonant Comparison]
\label{prop:log-pointwise-comparison}
Let $\omega\in\mathrm{DC}(\gamma,\tau)$, $q>\nu$ and
$p=q+2\tau<\min_j\beta_j$.
There exist a small enough constant $c>0$, a constant $C>0$
and $\vartheta\in(0,1)$, independent of
$\varepsilon$ and the sample point, such that, for sufficiently
small $\varepsilon$, on $G_{\delta_\varepsilon}$ the nonlinear
initial-value problem
\eqref{eq:log-normalized-equations}--\eqref{eq:log-normalized-initial-data}
has a unique solution
$\widetilde a\in C([0,\mathscr T_\varepsilon];\ell^\infty_p(\Z^\nu))$,
where
\begin{equation}\label{eq:log-observation-time}
 \mathscr T_\varepsilon:=c\mu_\varepsilon^{-1}|\log\mu_\varepsilon|.
\end{equation}
This solution agrees with $\sqrt\varepsilon\,a$ on their common
interval. It also satisfies the following estimates.
\begin{enumerate}[label=\textup{(\roman*)},leftmargin=*]
\item Low Weight Bound:
\begin{equation}\label{eq:log-low-weight-bound}
 \sup_{0\le t\le\mathscr T_\varepsilon}
 \|\widetilde a(t)\|_{\ell^\infty_q(\Z^\nu)}\le C.
\end{equation}
\item High Weight Bound:
\begin{equation}\label{eq:log-pointwise-high-decay}
  \sup_{0\le r\le t}\|\widetilde a(r)\|_{\ell^\infty_p(\Z^\nu)}
 \le C\exp(C\mu_\varepsilon t),\qquad 0\le t\le\mathscr T_\varepsilon.
\end{equation}
\item Fourier-Space Comparison:
\begin{equation}\label{eq:log-pointwise-comparison}
 \sup_{0\le t\le\mathscr T_\varepsilon}
 \|\widetilde a(t)-\widetilde b(t)\|_{\ell^\infty_q(\Z^\nu)}
 \le C\mu_\varepsilon^\vartheta.
\end{equation}
\item Physical-Space Comparison: 
\begin{equation}\label{eq:log-physical-comparison}
 \sup_{0\le t\le\mathscr T_\varepsilon}
 \sqrt\varepsilon\,
 \|u(t)-u_{\mathrm{res}}(t)\|_\infty
 \le C\mu_\varepsilon^\vartheta.
\end{equation}
\end{enumerate}
\end{proposition}

\subsubsection{Fourier-Space Comparison}
\label{subsubsec:log-bootstrap-comparison}
Define
\begin{equation}\label{eq:log-comparison-envelope}
 E_{\log}(t):=\|\widetilde a(t)-\widetilde b(t)\|_{\ell^\infty_q(\Z^\nu)}.
\end{equation}
\begin{lemma}[High Weight Growth and Fourier-Space Comparison]
\label{lem:log-bootstrap-comparison}
Let $\widetilde a$ be the solution of
\autoref{prop:log-local-wellposedness}, extended to $[0,T]$,
and let $\widetilde b$ be the resonant solution with the same initial data $\widetilde d$.
Suppose $\|\widetilde d\|_{\ell^\infty_p(\Z^\nu)}\le R$, $R\ge1$ and
\begin{equation}\label{eq:log-bootstrap-bound}
 \|\widetilde a(t)\|_{\ell^\infty_q(\Z^\nu)}\le2R,\quad 0\le t\le T.
\end{equation}
Then, for $0\le t\le T$,
\begin{equation}\label{eq:log-high-envelope}
 \sup_{0\le r\le t}\|\widetilde a(r)\|_{\ell^\infty_p(\Z^\nu)}
 \le R\exp(C_{p,q}R^2\mu_\varepsilon t),
\end{equation}
and there are constants $C>0$ and
$\Lambda>0$, depending only on
$R$ and the fixed parameters, such that
\begin{equation}\label{eq:log-error-growth}
 E_{\log}(t)\le C\mu_\varepsilon e^{\Lambda\mu_\varepsilon t},
 \qquad 0\le t\le T.
\end{equation}
The constants are independent of $\varepsilon$, $T$ and the sample point.
\end{lemma}
\begin{remark}[Explicit Growth Constant]
\label{rem:log-explicit-growth-constant}
We choose a common exponential rate for the high weight bound
\eqref{eq:log-high-envelope} and the error estimate
\eqref{eq:log-error-growth}. The cubic constant from
\eqref{eq:log-momentum-weight}--\eqref{eq:log-weighted-convolution-bound},
together with the low weight bound \eqref{eq:log-bootstrap-bound}, gives
\begin{equation}\label{eq:log-explicit-growth-constant}
 \Lambda=(2R)^2\,3^p
       \|\langle\cdot\rangle^{-q}\|_{\ell^1(\Z^\nu)}^2.
\end{equation}
This choice bounds the high weight norm by $Re^{\Lambda\mu_\varepsilon t}$.
The boundary and differentiated terms in Steps~4--5 below contain only
one high weight factor, so
\begin{equation}\label{eq:log-explicit-error-integral}
 E_{\log}(t)
 \le C\mu_\varepsilon e^{\Lambda\mu_\varepsilon t}
    +12R^2\mu_\varepsilon\int_0^t E_{\log}(r)\,dr.
\end{equation}
Moreover, $p>1$ and
$\|\langle\cdot\rangle^{-q}\|_{\ell^1(\Z^\nu)}\ge1$
imply $\Lambda>12R^2$. Thus Gronwall's inequality applied to
\eqref{eq:log-explicit-error-integral} gives
\begin{equation}\label{eq:log-explicit-gronwall}
 E_{\log}(t)
 \le C\frac{\Lambda}{\Lambda-12R^2}
          \mu_\varepsilon e^{\Lambda\mu_\varepsilon t},
\end{equation}
which proves \eqref{eq:log-error-growth} with the same $\Lambda$.
In particular, $0<c<\Lambda^{-1}$ ensures
$E_{\log}(t)\lesssim\mu_\varepsilon^{1-\Lambda c}\to0$
on $0\le t\le c\mu_\varepsilon^{-1}|\log\mu_\varepsilon|$,
as required to close the bootstrap.
\end{remark}

\begin{proof}[Proof of \autoref{lem:log-bootstrap-comparison}]
\par\medskip\noindent\textbf{Step 1. High Weight Bound.}
Taking the high weight norm in \eqref{eq:log-local-integral} and
using \eqref{eq:log-high-cubic-tame} give
\begin{align*}
 \|\widetilde a(t)\|_{\ell^\infty_p(\Z^\nu)}
 &\le R+C_{p,q}\mu_\varepsilon\int_0^t
       \|\widetilde a(r)\|_{\ell^\infty_p(\Z^\nu)}
       \|\widetilde a(r)\|_{\ell^\infty_q(\Z^\nu)}^2\,dr\\
 &\le R+4C_{p,q}R^2\mu_\varepsilon\int_0^t
       \|\widetilde a(r)\|_{\ell^\infty_p(\Z^\nu)}\,dr.
\end{align*}
Gronwall's inequality yields \eqref{eq:log-high-envelope}.

\par\medskip\noindent\textbf{Step 2. Integration by Parts.}
Subtracting \eqref{eq:log-normalized-equations} and using
\eqref{eq:nf-gauged-ode} and the common initial data give
\begin{align}
 \widetilde a_n(t)-\widetilde b_n(t)
 ={}&-\iu\mu_\varepsilon\int_0^t
 \bigl(|\widetilde a_n|^2\widetilde a_n-
       |\widetilde b_n|^2\widetilde b_n\bigr)(r)\,dr\nonumber\\
 &+\iu\mu_\varepsilon\int_0^t
        \mathcal N_1^3[\widetilde a](r,n)\,dr.\label{eq:log-difference-integral}
\end{align}
For every nonresonant triple,
$\partial_r e^{-\iu r\Sigma_3}=-\iu\Sigma_3e^{-\iu r\Sigma_3}$.
\emph{Integration by parts} gives
\begin{align}
 \widetilde a_n(t)-\widetilde b_n(t)
 ={}&-\iu\mu_\varepsilon\int_0^t
 \bigl(|\widetilde a_n|^2\widetilde a_n-
       |\widetilde b_n|^2\widetilde b_n\bigr)(r)\,dr\nonumber\\
 &-\mu_\varepsilon\mathcal B_{t,n}[\widetilde a(t)]
  +\mu_\varepsilon\mathcal B_{0,n}[\widetilde d]\nonumber\\
 &+\mu_\varepsilon\int_0^t
 \sum_{\substack{n_1-n_2+n_3=n\\\Sigma_3\ne0}}
 \frac{e^{-\iu r\Sigma_3}}{\Sigma_3}
 \frac{d}{dr}
 \bigl(\widetilde a_{n_1}\overline{\widetilde a_{n_2}}
                    \widetilde a_{n_3}\bigr)(r)\,dr.\label{eq:log-normal-form-identity}
\end{align}

\par\medskip\noindent\textbf{Step 3. Diagonal Cubic Difference.}
The explicit resonant solution in
\autoref{subsec:log-normalization} satisfies
\begin{equation}\label{eq:log-resonant-norm}
 \|\widetilde b(t)\|_{\ell^\infty_q(\Z^\nu)}
 =\|\widetilde d\|_{\ell^\infty_q(\Z^\nu)}\le R,
 \qquad t\ge0.
\end{equation}

The scalar cubic difference estimate used in
\eqref{eq:resonant-diagonal-pointwise-difference}, together with
\eqref{eq:log-resonant-norm} and \eqref{eq:log-bootstrap-bound}, gives
\begin{align}
 \bigl||\widetilde a_n|^2\widetilde a_n-
              |\widetilde b_n|^2\widetilde b_n\bigr|
 &\le3\max\{|\widetilde a_n|,|\widetilde b_n|\}^2
               |\widetilde a_n-\widetilde b_n|\nonumber\\
 &\le12R^2E_{\log}(r)\langle n\rangle^{-q},\label{eq:log-diagonal-difference}
\end{align}
where all coefficients are evaluated at $r$.
\par\medskip\noindent\textbf{Step 4. Boundary Terms.}
The two boundary terms in \eqref{eq:log-normal-form-identity}
are controlled by \eqref{eq:log-primitive-tame}:
\begin{align}
 &\|\mathcal B_t[\widetilde a(t)]\|_{\ell^\infty_q(\Z^\nu)}
  +\|\mathcal B_0[\widetilde d]\|_{\ell^\infty_q(\Z^\nu)}\nonumber\\
 &\qquad\le C\|\widetilde a(t)\|_{\ell^\infty_p(\Z^\nu)}
              \|\widetilde a(t)\|_{\ell^\infty_q(\Z^\nu)}^2
       +C\|\widetilde d\|_{\ell^\infty_p(\Z^\nu)}
          \|\widetilde d\|_{\ell^\infty_q(\Z^\nu)}^2\nonumber\\
 &\qquad\le CR^3e^{CR^2\mu_\varepsilon t}.\label{eq:log-boundary-bound}
\end{align}
\par\medskip\noindent\textbf{Step 5. Differentiated Terms.}
\par\medskip\noindent\emph{Derivative Bounds.}
Applying
\eqref{eq:log-high-cubic-tame} and \eqref{eq:log-cubic-tame} to
\eqref{eq:log-normalized-equations} also gives
\begin{align}
 \|\dot{\widetilde a}(r)\|_{\ell^\infty_p(\Z^\nu)}
 &\le C\mu_\varepsilon
       \|\widetilde a(r)\|_{\ell^\infty_p(\Z^\nu)}
       \|\widetilde a(r)\|_{\ell^\infty_q(\Z^\nu)}^2
 \le C\mu_\varepsilon R^3e^{CR^2\mu_\varepsilon r},
 \label{eq:log-high-derivative}\\
 \|\dot{\widetilde a}(r)\|_{\ell^\infty_q(\Z^\nu)}
 &\le C\mu_\varepsilon
       \|\widetilde a(r)\|_{\ell^\infty_q(\Z^\nu)}^3
 \le C\mu_\varepsilon R^3.
 \label{eq:log-low-derivative}
\end{align}

\par\medskip\noindent\emph{Differentiated Cubic Terms.}
For the differentiated term in \eqref{eq:log-normal-form-identity},
the \emph{product rule} gives
\begin{align*}
 \frac{d}{dr}\bigl(
 \widetilde a_{n_1}\overline{\widetilde a_{n_2}}\widetilde a_{n_3}\bigr)
 ={}&\dot{\widetilde a}_{n_1}
       \overline{\widetilde a_{n_2}}\widetilde a_{n_3}
 +\widetilde a_{n_1}
       \overline{\dot{\widetilde a}_{n_2}}\widetilde a_{n_3}+\widetilde a_{n_1}
       \overline{\widetilde a_{n_2}}\dot{\widetilde a}_{n_3}.
\end{align*}

For the first term,
\eqref{eq:log-pointwise-divisor} gives
\begin{align}
 &\langle n\rangle^q
 \sum_{\substack{n_1-n_2+n_3=n\\\Sigma_3\ne0}}
 \frac{|\dot{\widetilde a}_{n_1}|
       |\widetilde a_{n_2}|\,|\widetilde a_{n_3}|}{|\Sigma_3|}\le C\sum_{j=1}^3\sum_{n_1-n_2+n_3=n}
 \langle n_j\rangle^p
 |\dot{\widetilde a}_{n_1}|
 |\widetilde a_{n_2}|\,|\widetilde a_{n_3}|.\label{eq:log-differentiated-weight-split}
\end{align}
For $j=1$, solve
$n_1=n+n_2-n_3$ and bound the differentiated factor in the
high weight norm:
\begin{align*}
 &\sum_{n_1-n_2+n_3=n}\langle n_1\rangle^p
 |\dot{\widetilde a}_{n_1}|
 |\widetilde a_{n_2}|\,|\widetilde a_{n_3}|\le\|\dot{\widetilde a}(r)\|_{\ell^\infty_p(\Z^\nu)}
       \|\widetilde a(r)\|_{\ell^\infty_q(\Z^\nu)}^2
       \sum_{n_2,n_3}\langle n_2\rangle^{-q}\langle n_3\rangle^{-q}.
\end{align*}
For $j=2$, solve $n_2=n_1+n_3-n$ and place the high weight on
the second factor; the differentiated factor is then controlled in
$\ell^\infty_q$:
\begin{align*}
 \sum_{n_1-n_2+n_3=n}\langle n_2\rangle^p
 |\dot{\widetilde a}_{n_1}|
 |\widetilde a_{n_2}|\,|\widetilde a_{n_3}|\le\|\widetilde a(r)\|_{\ell^\infty_p}
       \|\dot{\widetilde a}(r)\|_{\ell^\infty_q}
       \|\widetilde a(r)\|_{\ell^\infty_q}
       \sum_{n_1,n_3}\langle n_1\rangle^{-q}\langle n_3\rangle^{-q}.
\end{align*}
For $j=3$, solving $n_3=n-n_1+n_2$ gives the same bound.
Each double sum equals
$\|\langle\cdot\rangle^{-q}\|_{\ell^1(\Z^\nu)}^2<\infty$ since $q>\nu$.
Substituting these three estimates into
\eqref{eq:log-differentiated-weight-split} therefore gives
\begin{equation}\label{eq:log-first-amplitude-derivative}
\begin{aligned}[b]
 &\langle n\rangle^q
 \sum_{\substack{n_1-n_2+n_3=n\\\Sigma_3\ne0}}
 \frac{|\dot{\widetilde a}_{n_1}|
       |\widetilde a_{n_2}|\,|\widetilde a_{n_3}|}{|\Sigma_3|}\\
&\quad\le C\|\langle\cdot\rangle^{-q}\|_{\ell^1(\Z^\nu)}^2
 \left[\begin{aligned}
 &\|\dot{\widetilde a}(r)\|_{\ell^\infty_p(\Z^\nu)}
       \|\widetilde a(r)\|_{\ell^\infty_q(\Z^\nu)}^2+{}\\
 &2\|\widetilde a(r)\|_{\ell^\infty_p(\Z^\nu)}
       \|\dot{\widetilde a}(r)\|_{\ell^\infty_q(\Z^\nu)}
       \|\widetilde a(r)\|_{\ell^\infty_q(\Z^\nu)}
 \end{aligned}\right].
\end{aligned}
\end{equation}

The other two derivative placements have the same bound.
Summing the three contributions and taking the weighted supremum yield
\begin{align}
 &\sup_n\langle n\rangle^q
 \sum_{\substack{n_1-n_2+n_3=n\\\Sigma_3\ne0}}
 \frac{\left|\dfrac{d}{dr}
 (\widetilde a_{n_1}\overline{\widetilde a_{n_2}}\widetilde a_{n_3})\right|}
      {|\Sigma_3|}\nonumber\\
 &\quad\le C\Bigl(
 \|\dot{\widetilde a}(r)\|_{\ell^\infty_p(\Z^\nu)}
 \|\widetilde a(r)\|_{\ell^\infty_q(\Z^\nu)}^2
 +\|\widetilde a(r)\|_{\ell^\infty_p(\Z^\nu)}
  \|\widetilde a(r)\|_{\ell^\infty_q(\Z^\nu)}
  \|\dot{\widetilde a}(r)\|_{\ell^\infty_q(\Z^\nu)}\Bigr)\nonumber\\
 &\quad\le C\Bigl[
 (\mu_\varepsilon R^3e^{CR^2\mu_\varepsilon r})(2R)^2
 +(Re^{CR^2\mu_\varepsilon r})(2R)(\mu_\varepsilon R^3)
 \Bigr]\nonumber\\
 &\quad\le C\mu_\varepsilon R^5e^{CR^2\mu_\varepsilon r}.\label{eq:log-amplitude-tame}
\end{align}
The last step uses \eqref{eq:log-bootstrap-bound},
\eqref{eq:log-high-envelope}, \eqref{eq:log-high-derivative} and
\eqref{eq:log-low-derivative}.

\par\medskip\noindent\textbf{Step 6. Comparison Estimate.}
Taking the weighted supremum in \eqref{eq:log-normal-form-identity}
and substituting \eqref{eq:log-diagonal-difference},
\eqref{eq:log-boundary-bound} and \eqref{eq:log-amplitude-tame}, we obtain
\begin{align}
 E_{\log}(t)
 &\le12R^2\mu_\varepsilon\int_0^tE_{\log}(r)\,dr
    +CR^3\mu_\varepsilon e^{CR^2\mu_\varepsilon t}
    +C\mu_\varepsilon^2R^5\int_0^t e^{CR^2\mu_\varepsilon r}\,dr\nonumber\\
 &\le CR^2\mu_\varepsilon\int_0^tE_{\log}(r)\,dr
    +CR^3\mu_\varepsilon(1+R^2\mu_\varepsilon t)
                         e^{CR^2\mu_\varepsilon t}.\label{eq:log-error-integral}
\end{align}
Since $R$ is fixed and $1+R^2\mu_\varepsilon t\le
 e^{R^2\mu_\varepsilon t}$, enlarging $C$ reduces this to
$E_{\log}(t)\le C\mu_\varepsilon e^{C\mu_\varepsilon t}
+C\mu_\varepsilon\int_0^t E_{\log}(r)\,dr$.
The \emph{integral form of Gronwall's inequality} therefore gives
\begin{align*}
 E_{\log}(t)
 &\le C\mu_\varepsilon e^{C\mu_\varepsilon t}
 +C\mu_\varepsilon\int_0^t e^{C\mu_\varepsilon(t-r)}
                  C\mu_\varepsilon e^{C\mu_\varepsilon r}\,dr\\
 &\le C\mu_\varepsilon(1+C\mu_\varepsilon t)
                         e^{C\mu_\varepsilon t}.
\end{align*}
Absorbing the remaining factor into the exponential proves
\eqref{eq:log-error-growth}.
\end{proof}

\subsubsection{Logarithmic Extension and Physical-Space Comparison}
\label{subsubsec:log-bootstrap-closure}
\begin{proof}[Proof of \autoref{prop:log-pointwise-comparison}]
\par\medskip\noindent\textbf{Step 1. Bootstrap Interval.}
Let $\widetilde a$ be the maximal extension of the local solution from
\autoref{prop:log-local-wellposedness}, with lifespan $T_{\max}$ as in
\autoref{rem:log-blowup-alternative}. Since $p\ge q$, continuity in
$\ell^\infty_p$ implies continuity in $\ell^\infty_q$. Together with
\[
 \|\widetilde a(0)\|_{\ell^\infty_q}
 \le\|\widetilde d\|_{\ell^\infty_p}\le R<2R,
\]
this ensures that the strict low weight bound holds on a positive interval.
Define 
\begin{equation}\label{eq:log-bootstrap-time}
 T_*:=\sup\left\{
 t\in\bigl(0,\min\{T_{\max},\mathscr T_\varepsilon\}\bigr):
 \sup_{0\le r\le t}
 \|\widetilde a(r)\|_{\ell^\infty_q}<2R
 \right\}.
\end{equation}
The defining set is nonempty, so
$0<T_*\le\min\{T_{\max},\mathscr T_\varepsilon\}$.

\par\medskip\noindent\textbf{Step 2. Improvement and High Weight Control.}
Choose $\vartheta\in(0,1)$ and set
$c=(1-\vartheta)/\Lambda$ in \eqref{eq:log-observation-time},
where $C,\Lambda$ are the constants in
\autoref{lem:log-bootstrap-comparison}.
The restriction $\Lambda c<1$ ensures that the comparison error vanishes
on the logarithmic interval; an explicit choice of $\Lambda$ is given
in \autoref{rem:log-explicit-growth-constant}.
For sufficiently small $\varepsilon$,
\begin{equation}\label{eq:log-bootstrap-smallness}
 0<\mu_\varepsilon<1,\qquad
 C\mu_\varepsilon^\vartheta\le\tfrac12R.
\end{equation}

Fix $0<t<T_*$. Applying
\autoref{lem:log-bootstrap-comparison} on $[0,t]$ and using
$t<\mathscr T_\varepsilon=c\mu_\varepsilon^{-1}|\log\mu_\varepsilon|$,
we obtain from \eqref{eq:log-error-growth}, for every $0\le r\le t$,
\begin{align}
 E_{\log}(r)
 &\le C\mu_\varepsilon e^{\Lambda\mu_\varepsilon r}
 \le C\mu_\varepsilon e^{\Lambda c|\log\mu_\varepsilon|}\nonumber\\
 &=C\mu_\varepsilon^{1-\Lambda c}
 =C\mu_\varepsilon^\vartheta.
 \label{eq:log-bootstrap-improvement}
\end{align}
By \eqref{eq:log-resonant-norm}, the resonant low weight norm is at most $R$.
The triangle inequality and \eqref{eq:log-bootstrap-smallness} therefore give
\begin{align}
 \sup_{0\le r\le t}\|\widetilde a(r)\|_{\ell^\infty_q}
 &\le\sup_{0\le r\le t}\|\widetilde b(r)\|_{\ell^\infty_q}
       +\sup_{0\le r\le t}E_{\log}(r)\nonumber\\
 &\le R+C\mu_\varepsilon^\vartheta
 \le\tfrac32R<2R.
 \label{eq:log-bootstrap-low-improvement}
\end{align}

Under the same bootstrap hypothesis, \eqref{eq:log-high-envelope} gives
\begin{align}
 \sup_{0\le r\le t}\|\widetilde a(r)\|_{\ell^\infty_p}
 &\le R\exp(C_{p,q}R^2\mu_\varepsilon t)\nonumber\\
 &\le R\exp(C_{p,q}R^2\mu_\varepsilon\mathscr T_\varepsilon)
 =R\mu_\varepsilon^{-C_{p,q}R^2c}<\infty.
 \label{eq:log-bootstrap-high-control}
\end{align}
This bound is uniform in $t<T_*$ and finite for each fixed $\varepsilon$,
as required by \autoref{lem:log-continuation}.

\par\medskip\noindent\textbf{Step 3. Closing the Bootstrap.}
If $T_{\max}=T_*$, the high weight bound
\eqref{eq:log-bootstrap-high-control} and
\autoref{lem:log-continuation} extend the solution beyond $T_*$,
contradicting maximality. Hence $T_*<T_{\max}$. Passing to the endpoint in
\eqref{eq:log-bootstrap-low-improvement} yields
\[
 \sup_{0\le r\le T_*}\|\widetilde a(r)\|_{\ell^\infty_q}
 \le\tfrac32R<2R.
\]

Suppose that $T_*<\mathscr T_\varepsilon$. Since $T_*<T_{\max}$,
the solution is defined on a neighborhood of $T_*$. The local construction
in \autoref{prop:log-local-wellposedness} and successive extensions by
\autoref{lem:log-continuation} give continuity in the high weight
norm. Since $p>q$, this implies
\[
 \|\widetilde a(r)-\widetilde a(T_*)\|_{\ell^\infty_q}
 \le \|\widetilde a(r)-\widetilde a(T_*)\|_{\ell^\infty_p}
 \longrightarrow0\qquad(r\to T_*).
\]
We may therefore choose $h>0$ sufficiently small that
\[
 T_*+h<\min\{T_{\max},\mathscr T_\varepsilon\},
 \qquad
 \sup_{T_*\le r\le T_*+h}
 \|\widetilde a(r)-\widetilde a(T_*)\|_{\ell^\infty_q}
 <\tfrac14R.
\]
For every $r\in[T_*,T_*+h]$, the triangle inequality gives
\[
 \|\widetilde a(r)\|_{\ell^\infty_q}
 \le \|\widetilde a(T_*)\|_{\ell^\infty_q}
    +\|\widetilde a(r)-\widetilde a(T_*)\|_{\ell^\infty_q}
 <\tfrac32R+\tfrac14R=\tfrac74R.
\]
Combining this with the bound on $[0,T_*]$, we obtain
\[
 \sup_{0\le r\le T_*+h}\|\widetilde a(r)\|_{\ell^\infty_q}
 \le\tfrac74R<2R.
\]
Thus $T_*+h$ belongs to the set whose supremum defines $T_*$ in
\eqref{eq:log-bootstrap-time}, contradicting the definition of $T_*$. Therefore
\begin{equation}\label{eq:log-bootstrap-closure}
 T_*=\mathscr T_\varepsilon<T_{\max}.
\end{equation}

The solution therefore exists on $[0,\mathscr T_\varepsilon]$.
Estimate \eqref{eq:log-high-envelope} gives
\eqref{eq:log-pointwise-high-decay} and continuity extends
\eqref{eq:log-bootstrap-low-improvement} and
\eqref{eq:log-bootstrap-improvement} to its endpoint, proving
\eqref{eq:log-low-weight-bound} and
\eqref{eq:log-pointwise-comparison}.

\par\medskip\noindent\textbf{Step 4. Physical-Space Comparison.}
\phantomsection\label{subsubsec:log-physical-comparison}
It remains to pass from the coefficient estimate to the physical field.
Using \eqref{eq:log-normalization}, the definition
\eqref{eq:log-comparison-envelope} and the bound
\eqref{eq:log-bootstrap-improvement} extended to $t=\mathscr T_\varepsilon$,
we obtain
\begin{align*}
 \sqrt\varepsilon\,\|u(t)-u_{\mathrm{res}}(t)\|_\infty
 &\le\sum_{n\in\Z^\nu}|\widetilde a_n(t)-\widetilde b_n(t)|\\
 &\le E_{\log}(t)\sum_{n\in\Z^\nu}\langle n\rangle^{-q}\\
 &\le C\mu_\varepsilon^\vartheta
 \|\langle\cdot\rangle^{-q}\|_{\ell^1(\Z^\nu)}.
\end{align*}
Taking the
supremum over $[0,\mathscr T_\varepsilon]$ proves
\eqref{eq:log-physical-comparison}.
\end{proof}

\subsection{LDP at the Critical Logarithmic Correction Time with a Small Enough Prefactor}
\label{subsec:log-ldp-transfer}
By \eqref{eq:critical-pointwise-lifespan},
\eqref{eq:nf-high-mode-delta} and \eqref{eq:log-observation-time}, we have
\begin{equation}\label{eq:log-time-interval-inclusion}
 \frac{\mathscr T_\varepsilon}{T_\varepsilon}
 \asymp|\log\mu_\varepsilon|\longrightarrow\infty, \qquad \varepsilon\to0^+.
 \qquad
\end{equation}
Thus $ [0,T_\varepsilon]\subset[0,\mathscr T_\varepsilon]$ for sufficiently small $\varepsilon$. We now extend \autoref{thm:first-critical-ldp-local} to this logarithmic interval.

\begin{theorem}[LDP at the Critical Logarithmic Correction Time with a Small Enough Prefactor]
\label{thm:log-critical-ldp}
Assume the hypotheses of \autoref{prop:log-pointwise-comparison}.
Let $u$ be the nonlinear field supplied by that proposition on
$G_{\delta_\varepsilon}$ for $0\le t\le\mathscr T_\varepsilon$. Then
\begin{equation}\label{eq:log-critical-ldp}
\begin{aligned}[t]
 \lim_{\varepsilon\to0^+}\sup_{0\le t\le\mathscr T_\varepsilon}
 \varepsilon\log\Pp\!\left(
 \mathscr A_{\varepsilon,z}^{\|u(t)\|_\infty}\right)
 =\lim_{\varepsilon\to0^+}\inf_{0\le t\le\mathscr T_\varepsilon}
 \varepsilon\log\Pp\!\left(
 \mathscr A_{\varepsilon,z}^{\|u(t)\|_\infty}\right)
 =-\frac{z^2}{S}.
\end{aligned}
\end{equation}

\end{theorem}

\begin{proof}[Proof of \autoref{thm:log-critical-ldp}]
\par\medskip\noindent\textbf{Step 1. Event Inclusions.}
Fix a deterministic $t\in[0,\mathscr T_\varepsilon]$. The estimate
\eqref{eq:log-physical-comparison} gives
\begin{equation}\label{eq:log-comparison-event-inclusion}
 G_{\delta_\varepsilon}\subseteq
 \mathscr B_{\varepsilon,C\mu_\varepsilon^\vartheta}
 ^{\|u(t)-u_{\mathrm{res}}(t)\|_\infty},
 \qquad
 \mathscr A_{\varepsilon,C\mu_\varepsilon^\vartheta}
 ^{\|u(t)-u_{\mathrm{res}}(t)\|_\infty}
 \subseteq\Omega_{\delta_\varepsilon}.
\end{equation}
On $G_{\delta_\varepsilon}$, the reverse triangle inequality yields
\[
 \sqrt\varepsilon\,\left|\|u(t)\|_\infty
       -\|u_{\mathrm{res}}(t)\|_\infty\right|
 \le C\mu_\varepsilon^\vartheta.
\]
For sufficiently small $\varepsilon$, $C\mu_\varepsilon^\vartheta<z$.
Using the event notation in \eqref{eq:nonlinear-spatial-exceedance-event}
and \eqref{eq:spatial-exceedance-complement}, we obtain
\begin{align}
 G_{\delta_\varepsilon}\cap
 \mathscr A_{\varepsilon,z+C\mu_\varepsilon^\vartheta}
 ^{\|u_{\mathrm{res}}(t)\|_\infty}
 \subseteq G_{\delta_\varepsilon}\cap
 \mathscr A_{\varepsilon,z}^{\|u(t)\|_\infty}\subseteq G_{\delta_\varepsilon}\cap
 \mathscr A_{\varepsilon,z-C\mu_\varepsilon^\vartheta}
 ^{\|u_{\mathrm{res}}(t)\|_\infty}.
 \label{eq:log-good-tail-inclusions}
\end{align}

\par\medskip\noindent\textbf{Step 2. Probability Bounds.}
As in \autoref{thm:first-critical-ldp-local}, we split the exceedance event
over $G_{\delta_\varepsilon}$ and $\Omega_{\delta_\varepsilon}$. The \emph{total
probability formula} gives
\begin{align}
 \Pp\!\left(\mathscr A_{\varepsilon,z}^{\|u(t)\|_\infty}\right)
 ={}&\Pp\!\left(\mathscr A_{\varepsilon,z}^{\|u(t)\|_\infty}
          \mid G_{\delta_\varepsilon}\right)\Pp(G_{\delta_\varepsilon})+\Pp\!\left(\mathscr A_{\varepsilon,z}^{\|u(t)\|_\infty}
          \mid\Omega_{\delta_\varepsilon}\right)
          \Pp(\Omega_{\delta_\varepsilon})\nonumber\\
 ={}&\Pp\!\left(\mathscr A_{\varepsilon,z}^{\|u(t)\|_\infty}
          \cap G_{\delta_\varepsilon}\right)
 +\Pp\!\left(\mathscr A_{\varepsilon,z}^{\|u(t)\|_\infty}
          \cap\Omega_{\delta_\varepsilon}\right).
 \label{eq:log-total-probability}
\end{align}
The second term in \eqref{eq:log-total-probability} satisfies
\[
 0\le\Pp\!\left(\mathscr A_{\varepsilon,z}^{\|u(t)\|_\infty}
       \cap\Omega_{\delta_\varepsilon}\right)
 \le\Pp(\Omega_{\delta_\varepsilon}).
\]
Thus \eqref{eq:log-total-probability} and
\eqref{eq:log-good-tail-inclusions} give the upper bound
\begin{align*}
 \Pp\!\left(\mathscr A_{\varepsilon,z}^{\|u(t)\|_\infty}\right)
 &\le\Pp\!\left(G_{\delta_\varepsilon}\cap
 \mathscr A_{\varepsilon,z-C\mu_\varepsilon^\vartheta}
 ^{\|u_{\mathrm{res}}(t)\|_\infty}\right)
 +\Pp(\Omega_{\delta_\varepsilon})\\
 &\le p_\varepsilon(z-C\mu_\varepsilon^\vartheta)
 +\Pp(\Omega_{\delta_\varepsilon}),
\end{align*}
where $p_\varepsilon$ is defined in \eqref{eq:gaussian-tail-probability};
we used \autoref{prop:nf-resonant-law} in the last step.
For the lower bound, the same event inclusions give
\begin{align*}
 \Pp\!\left(\mathscr A_{\varepsilon,z}^{\|u(t)\|_\infty}\right)
 &\ge\Pp\!\left(G_{\delta_\varepsilon}\cap
 \mathscr A_{\varepsilon,z+C\mu_\varepsilon^\vartheta}
 ^{\|u_{\mathrm{res}}(t)\|_\infty}\right)\\
 &=p_\varepsilon(z+C\mu_\varepsilon^\vartheta)
 -\Pp\!\left(\Omega_{\delta_\varepsilon}\cap
 \mathscr A_{\varepsilon,z+C\mu_\varepsilon^\vartheta}
 ^{\|u_{\mathrm{res}}(t)\|_\infty}\right)\\
 &\ge p_\varepsilon(z+C\mu_\varepsilon^\vartheta)
 -\Pp(\Omega_{\delta_\varepsilon}).
\end{align*}
 Consequently,
\begin{align}
 p_\varepsilon(z+C\mu_\varepsilon^\vartheta)
 -\Pp(\Omega_{\delta_\varepsilon})
 \le\Pp\!\left(\mathscr A_{\varepsilon,z}^{\|u(t)\|_\infty}\right)
 \le p_\varepsilon(z-C\mu_\varepsilon^\vartheta)
 +\Pp(\Omega_{\delta_\varepsilon}).
 \label{eq:log-tail-probability-sandwich}
\end{align}
Factoring out the resonant main terms in
\eqref{eq:log-tail-probability-sandwich}, we obtain
\begin{align*}
 p_\varepsilon(z+C\mu_\varepsilon^\vartheta)
 \left[1-\frac{\Pp(\Omega_{\delta_\varepsilon})}
 {p_\varepsilon(z+C\mu_\varepsilon^\vartheta)}\right]
 &\le\Pp\!\left(\mathscr A_{\varepsilon,z}^{\|u(t)\|_\infty}\right)\\
 &\le p_\varepsilon(z-C\mu_\varepsilon^\vartheta)
 \left[1+\frac{\Pp(\Omega_{\delta_\varepsilon})}
 {p_\varepsilon(z-C\mu_\varepsilon^\vartheta)}\right].
\end{align*}
For sufficiently small $\varepsilon$, the lower bracket is positive,
as follows from \eqref{eq:log-relative-exception-small} below.
Taking logarithms and multiplying by $\varepsilon$ give
\begin{align}
 &\varepsilon\log p_\varepsilon(z+C\mu_\varepsilon^\vartheta)
 +\varepsilon\log\left(1-
 \frac{\Pp(\Omega_{\delta_\varepsilon})}
 {p_\varepsilon(z+C\mu_\varepsilon^\vartheta)}\right)\nonumber\\
 &\quad\le\varepsilon\log\Pp\!\left(
 \mathscr A_{\varepsilon,z}^{\|u(t)\|_\infty}\right)\nonumber\\
 &\quad\le\varepsilon\log p_\varepsilon(z-C\mu_\varepsilon^\vartheta)
 +\varepsilon\log\left(1+
 \frac{\Pp(\Omega_{\delta_\varepsilon})}
 {p_\varepsilon(z-C\mu_\varepsilon^\vartheta)}\right).
 \label{eq:log-logarithmic-sandwich}
\end{align}

\par\medskip\noindent\textbf{Step 3. Limits of the Main and Correction Terms.}
Since $C\mu_\varepsilon^\vartheta\to0$, the moving-threshold argument in
\eqref{eq:critical-moving-threshold-ldp} applies and gives
\begin{equation}\label{eq:log-moving-threshold-ldp}
 \varepsilon\log p_\varepsilon(z\pm C\mu_\varepsilon^\vartheta)
 \longrightarrow-\frac{z^2}{S}.
\end{equation}
Indeed, for every fixed $0<h<z$ and sufficiently small $\varepsilon$,
\[
 p_\varepsilon(z+h)\le
 p_\varepsilon(z\pm C\mu_\varepsilon^\vartheta)
 \le p_\varepsilon(z-h),
\]
and \autoref{thm:nf-resonant-LDP}, followed by $h\downarrow0$, gives
\eqref{eq:log-moving-threshold-ldp}. Also,
\autoref{prop:exception-prob} and \eqref{eq:nf-high-mode-delta} imply
\[
 \Pp(\Omega_{\delta_\varepsilon})
 \le C_{\kappa,\nu}\exp\!\left(-\frac{1+z^2/S}{\varepsilon}\right).
\]
Thus, for either sign,
\begin{align*}
 \varepsilon\log
 \frac{\Pp(\Omega_{\delta_\varepsilon})}
 {p_\varepsilon(z\pm C\mu_\varepsilon^\vartheta)}
 &\le\varepsilon\log C_{\kappa,\nu}-1-\frac{z^2}{S}
 -\varepsilon\log p_\varepsilon(z\pm C\mu_\varepsilon^\vartheta)\\
 &=-1+o(1).
\end{align*}
For sufficiently small $\varepsilon$, it follows that
\begin{equation}\label{eq:log-relative-exception-small}
 0\le\frac{\Pp(\Omega_{\delta_\varepsilon})}
 {p_\varepsilon(z\pm C\mu_\varepsilon^\vartheta)}
 \le e^{-1/(2\varepsilon)}.
\end{equation}

This proves the positivity required above. It also controls both
logarithmic corrections in \eqref{eq:log-logarithmic-sandwich}.
Using $-2x\le\log(1-x)\le0$ and $0\le\log(1+x)\le x$ for
$0\le x\le1/2$, as in \eqref{eq:critical-logarithmic-corrections},
we obtain
\begin{align}
 -2\varepsilon e^{-1/(2\varepsilon)}
 &\le\varepsilon\log\left(1-
 \frac{\Pp(\Omega_{\delta_\varepsilon})}
 {p_\varepsilon(z+C\mu_\varepsilon^\vartheta)}\right)
 \le0,\nonumber\\
 0&\le\varepsilon\log\left(1+
 \frac{\Pp(\Omega_{\delta_\varepsilon})}
 {p_\varepsilon(z-C\mu_\varepsilon^\vartheta)}\right)
 \le\varepsilon e^{-1/(2\varepsilon)}.
 \label{eq:log-logarithmic-corrections}
\end{align}
Thus both logarithmic correction terms in \eqref{eq:log-logarithmic-sandwich} tend to zero,
while both main terms tend to $-z^2/S$ by
\eqref{eq:log-moving-threshold-ldp}. These bounds are independent of
$t\in[0,\mathscr T_\varepsilon]$. Taking the time infimum and supremum
in \eqref{eq:log-logarithmic-sandwich} proves \eqref{eq:log-critical-ldp}.
\end{proof}

\section{Failure of the Initial Gaussian Rate Function at the Critical Time Scale with a Large Logarithmic Correction Prefactor: Special Case}
\label{sec:longtime-obstruction}

The logarithmic extension in \autoref{thm:log-critical-ldp} requires a
small enough prefactor; see also \cite[Remark~1.2]{FY25} for the periodic
problem. We give an example in which the
initial Gaussian rate function fails along a sequence of times
$t_j=15\mu_j^{-1}|\log\mu_j|$.
To this end, we first construct rescaled Akhmediev breathers by
\emph{small frequency embedding into a fixed quasi-periodic module};
see \autoref{lem:embedded-breather-properties}.
Fourier coefficient estimates then show that the \emph{Gaussian cost}
of their initial data tends to one; see \eqref{eq:precursor-gaussian-cost}.
Next, \emph{deterministic tube stability} shows that all solutions with
initial data in the chosen neighborhoods reach normalized amplitudes above two;
see \autoref{lem:breather-tube}.
To estimate the probability of these neighborhoods, we combine
\emph{finite-dimensional Gaussian density calculations} with a
centered ball estimate; see \autoref{lem:gaussian-tube} and
\eqref{eq:counterexample-centered-tube}.
Finally, the \emph{cost comparison} $1<4/S$ proves failure of the initial
Gaussian rate function and rules out uniform persistence on any time scale
polynomially longer than the critical scale; see
\autoref{thm:longtime-rate-obstruction} and
\autoref{cor:counterexample-polynomial-obstruction}.

\subsection{Motivation: Akhmediev Breather}
Consider the following \emph{Akhmediev breather}, which is spatially periodic and temporally localized:
\begin{equation}\label{eq:explicit-breather}
 u_{\mathrm{AB}}(t,x)=e^{\iu t}
 \frac{h\cos x+\iu\sinh t}{h\cos x-\cosh t},
 \qquad (t,x)\in\R\times\R,
\end{equation}
where $h=2^{-1/2}$; see \cite[formula~(39)]{AEK87}, with the amplitude adjusted to the normalization below.

\begin{lemma}[Properties of the Akhmediev Breather]\label{lem:breather-identities}
Let $u_{\mathrm{AB}}$ be given by \eqref{eq:explicit-breather}, with Fourier coefficients
\[
 \widehat u_{\mathrm{AB}}(t,m):=\frac1{2\pi}\int_0^{2\pi}u_{\mathrm{AB}}(t,x)e^{-\iu mx}\,dx,
 \qquad m\in\Z.
\]
Then the following properties hold.
\begin{enumerate}[label=\textup{(\roman*)},leftmargin=*,itemsep=3pt]
\item The function $u_{\mathrm{AB}}$ is smooth on $\R^2$, $2\pi$-periodic in $x$ and solves the cubic NLS
\begin{equation}\label{eq:breather-NLS}
 \iu\partial_t u_{\mathrm{AB}}+\partial_{xx}u_{\mathrm{AB}}+|u_{\mathrm{AB}}|^2u_{\mathrm{AB}}=0.
\end{equation}
\item As $T\to+\infty$, the zero and
nonzero modes satisfy
\begin{align}
 |\widehat u_{\mathrm{AB}}(-T,0)|&=1+\mathcal{O}(e^{-2T}),
 \label{eq:breather-initial-zero-mode}\\
 |\widehat u_{\mathrm{AB}}(-T,m)|&\le\sqrt2(\sqrt2e^{-T})^{|m|},\qquad m\ne0.
 \label{eq:breather-initial-Fourier}
\end{align}
\item The nonzero modes and the full coefficient sequence satisfy
\begin{align}
 \sum_{m\ne0}|\widehat u_{\mathrm{AB}}(-T,m)|&=\mathcal{O}(e^{-T}),
 \label{eq:breather-initial-nonzero-sum}\\
 \|\widehat u_{\mathrm{AB}}(-T)\|_{\ell^1}&=1+\mathcal{O}(e^{-T}).
 \label{eq:breather-initial-fourier-norm}
\end{align}
\item The amplitude approaches the background uniformly, that is, 
\begin{equation}\label{eq:breather-uniform-background}
 \sup_{x\in\R}\bigl||u_{\mathrm{AB}}(t,x)|-1\bigr|
 =\mathcal{O}(e^{-|t|}),\qquad |t|\to\infty.
\end{equation}
\item At $(t,x)=(0,0)$, one has
\begin{equation}\label{eq:breather-size}
 |u_{\mathrm{AB}}(0,0)|=1+\sqrt2>2.
\end{equation}
\item The Fourier coefficients obey the uniform bound 
\begin{equation}\label{eq:breather-uniform-fourier-bound}
 \sup_{t\in\R}\|\widehat u_{\mathrm{AB}}(t)\|_{\ell^1(\Z)}\le3+\sqrt2.
\end{equation}
\end{enumerate}
\end{lemma}
\begin{proof}[Proof of \autoref{lem:breather-identities}]
\smallskip
\noindent\textup{(i)}
Write
\[
 u_{\mathrm{AB}}=e^{\iu t}F,\qquad F=\frac{y+\iu\sinh t}{D},\qquad
 y=h\cos x,\qquad D=y-\cosh t.
\]
Since $D\le h-1<0$ and the numerator and denominator are smooth,
$u_{\mathrm{AB}}$ is smooth on $\R^2$. Its spatial periodicity follows from that of $\cos x$.
For $F$, we have
\[
 \iu \partial_t u_{\mathrm{AB}}+\partial_{xx}u_{\mathrm{AB}}+|u_{\mathrm{AB}}|^2u_{\mathrm{AB}}
 =e^{\iu t}\bigl(\iu F_t+F_{xx}+(|F|^2-1)F\bigr).
\]
It remains to verify that the expression in parentheses vanishes.
Differentiating and using $y_x^2=h^2-y^2$, $y_{xx}=-y$, $2h^2=1$ and $\cosh^2t-\sinh^2t=1$ give
\begin{align*}
 \iu F_t
 &=\frac{1-y\cosh t+\iu y\sinh t}{D^2},\\
 F_x&=-\frac{(\cosh t+\iu\sinh t)y_x}{D^2},\\
 F_{xx}
 &
 =-\frac{(\cosh t+\iu\sinh t)(y^2+y\cosh t-1)}{D^3},\\
(|F|^2-1)F
  &=\frac{(2y\cosh t-1)(y+\iu\sinh t)}{D^3}.
\end{align*}
Consequently,
\begin{align*}
 D^3\bigl(\iu F_t+F_{xx}+(|F|^2-1)F\bigr)
 =0.
\end{align*}
Thus $u_{\mathrm{AB}}$ in \eqref{eq:explicit-breather} solves the cubic NLS \eqref{eq:breather-NLS}.

\smallskip
\textup{(ii)}
Separating the constant term in \eqref{eq:explicit-breather} gives
\[
 u_{\mathrm{AB}}(t,x)=e^{\iu t}\left(1-
 \frac{\cosh t+\iu\sinh t}{\cosh t-h\cos x}\right). 
\]
By the definition of the Fourier coefficients, for $m\in\Z$,
\begin{align}
 e^{-\iu t}\widehat u_{\mathrm{AB}}(t,m)
 ={}&\frac1{2\pi}\int_0^{2\pi}e^{-\iu mx}\,dx\nonumber\\
 &\quad-(\cosh t+\iu\sinh t)
 \frac1{2\pi}\int_0^{2\pi}
 \frac{e^{-\iu mx}}{\cosh t-h\cos x}\,dx.
 \label{eq:breather-fourier-integral}
\end{align}
The first integral equals $1$ if $m=0$ and $0$ otherwise.
To evaluate the second integral in \eqref{eq:breather-fourier-integral},
we first use $2\cos x=e^{\iu x}+e^{-\iu x}$ to rewrite its denominator:
\begin{align}
 \cosh t-h\cos x
 &=\cosh t-\frac h2\left(e^{\iu x}+e^{-\iu x}\right)\nonumber\\
 &=-\frac{h}{2e^{\iu x}}
   \left(e^{2\iu x}-\frac{2\cosh t}{h}e^{\iu x}+1\right).
 \label{eq:breather-denominator-quadratic}
\end{align}
For fixed $t$, the quadratic polynomial in $e^{\iu x}$ in
\eqref{eq:breather-denominator-quadratic} has roots
$(\cosh t\pm\sqrt{\cosh^2t-h^2})/h$, whose product is $1$.
Denote the smaller root by
\begin{equation}\label{eq:breather-smaller-root}
 r(t)=\frac{\cosh t-\sqrt{\cosh^2t-h^2}}h
 =\frac{h}{\cosh t+\sqrt{\cosh^2t-h^2}}\in(0,1).
\end{equation}
The other root is $r(t)^{-1}$. Factoring
\eqref{eq:breather-denominator-quadratic} therefore gives
\begin{align}
 \cosh t-h\cos x
 &=-\frac{h}{2e^{\iu x}}
   (e^{\iu x}-r(t))(e^{\iu x}-r(t)^{-1})\nonumber\\
 &=\frac{h}{2r(t)}
   (1-r(t)e^{\iu x})(1-r(t)e^{-\iu x})\nonumber\\
 &=\frac{h}{2r(t)}(1-2r(t)\cos x+r(t)^2).
 \label{eq:breather-denominator-factorization}
\end{align}
The coefficient in \eqref{eq:breather-denominator-factorization}
can be computed directly from \eqref{eq:breather-smaller-root}:
\begin{align}
 \frac{h}{2r(t)}
 &=\frac{hr(t)^{-1}-hr(t)}{2(1-r(t)^2)}\nonumber\\
 &=\frac{\left(\cosh t+\sqrt{\cosh^2t-h^2}\right)
         -\left(\cosh t-\sqrt{\cosh^2t-h^2}\right)}
        {2(1-r(t)^2)}\nonumber\\
 &=\frac{\sqrt{\cosh^2t-h^2}}{1-r(t)^2}.
 \label{eq:breather-factor-coefficient}
\end{align}
Substituting \eqref{eq:breather-factor-coefficient} into
\eqref{eq:breather-denominator-factorization} and taking reciprocals,
we express the reciprocal of the denominator as a constant
multiple of the \emph{Poisson kernel}
\begin{equation}\label{eq:breather-poisson-kernel}
 \frac1{\cosh t-h\cos x}
 =\frac1{\sqrt{\cosh^2t-h^2}}
   \frac{1-r(t)^2}{1-2r(t)\cos x+r(t)^2}.
\end{equation}
Using $2\cos x=e^{\iu x}+e^{-\iu x}$, we expand the Poisson kernel in \eqref{eq:breather-poisson-kernel} and obtain
\begin{align}
 \frac{1-r(t)^2}{1-2r(t)\cos x+r(t)^2}
 &=1+\frac{\bigl(r(t)e^{\iu x}-r(t)^2\bigr)
           +\bigl(r(t)e^{-\iu x}-r(t)^2\bigr)}
          {(1-r(t)e^{\iu x})(1-r(t)e^{-\iu x})}\nonumber\\
 &=1+\frac{r(t)e^{\iu x}(1-r(t)e^{-\iu x})
           +r(t)e^{-\iu x}(1-r(t)e^{\iu x})}
          {(1-r(t)e^{\iu x})(1-r(t)e^{-\iu x})}\nonumber\\
 &=1+\frac{r(t)e^{\iu x}}{1-r(t)e^{\iu x}}
     +\frac{r(t)e^{-\iu x}}{1-r(t)e^{-\iu x}}.
 \label{eq:breather-poisson-fraction-decomposition}
\end{align}
Since $0<r(t)<1$, each of the last two fractions in
\eqref{eq:breather-poisson-fraction-decomposition} has a convergent
geometric series expansion. Substituting this decomposition into
\eqref{eq:breather-poisson-kernel} gives
\begin{align}
 \frac1{\cosh t-h\cos x}
 &=\frac1{\sqrt{\cosh^2t-h^2}}
   \left(1+\frac{r(t)e^{\iu x}}{1-r(t)e^{\iu x}}
          +\frac{r(t)e^{-\iu x}}{1-r(t)e^{-\iu x}}\right)\nonumber\\
 &=\frac1{\sqrt{\cosh^2t-h^2}}
   \left(1+\sum_{j=1}^{\infty}r(t)^je^{\iu jx}
          +\sum_{j=1}^{\infty}r(t)^je^{-\iu jx}\right)\nonumber\\
 &=\frac1{\sqrt{\cosh^2t-h^2}}
   \sum_{j\in\Z}r(t)^{|j|}e^{\iu jx}.
 \label{eq:breather-poisson-fourier-series}
\end{align}
By \eqref{eq:breather-poisson-kernel}, the second integral in
\eqref{eq:breather-fourier-integral} is, up to the factor
$(\cosh^2t-h^2)^{-1/2}$, the Poisson integral of $e^{-\iu mx}$
at radius $r(t)$ and angle zero. Substituting
\eqref{eq:breather-poisson-fourier-series} and using orthogonality,
we evaluate the integral needed for the Fourier coefficients of $u_{\mathrm{AB}}$, that is, 
\begin{align}
 \frac1{2\pi}\int_0^{2\pi}
 \frac{e^{-\iu mx}}{\cosh t-h\cos x}\,dx
 &=\frac1{\sqrt{\cosh^2t-h^2}}
   \sum_{j\in\Z}r(t)^{|j|}
   \frac1{2\pi}\int_0^{2\pi}e^{\iu(j-m)x}\,dx\nonumber\\
 &=\frac{r(t)^{|m|}}{\sqrt{\cosh^2t-h^2}}.
 \label{eq:breather-poisson-mode-integral}
\end{align}
Substituting \eqref{eq:breather-poisson-mode-integral} into
\eqref{eq:breather-fourier-integral},
we obtain
\begin{empheq}[left={\widehat u_{\mathrm{AB}}(t,m)=\empheqlbrace}]{alignat=2}
 &e^{\iu t}\left(1-
   \frac{\cosh t+\iu\sinh t}{\sqrt{\cosh^2t-h^2}}\right),
 &\quad&m=0,
 \label{eq:breather-fourier-zero}\\
 &-e^{\iu t}\frac{\cosh t+\iu\sinh t}{\sqrt{\cosh^2t-h^2}}
   r(t)^{|m|},&\quad&m\ne0.
 \label{eq:breather-fourier-nonzero}
\end{empheq}

\smallskip
We now use the Fourier coefficient formulas
\eqref{eq:breather-fourier-zero} and \eqref{eq:breather-fourier-nonzero}
at $t=-T$ to estimate the zero and nonzero modes separately.

\noindent\emph{The zero mode.}
For every $t\in\R$, $\cosh^2t+\sinh^2t=2\cosh^2t-1$ and $h^2=1/2$ give
\begin{equation}\label{eq:breather-fourier-prefactor}
 \frac{|\cosh t+\iu\sinh t|}{\sqrt{\cosh^2t-h^2}}
=\sqrt2.
\end{equation}
By \eqref{eq:breather-fourier-zero}, \eqref{eq:breather-fourier-prefactor} at $t=-T$ and $|e^{-\iu T}|=1$,
\begin{align}
 |\widehat u_{\mathrm{AB}}(-T,0)|^2
 &=\left|1-\frac{\cosh T-\iu\sinh T}
                     {\sqrt{\cosh^2T-h^2}}\right|^2\nonumber\\
 &=\left(1-\frac{\cosh T}{\sqrt{\cosh^2T-h^2}}\right)^2
   +\frac{\sinh^2T}{\cosh^2T-h^2}\nonumber\\
 &=1-\frac{2\cosh T}{\sqrt{\cosh^2T-h^2}}
   +\frac{\cosh^2T+\sinh^2T}{\cosh^2T-h^2}\nonumber\\
&=3-2\left(1-\frac{h^2}{\cosh^2T}\right)^{-1/2}.
\label{eq:breather-zero-modulus-identity}
\end{align}
Since $\cosh T=(e^T+e^{-T})/2$, we have, as $T\to+\infty$,
\begin{align*}
 \frac1{\cosh T}
 &=\frac{2e^{-T}}{1+e^{-2T}}\\
 &=2e^{-T}\bigl(1-e^{-2T}+\mathcal{O}(e^{-4T})\bigr)\\
 &=2e^{-T}+\mathcal{O}(e^{-3T}),
\end{align*}
where we used $(1+s)^{-1}=1-s+\mathcal{O}(s^2)$ as $s\to0$.
In particular, $h^2/\cosh^2T\to0$ as $T\to+\infty$. Applying the \emph{Taylor expansion}
$(1-s)^{-1/2}=1+\tfrac12s+\mathcal{O}(s^2)$ as $s\to0$, we obtain
\begin{align*}
 \left(1-\frac{h^2}{\cosh^2T}\right)^{-1/2}
 &=1+\frac{h^2}{2\cosh^2T}+\mathcal{O}(\cosh^{-4}T)\\
 &=1+\frac{h^2}{2}\bigl(2e^{-T}+\mathcal{O}(e^{-3T})\bigr)^2+\mathcal{O}(e^{-4T})\\
&=1+2h^2e^{-2T}+\mathcal{O}(e^{-4T})\\
 &=1+e^{-2T}+\mathcal{O}(e^{-4T}),
\end{align*}
where the last equality uses $h^2=1/2$. Hence
\begin{equation}\label{eq:breather-zero-taylor-bound}
 \left(1-\frac{h^2}{\cosh^2T}\right)^{-1/2}
 =1+\mathcal{O}(e^{-2T}).
\end{equation}
Substituting \eqref{eq:breather-zero-taylor-bound} into \eqref{eq:breather-zero-modulus-identity} gives $|\widehat u_{\mathrm{AB}}(-T,0)|^2=1+\mathcal{O}(e^{-2T})$.
Taking the positive square root gives
\[
 |\widehat u_{\mathrm{AB}}(-T,0)|=1+\mathcal{O}(e^{-2T}),
\]
which proves \eqref{eq:breather-initial-zero-mode}.

\smallskip
\noindent\emph{The nonzero modes.}
For $T\ge0$, \eqref{eq:breather-smaller-root} implies
\begin{equation}\label{eq:breather-initial-ratio-bound}
 r(-T)\le\frac{h}{\cosh T}\le\sqrt2e^{-T}.
\end{equation}
Combining \eqref{eq:breather-fourier-prefactor} and \eqref{eq:breather-initial-ratio-bound} with \eqref{eq:breather-fourier-nonzero} gives, for $m\ne0$,
\begin{align*}
 |\widehat u_{\mathrm{AB}}(-T,m)|
 &=\sqrt2\,r(-T)^{|m|}\\
 &\le\sqrt2\,(\sqrt2e^{-T})^{|m|},
\end{align*}
which proves \eqref{eq:breather-initial-Fourier}.

\smallskip
\textup{(iii)}
For all sufficiently large $T$, one has $\sqrt2e^{-T}\le1/2$.
Summing \eqref{eq:breather-initial-Fourier} over the nonzero modes yields
\[
 \sum_{m\ne0}|\widehat u_{\mathrm{AB}}(-T,m)|
 \le2\sqrt2\sum_{m=1}^{\infty}(\sqrt2e^{-T})^m
 =\frac{4e^{-T}}{1-\sqrt2e^{-T}}
 \le8e^{-T}.
\]
This proves \eqref{eq:breather-initial-nonzero-sum}. Adding the zero
mode and using \eqref{eq:breather-initial-zero-mode}, we obtain
\[
 \|\widehat u_{\mathrm{AB}}(-T)\|_{\ell^1}
 =|\widehat u_{\mathrm{AB}}(-T,0)|+\sum_{m\ne0}|\widehat u_{\mathrm{AB}}(-T,m)|
 =1+\mathcal{O}(e^{-T}),
\]
as claimed in \eqref{eq:breather-initial-fourier-norm}.

\smallskip
\textup{(iv)}
To prove \eqref{eq:breather-uniform-background}, we first reconstruct the Fourier series using
\eqref{eq:breather-poisson-fourier-series}:
\begin{equation}\label{eq:breather-fourier-reconstruction}
\begin{aligned}
 u_{\mathrm{AB}}(t,x)
 &=\sum_{m\in\Z}\widehat u_{\mathrm{AB}}(t,m)e^{\iu mx}\\
 &=\widehat u_{\mathrm{AB}}(t,0)
   +\sum_{m\ne0}\widehat u_{\mathrm{AB}}(t,m)e^{\iu mx}.
\end{aligned}
\end{equation}
Substituting \eqref{eq:breather-fourier-reconstruction} and applying
the reverse triangle inequality, we obtain, for every $x\in\R$,
\begin{align*}
 \bigl||u_{\mathrm{AB}}(t,x)|-1\bigr|
 &\le\bigl||u_{\mathrm{AB}}(t,x)|
             -|\widehat u_{\mathrm{AB}}(t,0)|\bigr|
       +\bigl||\widehat u_{\mathrm{AB}}(t,0)|-1\bigr|\\
 &\le|u_{\mathrm{AB}}(t,x)-\widehat u_{\mathrm{AB}}(t,0)|
       +\bigl||\widehat u_{\mathrm{AB}}(t,0)|-1\bigr|\\
 &=\left|\sum_{m\ne0}\widehat u_{\mathrm{AB}}(t,m)e^{\iu mx}\right|
       +\bigl||\widehat u_{\mathrm{AB}}(t,0)|-1\bigr|\\
 &\le\sum_{m\ne0}|\widehat u_{\mathrm{AB}}(t,m)|
       +\bigl||\widehat u_{\mathrm{AB}}(t,0)|-1\bigr|.
\end{align*}
The right-hand side is independent of $x$. Taking the supremum and
using \eqref{eq:breather-initial-nonzero-sum} and
\eqref{eq:breather-initial-zero-mode} give
\[
 \sup_{x\in\R}\bigl||u_{\mathrm{AB}}(t,x)|-1\bigr|
 \lesssim e^t+e^{2t}\lesssim e^t\to0,\qquad t\to-\infty.
\]
Also, \eqref{eq:explicit-breather} gives
\[
 |u_{\mathrm{AB}}(t,x)|^2
 =\frac{h^2\cos^2x+\sinh^2t}{(\cosh t-h\cos x)^2}
 =|u_{\mathrm{AB}}(-t,x)|^2.
\]
By evenness, the corresponding bound is $\mathcal{O}(e^{-t})$ as $t\to+\infty$.
This proves \eqref{eq:breather-uniform-background}.

\smallskip
\textup{(v)}
Substituting $t=x=0$ into \eqref{eq:explicit-breather} gives
\[
 u_{\mathrm{AB}}(0,0)=\frac{h}{h-1}=-(1+\sqrt2).
\]
Taking the modulus proves \eqref{eq:breather-size}.

\smallskip
\textup{(vi)}
Substituting \eqref{eq:breather-fourier-zero} and
\eqref{eq:breather-fourier-nonzero} into the $\ell^1$ norm and using
\eqref{eq:breather-fourier-prefactor}, we pair the modes $m$ and $-m$
and sum the geometric series to obtain
\begin{align}
 \|\widehat u_{\mathrm{AB}}(t)\|_{\ell^1}
 &=|\widehat u_{\mathrm{AB}}(t,0)|
   +\sum_{m\ne0}|\widehat u_{\mathrm{AB}}(t,m)|\nonumber\\
 &=\left|1-\frac{\cosh t+\iu\sinh t}{\sqrt{\cosh^2t-h^2}}\right|
   +2\sqrt2\sum_{m=1}^{\infty}r(t)^m\nonumber\\
 &=\left|1-\frac{\cosh t+\iu\sinh t}{\sqrt{\cosh^2t-h^2}}\right|
   +\frac{2\sqrt2\,r(t)}{1-r(t)}.
 \label{eq:breather-fourier-norm-series}
\end{align}
For the first term, the triangle inequality and
\eqref{eq:breather-fourier-prefactor} give
\begin{equation}\label{eq:breather-zero-uniform-bound}
 \left|1-\frac{\cosh t+\iu\sinh t}{\sqrt{\cosh^2t-h^2}}\right|
 \le1+\frac{|\cosh t+\iu\sinh t|}{\sqrt{\cosh^2t-h^2}}
 =1+\sqrt2.
\end{equation}
For the second term, \eqref{eq:breather-smaller-root} and
$\cosh t\ge1$ yield
\begin{equation}\label{eq:breather-ratio-uniform-bound}
 0<r(t)=\frac{h}{\cosh t+\sqrt{\cosh^2t-h^2}}
 \le\frac{h}{1+\sqrt{1-h^2}}=\sqrt2-1.
\end{equation}
Since $r/(1-r)$ is increasing for $0<r<1$, substituting
\eqref{eq:breather-zero-uniform-bound} and
\eqref{eq:breather-ratio-uniform-bound} into
\eqref{eq:breather-fourier-norm-series} gives
\begin{align*}
 \|\widehat u_{\mathrm{AB}}(t)\|_{\ell^1}
 &\le1+\sqrt2+\frac{2\sqrt2\,r(t)}{1-r(t)}\\
 &\le1+\sqrt2+\frac{2\sqrt2(\sqrt2-1)}{2-\sqrt2}
 =3+\sqrt2.
\end{align*}
The bound is uniform in $t$, proving
\eqref{eq:breather-uniform-fourier-bound}.
\end{proof}

\subsection{Embedding into a Fixed Quasi-Periodic Frequency Module}
\label{subsec:counterexample-embedding}

Fix
\[
 \begin{gathered}
 \nu=2,\quad \omega=(1,\sqrt2),\quad z=2,\\
 c(n)=(1+|n_1|)^{-12}(1+|n_2|)^{-12},\quad A=1,\\
 \rho=(12,12),\quad\kappa=(2,2),\quad\beta=\rho-\kappa=(10,10).
 \end{gathered}
\]
For these coefficients, the variance sum satisfies
\begin{equation}\label{eq:counterexample-variance-cutoff}
 1\le S=\left(\sum_{m\in\Z}(1+|m|)^{-24}\right)^2
 \le\left(1+\frac2{23}\right)^2<\frac43,
 \qquad
 \delta_\varepsilon=\frac{\varepsilon}{1+4/S}.
\end{equation}

\medskip
\noindent\emph{Diophantine property of $\omega$.}
Under this parameter choice,
\[
 \omega\in\mathrm{DC}\!\left(\frac1{\sqrt2},\frac54\right).
\]
Indeed, let $(p,q)\in\Z^2\setminus\{(0,0)\}$.
If one coordinate is zero, then
\[
 |p+\sqrt2q|\ge1
 \ge\frac1{\sqrt2}(|p|+|q|)^{-5/4}.
\]
If both coordinates are nonzero, irrationality of $\sqrt2$
implies $p^2-2q^2\ne0$. Since this is an integer,
$|p^2-2q^2|\ge1$, and factorization gives
\begin{align*}
 |p+\sqrt2q|
 &=\frac{|p^2-2q^2|}{|p-\sqrt2q|}
 \ge\frac1{|p|+\sqrt2|q|}\\
 &\ge\frac1{\sqrt2(|p|+|q|)}
 \ge\frac1{\sqrt2}(|p|+|q|)^{-5/4},
\end{align*}
where the last inequality uses $|p|+|q|\ge1$.
Consequently, \eqref{eq:HDC} holds with $\gamma=1/\sqrt2$ and
$\tau=5/4>\nu-1$.

\medskip
\noindent\emph{Construction of a sequence of small frequency breathers.}
Choose a sequence of positive integer solutions $(p_j,q_j)$ of the \emph{Pell equation}
$p_j^2-2q_j^2=1$ such that $p_j,q_j\to\infty$ as $j\to\infty$. Set
\begin{align}
 n_j&=(-p_j,q_j),\label{eq:pell-lattice-index}\\
 k_j&=\langle n_j,\omega\rangle
 =-\frac1{p_j+\sqrt2q_j}\to0,\qquad j\to\infty,\label{eq:pell-physical-frequency}\\
 \mu_j&=k_j^2,\label{eq:pell-coupling}\\
 T_j&=15|\log\mu_j|,\label{eq:pell-scaled-waiting-time}\\
 t_j&=T_j\mu_j^{-1}.\label{eq:pell-focusing-time}
\end{align}
With these parameters, rescale and translate the Akhmediev breather as follows:
\begin{align}
 v_j(t,x):&=u_{\mathrm{AB}}\bigl(\mu_j(t-t_j),k_jx\bigr)\nonumber\\
          &=u_{\mathrm{AB}}(\mu_jt-T_j,k_jx).
 \label{eq:embedded-breather}
\end{align}

\begin{lemma}[Properties of the Embedded Breathers]
\label{lem:embedded-breather-properties}
For each $j$, the function $v_j$ in \eqref{eq:embedded-breather}
has the following properties.
\begin{enumerate}[label=\textup{(\roman*)},leftmargin=*,itemsep=3pt]
\item It is smooth on $\R^2$ and solves the weakly nonlinear NLS
\begin{equation}\label{eq:embedded-breather-NLS}
 \iu\partial_t v_j+\partial_{xx}v_j+\mu_j|v_j|^2v_j=0,
 \qquad (t,x)\in\R\times\R.
\end{equation}
\item For every $t\in\R$, it is $2\pi/|k_j|$-periodic in $x$
and belongs to $\mathrm{QP}_\omega(\R)$, with Fourier expansion
\begin{equation}\label{eq:embedded-breather-fourier}
 v_j(t,x)=\sum_{m\in\Z}\widehat u_{\mathrm{AB}}(\mu_jt-T_j,m)
                   e^{\iu\langle mn_j,\omega\rangle x}.
\end{equation}
Its initial coefficient sequence $h_j$ is supported on
$\{mn_j:m\in\Z\}$, with
\begin{equation}\label{eq:embedded-breather-initial-coefficients}
 (h_j)_{mn_j}=\widehat u_{\mathrm{AB}}(-T_j,m),
 \qquad m\in\Z.
\end{equation}
\item Its Fourier norm is uniformly bounded, and its amplitude at
the focusing time satisfies
\begin{align}
 \sup_{t\in\R}\|\widehat v_j(t)\|_{\ell^1(\Z^2)}&\le3+\sqrt2,
 \label{eq:embedded-breather-uniform-fourier}\\
 |v_j(t_j,2\pi m/k_j)|&=1+\sqrt2,\qquad m\in\Z.
 \label{eq:embedded-breather-focus}
\end{align}
\end{enumerate}
\end{lemma}

\begin{proof}[Proof of \autoref{lem:embedded-breather-properties}]
\textup{(i)} Smoothness follows from \autoref{lem:breather-identities}.
By the \emph{chain rule}, time differentiation contributes $\mu_j$ and
two spatial derivatives contribute $k_j^2=\mu_j$.
Thus \eqref{eq:breather-NLS} gives \eqref{eq:embedded-breather-NLS}.

\par\smallskip
\textup{(ii)} Spatial $2\pi$-periodicity of $u_{\mathrm{AB}}$
gives the period $2\pi/|k_j|$ for $v_j$.
Substituting $k_j=\langle n_j,\omega\rangle$ into its Fourier
series gives \eqref{eq:embedded-breather-fourier}, whose frequencies
belong to $\mathcal M_\omega$. This shows that $v_j\in\mathrm{QP}_\omega(\R)$.
Setting $t=0$ gives \eqref{eq:embedded-breather-initial-coefficients}.

\par\smallskip
\textup{(iii)} Since $m\mapsto mn_j$ is injective, \eqref{eq:breather-uniform-fourier-bound} gives
\[
 \|\widehat v_j(t)\|_{\ell^1(\Z^2)}
 =\sum_{m\in\Z}|\widehat u_{\mathrm{AB}}(\mu_jt-T_j,m)|
 \le3+\sqrt2.
\]
Substituting $(t,x)=(t_j,2\pi m/k_j)$ into
\eqref{eq:embedded-breather} and using \eqref{eq:breather-size}
and spatial periodicity give \eqref{eq:embedded-breather-focus}.
\end{proof}

\subsection{Gaussian Cost of the Embedded Breather}
\label{subsec:breather-initial-neighborhoods}

The embedded breather $v_j$ serves as a deterministic reference
for the normalized random solution. We therefore consider events
on which the normalized random initial coefficients $\widetilde d$
lie near $h_j$, the initial Fourier coefficients of $v_j$.
The initial bounds and Gaussian cost below prepare the neighborhood
probability estimate in the next subsection.

\begin{definition}[Gaussian Cost]
\label{def:gaussian-cost}
The deterministic coefficient sequences with finite
\emph{Gaussian cost} form
\begin{equation}\label{eq:gaussian-cost-space}
 \mathcal H
 :=\left\{h\in\C^{\Z^2}:
             \|h\|_{\mathcal H}<\infty\right\},
\end{equation}
where the cost is defined by
\begin{equation}\label{eq:gaussian-cost}
 \|h\|_{\mathcal H}^2
 :=\sum_{n\in\Z^2}\frac{|h_n|^2}{|c(n)|^2}.
\end{equation}

\end{definition}

\begin{lemma}[Initial Bounds]\label{lem:precursor-bounds}
Let $h_j$ be the initial Fourier coefficient sequence \eqref{eq:embedded-breather-initial-coefficients}. Then $h_j\in\mathcal H$ for
every $j$, and the following bounds hold, with asymptotics taken
as $j\to\infty$.
\begin{enumerate}[label=\textup{(\roman*)},leftmargin=*,itemsep=3pt]
\item \textup{Initial Fourier Sum.}
\begin{equation}\label{eq:precursor-one-norm}
 \|h_j\|_{\ell^1(\Z^2)}=1+o(1).
\end{equation}
\item \textup{Weighted Initial Sum.}
\begin{equation}\label{eq:precursor-weighted-norm}
 \|w_\beta h_j\|_{\ell^1(\Z^2)}=1+\mathcal{O}(\mu_j^5).
\end{equation}
\item \textup{Gaussian Cost.}
\begin{equation}\label{eq:precursor-gaussian-cost}
 \|h_j\|_{\mathcal H}^2=1+\mathcal{O}(\mu_j^6).
\end{equation}
\end{enumerate}
\end{lemma}

\begin{proof}[Proof of \autoref{lem:precursor-bounds}]
\textup{(i)}
By \eqref{eq:embedded-breather-initial-coefficients}, we separate the zero and nonzero modes as
\begin{align}
 \|h_j\|_{\ell^1(\Z^2)}
 &=|(h_j)_0|+\sum_{m\ne0}|(h_j)_{mn_j}|\nonumber\\
 &=|\widehat u_{\mathrm{AB}}(-T_j,0)|
  +\sum_{m\ne0}|\widehat u_{\mathrm{AB}}(-T_j,m)|.
 \label{eq:precursor-one-norm-decomposition}
\end{align}
For the first term, \eqref{eq:breather-initial-zero-mode} gives
\begin{equation}\label{eq:precursor-zero-mode-bound}
|\widehat u_{\mathrm{AB}}(-T_j,0)|
 =1+\mathcal{O}(e^{-2T_j})=1+\mathcal{O}(\mu_j^{30}),
\end{equation}
where $e^{-T_j}=\mu_j^{15}$; see 
\eqref{eq:pell-scaled-waiting-time}.
For the second term, \eqref{eq:breather-initial-Fourier} yields
\[
|\widehat u_{\mathrm{AB}}(-T_j,m)|
 \le\sqrt2(\sqrt2\mu_j^{15})^{|m|},
 \qquad m\ne0.
\]
For sufficiently large $j$, $\sqrt2\mu_j^{15}\le1/2$, so
\begin{equation}\label{eq:precursor-nonzero-mode-sum}
\sum_{m\ne0}|\widehat u_{\mathrm{AB}}(-T_j,m)|
 \le
2\sqrt2\sum_{m=1}^{\infty}(\sqrt2\mu_j^{15})^m
 =\mathcal{O}(\mu_j^{15}).
\end{equation}
Substituting \eqref{eq:precursor-zero-mode-bound} and
\eqref{eq:precursor-nonzero-mode-sum} into
\eqref{eq:precursor-one-norm-decomposition}, we obtain
\[
 \|h_j\|_{\ell^1(\Z^2)}
 =1+\mathcal{O}(\mu_j^{15}),
\]
which proves \eqref{eq:precursor-one-norm}.

\par\smallskip\indent\textup{(ii)}
Using \eqref{eq:embedded-breather-initial-coefficients} and
$w_\beta(0)=1$, we decompose the weighted norm as
\begin{align}
 \|w_\beta h_j\|_{\ell^1(\Z^2)}
 &=w_\beta(0)|(h_j)_0|
   +\sum_{m\ne0}w_\beta(mn_j)|(h_j)_{mn_j}|\nonumber\\
 &=|\widehat u_{\mathrm{AB}}(-T_j,0)|
   +\sum_{m\ne0}w_\beta(mn_j)|\widehat u_{\mathrm{AB}}(-T_j,m)|.
 \label{eq:precursor-weighted-norm-decomposition}
\end{align}
The first term satisfies \eqref{eq:precursor-zero-mode-bound}.
For the second term, $\beta=(10,10)$, $n_j=(-p_j,q_j)$, and
the inequality $1+xy\le(1+x)(1+y)$ for $x,y\ge0$ give
\begin{align}
 w_\beta(mn_j)
 &=(1+|m|p_j)^{10}(1+|m|q_j)^{10}\nonumber\\
 &\le(1+p_j)^{10}(1+q_j)^{10}(1+|m|)^{20}\nonumber\\
 &\le \langle n_j\rangle^{20}(1+|m|)^{20}.
 \label{eq:precursor-product-weight-bound}
\end{align}
To express this bound in terms of $\mu_j$, we use
\eqref{eq:pell-lattice-index}--\eqref{eq:pell-coupling}, which give
$\langle n_j\rangle=1+p_j+q_j$ and
$\mu_j^{-1/2}=|k_j|^{-1}=p_j+\sqrt2q_j$.
Since $p_j,q_j\ge1$,
\[
 \frac1{\sqrt2}\mu_j^{-1/2}
 =\frac{p_j+\sqrt2q_j}{\sqrt2}
 \le p_j+q_j
 \le\langle n_j\rangle
 \le2(p_j+q_j)
 \le2\mu_j^{-1/2}.
\]
It follows that, with constants independent of $j$,
\begin{equation}\label{eq:precursor-parameter-comparison}
 \langle n_j\rangle\asymp\mu_j^{-1/2}.
\end{equation}
Substituting \eqref{eq:precursor-parameter-comparison} into
\eqref{eq:precursor-product-weight-bound}, then substituting
\eqref{eq:precursor-product-weight-bound} into the second term of
\eqref{eq:precursor-weighted-norm-decomposition}, and using
\eqref{eq:breather-initial-Fourier}, we obtain
\begin{align}
 \sum_{m\ne0}w_\beta(mn_j)|\widehat u_{\mathrm{AB}}(-T_j,m)|
 &\le\langle n_j\rangle^{20}
       \sum_{m\ne0}(1+|m|)^{20}
          |\widehat u_{\mathrm{AB}}(-T_j,m)|\nonumber\\
 &\le\sqrt2\langle n_j\rangle^{20}
       \sum_{m\ne0}(1+|m|)^{20}
                    (\sqrt2\mu_j^{15})^{|m|}\nonumber\\
 &=2\sqrt2\langle n_j\rangle^{20}
       \sum_{m=1}^{\infty}(1+m)^{20}
                    (\sqrt2\mu_j^{15})^m\nonumber\\
 &\le4\langle n_j\rangle^{20}\mu_j^{15}
       \sum_{m=1}^{\infty}(1+m)^{20}2^{1-m}\nonumber\\
 &\le C\langle n_j\rangle^{20}\mu_j^{15}
 \le C\mu_j^{-10}\mu_j^{15}
 =C\mu_j^5.
 \label{eq:precursor-nonzero-weighted}
\end{align}
Substituting \eqref{eq:precursor-zero-mode-bound} and
\eqref{eq:precursor-nonzero-weighted} into
\eqref{eq:precursor-weighted-norm-decomposition} yields
\eqref{eq:precursor-weighted-norm}.

\par\smallskip\indent\textup{(iii)}
By \eqref{eq:gaussian-cost} and
\eqref{eq:embedded-breather-initial-coefficients},
\begin{align}
 \|h_j\|_{\mathcal H}^2
 &=\frac{|(h_j)_0|^2}{|c(0)|^2}
   +\sum_{m\ne0}\frac{|(h_j)_{mn_j}|^2}{|c(mn_j)|^2}\nonumber\\*
 &=\frac{|\widehat u_{\mathrm{AB}}(-T_j,0)|^2}{|c(0)|^2}
   +\sum_{m\ne0}
       \frac{|\widehat u_{\mathrm{AB}}(-T_j,m)|^2}{|c(mn_j)|^2}.
 \label{eq:precursor-gaussian-cost-decomposition}
\end{align}
For the first term, $c(0)=1$ and
\eqref{eq:precursor-zero-mode-bound} give
\begin{equation}\label{eq:precursor-zero-gaussian-cost}
 \frac{|\widehat u_{\mathrm{AB}}(-T_j,0)|^2}{|c(0)|^2}
 =1+\mathcal{O}(\mu_j^{30}).
\end{equation}
For the second term, the definition of $c$ gives
\begin{align}
 |c(mn_j)|^{-2}
 &=(1+|m|p_j)^{24}(1+|m|q_j)^{24}\nonumber\\
 &\le\langle n_j\rangle^{48}(1+|m|)^{48}.
 \label{eq:precursor-inverse-coefficient-weight}
\end{align}
Substituting \eqref{eq:precursor-inverse-coefficient-weight} into
the second term of \eqref{eq:precursor-gaussian-cost-decomposition}
and using \eqref{eq:breather-initial-Fourier} and
\eqref{eq:precursor-parameter-comparison}, we obtain
\begin{align}
 \sum_{m\ne0}
     \frac{|\widehat u_{\mathrm{AB}}(-T_j,m)|^2}{|c(mn_j)|^2}
 &\le\langle n_j\rangle^{48}
       \sum_{m\ne0}(1+|m|)^{48}
          |\widehat u_{\mathrm{AB}}(-T_j,m)|^2\nonumber\\
 &\le2\langle n_j\rangle^{48}
       \sum_{m\ne0}(1+|m|)^{48}
                    (2\mu_j^{30})^{|m|}\nonumber\\
 &\le C\langle n_j\rangle^{48}\mu_j^{30}
 \le C\mu_j^{-24}\mu_j^{30}
 =C\mu_j^6.
 \label{eq:precursor-nonzero-gaussian-cost}
\end{align}
Substituting \eqref{eq:precursor-zero-gaussian-cost} and
\eqref{eq:precursor-nonzero-gaussian-cost} into
\eqref{eq:precursor-gaussian-cost-decomposition}, we obtain
\[
 \|h_j\|_{\mathcal H}^2
 =1+\mathcal{O}(\mu_j^{30})+\mathcal{O}(\mu_j^6)
 =1+\mathcal{O}(\mu_j^6),
\]
which proves \eqref{eq:precursor-gaussian-cost}. Thus $h_j\in\mathcal H$. 

\end{proof}

\subsection{Breather Tube and Failure of the Initial Gaussian Rate Function}

We first control solutions whose initial coefficients lie near $h_j$.
We then estimate the probability that $\widetilde d$ belongs to this
initial neighborhood. Set
\begin{equation}\label{eq:tube-radius}
 r_j=\tfrac18e^{-150T_j}=\tfrac18\mu_j^{2250}.
\end{equation}

\begin{lemma}[Stability in a Breather Tube]\label{lem:breather-tube}
For every $\mathsf e\in\ell^1(\Z^2)$ with $\|\mathsf e\|_{\ell^1(\Z^2)}<r_j$,
\eqref{eq:embedded-breather-NLS} with initial coefficients $h_j+\mathsf e$
has a unique solution $v$ that exists beyond
$t_j$. Moreover,
\begin{equation}\label{eq:breather-tube-stability}
 \sup_{0\le t\le t_j}
 \|\widehat v(t)-\widehat v_j(t)\|_{\ell^1(\Z^2)}\le\tfrac18,
 \qquad |v(t_j,0)|\ge1+\sqrt2-\tfrac18>2.
\end{equation}
\end{lemma}
\begin{remark}[Visualization of Tube Stability]
The tube estimate \eqref{eq:breather-tube-stability}
is illustrated in \autoref{fig:breather-solution-tube}.
\begin{figure}[htbp]
\centering
\resizebox{\dimexpr\linewidth-2\captionindent\relax}{!}{
\begingroup
\definecolor{breatherblue}{RGB}{35,87,150}\definecolor{sampleorange}{RGB}{177,82,30}\begin{tikzpicture}[x=1.3cm,y=1.8cm,font=\small,>=Stealth]
\fill[breatherblue!12]
 (.8,.70) .. controls (3.2,.70) and (4.1,.77) .. (5,.97)
 .. controls (6.4,1.27) and (6.9,1.80) .. (8.1,2.114)
 -- (8.1,2.714)
 .. controls (6.9,2.40) and (6.4,1.87) .. (5,1.57)
 .. controls (4.1,1.37) and (3.2,1.30) .. (.8,1.30) -- cycle;
\draw[breatherblue!40,thin]
 (.8,.70) .. controls (3.2,.70) and (4.1,.77) .. (5,.97)
 .. controls (6.4,1.27) and (6.9,1.80) .. (8.1,2.114);
\draw[breatherblue!40,thin]
 (.8,1.30) .. controls (3.2,1.30) and (4.1,1.37) .. (5,1.57)
 .. controls (6.4,1.87) and (6.9,2.40) .. (8.1,2.714);
\draw[->,black!75] (.8,0)--(8.8,0) node[right] {$t$};
\draw[->,black!75] (.8,0)--(.8,3.3);
\node[anchor=south west] at (.9,3.15) {Normalized Amplitude at $x=0$};
\draw[gray!60,densely dotted] (.8,1)--(2.7,1);
\draw[gray!65,dashed] (.8,2)--(8.55,2);
\draw[black!65] (.73,2)--(.87,2) node[left=3pt] {$2$};
\node[anchor=south,fill=white,inner sep=2pt,text=black!65]
 at (4.25,2.01) {Exceedance Threshold};

\draw[gray!60,densely dotted] (8.1,0)--(8.1,2.714);
\node[below=4pt] at (.8,0) {$0$};
\node[below=4pt] at (8.1,0) {$t_j$};
\draw[sampleorange,line width=.65pt] plot[smooth,tension=.3] coordinates {
(0.80000,1.18000)
(0.88930,1.17571)
(0.97721,1.17114)
(1.06375,1.16613)
(1.14894,1.16064)
(1.23279,1.15468)
(1.31532,1.14832)
(1.39655,1.14172)
(1.47650,1.13505)
(1.55518,1.12854)
(1.63262,1.12241)
(1.70882,1.11687)
(1.78381,1.11209)
(1.85761,1.10816)
(1.93023,1.10512)
(2.00168,1.10289)
(2.07200,1.10136)
(2.14119,1.10030)
(2.20927,1.09943)
(2.27627,1.09844)
(2.34219,1.09698)
(2.40705,1.09475)
(2.47088,1.09148)
(2.53369,1.08697)
(2.59550,1.08111)
(2.65632,1.07394)
(2.71618,1.06560)
(2.77509,1.05635)
(2.83306,1.04655)
(2.89012,1.03665)
(2.94629,1.02718)
(3.00158,1.01867)
(3.05600,1.01165)
(3.10958,1.00658)
(3.16234,1.00383)
(3.21428,1.00367)
(3.26544,1.00621)
(3.31582,1.01140)
(3.36545,1.01907)
(3.41433,1.02886)
(3.46250,1.04034)
(3.50996,1.05296)
(3.55674,1.06614)
(3.60285,1.07930)
(3.64831,1.09190)
(3.69314,1.10346)
(3.73735,1.11367)
(3.78097,1.12233)
(3.82400,1.12937)
(3.86647,1.13492)
(3.90840,1.13923)
(3.94980,1.14269)
(3.99069,1.14573)
(4.03108,1.14886)
(4.07101,1.15255)
(4.11047,1.15724)
(4.14950,1.16323)
(4.18810,1.17073)
(4.22630,1.17976)
(4.26412,1.19016)
(4.30156,1.20162)
(4.33866,1.21368)
(4.37541,1.22574)
(4.41186,1.23714)
(4.44800,1.24719)
(4.48386,1.25524)
(4.51946,1.26072)
(4.55481,1.26321)
(4.58994,1.26246)
(4.62485,1.25845)
(4.65957,1.25136)
(4.69411,1.24164)
(4.72850,1.22989)
(4.76275,1.21690)
(4.79687,1.20356)
(4.83088,1.19085)
(4.86481,1.17970)
(4.89867,1.17098)
(4.93248,1.16542)
(4.96625,1.16356)
(5.00000,1.16572)
(5.05208,1.17575)
(5.10334,1.18976)
(5.15379,1.20734)
(5.20345,1.22786)
(5.25234,1.25056)
(5.30049,1.27458)
(5.34790,1.29908)
(5.39460,1.32324)
(5.44061,1.34637)
(5.48594,1.36796)
(5.53061,1.38769)
(5.57465,1.40549)
(5.61807,1.42156)
(5.66089,1.43626)
(5.70313,1.45018)
(5.74480,1.46401)
(5.78593,1.47849)
(5.82654,1.49439)
(5.86664,1.51236)
(5.90625,1.53291)
(5.94539,1.55635)
(5.98409,1.58271)
(6.02235,1.61177)
(6.06020,1.64305)
(6.09766,1.67580)
(6.13474,1.70909)
(6.17146,1.74181)
(6.20785,1.77287)
(6.24392,1.80111)
(6.27969,1.82555)
(6.31518,1.84534)
(6.35040,1.85994)
(6.38538,1.86906)
(6.42014,1.87279)
(6.45469,1.87153)
(6.48905,1.86598)
(6.52324,1.85714)
(6.55729,1.84618)
(6.59120,1.83441)
(6.62500,1.82316)
(6.65871,1.81370)
(6.69234,1.80713)
(6.72591,1.80433)
(6.75945,1.80591)
(6.79297,1.81211)
(6.82649,1.82286)
(6.86003,1.83778)
(6.89360,1.85619)
(6.92723,1.87720)
(6.96094,1.89979)
(6.99474,1.92289)
(7.02865,1.94547)
(7.06269,1.96669)
(7.09689,1.98587)
(7.13125,2.00264)
(7.16580,2.01693)
(7.20056,2.02900)
(7.23554,2.03940)
(7.27076,2.04895)
(7.30625,2.05866)
(7.34202,2.06963)
(7.37809,2.08297)
(7.41448,2.09969)
(7.45120,2.12061)
(7.48828,2.14627)
(7.52574,2.17686)
(7.56359,2.21220)
(7.60185,2.25173)
(7.64054,2.29453)
(7.67969,2.33934)
(7.71930,2.38472)
(7.75940,2.42903)
(7.80001,2.47066)
(7.84114,2.50804)
(7.88281,2.53987)
(7.92505,2.56508)
(7.96787,2.58298)
(8.01129,2.59336)
(8.05533,2.59640)
(8.10000,2.59271)
};
\draw[sampleorange,line width=.65pt] plot[smooth,tension=.3] coordinates {
(0.80000,0.82000)
(0.88930,0.83109)
(0.97721,0.84107)
(1.06375,0.84992)
(1.14894,0.85771)
(1.23279,0.86465)
(1.31532,0.87104)
(1.39655,0.87720)
(1.47650,0.88356)
(1.55518,0.89047)
(1.63262,0.89825)
(1.70882,0.90709)
(1.78381,0.91706)
(1.85761,0.92810)
(1.93023,0.93994)
(2.00168,0.95223)
(2.07200,0.96448)
(2.14119,0.97618)
(2.20927,0.98678)
(2.27627,0.99580)
(2.34219,1.00289)
(2.40705,1.00783)
(2.47088,1.01059)
(2.53369,1.01134)
(2.59550,1.01044)
(2.65632,1.00835)
(2.71618,1.00567)
(2.77509,1.00302)
(2.83306,1.00100)
(2.89012,1.00008)
(2.94629,1.00061)
(3.00158,1.00272)
(3.05600,1.00634)
(3.10958,1.01117)
(3.16234,1.01677)
(3.21428,1.02252)
(3.26544,1.02777)
(3.31582,1.03191)
(3.36545,1.03435)
(3.41433,1.03473)
(3.46250,1.03287)
(3.50996,1.02887)
(3.55674,1.02306)
(3.60285,1.01600)
(3.64831,1.00844)
(3.69314,1.00126)
(3.73735,0.99532)
(3.78097,0.99147)
(3.82400,0.99034)
(3.86647,0.99240)
(3.90840,0.99780)
(3.94980,1.00642)
(3.99069,1.01783)
(4.03108,1.03139)
(4.07101,1.04625)
(4.11047,1.06149)
(4.14950,1.07620)
(4.18810,1.08960)
(4.22630,1.10110)
(4.26412,1.11037)
(4.30156,1.11742)
(4.33866,1.12256)
(4.37541,1.12636)
(4.41186,1.12964)
(4.44800,1.13333)
(4.48386,1.13841)
(4.51946,1.14573)
(4.55481,1.15599)
(4.58994,1.16957)
(4.62485,1.18653)
(4.65957,1.20656)
(4.69411,1.22902)
(4.72850,1.25297)
(4.76275,1.27728)
(4.79687,1.30072)
(4.83088,1.32209)
(4.86481,1.34038)
(4.89867,1.35482)
(4.93248,1.36501)
(4.96625,1.37090)
(5.00000,1.37292)
(5.05208,1.37556)
(5.10334,1.37617)
(5.15379,1.37593)
(5.20345,1.37606)
(5.25234,1.37771)
(5.30049,1.38180)
(5.34790,1.38892)
(5.39460,1.39927)
(5.44061,1.41265)
(5.48594,1.42842)
(5.53061,1.44565)
(5.57465,1.46317)
(5.61807,1.47969)
(5.66089,1.49398)
(5.70313,1.50503)
(5.74480,1.51210)
(5.78593,1.51489)
(5.82654,1.51353)
(5.86664,1.50865)
(5.90625,1.50128)
(5.94539,1.49275)
(5.98409,1.48459)
(6.02235,1.47839)
(6.06020,1.47556)
(6.09766,1.47727)
(6.13474,1.48431)
(6.17146,1.49694)
(6.20785,1.51496)
(6.24392,1.53767)
(6.27969,1.56396)
(6.31518,1.59240)
(6.35040,1.62149)
(6.38538,1.64969)
(6.42014,1.67571)
(6.45469,1.69859)
(6.48905,1.71780)
(6.52324,1.73335)
(6.55729,1.74576)
(6.59120,1.75599)
(6.62500,1.76539)
(6.65871,1.77548)
(6.69234,1.78785)
(6.72591,1.80391)
(6.75945,1.82477)
(6.79297,1.85109)
(6.82649,1.88299)
(6.86003,1.91999)
(6.89360,1.96111)
(6.92723,2.00486)
(6.96094,2.04948)
(6.99474,2.09303)
(7.02865,2.13363)
(7.06269,2.16966)
(7.09689,2.19987)
(7.13125,2.22359)
(7.16580,2.24071)
(7.20056,2.25176)
(7.23554,2.25780)
(7.27076,2.26033)
(7.30625,2.26111)
(7.34202,2.26196)
(7.37809,2.26457)
(7.41448,2.27034)
(7.45120,2.28013)
(7.48828,2.29424)
(7.52574,2.31237)
(7.56359,2.33359)
(7.60185,2.35649)
(7.64054,2.37933)
(7.67969,2.40018)
(7.71930,2.41722)
(7.75940,2.42884)
(7.80001,2.43397)
(7.84114,2.43206)
(7.88281,2.42328)
(7.92505,2.40846)
(7.96787,2.38903)
(8.01129,2.36691)
(8.05533,2.34430)
(8.10000,2.32345)
};
\draw[sampleorange,line width=.65pt] plot[smooth,tension=.3] coordinates {
(0.80000,1.09000)
(0.88930,1.08650)
(0.97721,1.08355)
(1.06375,1.08104)
(1.14894,1.07887)
(1.23279,1.07694)
(1.31532,1.07518)
(1.39655,1.07354)
(1.47650,1.07199)
(1.55518,1.07052)
(1.63262,1.06917)
(1.70882,1.06796)
(1.78381,1.06696)
(1.85761,1.06622)
(1.93023,1.06580)
(2.00168,1.06576)
(2.07200,1.06613)
(2.14119,1.06690)
(2.20927,1.06805)
(2.27627,1.06954)
(2.34219,1.07124)
(2.40705,1.07304)
(2.47088,1.07477)
(2.53369,1.07626)
(2.59550,1.07731)
(2.65632,1.07775)
(2.71618,1.07739)
(2.77509,1.07610)
(2.83306,1.07376)
(2.89012,1.07035)
(2.94629,1.06586)
(3.00158,1.06039)
(3.05600,1.05411)
(3.10958,1.04725)
(3.16234,1.04009)
(3.21428,1.03301)
(3.26544,1.02638)
(3.31582,1.02064)
(3.36545,1.01620)
(3.41433,1.01347)
(3.46250,1.01279)
(3.50996,1.01446)
(3.55674,1.01869)
(3.60285,1.02559)
(3.64831,1.03514)
(3.69314,1.04724)
(3.73735,1.06163)
(3.78097,1.07794)
(3.82400,1.09572)
(3.86647,1.11442)
(3.90840,1.13344)
(3.94980,1.15212)
(3.99069,1.16979)
(4.03108,1.18586)
(4.07101,1.19973)
(4.11047,1.21093)
(4.14950,1.21909)
(4.18810,1.22395)
(4.22630,1.22543)
(4.26412,1.22359)
(4.30156,1.21862)
(4.33866,1.21090)
(4.37541,1.20091)
(4.41186,1.18927)
(4.44800,1.17664)
(4.48386,1.16378)
(4.51946,1.15144)
(4.55481,1.14034)
(4.58994,1.13120)
(4.62485,1.12458)
(4.65957,1.12099)
(4.69411,1.12078)
(4.72850,1.12415)
(4.76275,1.13115)
(4.79687,1.14169)
(4.83088,1.15549)
(4.86481,1.17218)
(4.89867,1.19125)
(4.93248,1.21213)
(4.96625,1.23420)
(5.00000,1.25681)
(5.05208,1.28311)
(5.10334,1.30882)
(5.15379,1.33343)
(5.20345,1.35649)
(5.25234,1.37770)
(5.30049,1.39688)
(5.34790,1.41396)
(5.39460,1.42898)
(5.44061,1.44211)
(5.48594,1.45358)
(5.53061,1.46371)
(5.57465,1.47285)
(5.61807,1.48136)
(5.66089,1.48959)
(5.70313,1.49785)
(5.74480,1.50641)
(5.78593,1.51544)
(5.82654,1.52503)
(5.86664,1.53521)
(5.90625,1.54590)
(5.94539,1.55695)
(5.98409,1.56819)
(6.02235,1.57936)
(6.06020,1.59024)
(6.09766,1.60059)
(6.13474,1.61024)
(6.17146,1.61907)
(6.20785,1.62705)
(6.24392,1.63424)
(6.27969,1.64083)
(6.31518,1.64709)
(6.35040,1.65340)
(6.38538,1.66022)
(6.42014,1.66809)
(6.45469,1.67753)
(6.48905,1.68910)
(6.52324,1.70332)
(6.55729,1.72060)
(6.59120,1.74129)
(6.62500,1.76556)
(6.65871,1.79341)
(6.69234,1.82470)
(6.72591,1.85906)
(6.75945,1.89595)
(6.79297,1.93466)
(6.82649,1.97433)
(6.86003,2.01399)
(6.89360,2.05258)
(6.92723,2.08906)
(6.96094,2.12238)
(6.99474,2.15159)
(7.02865,2.17586)
(7.06269,2.19456)
(7.09689,2.20729)
(7.13125,2.21385)
(7.16580,2.21437)
(7.20056,2.20919)
(7.23554,2.19898)
(7.27076,2.18461)
(7.30625,2.16718)
(7.34202,2.14794)
(7.37809,2.12826)
(7.41448,2.10954)
(7.45120,2.09319)
(7.48828,2.08050)
(7.52574,2.07265)
(7.56359,2.07059)
(7.60185,2.07504)
(7.64054,2.08642)
(7.67969,2.10487)
(7.71930,2.13016)
(7.75940,2.16183)
(7.80001,2.19908)
(7.84114,2.24090)
(7.88281,2.28604)
(7.92505,2.33316)
(7.96787,2.38082)
(8.01129,2.42758)
(8.05533,2.47204)
(8.10000,2.51293)
};
\draw[breatherblue,very thick]
 (.8,1) .. controls (3.2,1) and (4.1,1.07) .. (5,1.27)
 .. controls (6.4,1.57) and (6.9,2.10) .. (8.1,2.414);
\node[text=breatherblue,inner sep=1pt]
 (blueinitiallabel) at (2.05,1.52) {$h_j$};
\draw[->,breatherblue!75,thin] (.8,1)
 .. controls (1.13,1.18) and (1.50,1.42)
 .. ([xshift=-2pt]blueinitiallabel.west);
\node[text=sampleorange,inner sep=1pt]
 (orangeinitiallabel) at (2.15,.51) {$\mathsf e+h_j$};
\draw[->,sampleorange!75,thin] (.8,.82)
 .. controls (.95,.61) and (1.24,.51)
 .. ([xshift=-2pt]orangeinitiallabel.west);
\fill[sampleorange] (.8,1.18) circle[radius=1.15pt];
\fill[sampleorange] (.8,.82) circle[radius=1.15pt];
\fill[sampleorange] (.8,1.09) circle[radius=1.15pt];
\fill[breatherblue] (.8,1) circle[radius=1.7pt];
\fill[breatherblue] (8.1,2.414) circle[radius=1.7pt];
\node[text=breatherblue] at (7.3,3.12) {$|v_j(t_j,0)|=1+\sqrt2$};
\draw[<-,breatherblue!75,thin] (8.35,2.99) .. controls (8.6,2.85) and (8.5,2.55) .. (8.12,2.43);
\draw[breatherblue,very thick] (1.15,2.90)--(1.75,2.90);
\node[anchor=west] at (1.88,2.90) {Breather $|v_j(t,0)|$};
\draw[sampleorange,line width=.8pt] (1.15,2.65)--(1.75,2.65);
\node[anchor=west] at (1.88,2.65)
 {Random Solutions $|v(t,0)|=\sqrt{\varepsilon_j}|u(t,0)|$};
\node[align=center] at (4.8,.36)
 {$\displaystyle\sup_{0\le t\le t_j}\|\widehat v(t)-\widehat v_j(t)\|_{\ell^1(\mathbb Z^2)}\le\frac18$};
\path[use as bounding box]
 ([xshift=-5pt,yshift=-5pt]current bounding box.south west)
 rectangle ([xshift=5pt,yshift=5pt]current bounding box.north east);
\end{tikzpicture}\endgroup

}
\caption{The breather tube at $x=0$.
The blue curve shows the amplitude $|v_j(t,0)|$ of the deterministic
breather with initial coefficient sequence $h_j$.
The orange curves show the amplitudes
$|v(t,0)|=\sqrt{\varepsilon_j}|u(t,0)|$ of normalized random solutions
with initial coefficient sequences $\mathsf e+h_j$, where
$\|w_\beta\mathsf e\|_{\ell^1(\Z^2)}<r_j$.}
\label{fig:breather-solution-tube}
\end{figure}
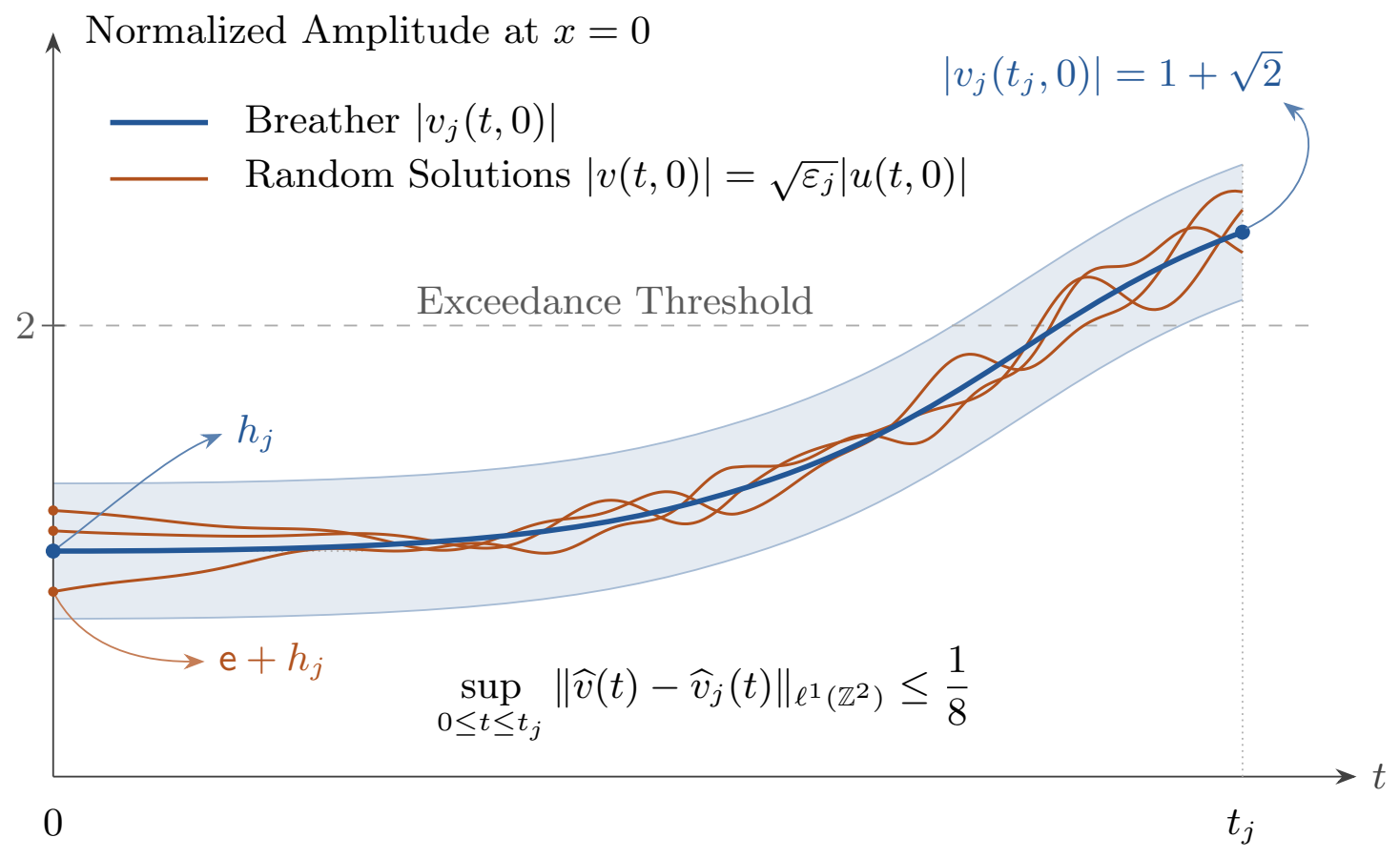
\end{remark}
\begin{proof}[Proof of \autoref{lem:breather-tube}]
\textbf{Step 1. Local Construction.}
Choose $t>0$ such that
\begin{equation}\label{eq:breather-tube-local-time}
 24\mu_jt(4+\sqrt2)^2\le1.
\end{equation}
\par\smallskip\noindent\emph{Picard Iteration.} 
Set
\begin{equation}\label{eq:breather-tube-initial-coefficients}
 \widehat v(0,n):=(h_j)_n+\mathsf e_n,\qquad n\in\Z^2.
\end{equation}
By \eqref{eq:embedded-breather-uniform-fourier} and
\eqref{eq:tube-radius},
\begin{align*}
 \|\widehat v(0)\|_{\ell^1(\Z^2)}
 &\le\|h_j\|_{\ell^1(\Z^2)}+\|\mathsf e\|_{\ell^1(\Z^2)}\\
 &<3+\sqrt2+\tfrac18<4+\sqrt2.
\end{align*}
Define the Picard iterates
\begin{align}
 \widehat v^{(0)}(s,n)
 &=e^{-\iu s\langle n,\omega\rangle^2}\widehat v(0,n),
 \label{eq:breather-tube-picard-base}\\
 \widehat v^{(k+1)}(s,n)
 &=e^{-\iu s\langle n,\omega\rangle^2}\widehat v(0,n)+\iu\mu_j\int_0^s e^{-\iu(s-r)\langle n,\omega\rangle^2}\nonumber\\*
 &\qquad\qquad\times\sum_{n_1-n_2+n_3=n}
    \widehat v^{(k)}(r,n_1)
    \overline{\widehat v^{(k)}(r,n_2)}
    \widehat v^{(k)}(r,n_3)\,dr.
 \label{eq:breather-tube-picard}
\end{align}
\par\smallskip\noindent\emph{Uniform Bounds.} For $k\ge0$, we have
\begin{equation}\label{eq:breather-tube-picard-bound}
 \sup_{0\le s\le t}\|\widehat v^{(k)}(s)\|_{\ell^1(\Z^2)}
 \le2(4+\sqrt2).
\end{equation}
Indeed, for $k=0$, \eqref{eq:breather-tube-picard-base} and the initial bound give
\[
 \sup_{0\le s\le t}\|\widehat v^{(0)}(s)\|_{\ell^1(\Z^2)}
 =\|\widehat v(0)\|_{\ell^1(\Z^2)}<4+\sqrt2.
\]
Suppose \eqref{eq:breather-tube-picard-bound} holds for some $k\ge0$.
By \eqref{eq:breather-tube-picard},
\autoref{prop:qp-fourier-calculus}, the triangle inequality
and \eqref{eq:breather-tube-local-time},
\begin{align*}
 \sup_{0\le s\le t}\|\widehat v^{(k+1)}(s)\|_{\ell^1(\Z^2)}
 &\le4+\sqrt2+\mu_jt[2(4+\sqrt2)]^3\\
 &\le\tfrac43(4+\sqrt2)<2(4+\sqrt2).
\end{align*}

\par\smallskip\noindent\emph{Cauchy Convergence.} Subtracting successive instances of \eqref{eq:breather-tube-picard}
and applying \eqref{eq:breather-tube-picard-bound} to the cubic
difference, we obtain, for $k\ge1$,
\begin{align}
 &\sup_{0\le s\le t}
 \|\widehat v^{(k+1)}(s)-\widehat v^{(k)}(s)\|_{\ell^1(\Z^2)}\nonumber\\*
 &\quad\le12\mu_jt(4+\sqrt2)^2
       \sup_{0\le s\le t}
       \|\widehat v^{(k)}(s)-\widehat v^{(k-1)}(s)\|_{\ell^1(\Z^2)}\nonumber\\*
 &\quad\le\tfrac12\sup_{0\le s\le t}
       \|\widehat v^{(k)}(s)-\widehat v^{(k-1)}(s)\|_{\ell^1(\Z^2)}.
 \label{eq:breather-tube-picard-cauchy}
\end{align}
Hence the iterates converge in $C([0,t];\ell^1(\Z^2))$ to $\widehat v$.
\par\smallskip\noindent\emph{Passage to the Limit.} The cubic difference estimate used in
\eqref{eq:breather-tube-picard-cauchy} gives uniform convergence of the
nonlinear terms in $\ell^1(\Z^2)$ on $[0,t]$, allowing passage to the
limit in \eqref{eq:breather-tube-picard} and yielding
\begin{align}
 \widehat v(s,n)
 &=e^{-\iu s\langle n,\omega\rangle^2}\widehat v(0,n)+\iu\mu_j\int_0^s e^{-\iu(s-r)\langle n,\omega\rangle^2}\nonumber\\*
 &\qquad\qquad\qquad\times\sum_{n_1-n_2+n_3=n}
    \widehat v(r,n_1)\overline{\widehat v(r,n_2)}
    \widehat v(r,n_3)\,dr.
 \label{eq:breather-tube-mild}
\end{align}
\par\smallskip\noindent\emph{Uniqueness.} Applying the cubic difference estimate to
\eqref{eq:breather-tube-mild} for two solutions with the same initial data,
Gronwall's inequality gives zero difference.
\par\medskip\noindent\textbf{Step 2. Bootstrap Setup.}
Let $T_{\max}$ be the maximal existence time of $v$, and define
\[
 \tau_*:=\sup\left\{t\in[0,t_j]\cap[0,T_{\max}):
 \sup_{0\le s\le t}\|\widehat v(s)-\widehat v_j(s)\|_{\ell^1(\Z^2)}<1\right\}.
\]
It follows from
\[
 \|\widehat v(0)-\widehat v_j(0)\|_{\ell^1(\Z^2)}
 =\|\mathsf e\|_{\ell^1(\Z^2)}<r_j<\tfrac18<1
\]
and the local construction that
$\tau_*>0$. By the definition of $\tau_*$, we have
\begin{equation}\label{eq:breather-tube-bootstrap-hypothesis}
 \sup_{0\le s\le t}\|\widehat v(s)-\widehat v_j(s)\|_{\ell^1(\Z^2)}<1,\qquad 0<t<\tau_*.
\end{equation}
Combining this hypothesis with
\eqref{eq:embedded-breather-uniform-fourier}, we obtain
\begin{equation}\label{eq:breather-tube-bootstrap-norm}
 \|\widehat v(s)\|_{\ell^1(\Z^2)}
 \le\|\widehat v_j(s)\|_{\ell^1(\Z^2)}+1
 \le4+\sqrt2,
 \qquad 0\le s\le t.
\end{equation}
These will be used to control the factors in the cubic difference.

\par\medskip\noindent\textbf{Step 3. Bootstrap Improvement and Extension.}
Using the cubic difference decomposition in
\eqref{eq:log-cubic-increment-expansion}, the convolution inequality,
invariance of the coefficient norm under conjugation and
\eqref{eq:breather-tube-bootstrap-norm}, we obtain, for $s\in[0,t]$,
\begin{align}
 \bigl\|\widehat{|v|^2v}-\widehat{|v_j|^2v_j}\bigr\|_{\ell^1(\Z^2)}
 &\le\bigl(\|\widehat v\|_{\ell^1(\Z^2)}^2
       +\|\widehat v\|_{\ell^1(\Z^2)}\|\widehat v_j\|_{\ell^1(\Z^2)}\nonumber\\
       &\qquad+\|\widehat v_j\|_{\ell^1(\Z^2)}^2\bigr)
       \|\widehat v-\widehat v_j\|_{\ell^1(\Z^2)}\nonumber\\
 &\le150\|\widehat v-\widehat v_j\|_{\ell^1(\Z^2)}.
 \label{eq:breather-tube-cubic-difference}
\end{align}
Subtracting the integral equations for $v$ and $v_j$, obtained from
\eqref{eq:breather-tube-mild}, gives
\begin{align}
 v(t)-v_j(t)
 &=e^{\iu t\partial_{xx}}\bigl(v(0)-v_j(0)\bigr)\nonumber\\
 &\quad+\iu\mu_j\int_0^t e^{\iu(t-s)\partial_{xx}}
             \bigl(|v(s)|^2v(s)-|v_j(s)|^2v_j(s)\bigr)\,ds.
 \label{eq:breather-tube-difference-integral}
\end{align}
The initial difference has coefficient sequence $\mathsf e$.
Taking coefficient norms and applying
\eqref{eq:breather-tube-cubic-difference}, we obtain
\begin{align}
 \|\widehat v(t)-\widehat v_j(t)\|_{\ell^1(\Z^2)}
 &\le\|\mathsf e\|_{\ell^1(\Z^2)}
   +150\mu_j\int_0^t
       \|\widehat v(s)-\widehat v_j(s)\|_{\ell^1(\Z^2)}\,ds.
 \label{eq:breather-tube-difference}
\end{align}
Applying Gronwall's inequality to
\eqref{eq:breather-tube-difference} on $[0,t]$ and using
$\mu_jt\le\mu_jt_j=T_j$ and \eqref{eq:tube-radius}, we obtain
\begin{equation}\label{eq:breather-tube-gronwall}
 \sup_{0\le s\le t}\|\widehat v(s)-\widehat v_j(s)\|_{\ell^1(\Z^2)}
 \le\|\mathsf e\|_{\ell^1(\Z^2)}e^{150\mu_jt}
 <r_je^{150T_j}=\tfrac18.
\end{equation}

Next, We show that $\tau_*=t_j<T_{\max}$.

Since $\tau_*\le T_{\max}$, it suffices to exclude equality.
If $\tau_*=T_{\max}$, \eqref{eq:breather-tube-bootstrap-norm} bounds
$\|\widehat v(s)\|_{\ell^1(\Z^2)}$ by $4+\sqrt2$ for all $s<T_{\max}$.
Step~1 therefore gives the same positive local lifespan
$[24\mu_j(4+\sqrt2)^2]^{-1}$ at every such time.
Restarting sufficiently close to $T_{\max}$ and using uniqueness
on the overlap extends $v$ beyond $T_{\max}$, contradicting maximality.
Hence $T_{\max}>\tau_*$.

If $\tau_*<t_j$, continuity and \eqref{eq:breather-tube-gronwall}
give $\sup_{0\le s\le\tau_*}\|\widehat v(s)-\widehat v_j(s)\|_{\ell^1(\Z^2)}\le1/8$.
Since $T_{\max}>\tau_*$, we may choose
$0<\Delta t<\min\{t_j,T_{\max}\}-\tau_*$ so that
$\|\widehat v(s)-\widehat v_j(s)\|_{\ell^1(\Z^2)}<1/2$ for
$\tau_*\le s\le\tau_*+\Delta t$.
The supremum over $[0,\tau_*+\Delta t]$ is therefore at most $1/2<1$,
contradicting the definition of $\tau_*$. 

Thus $\tau_*=t_j<T_{\max}$.
Taking $t\uparrow t_j$ in
\eqref{eq:breather-tube-gronwall} proves the first bound in
\eqref{eq:breather-tube-stability}.

\par\medskip\noindent\textbf{Step 4. Pointwise Exceedance.}
Absolute convergence of the Fourier series turns the coefficient
estimate into a pointwise bound 
\[
 \sup_{x\in\R}|v(t,x)-v_j(t,x)|
 \le\|\widehat v(t)-\widehat v_j(t)\|_{\ell^1(\Z^2)}\le\tfrac18,
 \qquad 0\le t\le t_j.
\]
By \eqref{eq:embedded-breather-focus} with $m=0$,
$|v_j(t_j,0)|=1+\sqrt2$. Hence
\[
 |v(t_j,0)|
 \ge|v_j(t_j,0)|-|v(t_j,0)-v_j(t_j,0)|
 \ge1+\sqrt2-\tfrac18>2.
 \qedhere
\]
\end{proof}

The tube stability in \autoref{lem:breather-tube} applies when the
normalized random initial coefficients $\widetilde d$ lie near $h_j$.
We therefore estimate the probability of this initial neighborhood.
Recall from \eqref{eq:log-normalized-initial-data} that
$\widetilde d=(\sqrt\varepsilon\,c(n)g_n)_n$.

By \autoref{def:gaussian}, $|g_0|^2$ has the \emph{exponential distribution}
with parameter $1$. Hence
\[
 \Ee|g_0|=\Ee\bigl[(|g_0|^2)^{1/2}\bigr]
 =\int_0^\infty s^{1/2}e^{-s}\,ds
 =\Gamma\!\left(\tfrac32\right)=\frac{\sqrt\pi}{2}.
\]
Tonelli's theorem gives
\begin{align}
 \Ee\|w_\beta\widetilde d\|_{\ell^1(\Z^2)}
 &=\sqrt\varepsilon\,\Ee|g_0|
     \sum_{n\in\Z^2}w_\beta(n)|c(n)|\nonumber\\
 &=\frac{\sqrt{\pi\varepsilon}}{2}
     \sum_{n\in\Z^2}(1+|n_1|)^{-2}(1+|n_2|)^{-2}
 <\infty.
 \label{eq:gaussian-weighted-mean}
\end{align}
In particular, $\|w_\beta\widetilde d\|_{\ell^1(\Z^2)}<\infty$ almost surely.
For a deterministic $h\in\mathcal H$, the Cauchy--Schwarz inequality gives
\begin{equation}\label{eq:gaussian-cost-weighted-embedding}
 \|w_\beta h\|_{\ell^1(\Z^2)}
 \le\|h\|_{\mathcal H}
     \left(\sum_{n\in\Z^2}w_\beta(n)^2|c(n)|^2\right)^{1/2}
 <\infty,\qquad h\in\mathcal H.
\end{equation}

\begin{lemma}[Gaussian Probability of a Translated Ball]
\label{lem:gaussian-tube}
Let $r>0$ and
\begin{equation}\label{eq:gaussian-weighted-ball}
 \mathsf B_r:=\left\{\mathsf e:
            \|w_\beta\mathsf e\|_{\ell^1(\Z^2)}<r\right\}.
\end{equation}
For every $\varepsilon>0$ and $h\in\mathcal H$, we have
\begin{equation}\label{eq:gaussian-tube}
 \Pp(\widetilde d\in h+\mathsf B_r)
 \ge e^{-\|h\|_{\mathcal H}^2/\varepsilon}
       \Pp(\widetilde d\in \mathsf B_r).
\end{equation}
\end{lemma}

\begin{remark}
\label{rem:gaussian-tube-cameron-martin}
The estimate \eqref{eq:gaussian-tube} also follows from the
\emph{Cameron--Martin formula} and the symmetry of $\mathsf B_r$;
see \cite{Bog98}. We give a direct proof below by \emph{finite-dimensional
truncation} and passage to the limit.
\end{remark}
\begin{proof}[Proof of \autoref{lem:gaussian-tube}]
Fix $\varepsilon>0$ and $h\in\mathcal H$.
By \eqref{eq:gaussian-weighted-mean} and
\eqref{eq:gaussian-cost-weighted-embedding},
$\|w_\beta\widetilde d\|_{\ell^1(\Z^2)}<\infty$ almost surely and
$\|w_\beta h\|_{\ell^1(\Z^2)}<\infty$.

\par\medskip\noindent\textbf{Step 1. Finite-Dimensional Density.}
Use $\Lambda_N$ from \eqref{eq:nf-truncation-set} with $\nu=2$.
The coefficients $c(n)$ are nonzero by the choice in
\autoref{subsec:counterexample-embedding}.
By \autoref{def:gaussian}, the real and imaginary parts of
$\widetilde d_n=\sqrt\varepsilon\,c(n)g_n$ have variance
$\varepsilon|c(n)|^2/2$ and are independent.
Thus $\widetilde d_n$ has density
\[
 \frac1{\pi\varepsilon|c(n)|^2}
       \exp\!\left(-\frac{|z_n|^2}{\varepsilon|c(n)|^2}\right)
\]
with respect to Lebesgue measure on $\C$.
Independence across $n$ gives the joint density
\begin{equation}\label{eq:gaussian-finite-density}
 p_N(z)=
 \left(\prod_{n\in\Lambda_N}\frac1{\pi\varepsilon|c(n)|^2}\right)
 \exp\!\left(-\frac1\varepsilon
               \sum_{n\in\Lambda_N}\frac{|z_n|^2}{|c(n)|^2}\right),
\end{equation}
where $z\in\C^{\Lambda_N}$ and $dz$ denotes
product Lebesgue measure in the real and imaginary coordinates.

\par\medskip\noindent\textbf{Step 2. Translation.}
Fix $0<s<r$ and define the closed ball
\begin{equation}\label{eq:gaussian-finite-ball}
 \overline{\mathsf B}_{s,N}:=\left\{y\in\C^{\Lambda_N}:
             \|w_\beta y\|_{\ell^1(\Lambda_N)}\le s\right\}.
\end{equation}
Let $f|_{\Lambda_N}:=(f_n)_{n\in\Lambda_N}$ be the restriction of $f$
to $\Lambda_N$.
By \eqref{eq:gaussian-finite-density} and the \emph{change
of variables} $z_n=y_n+h_n$ for $n\in\Lambda_N$, whose real
\emph{Jacobian determinant} is $1$, we have
\begin{align}
 \Pp\!\left(\widetilde d|_{\Lambda_N}\in h|_{\Lambda_N}+\overline{\mathsf B}_{s,N}\right)
 &=\int_{h|_{\Lambda_N}+\overline{\mathsf B}_{s,N}}p_N(z)\,dz\nonumber\\
 &=\int_{\overline{\mathsf B}_{s,N}}p_N(y+h|_{\Lambda_N})\,dy.
 \label{eq:gaussian-finite-shifted-integral}
\end{align}
For each $n\in\Lambda_N$,
\begin{equation}\label{eq:gaussian-square-expansion}
 |y_n+h_n|^2=|y_n|^2+|h_n|^2+2\Ree(\overline{h_n}y_n).
\end{equation}
Substituting \eqref{eq:gaussian-square-expansion} into
\eqref{eq:gaussian-finite-shifted-integral}, using
\eqref{eq:gaussian-finite-density}, and collecting terms, we obtain
\begin{align}
 \Pp\!\left(\widetilde d|_{\Lambda_N}\in h|_{\Lambda_N}+\overline{\mathsf B}_{s,N}\right)
 &=\exp\!\left(-\frac1\varepsilon
              \sum_{n\in\Lambda_N}\frac{|h_n|^2}{|c(n)|^2}\right)
\nonumber\\*
 &\quad\times           \int_{\overline{\mathsf B}_{s,N}}p_N(y)\exp\!\left(-\frac2\varepsilon\Ree
              \sum_{n\in\Lambda_N}
                    \frac{\overline{h_n}y_n}{|c(n)|^2}\right)\,dy.
 \label{eq:gaussian-finite-translation}
\end{align}

\par\medskip\noindent\textbf{Step 3. Symmetry.}
Under the substitution $y\mapsto-y$, the integration domain becomes
$-\overline{\mathsf B}_{s,N}=\overline{\mathsf B}_{s,N}$, Lebesgue measure is preserved, and
$p_N(-y)=p_N(y)$; only the linear exponent changes sign.
Relabeling the integration variable as $y$, we obtain
\begin{align*}
 &\int_{\overline{\mathsf B}_{s,N}}p_N(y)
       \exp\!\left(-\frac2\varepsilon\Ree
          \sum_{n\in\Lambda_N}\frac{\overline{h_n}y_n}{|c(n)|^2}\right)\,dy
\\
 &\quad=\int_{\overline{\mathsf B}_{s,N}}p_N(y)
       \exp\!\left(\frac2\varepsilon\Ree
          \sum_{n\in\Lambda_N}\frac{\overline{h_n}y_n}{|c(n)|^2}\right)\,dy.
\end{align*}
Averaging the two equal integrals and using
$\cosh x=(e^x+e^{-x})/2\ge1$ give
\begin{align}
 &\int_{\overline{\mathsf B}_{s,N}}p_N(y)
       \exp\!\left(-\frac2\varepsilon\Ree
          \sum_{n\in\Lambda_N}\frac{\overline{h_n}y_n}{|c(n)|^2}\right)\,dy
       \nonumber\\
 &\quad=\int_{\overline{\mathsf B}_{s,N}}p_N(y)
       \cosh\!\left(\frac2\varepsilon\Ree
          \sum_{n\in\Lambda_N}\frac{\overline{h_n}y_n}{|c(n)|^2}\right)\,dy
       \nonumber\\
 &\quad\ge\int_{\overline{\mathsf B}_{s,N}}p_N(y)\,dy   \nonumber\\
 &\quad=\Pp\!\left(\widetilde d|_{\Lambda_N}\in \overline{\mathsf B}_{s,N}\right).
 \label{eq:gaussian-finite-symmetry}
\end{align}
Combining \eqref{eq:gaussian-finite-translation} and
\eqref{eq:gaussian-finite-symmetry}, we obtain
\begin{align}
 \Pp\!\left(\widetilde d|_{\Lambda_N}\in h|_{\Lambda_N}+\overline{\mathsf B}_{s,N}\right)\ge
 \exp\!\left(-\frac1\varepsilon
             \sum_{n\in\Lambda_N}\frac{|h_n|^2}{|c(n)|^2}\right)
\Pp\!\left(\widetilde d|_{\Lambda_N}\in \overline{\mathsf B}_{s,N}\right).
 \label{eq:gaussian-finite-ball-bound}
\end{align}

\par\medskip\noindent\textbf{Step 4. Limit as $N\to\infty$.}
Fix $\varepsilon$, $h$ and $0<s<r$, and set
$\overline{\mathsf B}_s:=\{\mathsf e:\|w_\beta\mathsf e\|_{\ell^1(\Z^2)}\le s\}$.

\par\smallskip\noindent\emph{Events.}
As $\Lambda_N\uparrow\Z^2$, the partial weighted norms increase to
their full norms, so the events in \eqref{eq:gaussian-finite-ball-bound}
decrease and satisfy
\begin{align}
 &\bigcap_{N\ge1}\left\{
       \widetilde d|_{\Lambda_N}\in h|_{\Lambda_N}+\overline{\mathsf B}_{s,N}\right\}
=\{\widetilde d\in h+\overline{\mathsf B}_s\},
 \label{eq:gaussian-translated-closed-events}\\
 &\bigcap_{N\ge1}\left\{
       \widetilde d|_{\Lambda_N}\in \overline{\mathsf B}_{s,N}\right\}
=\{\widetilde d\in \overline{\mathsf B}_s\}.
 \label{eq:gaussian-centered-closed-events}
\end{align}

\par\smallskip\noindent\emph{Gaussian Factor.}
By \eqref{eq:gaussian-cost},
\begin{equation}\label{eq:gaussian-truncated-cost-limit}
 \sum_{n\in\Lambda_N}\frac{|h_n|^2}{|c(n)|^2}
 \uparrow\sum_{n\in\Z^2}\frac{|h_n|^2}{|c(n)|^2}
 =\|h\|_{\mathcal H}^2<\infty.
\end{equation}
Thus the exponential factor in \eqref{eq:gaussian-finite-ball-bound}
decreases to $e^{-\|h\|_{\mathcal H}^2/\varepsilon}$.

\par\smallskip\noindent\emph{Probability Limits.}
Since the events above decrease, the
continuity from above of the probability measure $\Pp$,
together with
\eqref{eq:gaussian-translated-closed-events} and
\eqref{eq:gaussian-centered-closed-events}, gives
\begin{align}
 \lim_{N\to\infty}\Pp\!\left(\widetilde d|_{\Lambda_N}
       \in h|_{\Lambda_N}+\overline{\mathsf B}_{s,N}\right)
 &=\Pp\!\left(\bigcap_{N\ge1}\left\{
       \widetilde d|_{\Lambda_N}\in
       h|_{\Lambda_N}+\overline{\mathsf B}_{s,N}\right\}\right)
       \nonumber\\*
 &=\Pp(\widetilde d\in h+\overline{\mathsf B}_s),
 \label{eq:gaussian-translated-probability-limit}\\
 \lim_{N\to\infty}\Pp\!\left(\widetilde d|_{\Lambda_N}
       \in\overline{\mathsf B}_{s,N}\right)
 &=\Pp\!\left(\bigcap_{N\ge1}\left\{
       \widetilde d|_{\Lambda_N}\in
       \overline{\mathsf B}_{s,N}\right\}\right)
       \nonumber\\*
 &=\Pp(\widetilde d\in\overline{\mathsf B}_s).
 \label{eq:gaussian-centered-probability-limit}
\end{align}
Using \eqref{eq:gaussian-truncated-cost-limit}--\eqref{eq:gaussian-centered-probability-limit},
we let $N\to\infty$ in \eqref{eq:gaussian-finite-ball-bound} to obtain
\begin{equation}\label{eq:gaussian-closed-ball-bound}
 \Pp(\widetilde d\in h+\overline{\mathsf B}_s)
 \ge e^{-\|h\|_{\mathcal H}^2/\varepsilon}
       \Pp(\widetilde d\in \overline{\mathsf B}_s).
\end{equation}

\par\medskip\noindent\textbf{Step 5. Limit as $s\uparrow r$.}
Having taken the limit in $N$, we now let $s\uparrow r$ through an
increasing sequence, with $\varepsilon$, $r$ and $h$ fixed.
By \eqref{eq:gaussian-weighted-ball} and the definition of $\overline{\mathsf B}_s$,
$\overline{\mathsf B}_s\uparrow \mathsf B_r$ and $h+\overline{\mathsf B}_s\uparrow h+\mathsf B_r$.
Thus both families of events increase monotonically to the corresponding
open-ball events. By the \emph{continuity from below} of the probability
measure $\Pp$, we obtain
\begin{align}
 \lim_{s\uparrow r}\Pp(\widetilde d\in h+\overline{\mathsf B}_s)
 &=\Pp(\widetilde d\in h+\mathsf B_r),
 \label{eq:gaussian-translated-open-limit}\\
 \lim_{s\uparrow r}\Pp(\widetilde d\in \overline{\mathsf B}_s)
 &=\Pp(\widetilde d\in \mathsf B_r).
 \label{eq:gaussian-centered-open-limit}
\end{align}
Taking $s\uparrow r$ in \eqref{eq:gaussian-closed-ball-bound}
proves \eqref{eq:gaussian-tube}.
\end{proof}

\begin{theorem}[Failure of the Gaussian Rate Function at Longer Logarithmic Time]
\label{thm:longtime-rate-obstruction}
Let $u$ be the maximal solution of
\eqref{eq:qp-cubic-NLS} with initial data
\eqref{eq:qp-randomized-initial-data}.
For the parameters in \autoref{subsec:counterexample-embedding}, let
\begin{equation}\label{eq:counterexample-alpha}
 1<\alpha<1+\frac1{4500},\qquad
 \varepsilon_j=\mu_j^{1/(\alpha-1)},\qquad
 t_j=15\mu_j^{-1}|\log\mu_j|.
\end{equation}
Then
\begin{equation}\label{eq:counterexample-LDP-bound}
 \liminf_{j\to\infty}\varepsilon_j\log\Pp
 \bigl(\mathscr A_{\varepsilon_j,2}^{\|u(t_j)\|_\infty}\bigr)\ge-1.
\end{equation}
Consequently, the Gaussian rate function $4/S$ cannot hold uniformly on
$[0,15\mu_\varepsilon^{-1}|\log\mu_\varepsilon|]$.
\end{theorem}

\begin{proof}[Proof of \autoref{thm:longtime-rate-obstruction}]
Throughout the proof, take $\varepsilon=\varepsilon_j$, so
$\widetilde d_n=\sqrt{\varepsilon_j}\,c(n)g_n$ and
$\varepsilon_j^{\alpha-1}=\mu_j$.
We claim that, for all sufficiently large $j$,
\begin{align}
 \mathscr A_{\varepsilon_j,2}^{\|u(t_j)\|_\infty}
 &\supseteq
 \mathscr A_{\varepsilon_j,2}^{\|u(t_j)\|_\infty}
 \cap G_{\delta_{\varepsilon_j}}\nonumber\\
 &\supseteq\{\widetilde d\in h_j+B_j\},
 \label{eq:counterexample-tube-inclusion}
\end{align}
where
\begin{align*}
 B_j&:=\{\mathsf e:
          \|w_\beta\mathsf e\|_{\ell^1(\Z^2)}<r_j\},\\
 h_j+B_j&:=\{\tilde d:
          \|w_\beta(\tilde d-h_j)\|_{\ell^1(\Z^2)}<r_j\};
\end{align*}
see \autoref{fig:initial-neighborhood-translation}. 
\begin{figure}[htbp]
\centering
\setlength{\abovecaptionskip}{4pt}
\scalebox{0.8}{
\begin{tikzpicture}[x=1cm,y=1cm,font=\small,line width=.55pt]
 \definecolor{neighborhoodblue}{RGB}{35,91,156}
 \definecolor{initialorange}{RGB}{177,82,30}
 \foreach \ballx in {-2.65,2.65} {
  \filldraw[draw=neighborhoodblue,fill=neighborhoodblue!6]
    (\ballx,0) circle[radius=1.05];
  \draw[dashed,black!55] (\ballx,0) -- ++(1.05,0)
    node[midway,below=3pt,text=black] {$r_j$};
  \fill (\ballx,0) circle[radius=1.5pt];
 }
 \node[below left=1pt,inner sep=1pt] at (-2.65,0) {$0$};
 \node[below left=1pt,inner sep=1pt] at (2.65,0) {$h_j$};
 \draw[->,>=stealth] (-1.3,0) -- (1.3,0)
    node[midway,above=5pt] {$+h_j$};
 \fill[initialorange] (2.45,.45) circle[radius=1.8pt];
 \node[left=3pt,inner sep=1pt,text=initialorange] at (2.45,.45)
    {$\widetilde d$};
 \node[anchor=south,font=\footnotesize,inner sep=0pt]
    (randominitiallabel) at (3.10,1.48) {\strut Random};
 \node[anchor=east,font=\footnotesize,inner sep=0pt]
    at ([xshift=-2.5pt]randominitiallabel.west) {\strut Normalized};
 \node[anchor=west,font=\footnotesize,inner sep=0pt]
    at ([xshift=2.5pt]randominitiallabel.east) {\strut Initial Data};
 \draw[dashed,->,>=stealth,black!65] (2.45,.45)
    -- ([yshift=-2pt]randominitiallabel.south);
 \node at (-2.65,-1.6) {$B_j$};
 \node at (2.65,-1.6) {$h_j+B_j$};
\end{tikzpicture}

}
\caption{
The ball $B_j$ is centered at $0$, and $h_j+B_j$ is centered at $h_j$.
Both have radius $r_j$ in the norm
$\|w_\beta\cdot\|_{\ell^1(\Z^2)}$.
}
\label{fig:initial-neighborhood-translation}
\end{figure}
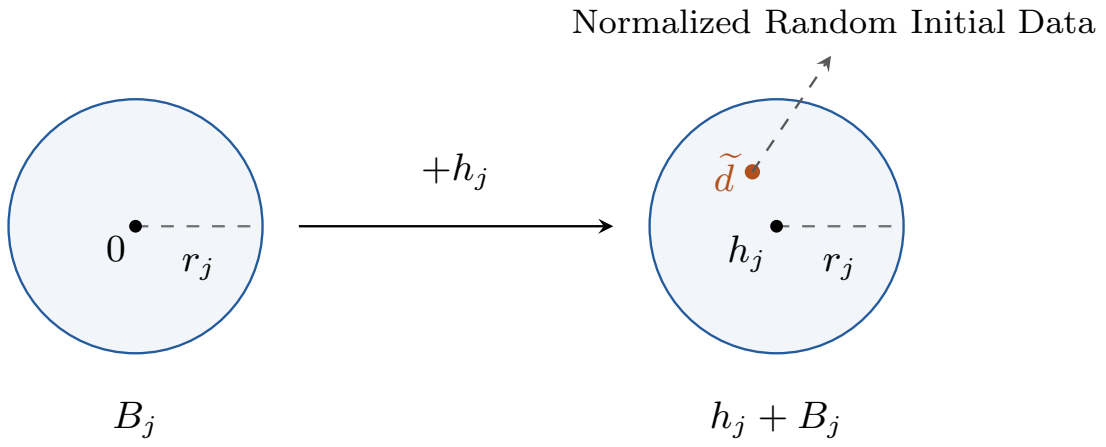
\noindent The corresponding events are
\begin{align}
 \{\widetilde d\in B_j\}
 &=\{\xi\in\Omega:
       \|w_\beta\widetilde d(\xi)\|_{\ell^1(\Z^2)}<r_j\},
 \label{eq:counterexample-centered-neighborhood-event}\\
 \{\widetilde d\in h_j+B_j\}
 &=\{\xi\in\Omega:
       \|w_\beta(\widetilde d(\xi)-h_j)\|_{\ell^1(\Z^2)}<r_j\}.
 \label{eq:counterexample-initial-neighborhood}
\end{align}
The first inclusion in \eqref{eq:counterexample-tube-inclusion}
is immediate. To prove the second, fix any $\xi\in\Omega$ such that
$\widetilde d(\xi)\in h_j+B_j$ and set
$\mathsf e(\xi)=\widetilde d(\xi)-h_j$.

First, we show that
$\xi\in\mathscr A_{\varepsilon_j,2}^{\|u(t_j)\|_\infty}$.
Since $w_\beta(n)\ge1$,
\[
 \|\mathsf e(\xi)\|_{\ell^1(\Z^2)}
 \le\|w_\beta\mathsf e(\xi)\|_{\ell^1(\Z^2)}<r_j.
\]
By \autoref{lem:breather-tube}, the solution $v(t,x;\xi)$ of
\eqref{eq:embedded-breather-NLS} with initial coefficients
$h_j+\mathsf e(\xi)=\widetilde d(\xi)$ exists beyond $t_j$.
Since $\varepsilon_j^{\alpha-1}=\mu_j$, the rescaled field
$\varepsilon_j^{-1/2}v(t,x;\xi)$ solves
\eqref{eq:qp-cubic-NLS}--\eqref{eq:qp-randomized-initial-data}.
Uniqueness identifies it with $u(t,x;\xi)$. In particular,
$u$ exists beyond $t_j$, and
\eqref{eq:breather-tube-stability} gives
\begin{align}
 \sqrt{\varepsilon_j}\|u(t_j,\cdot;\xi)\|_\infty
 &\ge\sqrt{\varepsilon_j}|u(t_j,0;\xi)|\nonumber\\*
 &=|v(t_j,0;\xi)|\ge1+\sqrt2-\tfrac18>2.
 \label{eq:counterexample-tube-amplitude}
\end{align}
Thus
$\xi\in\mathscr A_{\varepsilon_j,2}^{\|u(t_j)\|_\infty}$.

Next, we verify that $\xi\in G_{\delta_{\varepsilon_j}}$.
By \eqref{eq:precursor-weighted-norm}, \eqref{eq:tube-radius},
and \eqref{eq:counterexample-variance-cutoff}, for all sufficiently
large $j$,
\[
 \|w_\beta\widetilde d(\xi)\|_{\ell^1(\Z^2)}
 \le\|w_\beta h_j\|_{\ell^1(\Z^2)}
        +\|w_\beta\mathsf e(\xi)\|_{\ell^1(\Z^2)}
 <2<\sqrt{1+4/S}.
\]
Thus, for every $n\in\Z^2$,
\[
 w_\beta(n)|\widetilde d_n(\xi)|
 \le\|w_\beta\widetilde d(\xi)\|_{\ell^1(\Z^2)}
 <\sqrt{1+4/S}.
\]
Using $c(n)=(1+|n_1|)^{-12}(1+|n_2|)^{-12}$,
$\beta=(10,10)$ and $\kappa=(2,2)$, we obtain
\begin{align}
 |g_n(\xi)|
 &=\frac{|\widetilde d_n(\xi)|}{\sqrt{\varepsilon_j}|c(n)|}\nonumber\\
 &\le\sqrt{\frac{1+4/S}{\varepsilon_j}}\,
             \frac{w_{-\beta}(n)}{|c(n)|}
 =\delta_{\varepsilon_j}^{-1/2}w_\kappa(n).
 \label{eq:counterexample-good-event-bound}
\end{align}
The last equality uses the definition of $\delta_{\varepsilon_j}$ in
\eqref{eq:counterexample-variance-cutoff}.
Thus \eqref{eq:Omega-delta} gives
$\xi\in G_{\delta_{\varepsilon_j}}$, completing the proof of
\eqref{eq:counterexample-tube-inclusion}.

We now estimate $\Pp(\widetilde d\in h_j+B_j)$.
By \autoref{lem:precursor-bounds}, $h_j\in\mathcal H$.
Applying \autoref{lem:gaussian-tube} with $h=h_j$, $r=r_j$
and $\varepsilon=\varepsilon_j$ gives
\begin{equation}\label{eq:counterexample-tube-probability}
 \Pp(\widetilde d\in h_j+B_j)
 \ge e^{-\|h_j\|_{\mathcal H}^2/\varepsilon_j}
       \Pp(\widetilde d\in B_j).
\end{equation}
To estimate the last factor, \emph{Markov's inequality} and
\eqref{eq:gaussian-weighted-mean} give
\begin{align}
 \Pp(\widetilde d\notin B_j)
 &=\Pp\!\left(\|w_\beta\widetilde d\|_{\ell^1(\Z^2)}
                    \ge r_j\right)\nonumber\\*
 &\le\frac{\Ee\|w_\beta\widetilde d\|_{\ell^1(\Z^2)}}{r_j}
   \nonumber\\
 &=\frac{\sqrt{\varepsilon_j}}{r_j}\Ee|g_0|
      \sum_{n\in\Z^2}w_\beta(n)|c(n)|\nonumber\\
 &=\frac{\sqrt{\pi\varepsilon_j}}{2r_j}
      \left(\sum_{m\in\Z}(1+|m|)^{-2}\right)^2\nonumber\\
 &\lesssim\frac{\sqrt{\varepsilon_j}}{r_j}
 =8\mu_j^{1/(2(\alpha-1))-2250}\longrightarrow0.
 \label{eq:counterexample-centered-tube}
\end{align}
The last equality uses \eqref{eq:tube-radius}, and the exponent
$1/(2(\alpha-1))-2250$ is positive by
\eqref{eq:counterexample-alpha}. Hence
$\Pp(\widetilde d\in B_j)\to1$, so
\begin{equation}\label{eq:counterexample-centered-probability}
 \varepsilon_j\log\Pp(\widetilde d\in B_j)\longrightarrow0.
\end{equation}
Using \eqref{eq:counterexample-tube-inclusion} and
\eqref{eq:counterexample-tube-probability}, we take logarithms
and multiply by $\varepsilon_j$ to obtain
\begin{align}
 \varepsilon_j\log\Pp
   \bigl(\mathscr A_{\varepsilon_j,2}^{\|u(t_j)\|_\infty}\bigr)
 &\ge\varepsilon_j\log\Pp
   \bigl(G_{\delta_{\varepsilon_j}}\cap
    \mathscr A_{\varepsilon_j,2}^{\|u(t_j)\|_\infty}\bigr)
   \nonumber\\
 &\ge\varepsilon_j\log
                 \Pp(\widetilde d\in h_j+B_j)\nonumber\\
 &\ge-\|h_j\|_{\mathcal H}^2
           +\varepsilon_j\log\Pp(\widetilde d\in B_j)\nonumber\\
 &=-1+o(1).
 \label{eq:counterexample-log-probability}
\end{align}
The last equality follows from
\eqref{eq:precursor-gaussian-cost} and
\eqref{eq:counterexample-centered-probability}.
Taking the lower limit proves \eqref{eq:counterexample-LDP-bound}.
Since $\mu_{\varepsilon_j}=\mu_j$, the times $t_j$ lie in the
stated logarithmic interval. The strict inequality $-1>-4/S$ from
\eqref{eq:counterexample-variance-cutoff} therefore rules out the
Gaussian rate function uniformly on that interval.
\end{proof}

\section{Failure of the Initial Gaussian Rate Function at the Critical Time Scale with a Large Logarithmic Correction Prefactor: General Case}
\label{sec:general-failure}

The special construction in \autoref{sec:longtime-obstruction}
suggests how to extend the failure result to arbitrary finite frequency rank $\nu\ge2$
and every fixed $\omega\in\mathrm{DC}(\gamma,\tau)$.
Tildes distinguish the sequences below from those in the Pell construction.
The key is to replace the \emph{Pell construction} in
\eqref{eq:pell-lattice-index}--\eqref{eq:pell-physical-frequency}
by a \emph{Dirichlet-type approximation} of the linear form
$\langle n,\omega\rangle$; see \cite[subsection 11.3]{HW08}.
For this construction, we need both $\widetilde{k}_j=\langle \widetilde{n}_j,\omega\rangle\to0$, as in
\eqref{eq:pell-physical-frequency}, and control of the growth of
$\widetilde{n}_j$ through
$\langle \widetilde{n}_j\rangle\lesssim_\omega |\widetilde{k}_j|^{-1/(\nu-1)}$;
see \eqref{eq:general-small-frequencies}.
The first property permits small frequency embedding of the breather;
the second controls the resulting weight factors and ensures that
the initial Gaussian cost tends to one;
see \eqref{eq:general-initial-cost}.
The tube and probability arguments then apply as before.

\subsection{Construction of Small Frequencies and Embedded Breather}
\label{subsec:general-embedding}
\begin{lemma}[Small Frequencies in a Fixed Module]
\label{lem:general-small-frequencies}
There exist $\widetilde{n}_j\in\Z^\nu\setminus\{0\}$ such that, with
$\widetilde{k}_j=\langle \widetilde{n}_j,\omega\rangle$ and $\widetilde{\mu}_j=\widetilde{k}_j^2$,
\begin{equation}\label{eq:general-small-frequencies}
 0<|\widetilde{k}_j|\longrightarrow0,
 \qquad
 \langle \widetilde{n}_j\rangle\lesssim_\omega
       \widetilde{\mu}_j^{-1/[2(\nu-1)]}.
\end{equation}
\end{lemma}

\begin{proof}[Proof of \autoref{lem:general-small-frequencies}]
      Fix an integer $Q\ge1$. There are $(Q+1)^\nu$ integer vectors
$a\in\{0,\ldots,Q\}^\nu$.
\leavevmode\par\nobreak\smallskip
\noindent\emph{Distinct Values.}
 Their values $\langle a,\omega\rangle$
are distinct: equality for two vectors $a,b$ would give
$\langle a-b,\omega\rangle=0$, and rational independence would
force $a=b$.

\medskip
\noindent\emph{Total Span.}
For each coordinate, $0\le a_r\le Q$ gives
$0\le a_r\omega_r\le Q\omega_r$ if $\omega_r\ge0$, and
$Q\omega_r\le a_r\omega_r\le0$ if $\omega_r<0$.
Equivalently,
\[
 Q\min\{\omega_r,0\}\le a_r\omega_r
 \le Q\max\{\omega_r,0\}.
\]
Summing over $r=1,\ldots,\nu$ yields
\begin{equation}\label{eq:general-frequency-interval}
 Q\sum_{r=1}^\nu\min\{\omega_r,0\}
 \le\langle a,\omega\rangle
 \le Q\sum_{r=1}^\nu\max\{\omega_r,0\}.
\end{equation}
The difference between the two endpoints is
\begin{equation}\label{eq:general-frequency-interval-length}
 Q\sum_{r=1}^\nu\max\{\omega_r,0\}
 -Q\sum_{r=1}^\nu\min\{\omega_r,0\}
 =Q\sum_{r=1}^\nu|\omega_r|,
\end{equation}
since $\max\{\omega_r,0\}-\min\{\omega_r,0\}=|\omega_r|$. This implies that the $(Q+1)^\nu$ values $\langle a,\omega\rangle$ lie in an interval of length
$Q\sum_{r=1}^\nu|\omega_r|$.

\noindent\emph{A Small Consecutive Gap.}
Order the values as $x_1<\cdots<x_{(Q+1)^\nu}$, where
$x_s=\langle a_s,\omega\rangle$ and
$a_s\in\{0,\ldots,Q\}^\nu$.
By \eqref{eq:general-frequency-interval} and
\eqref{eq:general-frequency-interval-length}, we have
\begin{equation}\label{eq:general-frequency-span}
 x_{(Q+1)^\nu}-x_1\le Q\sum_{r=1}^\nu|\omega_r|.
\end{equation}
Notice that there are $(Q+1)^\nu-1$ consecutive gaps.
By the \emph{pigeonhole principle}, we claim that there exists an index $1\le s_\ast<(Q+1)^\nu$ such that
\begin{equation}\label{eq:general-consecutive-gap}
 0<x_{s_\ast+1}-x_{s_\ast}
 \le\frac{Q\sum_{r=1}^\nu|\omega_r|}{(Q+1)^\nu-1}
 \le\left(\sum_{r=1}^\nu|\omega_r|\right)Q^{-(\nu-1)}.
\end{equation}
If not, then for every $s$,
\[
 x_{s+1}-x_s>
 \frac{Q}{(Q+1)^\nu-1}\sum_{r=1}^\nu|\omega_r|.
\]
Summing over $s$ would give
\[
 x_{(Q+1)^\nu}-x_1
 =\sum_{s=1}^{(Q+1)^\nu-1}(x_{s+1}-x_s)
 >Q\sum_{r=1}^\nu|\omega_r|,
\]
contrary to \eqref{eq:general-frequency-span}.
The last inequality in \eqref{eq:general-consecutive-gap} follows from
$(Q+1)^\nu-1\ge Q^\nu$.

\medskip
\noindent\emph{Construction of the Integer Vectors.}
Set
\begin{equation}\label{eq:general-integer-vector}
 \widetilde{n}_Q=a_{s_\ast+1}-a_{s_\ast}\ne0.
\end{equation}
By \eqref{eq:general-integer-vector}, every coordinate of
$\widetilde{n}_Q$ lies in $[-Q,Q]$, and
$\langle \widetilde{n}_Q,\omega\rangle=x_{s_\ast+1}-x_{s_\ast}$. Therefore
\begin{equation}\label{eq:general-pigeonhole-frequency}
 |\widetilde{n}_Q|_\infty\le Q,
 \qquad
 0<|\langle \widetilde{n}_Q,\omega\rangle|
 \le C_\omega Q^{-(\nu-1)}.
\end{equation}
Here $C_\omega=\sum_{r=1}^\nu|\omega_r|$. Since $\nu\ge2$, these positive frequencies tend to zero as
$Q\to\infty$. We may therefore choose $Q_j\to\infty$ so that
$|\langle \widetilde{n}_{Q_j},\omega\rangle|$ strictly decreases to zero.
Define
\begin{equation}\label{eq:general-integer-subsequence}
 \widetilde{n}_j:=\widetilde{n}_{Q_j}.
\end{equation}
The vectors in \eqref{eq:general-integer-subsequence} are distinct.
Set $\widetilde{k}_j=\langle \widetilde{n}_j,\omega\rangle$ and $\widetilde{\mu}_j=\widetilde{k}_j^2$.
By \eqref{eq:general-pigeonhole-frequency} and
\eqref{eq:general-integer-subsequence},
$Q_j\le C_\omega^{1/(\nu-1)}|\widetilde{k}_j|^{-1/(\nu-1)}$.
Combining this with $|\widetilde{n}_j|_1\le\nu|\widetilde{n}_j|_\infty\le\nu Q_j$
and $Q_j\ge1$, we obtain
\begin{align*}
 \langle \widetilde{n}_j\rangle
 &=1+|\widetilde{n}_j|_1\\
 &\le1+\nu Q_j\\
 &\le(\nu+1)Q_j\\
 &\le(\nu+1)C_\omega^{1/(\nu-1)}|\widetilde{k}_j|^{-1/(\nu-1)}\\
 &\lesssim_\omega\widetilde{\mu}_j^{-1/[2(\nu-1)]}.
\end{align*}
This proves \eqref{eq:general-small-frequencies}.
\end{proof}

We now choose the decay and time parameters so that the Gaussian cost tends to one. Let $\rho=(\rho_1,\ldots,\rho_\nu)\in(0,\infty)^\nu$ and set
\begin{equation}\label{eq:general-coefficients}
 c(n)=\prod_{r=1}^\nu(1+|n_r|)^{-\rho_r}.
\end{equation}
Let $\kappa\in(0,\infty)^\nu$ be the parameter in the exceptional set; see
\autoref{def:exceptional}. In this section, assume $\kappa_r>1$ for every $r$ and let
\begin{equation}\label{eq:general-decay-parameters}
 \beta_r:=\rho_r-\kappa_r>\max\{4\tau+2,\nu+2\tau\},
 \qquad r=1,\ldots,\nu,
\end{equation}
consistent with the notation in \autoref{cor:randomized-decay}.
In this setting, we have
\begin{equation}\label{eq:general-variance}
 S=\sum_{n\in\Z^\nu}|c(n)|^2<\frac43.
\end{equation}
Choose a logarithmic time prefactor $\mathsf{C}_{\mathrm{fail}}>0$
and $\alpha>1$ satisfying
\begin{equation}\label{eq:general-parameter-compatibility}
 \frac{\sum_r\rho_r}{2(\nu-1)}<\mathsf{C}_{\mathrm{fail}}<\frac1{300(\alpha-1)},
 \qquad \alpha>1.
\end{equation}
For every fixed $\rho$, this interval is nonempty if
\begin{equation}\label{eq:general-alpha-range}
 1<\alpha<1+\frac{\nu-1}{150\sum_r\rho_r}.
\end{equation}
Discarding finitely many terms, assume $0<\widetilde{\mu}_j<1$ and define
\begin{equation}\label{eq:general-embedding-parameters}
 \widetilde{T}_j=\mathsf{C}_{\mathrm{fail}}|\log\widetilde{\mu}_j|,\qquad
 \widetilde{t}_j=\widetilde{\mu}_j^{-1}\widetilde{T}_j,\qquad
 \widetilde{\varepsilon}_j=\widetilde{\mu}_j^{1/(\alpha-1)}.
\end{equation}
Thus $\mu_{\widetilde\varepsilon_j}=\widetilde\mu_j$, with
$\mu_\varepsilon=\varepsilon^{\alpha-1}$ as before.

In \eqref{eq:embedded-breather} and
\eqref{eq:embedded-breather-initial-coefficients}, replace
$(n_j,k_j,\mu_j)$ from
\eqref{eq:pell-lattice-index}--\eqref{eq:pell-coupling}
by $(\widetilde n_j,\widetilde k_j,\widetilde\mu_j)$ from
\autoref{lem:general-small-frequencies}, and replace $(T_j,t_j)$ from
\eqref{eq:pell-scaled-waiting-time}--\eqref{eq:pell-focusing-time}
by $(\widetilde T_j,\widetilde t_j)$ from
\eqref{eq:general-embedding-parameters}.
Denote the resulting breather and initial coefficients by
$\widetilde v_j$ and $\widetilde h_j$.
Thus $15$ is replaced by $\mathsf C_{\mathrm{fail}}$;
we also replace $\Z^2$ by $\Z^\nu$, set
$\varepsilon=\widetilde\varepsilon_j$, and use $c$, $\kappa$ and $\beta$
from \eqref{eq:general-coefficients}--\eqref{eq:general-decay-parameters}.
Likewise, write $\mathsf B_j=\{\mathsf e:
\|w_\beta\mathsf e\|_{\ell^1(\Z^\nu)}<r_j\}$ for the initial
neighborhood, with radius $r_j$ specified in \eqref{eq:general-tube-radius}.
The properties in \autoref{lem:embedded-breather-properties} carry over:
$\widetilde k_j^2=\widetilde\mu_j$ gives the rescaled equation, and
$m\mapsto m\widetilde n_j$ is injective.
The weighted initial bounds, Gaussian cost and neighborhood
probability estimates are verified below with these parameters.

\subsection{Failure of the Initial Gaussian Rate Function}
\label{subsec:general-tube-failure}
\label{subsec:general-rate-failure}

\begin{theorem}[Failure for General Diophantine Frequencies]
\label{thm:general-failure}
Fix any $\nu\ge2$ and $\omega\in\mathrm{DC}(\gamma,\tau)$.
Assume \eqref{eq:general-coefficients}--\eqref{eq:general-parameter-compatibility}.
Let $u$ be the maximal solution of
\eqref{eq:qp-cubic-NLS}--\eqref{eq:qp-randomized-initial-data}.
For the sequences in \eqref{eq:general-embedding-parameters},
\begin{equation}\label{eq:general-failure-lower-bound}
 \liminf_{j\to\infty}\widetilde{\varepsilon}_j\log\Pp\!\left(
 \mathscr A_{\widetilde{\varepsilon}_j,2}^{\|u(\widetilde{t}_j)\|_\infty}
 \right)\ge-1>-\frac4S.
\end{equation}
Thus uniform persistence fails on
$[0,\mathsf{C}_{\mathrm{fail}}\mu_\varepsilon^{-1}|\log\mu_\varepsilon|]$.
\end{theorem}

\begin{proof}[Proof of \autoref{thm:general-failure}]
We prove the general case by following the proof strategy for the special
case in \autoref{sec:longtime-obstruction}, outlined in
\autoref{fig:proof-routes}, with the replacements specified after
\eqref{eq:general-embedding-parameters}.
Throughout, $\varepsilon=\widetilde{\varepsilon}_j$,
$\widetilde d_n=\sqrt{\widetilde{\varepsilon}_j}\,c(n)g_n$ and
$\widetilde{\varepsilon}_j^{\alpha-1}=\widetilde{\mu}_j$.

\medskip
\noindent\emph{Breather Embedding.}
\autoref{lem:embedded-breather-properties} gives
\begin{equation}\label{eq:general-reference-breather}
 \widetilde v_j(t,x)=u_{\mathrm{AB}}(\widetilde{\mu}_jt-\widetilde{T}_j,\widetilde{k}_jx),
 \qquad \widetilde h_j(n)=\widehat{\widetilde v}_j(0,n),
\end{equation}
and the counterparts of \eqref{eq:embedded-breather-uniform-fourier}
and \eqref{eq:embedded-breather-focus} are
\begin{equation}\label{eq:general-reference-bounds}
 \sup_{t\in\R}\|\widehat{\widetilde v}_j(t)\|_{\ell^1(\Z^\nu)}\le3+\sqrt2,
 \qquad |\widetilde v_j(\widetilde{t}_j,0)|=1+\sqrt2.
\end{equation}
Here the only change in the embedding is the small frequency
sequence supplied by \autoref{lem:general-small-frequencies}.

\medskip
\noindent\emph{Weighted Initial Bound.}
We first show that
\begin{align}
 \|w_\beta \widetilde h_j\|_{\ell^1(\Z^\nu)}
 &=1+\mathcal{O}\!\left(\widetilde{\mu}_j^{\mathsf{C}_{\mathrm{fail}}-\sum_r\beta_r/[2(\nu-1)]}\right),
 \label{eq:general-weighted-initial-bound}
\end{align}
By \eqref{eq:embedded-breather-fourier}, split the weighted sum as
\begin{equation}\label{eq:general-weighted-decomposition}
 \|w_\beta \widetilde h_j\|_{\ell^1(\Z^\nu)}
 =|\widehat u_{\mathrm{AB}}(-\widetilde{T}_j,0)|
 +\sum_{m\ne0}w_\beta(m\widetilde{n}_j)
       |\widehat u_{\mathrm{AB}}(-\widetilde{T}_j,m)|.
\end{equation}
The zero mode is $1+\mathcal{O}(\widetilde{\mu}_j^{2\mathsf{C}_{\mathrm{fail}}})$ by
\eqref{eq:breather-initial-zero-mode}.
For $|m|\ge1$,
\[
 w_\beta(m\widetilde{n}_j)
 \le |m|^{\sum_r\beta_r}\langle \widetilde{n}_j\rangle^{\sum_r\beta_r}.
\]
For any fixed $a\ge0$, $\sum_{m\ge1}m^a r^m\lesssim_a r$
when $0<r\le1/2$, since
$\sum_{m\ge1}m^a r^{m-1}\le\sum_{m\ge1}m^a2^{1-m}<\infty$.
Using \eqref{eq:breather-initial-Fourier} and then
\eqref{eq:general-small-frequencies}, we obtain
\begin{align}
 \sum_{m\ne0}w_\beta(m\widetilde{n}_j)|\widehat u_{\mathrm{AB}}(-\widetilde{T}_j,m)|
 &\lesssim\langle \widetilde{n}_j\rangle^{\sum_r\beta_r}
       \sum_{m\ge1}m^{\sum_r\beta_r}(\sqrt2\widetilde{\mu}_j^{\mathsf{C}_{\mathrm{fail}}})^m
 \nonumber\\
 &\lesssim\langle \widetilde{n}_j\rangle^{\sum_r\beta_r}\widetilde{\mu}_j^{\mathsf{C}_{\mathrm{fail}}}
 \lesssim\widetilde{\mu}_j^{\mathsf{C}_{\mathrm{fail}}-\sum_r\beta_r/[2(\nu-1)]}.
 \label{eq:general-weighted-nonzero-bound}
\end{align}
Substitution into \eqref{eq:general-weighted-decomposition} proves
\eqref{eq:general-weighted-initial-bound}.

\medskip
\noindent\emph{Gaussian Cost.}
The counterpart of \eqref{eq:precursor-gaussian-cost} is
\begin{align}
 \|\widetilde h_j\|_{\mathcal H}^2
 &=1+\mathcal{O}\!\left(\widetilde{\mu}_j^{2\mathsf{C}_{\mathrm{fail}}-\sum_r\rho_r/(\nu-1)}\right).
 \label{eq:general-initial-cost}
\end{align}
Since $c(0)=1$, the cost decomposes as
\begin{equation}\label{eq:general-cost-decomposition}
 \|\widetilde h_j\|_{\mathcal H}^2
 =|\widehat u_{\mathrm{AB}}(-\widetilde{T}_j,0)|^2
 +\sum_{m\ne0}
       \frac{|\widehat u_{\mathrm{AB}}(-\widetilde{T}_j,m)|^2}{|c(m\widetilde{n}_j)|^2}.
\end{equation}
The first term is $1+\mathcal{O}(\widetilde{\mu}_j^{2\mathsf{C}_{\mathrm{fail}}})$.
For the second, \eqref{eq:general-coefficients} gives
$|c(m\widetilde{n}_j)|^{-2}\le |m|^{2\sum_r\rho_r}\langle \widetilde{n}_j\rangle^{2\sum_r\rho_r}$.
The same series estimate, now with ratio $2\widetilde{\mu}_j^{2\mathsf{C}_{\mathrm{fail}}}$, yields
\begin{align}
 \sum_{m\ne0}
       \frac{|\widehat u_{\mathrm{AB}}(-\widetilde{T}_j,m)|^2}{|c(m\widetilde{n}_j)|^2}
 &\lesssim\langle \widetilde{n}_j\rangle^{2\sum_r\rho_r}
       \sum_{m\ge1}m^{2\sum_r\rho_r}(2\widetilde{\mu}_j^{2\mathsf{C}_{\mathrm{fail}}})^m
 \nonumber\\
 &\lesssim\langle \widetilde{n}_j\rangle^{2\sum_r\rho_r}\widetilde{\mu}_j^{2\mathsf{C}_{\mathrm{fail}}}
 \lesssim\widetilde{\mu}_j^{2\mathsf{C}_{\mathrm{fail}}-\sum_r\rho_r/(\nu-1)}.
 \label{eq:general-nonzero-cost}
\end{align}
Substituting into \eqref{eq:general-cost-decomposition} proves
\eqref{eq:general-initial-cost}. Since $\sum_r\beta_r<\sum_r\rho_r$, both exponents in \eqref{eq:general-weighted-initial-bound}
and \eqref{eq:general-initial-cost} are positive by \eqref{eq:general-parameter-compatibility}.
In particular, $\widetilde h_j\in\mathcal H$ for all sufficiently large $j$,
and $\|\widetilde h_j\|_{\mathcal H}^2\to1$.

\medskip
\noindent\emph{Tube Stability.}
The radius corresponding to \eqref{eq:tube-radius} is
\begin{equation}\label{eq:general-tube-radius}
 r_j=\tfrac18e^{-150\widetilde{T}_j}
     =\tfrac18\widetilde{\mu}_j^{150\mathsf C_{\mathrm{fail}}}.
\end{equation}
For $\mathsf e\in\mathsf B_j$, the Picard construction of
\autoref{lem:breather-tube} gives a local solution $v$ of \eqref{eq:embedded-breather-NLS} with initial
coefficients $\widetilde h_j+\mathsf e$. The linear isometry and cubic
$\ell^1$ convolution estimate have the same constants on $\Z^\nu$.
Consequently, on a bootstrap interval where the difference norm is
at most $1/8$, \eqref{eq:general-reference-bounds} bounds both
coefficient norms by $4+\sqrt2$. The proof of
\eqref{eq:breather-tube-gronwall} then gives, for $t\le \widetilde{t}_j$,
\begin{align}
 \sup_{0\le s\le t}
 \|\widehat v(s)-\widehat{\widetilde v}_j(s)\|_{\ell^1(\Z^\nu)}
 &\le\|\mathsf e\|_{\ell^1(\Z^\nu)}e^{150\widetilde{\mu}_jt}
 \nonumber\\
 &\le\|\mathsf e\|_{\ell^1(\Z^\nu)}e^{150\widetilde{T}_j}
 <r_je^{150\widetilde{T}_j}=\tfrac18.
 \label{eq:general-tube-stability}
\end{align}
For each fixed $\mathsf e$, the strict inequality has a positive
margin independent of $t\le \widetilde{t}_j$. Continuity therefore closes the
bootstrap, and the uniform coefficient bound permits extension
beyond $\widetilde{t}_j$, exactly as in \autoref{lem:breather-tube}.
In particular,
\begin{equation}\label{eq:general-tube-amplitude}
 |v(\widetilde{t}_j,0)|\ge1+\sqrt2-\tfrac18>2.
\end{equation}

\medskip
\noindent\emph{Neighborhood Inclusion.}
Set $\delta_{\widetilde{\varepsilon}_j}=\widetilde{\varepsilon}_j/(1+4/S)$.
The counterpart of \eqref{eq:counterexample-tube-inclusion} is
\begin{equation}\label{eq:general-event-inclusion}
 \{\widetilde d\in \widetilde h_j+\mathsf B_j\}
 \subseteq G_{\delta_{\widetilde{\varepsilon}_j}}
 \cap\mathscr A_{\widetilde{\varepsilon}_j,2}^{\|u(\widetilde{t}_j)\|_\infty}.
\end{equation}
Indeed, fix a sample point $\xi$ in the event on the left, and put
$\mathsf e=\widetilde d(\xi)-\widetilde h_j\in\mathsf B_j$.
Uniqueness identifies the solution above with
$v=\sqrt{\widetilde{\varepsilon}_j}u(\cdot;\xi)$, since
$\widetilde{\varepsilon}_j^{\alpha-1}=\widetilde{\mu}_j$.
Thus \eqref{eq:general-tube-amplitude} implies
$\sqrt{\widetilde{\varepsilon}_j}\|u(\widetilde{t}_j;\xi)\|_\infty>2$.
Moreover, \eqref{eq:general-weighted-initial-bound} gives, for large $j$,
\begin{equation}\label{eq:general-good-event-weight}
 \|w_\beta\widetilde d(\xi)\|_{\ell^1(\Z^\nu)}
 \le\|w_\beta \widetilde h_j\|_{\ell^1(\Z^\nu)}+r_j
 <2<\sqrt{1+4/S}.
\end{equation}
Using $w_\beta(n)c(n)=w_{-\kappa}(n)$, we obtain
\[
 |g_n(\xi)|
 =\frac{|\widetilde d_n(\xi)|}{\sqrt{\widetilde{\varepsilon}_j}c(n)}
 \le\delta_{\widetilde{\varepsilon}_j}^{-1/2}w_\kappa(n),
 \qquad n\in\Z^\nu.
\]
Hence $\xi\in G_{\delta_{\widetilde{\varepsilon}_j}}$ by \eqref{eq:Omega-delta},
proving \eqref{eq:general-event-inclusion}.

\medskip
\noindent\emph{Translated Ball.}
The estimate \eqref{eq:counterexample-tube-probability} becomes
\begin{equation}\label{eq:general-translated-probability}
 \Pp(\widetilde d\in \widetilde h_j+\mathsf B_j)
 \ge e^{-\|\widetilde h_j\|_{\mathcal H}^2/\widetilde{\varepsilon}_j}
       \Pp(\widetilde d\in\mathsf B_j).
\end{equation}
To apply \autoref{lem:gaussian-tube}, replace its index set $\Z^2$
by $\Z^\nu$ and its weight by $w_\beta$ from
\eqref{eq:general-decay-parameters}.
The required integrability follows from
\begin{equation}\label{eq:general-gaussian-weighted-mean}
 \Ee\|w_\beta\widetilde d\|_{\ell^1(\Z^\nu)}
 =\frac{\sqrt{\pi\widetilde{\varepsilon}_j}}2
   \prod_{r=1}^\nu\sum_{m\in\Z}(1+|m|)^{-\kappa_r}<\infty,
\end{equation}
and
$\sum_{n\in\Z^\nu}w_\beta(n)^2|c(n)|^2
=\prod_{r=1}^\nu\sum_{m\in\Z}(1+|m|)^{-2\kappa_r}<\infty$.
The finite-dimensional density and symmetry calculation is unchanged.
For each fixed $j$, its truncation limits follow from probability
continuity from above and then from below, as in that lemma;
no interchange with $j\to\infty$ is needed.

\medskip
\noindent\emph{Centered Ball.}
Applying \emph{Markov's inequality} to
\eqref{eq:general-gaussian-weighted-mean} gives the counterpart of
\eqref{eq:counterexample-centered-tube}, that is, 
\begin{equation}\label{eq:general-centered-probability}
 \Pp(\widetilde d\notin\mathsf B_j)
 \lesssim_{\kappa,\nu}\frac{\sqrt{\widetilde{\varepsilon}_j}}{r_j}
 =8\widetilde{\mu}_j^{1/[2(\alpha-1)]-150\mathsf C_{\mathrm{fail}}}
 \longrightarrow0.
\end{equation}
Here the exponent is positive by
\eqref{eq:general-parameter-compatibility}.
Consequently, just as in \eqref{eq:counterexample-centered-probability},
\begin{equation}\label{eq:general-centered-log-probability}
 \Pp(\widetilde d\in\mathsf B_j)\longrightarrow1,
 \qquad
 \widetilde{\varepsilon}_j\log\Pp(\widetilde d\in\mathsf B_j)\longrightarrow0.
\end{equation}

\medskip
\noindent\emph{Lower Bound with a Strictly Positive Gap.}
Combining \eqref{eq:general-event-inclusion} and
\eqref{eq:general-translated-probability}, we obtain
\begin{align}
 \widetilde{\varepsilon}_j\log\Pp\!\left(
 \mathscr A_{\widetilde{\varepsilon}_j,2}^{\|u(\widetilde{t}_j)\|_\infty}\right)
 &\ge\widetilde{\varepsilon}_j\log\Pp\!\left(
 G_{\delta_{\widetilde{\varepsilon}_j}}\cap
 \mathscr A_{\widetilde{\varepsilon}_j,2}^{\|u(\widetilde{t}_j)\|_\infty}\right)
 \nonumber\\
 &\ge\widetilde{\varepsilon}_j\log\Pp(\widetilde d\in \widetilde h_j+\mathsf B_j)
 \nonumber\\
 &\ge-\|\widetilde h_j\|_{\mathcal H}^2
       +\widetilde{\varepsilon}_j\log\Pp(\widetilde d\in\mathsf B_j)
 =-1+o(1).
 \label{eq:general-log-probability}
\end{align}
The last equality follows from \eqref{eq:general-initial-cost} and
\eqref{eq:general-centered-log-probability}.
Taking the lower limit and using $S<4/3$ from
\eqref{eq:general-variance} proves
\eqref{eq:general-failure-lower-bound}.
\end{proof}

\begin{remark}[Logarithmic Time Coefficients]
\label{rem:counterexample-time-range}
We compare persistence and failure for the same data. Fix
\begin{equation}\label{eq:general-example-parameters}
 \rho_r=4\tau+\nu+6+r,\qquad 1<\kappa_r\le2,
 \qquad r=1,\ldots,\nu.
\end{equation}
Then $\beta_r=\rho_r-\kappa_r>\max\{4\tau+2,\nu+2\tau\}$
and $1\le S<4/3$.

For persistence at threshold $z=2$, choose
\[
 q=\nu+1,\qquad p=q+2\tau<\min_r\beta_r,\qquad R=\sqrt5.
\]
Indeed, \eqref{eq:log-normalized-initial-decay}, $A=1$ and $S\ge1$
give $\|\widetilde d\|_{\ell^\infty_p(\Z^\nu)}
\le\sqrt{1+4/S}\le R$ on $G_{\delta_\varepsilon}$.
Counting nonnegative coordinate sums and their sign choices gives
$\#\{n:|n|_1=m\}\le2^\nu\binom{m+\nu-1}{\nu-1}
\le2^\nu(m+1)^{\nu-1}$; see \cite[Theorem~13.1]{vLW01}.
Consequently,
\[
 \|\langle\cdot\rangle^{-q}\|_{\ell^1(\Z^\nu)}
 \le1+2^\nu\sum_{m\ge1}(m+1)^{-2}\le1+2^\nu.
\]
By \eqref{eq:log-explicit-growth-constant}, we may therefore take
\begin{equation}\label{eq:counterexample-explicit-log-constant}
 \Lambda=20\cdot3^{\nu+1+2\tau}(1+2^\nu)^2,
 \qquad 0<c<\Lambda^{-1}.
\end{equation}
This gives a sufficient range of persistence coefficients.

Next choose
\[
 \mathsf C_{\mathrm{fail}}=\frac{\sum_r\rho_r}{2(\nu-1)}+1.
\]
It satisfies the lower bound in
\eqref{eq:general-parameter-compatibility} and
$\mathsf C_{\mathrm{fail}}>1>\Lambda^{-1}>c$.
Thus \autoref{thm:general-failure} gives failure at this larger
coefficient, provided that
$1<\alpha<1+(300\mathsf C_{\mathrm{fail}})^{-1}$.

\end{remark}

\begin{corollary}[Impossibility of Polynomial Extension Preserving the Initial Gaussian Rate Function]
\label{cor:counterexample-polynomial-obstruction}
\label{cor:general-polynomial-failure}
Under the assumptions of \autoref{thm:general-failure}, every fixed
$T>0$ and $\vartheta>1$ satisfy
\begin{align}
 \limsup_{\varepsilon\to0^+}\sup_{0\le t\le T\mu_\varepsilon^{-\vartheta}}
 \varepsilon\log\Pp\!\left(
 \mathscr A_{\varepsilon,2}^{\|u(t)\|_\infty}\right)
 \ge-1>-\frac4S.
 \label{eq:counterexample-polynomial-rate-obstruction}
\end{align}
Thus the initial Gaussian rate function cannot persist uniformly on
$[0,T\mu_\varepsilon^{-\vartheta}]
=[0,T\varepsilon^{-\vartheta(\alpha-1)}]$.
\end{corollary}

\begin{proof}[Proof of \autoref{cor:counterexample-polynomial-obstruction}]
Fix $T>0$ and $\vartheta>1$.
For the sequence in \eqref{eq:general-embedding-parameters},
$\widetilde{\varepsilon}_j\to0$ and $\mu_{\widetilde{\varepsilon}_j}=\widetilde{\mu}_j$. Moreover,
\begin{equation}\label{eq:counterexample-polynomial-inclusion}
 \frac{\widetilde{t}_j}{T\widetilde{\mu}_j^{-\vartheta}}
 =\frac{\mathsf{C}_{\mathrm{fail}}}{T}\widetilde{\mu}_j^{\vartheta-1}|\log\widetilde{\mu}_j|
 \longrightarrow0.
\end{equation}
Thus $[0,\widetilde{t}_j]\subset[0,T\mu_{\widetilde{\varepsilon}_j}^{-\vartheta}]$
for all sufficiently large $j$. In particular,
\begin{equation}
 \sup_{0\le t\le T\mu_{\widetilde{\varepsilon}_j}^{-\vartheta}}
 \widetilde{\varepsilon}_j\log\Pp\!\left(
 \mathscr A_{\widetilde{\varepsilon}_j,2}^{\|u(t)\|_\infty}\right)
 \ge\widetilde{\varepsilon}_j\log\Pp\!\left(
 \mathscr A_{\widetilde{\varepsilon}_j,2}^{\|u(\widetilde{t}_j)\|_\infty}\right).
 \label{eq:counterexample-polynomial-supremum}
\end{equation}
Taking the upper limit along $\widetilde{\varepsilon}_j\to0$ and using
\eqref{eq:general-failure-lower-bound}, we obtain
\begin{align*}
 &\limsup_{\varepsilon\to0^+}\sup_{0\le t\le T\mu_\varepsilon^{-\vartheta}}
 \varepsilon\log\Pp\!\left(
 \mathscr A_{\varepsilon,2}^{\|u(t)\|_\infty}\right)\\
 &\qquad\ge\liminf_{j\to\infty}\widetilde{\varepsilon}_j\log\Pp\!\left(
 \mathscr A_{\widetilde{\varepsilon}_j,2}^{\|u(\widetilde{t}_j)\|_\infty}\right)
 \ge-1.
\end{align*}
Finally, $S<4/3$ by \eqref{eq:general-variance}, so $-1>-4/S$.
This proves \eqref{eq:counterexample-polynomial-rate-obstruction}.
\end{proof}

\section*{Acknowledgements}
The authors acknowledge the use of ChatGPT for exploratory discussions and
language refinement during the preparation of this work. All mathematical statements, proofs and references were independently verified and finalized by the authors.

\begingroup

\newcommand{\etalchar}[1]{$^{#1}$}
\providecommand{\bysame}{\leavevmode\hbox to3em{\hrulefill}\thinspace}
\providecommand{\MR}{\relax\ifhmode\unskip\space\fi MR }
\providecommand{\MRhref}[2]{  \href{http://www.ams.org/mathscinet-getitem?mr=#1}{#2}
}
\providecommand{\href}[2]{#2}

\endgroup

\end{document}